\documentclass{amsbook}

\usepackage{amssymb}
\usepackage{graphicx}
\allowdisplaybreaks

\newtheorem{theorem}{Theorem}[chapter]
\newtheorem{lemma}[theorem]{Lemma}
\newtheorem{proposition}[theorem]{Proposition}
\newtheorem{corollary}[theorem]{Corollary}
\newtheorem{remark}[theorem]{Remark}
\newtheorem{example}[theorem]{Example}
\newtheorem{observation}[theorem]{Observation}

\newtheorem{theoA}{Theorem}
\newtheorem{theoB}{Theorem}
\newtheorem{theoC}{Theorem}
\newtheorem{theoD}{Theorem}

\numberwithin{figure}{chapter}
\numberwithin{section}{chapter}
\numberwithin{equation}{chapter}

\newcommand{\N}{\mathbb{N}}
\newcommand{\Z}{\mathbb{Z}}
\newcommand{\R}{\mathbb{R}}
\newcommand{\C}{\mathbb{C}}

\newcommand{\EME}{\mathcal{M}}
\newcommand{\ENE}{\mathcal{N}}
\newcommand{\Ha}{\mathcal{H}}
\newcommand{\oM}{\otimes_{\mathcal{M}}}

\newcommand{\al}{\mathbb{\alpha}}

\newcommand{\ten}{\otimes}

\newcommand{\fin}{\hspace*{\fill} $\square$ \vskip0.2cm}
\newcommand{\prodd}{\prod\nolimits}
\newcommand{\limm}{\lim\nolimits}

\def\bubl{{\displaystyle \mathop{\mathsf{A}}^{\circ}}}
\def\bubla{{\displaystyle \mathop{a}^{\circ}}}

\newcommand{\summ}{\sum\nolimits}

\makeindex

\begin{document}

%\frontmatter
%\title{}
%\address{}
%\email{}
%\thanks{}
%\author{}
%\address{}
%\email{}
%\date{}
%\subjclass{}
%\maketitle

\pagestyle{empty}

\pagenumbering{roman}

\null

\vskip2cm

\begin{center} \LARGE \textbf{Theory of Amalgamated $L_p$
Spaces \\ in Noncommutative Probability } \\
\vskip2cm Marius Junge \\ \vskip-2pt {\large Department of
Mathematics} \\ \vskip-2pt {\large University of Illinois at
Urbana-Champaign} \\ \vskip-2pt {\large junge@math.uiuc.edu} \\
\vskip10pt {\Large and} \\ \vskip10pt Javier Parcet \\ \vskip-2pt
{\large Instituto de Matem{\'a}ticas y F{\'\i}sica Fundamental} \\
\vskip-2pt {\large Consejo Superior de Investigaciones
Cient{\'\i}ficas} \vskip-2pt {\large javier.parcet@uam.es}
\end{center}

\newpage

\null

\newpage

\null

\vskip3cm

\begin{center}
{\LARGE \textbf{Abstract}}
\end{center}

\

\renewcommand{\theequation}{$\Sigma_{pq}$}
\addtocounter{equation}{-1}

Let $f_1, f_2, \ldots, f_n$ be a family of independent copies of a
given random variable $f$ in a probability space $(\Omega,
\mathcal{F}, \mu)$. Then, the following equivalence of norms holds
whenever $1 \le q \le p < \infty$
\begin{equation}
\Big( \int_{\Omega} \Big[ \sum_{k=1}^n |f_k|^q \Big]^{\frac{p}{q}}
d \mu \Big)^{\frac1p} \sim \max_{r \in \{p,q\}} \left\{
n^{\frac1r} \Big( \int_\Omega |f|^r d\mu \Big)^{\frac1r} \right\}.
\end{equation}
We prove a noncommutative analogue of this inequality for sums of
free random variables over a given von Neumann subalgebra. This
formulation leads to new classes of noncommutative function spaces
which appear in quantum probability as square functions,
conditioned square functions and maximal functions. Our main tools
are Rosenthal type inequalities for free random variables,
noncommutative martingale theory and factorization of
operator-valued analytic functions. This allows us to generalize
$(\Sigma_{pq})$ as a result for noncommutative $L_p$ in the
category of operator spaces. Moreover, the use of free random
variables produces the right formulation of $(\Sigma_{\infty q})$,
which has not a commutative counterpart.

\renewcommand{\theequation}{\arabic{equation}}

\footnote{Junge is partially supported by the NSF DMS-0301116.}
\footnote{Parcet is partially supported by \lq Programa Ram{\'o}n y
Cajal, 2005\rq${}$ \\ \indent and also by Grants MTM2004-00678 and
CCG06-UAM/ESP-0286, Spain.} \footnote{2000 Mathematics Subject
Classification: 46L07, 46L09,46L51, 46L52, 46L53, 46L54.}

\pagestyle{empty}

\tableofcontents

%\mainmatter

\setcounter{page}{0}
\renewcommand{\thepage}{\arabic{page}}

\pagestyle{headings}

\chapter*{Introduction}

Probabilistic methods play an important role in harmonic analysis
and Banach space theory. Let us just mention the relevance of sums
of independent random variables, $p$-stable processes or
martingale inequalities in both fields. The analysis of subspaces
of the classical $L_p$ spaces is specially benefited from such
probabilistic notions. Viceversa, Burkholder's martingale
inequality for the conditional square function has been discovered
in view of Rosenthal's inequality for the norm in $L_p$ of sums of
independent random variables. This is only one example of the
fruitful interplay between harmonic analysis, probability theory
and Banach space geometry carried out mostly in the 70's by
Burkholder, Gundy, Kwapie\'n, Maurey, Pisier, Rosenthal and many
others.

\vskip5pt

More recently it became clear that a similar endeavor for
noncommutative $L_p$ spaces requires an additional insight from
quantum probability and operator space theory \cite{J0,J2,PS}. A
noncommutative theory of martingale inequalities finds its
beginnings in the work of Lust-Piquard \cite{Lu} and
Lust-Piquard/Pisier \cite{LuP} on the noncommutative Khintchine
inequality. The seminal paper of Pisier and Xu on the
noncommutative analogue of Burkholder-Gundy inequality \cite{PX1}
started a new trend in quantum probability. Nowadays, most
classical martingale inequalities have a satisfactory
noncommutative analogue, see \cite{J1,JX,PR,Ran1}. In the proof of
these results the classical stopping time arguments are no longer
available, essentially because point sets disappear after
quantization. These arguments are replaced by functional analytic
or combinatorial arguments. In the functional analytic approach we
often encounter new spaces. Indeed, maximal functions in the
noncommutative context can only be understood and defined through
analogy with vector-valued $L_p$ spaces. A careful analysis of
these spaces is crucial in establishing basic results such as
Doob's inequality \cite{J1} for noncommutative martingales and the
noncommutative maximal theorem behind Birkhoff's ergodic theorem
\cite{JX2}. The proof of maximal theorems and noncommutative
versions of Rosenthal's inequality often uses square function and
conditioned square function estimates, see \cite{JPX} and the
references therein. These are examples of more general classes of
noncommutative function spaces to be defined below. However, all
of them illustrate our main motto in this paper. Namely,
\emph{certain problems can be solved by finding and analyzing the
appropriate class of Banach spaces}. We shall develop in this
paper a new theory of \emph{generalized noncommutative $L_p$
spaces} with three problems in mind for a given von Neumann
algebra $\mathcal{A}$.

\vskip5pt

\noindent {\bf Problem 1.} \emph{Calculate the
$L_p(\mathcal{A};\ell_q)$ norm for sums of free random variables}.

\vskip5pt

\noindent {\bf Problem 2.} \emph{Any reflexive subspace of
$\mathcal{A}_*$ embeds in some $L_p$ for certain $p > 1$}.

\vskip5pt

\noindent {\bf Problem 3.} \emph{If $1 < p \le 2$, find a complete
embedding of $L_p(\mathcal{A})$ into some $L_1$ space}.

\vskip5pt

The main contribution of this paper is a detailed solution of
Problem 1, see below for more information. Unfortunately, Problems
2 and 3 are beyond the scope of this paper. However, we should
note that the solutions to both problems are deeply related to the
main results in this paper. Problem 2 is the noncommutative
version of Rosenthal's famous theorem \cite{Ro}. A certain version
of this result was obtained by Pisier in \cite{Pis}. However, he
uses a nonconventional interpolation space which does not seem to
fit in the usual scale of noncommutative $L_p$ spaces. More
precisely, one has to understand the interpolation space between a
von Neumann algebra and an intersection of two Hilbert spaces. The
interplay of interpolation and intersection is a the heart of this
article. In the semifinite case, the problem is solved by
Randrianantoanina \cite{Narcisse}. On the other hand, Problem 3 is
motivated by the classical notion of $p$-stable variables for $0 <
p \le 2$. That is, a sequence $(\xi_k)$ of independent random
variables such that
$$\mathbb{E} \, \exp \Big( \summ_k i t a_k \xi_k \Big) = \exp
\Big( -c \summ_k |a_k|^p \Big).$$ Indeed, a positive solution of
Problem 3 for $p = 2$ has been obtained in \cite{J2} using norm
estimates for sums of independent random variables. Let us briefly
explain this. The simplest model of $2$-stable variables is
provided by normalized gaussians $(g_k)$. In this particular case
and after taking operator coefficients $(a_k)$ in some
noncommutative $L_1$ space, the noncommutative Khintchine
inequality \cite{LuP} tells us that
\begin{equation} \label{intro-k1}
\mathbb{E} \, \Big\| \summ_k a_k g_k \Big\|_1 \sim \inf_{a_k =
r_k+c_k} \Big\| \big( \summ_k r_k r_k^* \big)^{\frac12} \Big\|_1 +
\Big\| \big( \summ_k c_k^* c_k \big)^{\frac12} \Big\|_1.
\end{equation}
Let us point out that operator space theory provides a very
appropriate framework for analyzing noncommutative $L_p$ spaces
and linear maps between them. Indeed, the inequality
\eqref{intro-k1} describes the operator space structure of the
subspace spanned by the $g_k$'s in $L_1$ as the sum $R+C$ of row
and column subspaces of $\mathcal{B}(\ell_2)$. We refer to
\cite{ER} and \cite{P3} for background information on operator
spaces. In the language of noncommutative probability many
operator space inequalities translate into module valued versions
of scalar inequalities, this will be further explained below. The
only drawback of \eqref{intro-k1} is that it does not coincide
with Pisier's definition of the operator space $\ell_2$
\begin{equation} \label{intro-oh}
\Big\| \summ_k a_k \ten \delta_k \Big\|_{L_1(\mathcal{M}; \ell_2)}
= \inf_{a_k = \alpha \gamma_k \beta} \|\alpha\|_{L_4(\mathcal{M})}
\Big( \summ_k \|\gamma_k\|_{L_2(\mathcal{M})}^2 \Big)^{\frac12}
\|\beta\|_{L_4(\mathcal{M})}.
\end{equation}
However, it was proved by Pisier that the right side in
\eqref{intro-oh} is obtained by complex interpolation between the
row and the column square functions appearing on the right of
\eqref{intro-k1}. One the main results in this article is a far
reaching generalization of this observation. In fact, the solution
of Problem 3 in full generality is closely related to this
analysis.

\vskip5pt

Following our guideline we will now introduce and discuss the new
class of spaces relevant for  these problems and martingale
theory. These generalize Pisier's theory of $L_p(L_q)$ spaces over
hyperfinite von Neumann algebras. We begin with a brief review of
some noncommutative function spaces which have lately appeared in
the literature, mainly in noncommutative martingale theory. We
refer to Chapter \ref{Section1} below for a more detailed
exposition.

\vskip19pt

\noindent \textbf{1. Noncommutative function spaces.}

\vskip5pt

Inspired by Pisier's theory \cite{P2}, several noncommutative
function spaces have been recently introduced in quantum
probability. The first motivation comes from some of Pisier's
fundamental equalities, which we briefly review. Let
$\mathcal{N}_1$ and $\mathcal{N}_2$ be two hyperfinite von Neumann
algebras. Then, given $1 \le p,q \le \infty$ and defining $1/r =
|1/p - 1/q|$, we have

\vskip5pt

\begin{itemize}
\item[i)] If $p \le q$, the norm of $x \in L_p \big(
\mathcal{N}_1; L_q(\mathcal{N}_2) \big)$ is given by
$$\inf \Big\{ \|\alpha\|_{L_{2r}(\mathcal{N}_1)}
\|y\|_{L_q(\mathcal{N}_1 \bar\otimes \mathcal{N}_2)}
\|\beta\|_{L_{2r}(\mathcal{N}_1)} \, \big| \ x = \alpha y \beta
\Big\}.$$ \item[ii)] If $p \ge q$, the norm of $x \in L_p \big(
\mathcal{N}_1; L_q(\mathcal{N}_2) \big)$ is given by
$$\sup \Big\{ \|\alpha x \beta \|_{L_q(\mathcal{N}_1
\bar\otimes \mathcal{N}_2)} \, \big| \ \alpha, \beta \in
\mathsf{B}_{L_{2r}(\mathcal{N}_1)} \Big\}.$$
\end{itemize}
On the other hand, the row and column subspaces of $L_p$ are
defined as follows
$$L_p(\mathcal{M}; R_p^n) = \Big\{ \sum_{k=1}^n x_k \ten e_{1k} \,
\big| \ x_k \in L_p(\mathcal{M}) \Big\} \subset L_p \big(
\mathcal{M} \bar\ten \mathcal{B}(\ell_2) \big),$$
$$L_p(\mathcal{M}; C_p^n) = \Big\{ \sum_{k=1}^n x_k \ten
e_{k1} \, \big| \ x_k \in L_p(\mathcal{M}) \Big\} \subset L_p
\big( \mathcal{M} \bar\ten \mathcal{B}(\ell_2) \big),$$ where
$(e_{ij})$ denotes the unit vector basis of $\mathcal{B}(\ell_2)$.
These spaces are crucial in the noncommutative
Khintchine/Rosenthal type inequalities \cite{JPX,LuP,PP} and in
noncommutative martingale inequalities \cite{JX,PX1,Ran1}, where
the row and column spaces are traditionally denoted by
$L_p(\mathcal{M}; \ell_2^r)$ and $L_p(\mathcal{M}; \ell_2^c)$.
Now, considering a von Neumann subalgebra $\mathcal{N}$ of
$\mathcal{M}$ with a normal faithful conditional expectation
$\mathsf{E}: \mathcal{M} \to \mathcal{N}$, we may define $L_p$
norms of the conditional square functions
$$\Big( \sum_{k=1}^n \mathsf{E}(x_k x_k^*) \Big)^{\frac12} \quad
\mbox{and} \quad \Big( \sum_{k=1}^n \mathsf{E}(x_k^* x_k)
\Big)^{\frac12}.$$ The expressions $\mathsf{E}(x_k x_k^*)$ and
$\mathsf{E}(x_k^* x_k)$ have to be defined properly for $1 \le p
\le 2$, see \cite{J1} or Chapter \ref{Section1} below. Note that
the resulting spaces coincide with the row and column spaces
defined above when $\mathcal{N}$ is $\mathcal{M}$ itself. When
$n=1$ we recover the spaces $L_p^r(\mathcal{M}, \mathsf{E})$ and
$L_p^c(\mathcal{M}, \mathsf{E})$, which have been instrumental in
proving Doob's inequality \cite{J1}, see also \cite{JM,JX2} for
more applications.

\vskip19pt

\noindent \textbf{2. Amalgamated $L_p$ spaces}

\vskip5pt

The definition of amalgamated $L_p$ spaces is algebraic. We recall
that by H\"{o}lder's inequality $L_u(\mathcal{M}) L_q(\mathcal{M})
L_v(\mathcal{M})$ is contractively included in $L_p(\mathcal{M})$
when $1/p = 1/u + 1/q + 1/v$. Let us now assume that $\mathcal{N}$
is a von Neumann subalgebra of $\mathcal{M}$ with a normal
faithful conditional expectation $\mathsf{E}:\mathcal{M}\to
\mathcal{N}$. Then we have natural isometric inclusions
$L_s(\mathcal{N}) \subset L_s(\mathcal{M})$ for $0<s\le \infty$
and we may consider the \emph{amalgamated} $L_p$ space
\[ L_u(\mathcal{N})L_q(\mathcal{M})L_v(\mathcal{N})  \]
as the subset of elements $x$ in $L_p(\mathcal{M})$ which
factorize as $x = \alpha y \beta$ with $\alpha \in
L_u(\mathcal{N})$, $y \in L_q(\mathcal{M})$ and $\beta \in
L_v(\mathcal{N})$. The natural ``norm'' is then given by the
following expression $$\|x\|_{u \cdot q \cdot v} = \inf \Big\{
\|\alpha\|_{L_u(\mathcal{N})} \|y\|_{L_q(\mathcal{M})}
\|\beta\|_{L_v(\mathcal{N})} \, \big| \ x = \alpha y \beta
\Big\}.$$ However, the triangle inequality for the homogeneous
expression $\| \ \|_{u \cdot q \cdot v}$ is by no means trivial.
Moreover, it is not clear a priori that this subset of
$L_p(\mathcal{M})$ is indeed a linear space. Before explaining
these difficulties in some detail, let us consider some examples.
We fix an integer $n \ge 1$ and the subalgebra $\mathcal{N}$
embedded in the diagonal of the direct sum $\mathcal{M} =
\mathcal{N} \oplus_\infty \mathcal{N} \oplus_\infty \cdots
\oplus_\infty \mathcal{N}$ with $n$ terms. The natural conditional
expectation is $$\mathsf{E}_n \big(x_1, \ldots, x_n \big) =
\frac1n \sum_{k=1}^n x_k.$$ Then it is easy to see that for $1 \le
p \le 2$ and $1/p = 1/2 + 1/w$ we have
\begin{eqnarray*}
L_p(\mathcal{N}; R_p^n) & = & \sqrt{n} \, L_w(\mathcal{N})
L_2(\mathcal{M}) L_{\infty}(\mathcal{N}), \\ L_p(\mathcal{N};
C_p^n) & = & \sqrt{n} \, L_{\infty}(\mathcal{N}) L_2(\mathcal{M})
L_w(\mathcal{N}),
\end{eqnarray*}
isometrically. Here we use the notation $\gamma \mathrm{X}$ to
denote the space $\mathrm{X}$ equipped with the norm $\gamma \|
\cdot \|_{\mathrm{X}}$. At the time of this writing and with
independence of this paper, a result of Pisier \cite{P0} on
interpolation of these spaces for $p = \infty$ was generalized by
Xu \cite{X} for arbitrary $p$'s
$$\Big\| \sum_{k=1}^n x_k \otimes \delta_k \Big\|_{\left[
L_p(\mathcal{N}; R_p^n), L_p(\mathcal{N};C_p^n) \right]_{\theta}}
\, =\,  \inf_{x_k = \alpha y_k \beta} \|\alpha\|_{u_\theta} \Big(
\sum_{k=1}^n \|y_k\|_2^2 \Big)^{\frac12} \|\beta\|_{v_{\theta}}.$$
Here $(1/u_{\theta}, 1/v_\theta) = ((1-\theta)/w, \theta/w)$ and
for $\theta = 1/2$ we find Pisier's definition of
$L_p(\mathcal{N};\ell_2)$. That is, we obtain the space $\sqrt{n}
\, L_{u_\theta}(\mathcal{N}) L_2(\mathcal{M})
L_{v_\theta}(\mathcal{N})$. Our definition is flexible enough to
accommodate the conditional square function. Indeed, given a von
Neumann subalgebra $\mathcal{N}_0$ of $\mathcal{N}$ with a normal
faithful conditional expectation $\mathcal{E}_0: \mathcal{N} \to
\mathcal{N}_0$, we find
\begin{eqnarray*}
\Big\| \sum_{k=1}^n x_k \otimes \delta_k
\Big\|_{L_w(\mathcal{N}_0) L_2(\mathcal{M})
L_\infty(\mathcal{N}_0)} & = & n^{-\frac12} \Big\| \big(
\sum_{k=1}^n \mathcal{E}_0(x_k x_k^*)
\big)^{\frac12} \Big\|_{L_p(\mathcal{N}_0)}, \\
\Big\| \sum_{k=1}^n x_k \otimes \delta_k
\Big\|_{L_\infty(\mathcal{N}_0) L_2(\mathcal{M})
L_w(\mathcal{N}_0)} & = & n^{-\frac12} \Big\| \big( \sum_{k=1}^n
\mathcal{E}_0(x_k^* x_k) \big)^{\frac12}
\Big\|_{L_p(\mathcal{N}_0)}.
\end{eqnarray*}
Xu's interpolation does not apply in this more general setting,
which appears in the context of the noncommutative Rosenthal
inequality. This illustrates how certain amalgamated $L_p$ spaces
occur naturally in quantum probability. Now we want to understand
for which range of indices $(u,q,v)$ we have the triangle
inequality. In fact, our proof intertwines with the proof of our
main interpolation result which can be stated as follows. Let us
consider the solid $\mathsf{K}$ in $\R^3$ defined by
$$\mathsf{K} = \Big\{(1/u,1/v,1/q) \, \big| \ 2 \le u,v \le
\infty, \ 1 \le q \le \infty, \ 1/u + 1/q + 1/v \le 1 \Big\}.$$

\vskip2pt

\begin{theoA} \label{theoA} The amalgamated space
$L_u(\mathcal{N}) L_q(\mathcal{M}) L_v(\mathcal{N})$ is a Banach
space for any $(1/u,1/v,1/q) \in \mathsf{K}$. Moreover, if
$(1/u_j,1/v_j,1/q_j) \in \mathsf{K}$ for $j=0,1$, the space
$L_{u_{\theta}}(\mathcal{N}) L_{q_{\theta}}(\mathcal{M})
L_{v_{\theta}}(\mathcal{N})$ is isometrically isomorphic to
$$\Big[L_{u_0}(\mathcal{N})
L_{q_0}(\mathcal{M}) L_{v_0}(\mathcal{N}), L_{u_1}(\mathcal{N})
L_{q_1}(\mathcal{M}) L_{v_1}(\mathcal{N})
\Big]_{\theta}^{\null}.$$
\end{theoA}

The triangle inequality follows from Theorem A. On the other hand
we will need the triangle inequality in order to apply
interpolation and factorization techniques in proving Theorem A.
This intriguing interplay makes our proof quite involved. Our
\emph{first step} is showing that the triangle inequality holds in
the boundary region $$\partial_{\infty} \mathsf{K} = \Big\{
(1/u,1/v,1/q) \in \mathsf{K} \, \big| \ \min \big\{ 1/u, 1/q, 1/v
\big\} = 0 \Big\}.$$ Our argument uses the operator-valued
analogue of Szeg\"{o}'s factorization theorem, a technique which
will be used repeatedly throughout this paper. The triangle
inequality for other indices follows by convexity since
$\mathsf{K}$ is the convex hull of $\partial_{\infty} \mathsf{K}$,
so that any other point in $\mathsf{K}$ is associated to an
interpolation space between two spaces living in
$\partial_{\infty} \mathsf{K}$. Another technical difficulty is
the fact that the intersection of two amalgamated $L_p$ spaces is
in general quite difficult to describe. Thus, any attempt to use a
density argument meets this obstacle. The \emph{second step} is to
prove Theorem A for finite von Neumann algebras, where the
intersections are easier to handle. Moreover, most of the
factorization arguments (as Szeg\"{o}'s theorem) a priori only
apply in the finite setting. In the \emph{third step} we consider
general von Neumann algebras using Haagerup's crossed product
construction \cite{H2} to approximate $\sigma$-finite von Neumann
algebras by direct limits of finite von Neumann algebras. Finally,
we need a different argument for the case $\min(q_0,q_1) =
\infty$, which is out of the scope of Haagerup's construction. The
main technique here is a Grothendieck-Pietsch version of the
Hahn-Banach theorem.

\vskip5pt

Let us observe that in the hyperfinite case Pisier was able to
establish many of his results using the Haagerup tensor product.
Though similar in nature, we can not directly use tensor product
formulas for our interpolation results due to the complicated
structure of general von Neumann subalgebras. Theorem \ref{theoA}
will also be useful in understanding certain interpolation spaces
in martingale theory. Let us mention some open problems, for
partial results see Chapter \ref{Section5} below.

\vskip5pt

\noindent \textbf{Problem 4.} \emph{Let $\mathcal{M}$ be a von
Neumann algebra and denote by $\mathcal{H}_p^r(\mathcal{M})$ and
$\mathcal{H}_p^c(\mathcal{M})$ the row and column Hardy spaces of
noncommutative martingales over $\mathcal{M}$. Let us consider an
interpolation parameter $0 < \theta < 1$.}
\begin{itemize}
\item[\textbf{(a)}] \emph{Calculate the interpolation norms
$[\mathcal{H}_p^r(\mathcal{M}),
\mathcal{H}_p^c(\mathcal{M})]_{\theta}$.}

\item[\textbf{(b)}] \emph{If $x \in [\mathcal{H}_1^r(\mathcal{M}),
\mathcal{H}_1^c(\mathcal{M})]_{\theta}$, the maximal function is
in $L_1$}.
\end{itemize}

\vskip19pt

\noindent \textbf{3. Conditional \boldmath $L_p$ \unboldmath
spaces}

\vskip5pt

Once we know which amalgamated $L_p$ spaces are Banach spaces it
is natural to investigate their dual spaces. We assume as above
that $\mathcal{N}$ is a von Neumann subalgebra of $\mathcal{M}$
and $\mathsf{E}:\mathcal{M} \to \mathcal{N}$ is a normal faithful
conditional expectation. Let
$$1/s = 1/u+1/p+1/v \le 1.$$ The \emph{conditional} $L_p$ space
$$L_{(u,v)}^p(\mathcal{M}, \mathsf{E})$$ is defined as the
completion of $L_p(\mathcal{M})$ with respect to the norm
$$\|x\|_{L_{(u,v)}^p(\mathcal{M}, \mathsf{E})} = \sup \Big\{
\|axb\|_{L_s(\mathcal{M})} \, \big| \ \|a\|_{L_u(\mathcal{N})},
\|b\|_{L_v(\mathcal{N})} \le 1 \Big\}.$$ In our next result we
show that amalgamated and conditional $L_p$ are related by
anti-linear duality. This will allow us to translate the
interpolation identities in Theorem A in terms of conditional
$L_p$ spaces. In this context the correct set of parameters is
given by
$$\mathsf{K}_0 = \Big\{ (1/u,1/v,1/q) \in \mathsf{K} \, \big| \ 2
< u,v \le \infty, \ 1 < q < \infty, \ 1/u + 1/q + 1/v < 1
\Big\}.$$ %\pagebreak[3]

\begin{theoB}
Let $1 < p < \infty$ given by $1/q' = 1/u + 1/p + 1/v$, where the
indices $(u,q,v)$ belong to the solid $\mathsf{K}_0$ and $q'$ is
conjugate to $q$. Then, the following isometric isomorphisms hold
via the anti-linear duality bracket $\langle x,y \rangle =
\mathrm{tr}(x^* y)$ $$\big( L_u(\mathcal{N}) L_q(\mathcal{M})
L_v(\mathcal{N}) \big)^* = L_{(u,v)}^p(\mathcal{M}, \mathsf{E}),$$
$$\big( L_{(u,v)}^p(\mathcal{M}, \mathsf{E}) \big)^* =
L_u(\mathcal{N}) L_q(\mathcal{M}) L_v(\mathcal{N}).$$ In
particular, we obtain isometric isomorphisms $$\Big[
L_{(u_0,v_0)}^{p_0}(\mathcal{M}, \mathsf{E}),
L_{(u_1,v_1)}^{p_1}(\mathcal{M}, \mathsf{E}) \Big]_{\theta} =
L_{(u_\theta,v_\theta)}^{p_\theta}(\mathcal{M}, \mathsf{E}).$$
\end{theoB}

\vskip5pt

As we shall see, Theorem B generalizes the interpolation results
obtained by Pisier \cite{P0} and Xu \cite{X} mentioned above.
Pisier and Xu's results provide an explicit expression for the
operator space structure of $[C_p^n,R_p^n]_{\theta}$ with $0 <
\theta < 1$. Theorem B also provides explicit formulas for
$[\ell_p^n, C_p^n]_{\theta}$ and $[\ell_p^n, R_p^n]_{\theta}$. In
fact, a large variety of interesting formulas of this kind arise
from Theorem B. A detailed analysis of these applications is out
of the scope of this paper. On the other hand, the analogue of
Theorem B for $p=\infty$ (which we will investigate separately)
has been already applied to study the noncommutative
John-Nirenberg theorem \cite{JM}.

\vskip5pt

Exactly as it happens with amalgamated $L_p$ spaces, several
noncommutative function spaces arise as particular forms of
conditional $L_p$ spaces. Let us review the basic examples in both
cases.

\vskip5pt

\begin{itemize} \item[\textbf{(a)}] The spaces
$L_p(\mathcal{M})$ satisfy
$$\quad L_p(\mathcal{M}) = L_{\infty}(\mathcal{N}) L_p(\mathcal{M})
L_{\infty}(\mathcal{N}) \quad \mbox{and} \quad L_p(\mathcal{M}) =
L_{(\infty,\infty)}^p(\mathcal{M}, \mathsf{E}).$$

\item[\textbf{(b)}] The spaces $L_p(\mathcal{N}_1;
L_q(\mathcal{N}_2))$:

\vskip3pt

\begin{itemize} \item[$\bullet$] Let  $p \le q$ and $1/r =
1/p - 1/q$. Then  $$L_p(\mathcal{N}_1; L_q(\mathcal{N}_2)) =
L_{2r}(\mathcal{N}_1) L_q(\mathcal{N}_1 \bar\otimes \mathcal{N}_2)
L_{2r}(\mathcal{N}_1).$$

\item[$\bullet$] Let $p \ge q$ and $1/r = 1/q - 1/p$. Then
$$L_p(\mathcal{N}_1; L_q(\mathcal{N}_2)) =
L_{(2r,2r)}^p(\mathcal{N}_1 \bar\otimes \mathcal{N}_2,
\mathsf{E}),$$ where $\mathsf{E}: \mathcal{N}_1 \bar\otimes
\mathcal{N}_2 \to \mathcal{N}_1$ is given by $\mathsf{E} =
1_{\mathcal{N}_1} \otimes \varphi_{\mathcal{N}_2}$.
\end{itemize}

\vskip5pt

\item[\textbf{(c)}] The spaces $L_p^r(\mathcal{M}, \mathsf{E})$
and $L_p^c(\mathcal{M}, \mathsf{E})$:

\vskip3pt

\begin{itemize}
\item[$\bullet$] Let  $1 \le p \le 2$ and $1/p = 1/2 + 1/s$. Then
\begin{align*}
L_p^r(\mathcal{M}, \mathsf{E}) & = \, L_s(\mathcal{N})
L_2(\mathcal{M}) L_{\infty}(\mathcal{N}) ,\\
% \quad ,\quad
 L_p^c(\mathcal{M}, \mathsf{E}) & =   L_{\infty}(\mathcal{N})
L_2(\mathcal{M}) L_s(\mathcal{N}).
\end{align*}

\item[$\bullet$] Let  $2 \le p \le \infty$ and $1/p + 1/s = 1/2$.
Then
\begin{align*} L_p^r(\mathcal{M},
 \mathsf{E}) & = \, L_{(s,\infty)}^p(\mathcal{M}, \mathsf{E}) \\
 %\quad,\quad
 L_p^c(\mathcal{M}, \mathsf{E}) & = \,
L_{(\infty,s)}^p(\mathcal{M}, \mathsf{E}).
\end{align*}
\end{itemize}
In particular, taking $\mathsf{E}_n(x_1,...,x_n) = \frac1n \sum_k
x_k$ we find
\begin{eqnarray*}
L_p(\mathcal{M}; R_p^n) & = & \sqrt{n} \,
L_p^r(\ell_{\infty}^n(\mathcal{M}), \mathsf{E}_n),
\\ L_p(\mathcal{M}; C_p^n) & = &
\sqrt{n} \, L_p^c(\ell_{\infty}^n(\mathcal{M}), \mathsf{E}_n).
\end{eqnarray*}

\item[\textbf{(d)}] As we shall see through the text, asymmetric
$L_p$ spaces (a non-standard operator space structure on $L_p$
which will be crucial in this paper) also have representations in
terms of amalgamated or conditional $L_p$ spaces.
\end{itemize}

\vskip19pt

\noindent \textbf{4. Intersection spaces}

\vskip5pt

Intersection of $L_p$ spaces appear naturally in the theory of
noncommutative Hardy spaces. These spaces are also natural
byproducts of Rosenthal's inequality for sums of independent
random variables. Let us first illustrate this point in the
commutative setting and then provide the link to the spaces
defined above. Let us consider a finite collection $f_1, f_2,
\ldots, f_n$ of independent random variables on a probability
space $(\Omega, \mathcal{F}, \mu)$. The Khintchine inequality
implies for $0 < p < \infty$ \[ \Big( \int_\Omega \Big[
\sum_{k=1}^n |f_k|^2 \Big]^{\frac{p}{2}}  d\mu \Big)^{\frac{1}{p}}
\sim_{c_p} \mathbb{E} \, \Big\| \sum_{k=1}^n \varepsilon_k f_k
\Big\|_p\, .\] Therefore, Rosenthal's inequality \cite{Ro0} gives
for $2 \le p < \infty$
\renewcommand{\theequation}{$\Sigma_{p2}$}
\addtocounter{equation}{-1}
 \begin{equation}
 \Big( \int_\Omega \Big[ \sum_{k=1}^n |f_k|^2
 \Big]^{\frac{p}{2}}  d\mu \Big)^{\frac{1}{p}} \sim_{c_p}
 \Big( \sum_{k=1}^n \|f_k\|_p^p \Big)^{\frac{1}{p}} +
 \Big( \sum_{k=1}^n \|f_k\|_2^2 \Big)^{\frac{1}{2}}. \label{int-ros}
\end{equation}\renewcommand{\theequation}{\arabic{equation}}

\vspace{-.2cm}\noindent Here $(\varepsilon_k)$ is an independent
sequence of Bernoulli random variables equidistributed on $\pm 1$.
We can easily generalize this result for calculating $\ell_q$ sums
of independent random variables. Indeed, consider $1 \le q \le p <
\infty$ and define $g_1, g_2, \ldots, g_n$ by the relation $g_k =
|f_k|^{q/2}$ for $1 \le k \le n$. Then we have the following
identity for the index $s = 2p/q$
$$\Big( \int_\Omega \Big[ \sum_{k=1}^n |f_k|^q \Big]^{\frac{p}{q}}
d\mu \Big)^{\frac1p} = \Big( \int_\Omega \Big[ \sum_{k=1}^n
|g_k|^2 \Big]^{\frac{s}{2}} d\mu \Big)^{\frac{2}{qs}}.$$ Therefore
\eqref{int-ros} implies
\renewcommand{\theequation}{$\Sigma_{pq}$}
\addtocounter{equation}{-1}
\begin{equation}
\Big( \int_{\Omega} \Big[ \sum_{k=1}^n |f_k|^q \Big]^{\frac{p}{q}}
d \mu \Big)^{\frac1p} \sim_{c_p} \max_{r \in \{p,q\}} \left\{
\Big( \sum_{k=1}^n \int_{\Omega} |f_k|^r d \mu \Big)^{\frac1r}
\right\}.
\end{equation}
\renewcommand{\theequation}{\arabic{equation}}

\vskip-5pt

\noindent In particular, Rosenthal's inequality provides a natural
realization of $$\mathcal{J}_{p,q}^n(\Omega) = n^{\frac1p}
L_p(\Omega) \cap n^{\frac1q} L_q(\Omega)$$ into
$L_p(\Omega;\ell_q^n)$. More precisely, if $f_1, f_2, \ldots, f_n$
are taken to be independent copies of a given random variable $f$,
the right hand side of $(\Sigma_{pq})$ is the norm of $f$ in the
intersection space $\mathcal{J}_{p,q}^n(\Omega)$ and inequality
($\Sigma_{pq}$) provides an isomorphic embedding of
$\mathcal{J}_{p,q}^n(\Omega)$ into the space $L_p(\Omega;
\ell_q^n)$.

\vskip5pt

Quite surprisingly, replacing independent variables by matrices of
independent variables in $(\Sigma_{pq})$ requires to intersect
\emph{four} spaces using the so-called \emph{asymmetric} $L_p$
spaces. In other words, the \emph{natural} operator space
structure of $\mathcal{J}_{p,q}^n$ comes from a $4$-term
intersection space. We have already encountered such a phenomenon
in \cite{JP} for the case $q=1$. To justify this point, instead of
giving precise definitions we note that H\"{o}lder inequality
gives $L_p = L_{2p} L_{2p}$, meaning that the $p$-norm of $f$ is
the infimum of $\|g\|_{2p} \|h\|_{2p}$ over all factorizations $f
= g h$. If $L_p^r$ and $L_p^c$ denote the row and column
quantizations of $L_p$ (see Chapter \ref{Section1} for the
definition), the operator space analogue of the isometry above is
given by the complete isometry $L_p = L_{2p}^r L_{2p}^c$, see
Chapter \ref{Section7} for more details. In particular, according
to the algebraic definition of $L_p(\ell_q)$, the space
$\mathcal{J}_{p,q}^n$ has to be redefined as the product
$$\mathcal{J}_{p,q}^n = \Big( n^{\frac{1}{2p}} L_{2p}^r \cap
n^{\frac{1}{2q}} L_{2q}^r \Big) \Big( n^{\frac{1}{2p}} L_{2p}^c
\cap n^{\frac{1}{2q}} L_{2q}^c \Big).$$ We shall see in this paper
that $$\mathcal{J}_{p,q}^n = n^{\frac{1}{p}} L_{2p}^rL_{2p}^c \cap
n^{\frac{1}{2q} + \frac{1}{2p}} L_{2q}^r L_{2p}^c \cap
n^{\frac{1}{2p}+\frac{1}{2q}} L_{2p}^r L_{2q}^c \cap
n^{\frac{1}{q}} L_{2q}^rL_{2q}^c.$$ On the Banach space level we
have the isometries $$L_{2p}^r L_{2q}^c = L_s = L_{2q}^r L_{2p}^c
\quad \mbox{with} \quad 1/s = 1/2p + 1/2q.$$ Moreover, again by
H\"{o}lder inequality it is clear that
$$n^{\frac1s} \|f\|_s \le \max \Big\{ n^{\frac1p}
\|f\|_p, n^{\frac1q} \|f\|_q \Big\}.$$ Therefore, the two cross
terms in the middle disappear in the Banach space level. However,
replacing scalars by operators in the context of
independence/freness over a given von Neumann subalgebra, the
Banach space estimates from above are no longer valid and all four
terms may have a significant contribution.

\vskip5pt

It is the operator space structure of $\mathcal{J}_{p,q}^n$ what
originally led us to introduce amalgamated and conditional $L_p$
spaces. To be more precise, we consider a von Neumann algebra
$\mathcal{M}$ equipped with a normal faithful state $\varphi$ and
a von Neumann subalgebra $\mathcal{N}$ with associated normal
faithful conditional expectation $\mathsf{E}$. Then, if we fix $1
\le q \le p \le \infty$, we define
$$\mathcal{J}_{p,q}^n(\mathcal{M}, \mathsf{E}) = \bigcap_{u,v \in
\{2r, \infty\}} n^{\frac1u + \frac1p + \frac1v} \,
L_{(u,v)}^p(\mathcal{M}, \mathsf{E}) \qquad \mbox{with} \qquad 1/r
= 1/q - 1/p.$$ This definition is motivated by the fundamental
isometry
\begin{equation} \label{Eq-ossss}
S_p^m \big( \mathcal{J}_{p,q}^n(\mathcal{M}) \big) =
\mathcal{J}_{p,q}^n \big( \mathrm{M}_m \otimes \mathcal{M},
1_{\mathrm{M}_m} \otimes \varphi \big),
\end{equation}
which will be proved in Chapter \ref{Section7}. Some preliminary
results on $\mathcal{J}_{p,q}^n(\mathcal{M})$ (and vector-valued
generalizations) are already contained in the recent paper
\cite{JP}. We extend many results from \cite{JP} to the realm of
free random variables including the limit case $p=\infty$. Our
main result for the spaces $\mathcal{J}_{p,q}^n(\mathcal{M},
\mathsf{E})$ shows that we have an interpolation scale with
respect to the index $1 \le q \le p$.

\begin{theoC} \label{theoC}
If $1 \le p \le \infty$, then $$\big[
\mathcal{J}_{p,1}^n(\mathcal{M}, \mathsf{E}),
\mathcal{J}_{p,p}^n(\mathcal{M}, \mathsf{E}) \big]_{\theta} \simeq
\mathcal{J}_{p,q}^n(\mathcal{M}, \mathsf{E})$$ with $1/q =
1-\theta + \theta/p$ and with relevant constants independent of
$n$.
\end{theoC}

There seems to be no general argument to make intersections
commute with complex interpolation. For commutative $L_p$ spaces
or rearrangement invariant spaces one can often find concrete
formulas of the resulting interpolation norms, see \cite{LT}.
However, in the noncommutative context these arguments are no
longer valid and we need genuinely new tools. Theorem \ref{theoC}
is the new key ingredient in the positive solution of Problem 2.

\vskip19pt

\noindent \textbf{5. Main embedding theorem}

\vskip5pt

The central role played by Rosenthal inequality partially
justifies why the index $p$ must be finite in the commutative form
of ($\Sigma_{pq}$). This also happens in \cite{JP}, where we used
the noncommutative Rosenthal inequality from \cite{JX}. However,
mainly motivated by Problems 2 and 3, one of our main goals in
this paper is to obtain a \emph{right} formulation of
($\Sigma_{\infty q}$). As in several other inequalities involving
independent random variables, such as the noncommutative
Khintchine inequalities, the limit case as $p \to \infty$ holds
when replacing classical independence by Voiculescu's concept of
freeness \cite{VDN}. Therefore, it is not surprising that we shall
use in our proof the free analogue of Rosenthal inequality
\cite{JPX}. Following $(\Sigma_{pq})$ we have a natural candidate
for a complemented  embedding of $\mathcal{J}_{p,q}^n$ in
$L_p(\ell_q^n)$ using free probability.

\vskip5pt

We define $\mathsf{A}_k$ to be the direct sum $\mathcal{M} \oplus
\mathcal{M}$ for $1 \le k \le n$. Then we consider the reduced
amalgamated free product $\mathcal{A} = *_{\mathcal{N}}
\mathsf{A}_k$ where the conditional expectation
$\mathsf{E}_{\mathcal{N}}: \mathcal{A} \rightarrow \mathcal{N}$
has the form $\mathsf{E}_{\mathcal{N}}(x_1,x_2) = \frac{1}{2}
\big( \mathsf{E}(x_1) + \mathsf{E}(x_2) \big)$ when restricted to
the algebra $\mathsf{A}_k$.  Let $\pi_k: \mathsf{A}_k \to
\mathcal{A}$ denote the natural embedding of $\mathsf{A}_k$ into
$\mathcal{A}$. Moreover, given $x \in \mathcal{M}$ we shall write
$x_k$ as an abbreviation of $\pi_k(x,-x)$. Note that $x_k$ is a
mean-zero element for $1 \le k \le n$. Our main embedding result
reads as follows.

\begin{theoD}\label{theoD} Let  $1 \le q \le p \le \infty$. The map
$$u: x \in \mathcal{J}_{p,q}^n (\mathcal{M}, \mathsf{E}) \mapsto
\sum_{k=1}^n x_k \otimes \delta_k \in L_p(\mathcal{A}; \ell_q^n)$$
is an isomorphism with complemented image. The  constants are
independent of $n$.
\end{theoD}

The map $u$ is of course reminiscent of the fundamental mappings
employed in \cite{J0,J2,JP} constructing certain embeddings of
$L_p$ spaces. Theorem \ref{theoC} follows as a consequence of
Theorem \ref{theoD} using the fact that
\[ L_p(\mathcal{A},\ell_q^n) = \big[ L_p(\mathcal{A},\ell_1^n),
L_p(\mathcal{A},\ell_p^n) \big]_{\theta}.\] We have tried in vain
to prove Theorem \ref{theoD} directly using uniquely tools from
free probability. The methods around Voiculescu's inequality seem
to work perfectly fine in the limit case $p=\infty$, for which
there is no commutative version. However, basic tools from free
harmonic analysis are still missing for a direct proof of Theorem
\ref{theoD}. Instead, apart from the free analogue \cite{JPX} of
Rosenthal inequality, we also use factorization techniques and
interpolation results for noncommutative Hardy and BMO spaces, see
chapter \ref{Section5} for further details.

\vskip5pt

The interested reader might be surprised that we have not
formulated our results in category of operator spaces. However, as
so often in martingale theory these results are automatic,
provided the spaces carry the correct operator space structure,
given in this case by \eqref{Eq-ossss} and $$S_p^m \Big( L_p \big(
*_k \mathsf{A}_k ; \ell_q^n \big) \Big) = L_p \Big( \mathrm{M}_m
\otimes (*_k \mathsf{A}_k) ; \ell_q^n \Big) = L_p(\mathcal{A};
\ell_q^n),$$ where the von Neumann algebra $\mathcal{A}$ is now
given by $$\mathcal{A} = *_{\mathrm{M}_m} \widehat{\mathsf{A}}_k
\quad \mbox{with} \quad \widehat{\mathsf{A}}_k = \mathrm{M}_m
\otimes \mathsf{A}_k = \big( \mathrm{M}_m \otimes \mathcal{M}
\big) \oplus \big( \mathrm{M}_m \otimes \mathcal{M} \big).$$ The
trick is to write matrix norms as modular versions of scalar
norms. Theorem \ref{theoD} is already formulated in its modular
version. Therefore, we may replace bounded (complemented) by
cb-bounded (cb-complemented) for free. Theorem \ref{theoD} shall
be the starting point for a positive solution of Problem 3 in a
forthcoming paper.

\vskip19pt

\noindent \textbf{Background and notation.} We shall assume
certain familiarity with the theory of von Neumann algebras. Other
branches of operator algebra that are central for us are operator
space theory, noncommutative integration and free probability. We
will define some basic concepts of these theories along the paper.
We shall also use quite frequently Calder{\'o}n's complex
interpolation method \cite{BL}, Haagerup's approximation theorem
\cite{H2}, the Grothendieck-Pietsch separation argument \cite{P}
or Raynaud's results on ultraproducts of noncommutative $L_p$
spaces \cite{Ra}. Some of these techniques will be introduced
along the text. The inner products and duality brackets $\langle \
, \, \rangle$ will be anti-linear on the first component and
linear on the second one. Apart from this and the terminology
introduced along the paper, we shall use standard notation from
the literature such as e.g. \cite{P3,Ta3}.

\vskip19pt

\noindent \textbf{Acknowledgements.} The first-named author was
partially supported by the NSF DMS-0301116. The second-named
author was partially supported by the Project MTM2004-00678,
Spain. This work was mostly carried out in a one-year visit of the
second-named author to the University of Illinois at
Urbana-Champaign. The second-named author would like to thank the
Math Department of the University of Illinois for its support and
hospitality.

\chapter{Noncommutative integration}
\label{Section1}

\numberwithin{equation}{chapter}

In this chapter we review some basic notions on noncommutative
integration that will be frequently used through out this paper.
We begin by recalling Haagerup and Kosaki's constructions of
noncommutative $L_p$ spaces. Then we briefly introduce Pisier's
theory of vector-valued noncommutative $L_p$ spaces, giving some
emphasis to those aspects which are relevant in this work.
Finally, we analyze some basic properties of certain $L_p$ spaces
associated to a conditional expectation, which were recently
introduced in the literature and are basic for our further
purposes. We shall assume some familiarity with von Neumann
algebras. Basic concepts such as trace, state, commutant,
affiliated operator, crossed product... can be found in \cite{KR}
or \cite{Ta3} and will be freely used along the text.

\section{Noncommutative $L_p$ spaces}

Noncommutative $L_p$ spaces over non-semifinite von Neumann
algebras will be used quite frequently in this paper. In the
literature there exist two compatible constructions of $L_p$ in
such a general setting: Haagerup $L_p$ spaces and Kosaki's
interpolation spaces. These constructions and the associated
notion of conditional expectation will be considered in this
section.

\subsection{Haagerup $L_p$ spaces}

A full-detailed exposition of this theory is given in Haagerup
\cite{H} and Terp \cite{T1} papers. We just present the main
notions according to our purposes with an exposition similar to
\cite{JX}. A preliminary restriction is that, in view of our aims,
we can work in what follows with normal faithful (\emph{n.f.} in
short) states instead of normal semifinite faithful (\emph{n.s.f.}
in short) weights.

\vskip5pt

Let $\mathcal{M}$ be a von Neumann algebra equipped with a
distinguished \emph{n.f.} state $\varphi$. The GNS construction
applied to $\varphi$ yields a faithful representation $\rho$ of
$\mathcal{M}$ into $\mathcal{B(H)}$ so that $\rho(\mathcal{M})$ is
a von Neumann algebra acting on $\mathcal{H}$ with a separating
and generating unit vector $u$ satisfying $\varphi(x) = \langle u,
\rho(x) u \rangle$ for all $x \in \mathcal{M}$. Let us agree to
identify $\mathcal{M}$ with $\rho(\mathcal{M})$ in the following.
Then, the modular operator $\Delta$ is the (generally unbounded)
operator obtained from the polar decomposition $S = J
\Delta^{1/2}$ of the anti-linear map $S: \mathcal{M} u \rightarrow
\mathcal{M} u$ given by $S(xu) = x^* u$, see Section 9.2 in
\cite{KR} or \cite{Ta1}. We denote by $\sigma_t: \mathcal{M}
\rightarrow \mathcal{M}$ the one-parameter modular automorphism
group associated to the separating and generating unit vector $u$.
That is, for any $t \in \R$ we have an automorphism of
$\mathcal{M}$ given by \label{sigmat}
$$\sigma_t(x) = \Delta^{it} x \Delta^{-it}.$$

Then we consider the crossed product $\mathcal{R} = \mathcal{M}
\rtimes_{\sigma} \R$, which is defined as the von Neumann algebra
acting on $L_2(\R; \mathcal{H})$ and generated by the
representations $\pi: \mathcal{M} \rightarrow \mathcal{B}(L_2(\R;
\mathcal{H}))$ and $\lambda: \R \rightarrow \mathcal{B}(L_2(\R;
\mathcal{H}))$ with
$$\big( \pi(x) \xi \big) (t) = \sigma_{-t}(x) \xi(t) \qquad
\mbox{and} \qquad \big( \lambda(s) \xi \big) (t) = \xi(t-s)$$ for
$t \in \R$ and $\xi \in L_2(\R; \mathcal{H})$. Note that the
representation $\pi$ is faithful so that we can identify
$\mathcal{M}$ with $\pi(\mathcal{M})$. The dual action of $\R$ on
$\mathcal{R}$ is defined as follows. Let $\mathrm{W}: \R
\rightarrow \mathcal{B}(L_2(\R;\mathcal{H}))$ be the unitary
representation $\big( \mathrm{W}(t) \xi \big) (s) = e^{-its}
\xi(s)$. Then we define the one-parameter automorphism group
$\hat{\sigma}_t: \mathcal{R} \rightarrow \mathcal{R}$ by
$$\hat{\sigma}_t(x) = \mathrm{W}(t) x \mathrm{W}(t)^*.$$ It turns
out that $\mathcal{M}$ is the space of fixed points of the dual
action $$\mathcal{M} = \Big\{ x \in \mathcal{R} \, \big| \
\hat{\sigma}_t(x) = x \ \mbox{for all} \ t \in \R \Big\}.$$

Following \cite{PT}, the crossed product $\mathcal{R}$ is a
semifinite von Neumann algebra and admits a unique \emph{n.s.f.}
trace $\tau$ satisfying $\tau \circ \hat{\sigma}_t = e^{-t} \tau$
for all $t \in \R$. Let $L_0(\mathcal{R},\tau)$ denote the
topological $*$-algebra of $\tau$-measurable operators affiliated
with $\mathcal{R}$ and let $0 < p \le \infty$. The \emph{Haagerup
noncommutative $L_p$ space} over $\mathcal{M}$ is defined as
\label{Uffe} $$L_p(\mathcal{M},\varphi) = \Big\{ x \in
L_0(\mathcal{R},\tau) \, \big| \ \hat{\sigma}_t(x) = e^{-t/p} x \
\mbox{for all} \ t \in \R \Big\}.$$ It is clear from the
definition that $L_{\infty}(\mathcal{M},\varphi)$ coincides with
$\mathcal{M}$. Moreover, as it is to be expected,
$L_1(\mathcal{M},\varphi)$ can be canonically identified with the
predual $\mathcal{M}_*$ of the von Neumann algebra $\mathcal{M}$.
This requires a short explanation. Given a \emph{n.s.f.} weight
$\omega \in \mathcal{M}_*^+$, the dual weight $\tilde{\omega}:
\mathcal{R}_+ \rightarrow [0,\infty]$ is defined by
$$\tilde{\omega}(x) = \omega \Big( \int_{\R} \hat{\sigma}_s(x) \,
ds \Big).$$ Note that, by the translation invariance of the
Lebesgue measure, the operator valued integral above is invariant
under the dual action. In particular, it can be regarded as an
element of $\mathcal{M}$. As a \emph{n.s.f.} weight on
$\mathcal{R}$ and according to \cite{PT}, the dual weight
$\tilde{\omega}$ has a Radon-Nikodym derivative $h_{\omega}$ with
respect to $\tau$ so that
$$\tilde{\omega}(x) = \tau (h_{\omega} x)$$ for any $x \in
\mathcal{R}_+$. The operator $h_{\omega}$ so defined belongs to
$L_1(\mathcal{M}, \varphi)_+$. Indeed, $$\tau(h_{\omega}
\hat{\sigma}_t(x)) = \omega \Big( \int_{\R} \hat{\sigma}_s
(\hat{\sigma}_t(x)) \, ds \Big) = \omega \Big( \int_{\R}
\hat{\sigma}_s (x) \, ds \Big) = \tau (h_{\omega} x).$$ In
particular, $$\tau(\hat{\sigma}_t(h_{\omega}) \hat{\sigma}_t(x)) =
e^{-t} \tau (h_{\omega} x) = e^{-t} \tau(h_{\omega}
\hat{\sigma}_t(x)) \quad \mbox{for all} \qquad x \in
\mathcal{R},$$ which implies that $\hat{\sigma}_t(h_{\omega}) =
e^{-t} h_{\omega}$. Therefore, there exists a bijection between
$\mathcal{M}_*^+$ and $L_1(\mathcal{M},\varphi)_+$ which extends
to a bijection between the predual $\mathcal{M}_*$ and
$L_1(\mathcal{M},\varphi)$ by polar decomposition $$\omega = u
|\omega| \in \mathcal{M}_* \mapsto u h_{|\omega|} = h_{\omega} \in
L_1(\mathcal{M}, \varphi).$$ In fact, after imposing on
$L_1(\mathcal{M}, \varphi)$ the norm $$\|h_{\omega}\|_1 =
|\omega|(1) = \|\omega\|_{\mathcal{M}_*},$$ we obtain an isometry
between $\mathcal{M}_*$ and $L_1(\mathcal{M}, \varphi)$. There is
however a nicer way to describe this norm. As we have seen, for
any $x \in L_1(\mathcal{M}, \varphi)$ there exists a unique
$\omega_x \in \mathcal{M}_*$ such that $h_{\omega_x} = x$. This
gives rise to the functional $\mbox{tr}: L_1(\mathcal{M}, \varphi)
\rightarrow \C$ \label{TrTr} called \emph{trace} and defined by
$$\mbox{tr}(x) = \omega_x(1).$$ The functional $\mbox{tr}$ is
continuous since $|\mbox{tr}(x)| \le \mbox{tr} \big( |x| \big) =
\|x\|_1$ and satisfies the tracial property
$$\mathrm{tr}(xy) = \mathrm{tr}(yx).$$ Our distinguished state
$\varphi$ can be recovered from $\mbox{tr}$ as follows. First we
note as above that its dual weight $\tilde{\varphi}$ admits a
Radon-Nikodym derivative $\mathrm{D}_{\varphi}$ with respect to
$\tau$, so that $\tilde{\varphi}(x) = \tau (\mathrm{D}_{\varphi}
x)$ for $x \in \mathcal{R}_+$. Then, it turns out that
$$\varphi(x) = \mbox{tr}(\mathrm{D}_{\varphi}x) \quad \mbox{for}
\quad x \in \mathcal{M}.$$ According to this, we will refer in
what follows to $\mathrm{D}_{\varphi}$ as the \emph{density}
\label{densidad} of $\varphi$. Moreover, we shall write
$\mathrm{D}$ instead of $\mathrm{D}_{\varphi}$ whenever the
dependence on $\varphi$ is clear from the context. Given $0 < p <
\infty$ and $x \in L_p(\mathcal{M}, \varphi)$, we define $$\|x\|_p
= \big( \mbox{tr} |x|^p \big)^{1/p} \quad \mbox{and} \quad
\|x\|_{\infty} = \|x\|_{\mathcal{M}}.$$ $\| \ \|_p$ is a norm
(resp. $p$-norm) on $L_p(\mathcal{M}, \varphi)$ for $1 \le p \le
\infty$ (resp. $0 < p < 1$).

\begin{lemma} The Haagerup $L_p$ spaces satisfy the following
properties:
\begin{itemize}
\item[i)] If $0 < p,q,r \le \infty$ with $1/r = 1/p + 1/q$, we
have $$\|xy\|_r \le \|x\|_p \|y\|_q \qquad \mbox{for all} \qquad x
\in L_p(\mathcal{M}, \varphi), \ y \in L_q(\mathcal{M},
\varphi).$$ \item[ii)] If $1 \le p < \infty$, $L_p(\mathcal{M},
\varphi)^*$ is isometrically isomorphic to $L_{p'}(\mathcal{M},
\varphi)$ via $$x \in L_{p'}(\mathcal{M}, \varphi) \mapsto
\mathrm{tr}(x^* \, \cdot) \in L_p(\mathcal{M}, \varphi)^*.$$
\end{itemize}
\end{lemma}

An element $x \in \mathcal{M}$ is called analytic if the function
$t \in \R \mapsto \sigma_t(x) \in \mathcal{M}$ extends to an
analytic function $z \in \C \rightarrow \sigma_z(x) \in
\mathcal{M}$. By \cite{PT} we know that the subspace
$\mathcal{M}_a$ of analytic elements in $\mathcal{M}$ \label{Ma}
is a weak$^*$ dense $*$-subalgebra of $\mathcal{M}$. The proof of
the following result can be found in \cite{JX}. It will be useful
in the sequel.

\begin{lemma} \label{Lemma-Density-Analytic}
If $0 <  p < \infty$, we have
\begin{itemize}
\item[i)] $\mathcal{M}_a \mathrm{D}^{1/p}$ is dense in
$L_p(\mathcal{M}, \varphi)$. \item[ii)] $\mathrm{D}^{(1-\eta)/p}
\mathcal{M}_a \mathrm{D}^{\eta/p} = \mathcal{M}_a
\mathrm{D}^{1/p}$ for all $0 \le \eta \le 1$.
\end{itemize}
\end{lemma}

\subsection{Kosaki's interpolation}

The given definition of Haagerup $L_p$ space has the disadvantage
that the intersection of $L_p(\mathcal{M}, \varphi)$ and
$L_q(\mathcal{M}, \varphi)$ is trivial for $p \neq q$. In
particular, these spaces do not form an interpolation scale. All
these difficulties disappear with Kosaki's construction. As above,
we only consider von Neumann algebras equipped with \emph{n.f.}
states. The general construction for any von Neumann algebra can
be found in \cite{Ko} and \cite{T2}. Let us consider a von Neumann
algebra $\mathcal{M}$ equipped with a \emph{n.f.} state $\varphi$.
First we define
$$L_1(\mathcal{M}) = \mathcal{M}_*^{\mathrm{op}}.$$ Note that the
natural operator space structure for $L_1(\mathcal{M})$ requires
to consider $\mathcal{M}_*^{\mathrm{op}}$ instead of
$\mathcal{M}_*$, we refer the reader to \cite{P3} for a detailed
explanation. Then, given any real number $t$, we consider the map
$$j_t: x \in \mathcal{M} \mapsto \sigma_t(x) \varphi \in
L_1(\mathcal{M}) \qquad \mbox{with} \qquad \big( \sigma_t(x)
\varphi \big) (y) = \varphi(\sigma_t(x) y).$$ According to
\cite{Ko} there exists a unique extension $j_z: \mathcal{M}
\rightarrow L_1(\mathcal{M})$ such that, for any $0 \le \eta \le
1$, the map $j_{-i \eta}$ is injective. In particular, $(j_{-i
\eta}(\mathcal{M}), L_1(\mathcal{M}))$ is compatible for complex
interpolation and we define the \emph{Kosaki noncommutative $L_p$
spaces} as follows \label{KosLp} $$\mathcal{L}_p(\mathcal{M},
\varphi, \eta) = \big[ j_{-i \eta}(\mathcal{M}), L_1(\mathcal{M})
\big]_{\frac{1}{p}}$$ by specifying $$\|x\|_0 = \|j_{-i
\eta}^{-1}(x)\|_{\mathcal{M}} \quad \mbox{and} \quad \|x\|_1 =
\|x\|_{L_1(\mathcal{M})}.$$ The following result was proved in
\cite{H} except for the last isometry \cite{Ko}.

\begin{theorem} \label{Theorem-Haagerup-Kosaki}
Let $1 \le p \le \infty$ and let $\mathcal{M}$ be any von Neumann
algebra:
\begin{itemize}
\item[i)] If $\varphi_1$ and $\varphi_2$ are two n.f. states on
$\mathcal{M}$, we have an isometric isomorphism $$L_p(\mathcal{M},
\varphi_1) = L_p(\mathcal{M}, \varphi_2).$$ \item[ii)] If
$\varphi$ is a n.f. state and $0 \le \eta \le 1$, we have an
isometric isomorphism
$$L_p(\mathcal{M}, \varphi) = \mathcal{L}_p(\mathcal{M}, \varphi,
\eta).$$ More concretely, given $x \in \mathcal{M}$ $$\|j_{-i
\eta}(x)\|_{\mathcal{L}_p(\mathcal{M}, \varphi, \eta)} = \big\|
\mathrm{D}^{\eta/p} x \mathrm{D}^{(1-\eta)/p}
\big\|_{L_p(\mathcal{M},\varphi)}.$$
\end{itemize}
\end{theorem}

According to Theorem \ref{Theorem-Haagerup-Kosaki}, Haagerup and
Kosaki noncommutative $L_p$ spaces can be identified. We shall
write $L_p(\mathcal{M})$ \label{HaaLp} to denote in what follows
\emph{any} of the spaces defined above. In particular, after the
corresponding identifications, we may use the complex
interpolation method for Haagerup $L_p$ spaces. We should also
note that Kosaki's definition of $L_p$ presents some other
disadvantages. The main lacks in this paper will be the absence of
positive cones and the fact that the case $0 < p < 1$ is excluded
from the definition. In particular, we will use Haagerup $L_p$
spaces and we will apply Theorem \ref{Theorem-Haagerup-Kosaki}
when needed.

\subsection{Conditional expectations}

Let us consider a von Neumann algebra $\mathcal{M}$ equipped with
a \emph{n.f.} state $\varphi$ and a von Neumann subalgebra
$\mathcal{N}$ of $\mathcal{M}$. A \emph{conditional expectation}
$\mathsf{E}: \mathcal{M} \to \mathcal{N}$ is a positive
contractive projection. $\mathsf{E}$ is called faithful if
$\mathsf{E}(x^*x) \neq 0$ for any $x \in \mathcal{M}$ and normal
when it has a predual $\mathsf{E}_*: \mathcal{M}_* \to
\mathcal{N}_*$. According to Takesaki \cite{Ta2}, if $\mathcal{N}$
is invariant under the action of the modular automorphism group
(i.e. $\sigma_t(\mathcal{N}) \subset \mathcal{N}$ for all $t \in
\R$), there exists a unique faithful normal conditional
expectation $\mathsf{E}: \mathcal{M} \rightarrow \mathcal{N}$
satisfying $\varphi \circ \mathsf{E} = \varphi$. Moreover, by
Connes \cite{Co} it commutes with the modular automorphism group
$$\mathsf{E} \circ \sigma_t = \sigma_t \circ \mathsf{E}.$$ The required
invariance of $\mathcal{N}$ under the action of $\sigma_t$ implies
that the modular automorphism group associated to $\mathcal{N}$
coincides with the restriction of $\sigma$ to $\mathcal{N}$. It
follows that $\mathcal{N} \rtimes_{\sigma} \R$ is a von Neumann
subalgebra of $\mathcal{M} \rtimes_{\sigma} \R$. In particular,
the space $L_p(\mathcal{N})$ can be identified isometrically with
a subspace of $L_p(\mathcal{M})$, see \cite{JX} for details. In
this paper we shall permanently assume the existence of a normal
faithful conditional expectation $\mathsf{E}: \mathcal{M} \to
\mathcal{N}$.

\vskip5pt

It is well-known that in the tracial case, the conditional
expectation $\mathsf{E}$ extends to a contractive projection from
$L_p(\mathcal{M})$ onto $L_p(\mathcal{N})$ for any $1 \le p \le
\infty$, which is still positive and has the modular property
$\mathsf{E}(axb) = a \mathsf{E}(x) b$ for all $a,b \in
\mathcal{N}$ and $x \in L_p(\mathcal{M})$. These properties remain
valid in this context. We summarize in the following lemma the
main properties of $\mathsf{E}$ that will be used in the sequel
and refer the reader to \cite{JX} for a proof of these facts.

\begin{lemma} \label{Lemma-Properties-Expectation}
Let $\mathcal{M}$ and $\mathcal{N}$ be as above:
\begin{itemize}
\item[i)] If $2 \le p \le \infty$ and $x \in L_p(\mathcal{M})$, we
have $\mathsf{E}(x)^* \, \mathsf{E}(x) \le \mathsf{E}(x^*x)$.
\item[ii)] If $1 \le p \le \infty$, $\mathsf{E}: L_p(\mathcal{M})
\rightarrow L_p(\mathcal{N})$ is a positive contractive
projection. \item[iii)] If $1 \le p \le \infty$ and $x \in
\mathcal{M}$,
$$\mathsf{E}(\mathrm{D}^{1/p} x) = \mathrm{D}^{1/p} \mathsf{E}(x)
\quad \mbox{and} \quad \mathsf{E}(x \mathrm{D}^{1/p}) =
\mathsf{E}(x) \mathrm{D}^{1/p}.$$ \item[iv)] If $a \in
L_p(\mathcal{N}), x \in L_q(\mathcal{M}), b \in L_r(\mathcal{N})$
and $\frac{1}{p} + \frac{1}{q} + \frac{1}{r} \le 1$,
$\mathsf{E}(axb) = a \mathsf{E}(x) b$.
\end{itemize}
\end{lemma}

\section{Pisier's vector-valued $L_p$ spaces}
\label{Subsetion1.2}

Noncommutative vector-valued $L_p$ spaces
$L_p(\mathcal{M};\mathrm{X})$ appeared quite recently with
Pisier's work \cite{P2}. The reason for the novelty of such a
natural concept lies in the fact that $\mathrm{X}$ must be
equipped with an operator space structure rather than a Banach
space one. Moreover, many natural properties such as duality,
complex interpolation, etc... must be formulated in the category
of operator spaces. Let us begin by recalling the concept of
operator space.

\subsection{Operator spaces}

Operator space theory plays a central role in this paper. It can
be regarded as a noncommutative generalization of Banach space
theory and it has proved to be an essential tool in operator
algebra as well as in noncommutative harmonic analysis. An
\emph{operator space} $\mathrm{X}$ is defined as a closed subspace
of the space $\mathcal{B}(\mathcal{H})$ of bounded operators on
some Hilbert space $\mathcal{H}$. Let $\mathrm{M}_n(\mathrm{X})$
denote the vector space of $n \times n$ matrices with entries in
$\mathrm{X}$. According to Ruan's theorem \cite{R}, an operator
space $\mathrm{X}$ comes either with a concrete isometric
embedding $j: \mathrm{X} \to \mathcal{B}(\mathcal{H})$ of with a
sequence of norms on $\mathrm{M}_n(\mathrm{X})$ for $n \ge 1$
satisfying $$\big\| \big( x_{ij} \big)
\big\|_{\mathrm{M}_n(\mathrm{X})} = \big\| \big( j(x_{ij}) \big)
\big\|_{\mathcal{B}(\mathcal{H}^n)}.$$ Ruan's axioms
\cite{ER,P3,R} describe axiomatically those sequences of matrix
norms which can occur from an isometric embedding into
$\mathcal{B}(\mathcal{H})$. Any such sequence of norms provides
$\mathrm{X}$ with a so-called operator space structure. Every
Banach space can be equipped with several operator space
structures. In particular, the most important information carried
by an operator space is not the space itself but the way in which
it embeds isometrically into $\mathcal{B}(\mathcal{H})$. For this
reason, the main difference between Banach and operator space
theory lies on the morphisms rather than on the spaces. The
morphisms in the category of operator spaces are \emph{completely
bounded} linear maps $u: \mathrm{X} \to \mathrm{Y}$. That is,
linear maps satisfying that the quantity $$\|u\|_{cb(\mathrm{X},
\mathrm{Y})} = \sup_{n \ge 1} \big\| id_{\mathrm{M}_n} \ten u
\big\|_{\mathcal{B}(\mathrm{M}_n(\mathrm{X}),
\mathrm{M}_n(\mathrm{Y}))}$$ is finite. In this paper we shall
assume some familiarity with basic notions such as duality,
Haagerup tensor products, the OH spaces... that can be found in
the recent books \cite{ER,P3}.

\vskip5pt

\subsection{The hyperfinite case}

Before any other consideration, let us recall the natural operator
space structure (\emph{o.s.s.} in short) on $L_p(\mathcal{M})$.
Proceeding as in Chapter 3 of \cite{P2}, we regard $\mathcal{M}$
as a subspace of $\mathcal{B(H)}$ with $\mathcal{H}$ being the
Hilbert space arising from the GNS construction. Similarly, by
embedding the predual von Neumann algebra $\mathcal{M}_*$ on its
bidual $\mathcal{M}^*$, we obtain an \emph{o.s.s.} for
$\mathcal{M}_*$. The \emph{o.s.s.} on $L_1(\mathcal{M})$ is then
given by that of $\mathcal{M}_*^\mathrm{op}$, see page 139 in
\cite{P3} for a detailed justification of this definition. Then,
the complex interpolation method for operator spaces developed in
\cite{P1} provides a natural \emph{o.s.s.} on $L_p(\mathcal{M})$
$$L_p(\mathcal{M}) = \big[ L_\infty(\mathcal{M}), L_1(\mathcal{M})
\big]_{1/p} = \big[ \mathcal{M}, \mathcal{M}_*^\mathrm{op}
\big]_{1/p}.$$

Let us now describe the natural operator space structure for
vector-valued $L_p$ spaces. Let $\mathcal{M}$ be an hyperfinite
von Neumann algebra and let $\mathrm{X}$ be any operator space.
Then, recalling the definitions of the projective and the minimal
tensor product in the category of operator spaces, we define
$$L_1(\mathcal{M}; \mathrm{X}) = L_1(\mathcal{M})
\widehat{\otimes} \, \mathrm{X} \quad \mbox{and} \quad
L_{\infty}(\mathcal{M}; \mathrm{X}) = L_{\infty}(\mathcal{M})
\otimes_{\mathrm{min}} \mathrm{X}.$$ Then, the \emph{space}
\label{VVLp} $L_p(\mathcal{M}; \mathrm{X})$ is defined by complex
interpolation
$$L_p(\mathcal{M};\mathrm{X}) = \big[ L_\infty(\mathcal{M};
\mathrm{X}), L_1(\mathcal{M}; \mathrm{X}) \big]_{1/p}.$$ As it was
explained in \cite{P2}, the hyperfiniteness of $\mathcal{M}$ leads
to obtain some expected properties of these spaces. We shall
discuss the non-hyperfinite case in the next paragraph. We are not
reviewing here the basic results from Pisier's theory, for which
we refer the reader to Chapters 1,2 and 3 in \cite{P2}.

\vskip5pt

Let us fix some notation. $R$ and $C$ denote the row and column
operator spaces (\emph{c.f.} Chapter 1 of \cite{P3}) constructed
over $\ell_2$. Identifying $\mathcal{B}(\ell_2)$ (via the
canonical basis of $\ell_2$) with a space of matrices with
infinitely many rows and columns, $R$ and $C$ are the first row
and first column subspaces of $\mathcal{B}(\ell_2)$. Similarly, if
$1 \le p \le \infty$ and $S_p$ denotes the Schatten $p$-class over
$\ell_2$, the spaces $R_p$ and $C_p$ denote the row and column
subspaces of $S_p$. The finite-dimensional versions over
$\ell_2^n$ will be denoted by $R_p^n$ and $C_p^n$ respectively.
Note that, as in the infinite-dimensional case, we have $R_n =
R_\infty^n$ and $C_n = C_\infty^n$. Given an operator space
$\mathrm{X}$, the vector-valued Schatten classes with values in
$\mathrm{X}$ will be denoted by $S_p(\mathrm{X})$ and
$S_p^n(\mathrm{X})$ respectively. Among the several
characterizations of these spaces given in \cite{P2}, we have
$$S_p(\mathrm{X}) = C_p \otimes_h \mathrm{X} \ten_h R_p \quad
\mbox{and} \quad S_p^n(\mathrm{X}) = C_p^n \ten_h \mathrm{X}
\ten_h R_p^n.$$

\subsection{The non-hyperfinite case}

One of the main restrictions in Pisier's theory \cite{P2} comes
from the fact that the construction of $L_p(\mathcal{M};
\mathrm{X})$ requires $\mathcal{M}$ to be hyperfinite. This
excludes for instance free products of von Neumann algebras, a
very relevant tool in this paper. There exists however a very
recent construction in \cite{J4} of $L_p(\mathcal{M}; \mathrm{X})$
which is valid for any QWEP von Neumann algebra. Nevertheless,
since we only deal with very specific operator spaces, we briefly
discuss them.

\vskip5pt

\textsc{The spaces $L_p(\mathcal{M};R_p^n)$ and
$L_p(\mathcal{M};C_p^n)$}. As we explained in the Introduction,
these spaces are of capital importance. Let us recall that the
spaces $R_p^n(\mathrm{X})$ and $C_p^n(\mathrm{X})$ are defined as
subspaces of $S_p^n(\mathrm{X})$ for any operator space
$\mathrm{X}$. Therefore, motivated by Fubini's theorem for
noncommutative $L_p$ spaces \cite{P2}, we obtain the following
definition valid for arbitrary von Neumann algebras
$$L_p(\mathcal{M};R_p^n) = R_p^n(L_p(\mathcal{M})) \quad
\mbox{and} \quad L_p(\mathcal{M};C_p^n) =
C_p^n(L_p(\mathcal{M})).$$ These spaces appear quite frequently in
the theory of noncommutative martingale inequalities. However, in
that context they are respectively denoted by $L_p(\mathcal{M};
\ell_2^r)$ and $L_p(\mathcal{M}; \ell_2^c)$, see e.g. \cite{PX1}
or \cite{Ran1} for more details.

\vskip5pt

\textsc{The spaces $L_p(\mathcal{M}; \ell_q^n)$}. The spaces
$L_p(\mathcal{M}; \ell_1^n)$ and $L_p(\mathcal{M}; \ell_\infty^n)$
where defined in \cite{J1} for general von Neumann algebras to
study Doob's maximal inequality for noncommutative martingales.
Indeed, since the notion of maximal function does not make any
sense in the noncommutative setting, the norm in $L_p(\Omega)$ of
Doob's maximal function
$$f^*(w) = \sup_{n \ge 1} |f_n(w)|$$ was reinterpreted as the norm
in $L_p(\Omega;\ell_\infty)$ of the sequence $(f_1,f_2, \ldots)$.
This was the original motivation to study the spaces
$L_p(\mathcal{M};\ell_1)$ and $L_p(\mathcal{M};\ell_\infty)$. A
detailed exposition of these spaces can also be found in
\cite{JX2}. The spaces $L_p(\mathcal{M}; \ell_q^n)$ will be of
capital importance in the last chapter of this paper, where we
shall obtain our main results. These spaces where recently defined
in \cite{JX4} over non-hyperfinite von Neumann algebras as follows
$$L_p(\mathcal{M}; \ell_q^n) = \big[ L_p(\mathcal{M};
\ell_\infty^n), L_p(\mathcal{M}; \ell_1^n) \big]_{\frac1q}.$$ This
formula trivially generalizes by the reiteration theorem
\cite{BL}. A relevant property of these spaces proved in
\cite{JX4} is Fubini's isometry $L_p(\mathcal{M}; \ell_p^n) =
\ell_p^n(L_p(\mathcal{M}))$. It is also proved in \cite{JX4} that
Pisier's identities hold. In other words, defining the auxiliary
index $1/r = |1/p - 1/q|$ we find
\begin{itemize}
\item[i)] If $p \le q$, we have $$\|x\|_{pq} = \inf \Big\{
\|\alpha\|_{L_{2r}(\mathcal{M})} \|y\|_{L_q(\mathcal{M};
\ell_q^n)} \|\beta\|_{L_{2r}(\mathcal{M})} \, \big| \ x = \alpha y
\beta \Big\}.$$

\item[ii)] If $p \ge q$, we have $$\|x\|_{pq} = \sup \Big\{
\|\alpha y \beta \|_{L_q(\mathcal{M}; \ell_q^n)} \, \big| \
\alpha, \beta \in \mathsf{B}_{L_{2r}(\mathcal{M})} \Big\}.$$
\end{itemize}

\begin{remark}
\emph{Let $\mathcal{M}$ be an hyperfinite von Neumann algebra and
\begin{eqnarray*}
\mathsf{A}_p(\mathcal{M}; n) & = & \big[
L_p(\mathcal{M};\ell_\infty^n), L_p(\mathcal{M};\ell_1^n)
\big]_{\frac12}, \\ \mathsf{B}_p(\mathcal{M};n) & = & \big[
L_p(\mathcal{M};C_p^n), L_p(\mathcal{M}; R_p^n) \big]_{\frac12}.
\end{eqnarray*}
According to \cite{P2}, we have
$$\mathsf{A}_p(\mathcal{M};n) = L_p(\mathcal{M};\mathrm{OH}_n) =
\mathsf{B}_p(\mathcal{M};n).$$ It is therefore natural to ask
whether these identities remain valid with our definition of
$L_p(\mathcal{M}; \mathrm{OH}_n)$ for non-hyperfinite
$\mathcal{M}$. This is indeed the case since, following the same
terminology already defined in the Introduction, we have
\begin{itemize}
\item[i)] If $1 \le p \le 2 \, $ and $1/p = 1/2 + 1/q$, we have
$$\mathsf{A}_p(\mathcal{M};n) = L_{2q}(\mathcal{M}) \ell_2^n \big(
L_2(\mathcal{M}) \big) L_{2q}(\mathcal{M}) =
\mathsf{B}_p(\mathcal{M};n).$$ \item[ii)] If $2 \le p \le \infty$
and $1/2 = 1/p + 1/q$, we have
$$\mathsf{A}_p(\mathcal{M};n) = \sqrt{n} \, L_{(2q,2q)}^p
(\ell_\infty^n(\mathcal{M}), \mathsf{E}_n) =
\mathsf{B}_p(\mathcal{M};n),$$ where $\mathsf{E}_n$ is given by
$\mathsf{E}_n: \sum_{k=1}^n x_k \otimes \delta_k \in
\ell_\infty^n(\mathcal{M}) \mapsto \frac{1}{n} \sum_{k=1}^n x_k
\in \mathcal{M}.$
\end{itemize}
The isometric identities for the $\mathsf{A}_p$'s follow from
Pisier's identities above while the identities for the
$\mathsf{B}_p$'s follow from \cite{X} or Theorem \ref{TheoremA1}
below. Thus, we agree to define the \emph{o.s.s.} on
$L_p(\mathcal{M};\mathrm{OH}_n)$ as any of the interpolation
spaces above.}
\end{remark}

\begin{remark}
\emph{As we shall see through this paper, all the spaces mentioned
so far arise as particular cases of our definition below of
amalgamated and conditional $L_p$ spaces.}
\end{remark}

\section {The spaces $L_p^r(\EME, \mathsf{E})$ and $L_p^c(\mathcal{M},
\mathsf{E})$}

The spaces $L_p^r(\mathcal{M}, \mathsf{E})$ and
$L_p^c(\mathcal{M}, \mathsf{E})$ were introduced in \cite{J1},
where they turned out to be quite useful in the context of
noncommutative probability. Both spaces will play a central role
in this paper since they are the most significant examples of the
so-called amalgamated and conditional $L_p$ spaces, to be defined
below. Let $\mathcal{M}$ be a von Neumann algebra equipped with a
\emph{n.f.} state $\varphi$ and let $\mathcal{N}$ be a von Neumann
subalgebra of $\mathcal{M}$. Let $\mathsf{E}: \mathcal{M}
\rightarrow \mathcal{N}$ denote the conditional expectation of
$\mathcal{M}$ onto $\mathcal{N}$. Then, given $1 \le p \le \infty$
and $(\alpha,a) \in \mathcal{N} \times \mathcal{M}$, we define
\begin{equation} \label{Equation-Conditional-Dual}
\begin{array}{rcl} \big\| \alpha \mathrm{D}^{\frac{1}{p}} a
\big\|_{L_p^r(\mathcal{M}, \mathsf{E})} & = & \big\| \alpha
\mathrm{D}^{\frac{1}{p}} \mathsf{E}(aa^*) \mathrm{D}^{\frac{1}{p}}
\alpha^* \big\|_{L_{p/2}(\mathcal{N})}^{1/2},
\\ \\ \big\| a \mathrm{D}^{\frac{1}{p}} \alpha
\big\|_{L_p^c(\mathcal{M}, \mathsf{E})} & = & \big\| \alpha^*
\mathrm{D}^{\frac{1}{p}} \mathsf{E}(a^*a) \mathrm{D}^{\frac{1}{p}}
\alpha \big\|_{L_{p/2}(\mathcal{N})}^{1/2}.
\end{array}
\end{equation}
$L_p^r(\mathcal{M}, \mathsf{E})$ and $L_p^c(\mathcal{M},
\mathsf{E})$ will stand for the completions with respect to these
norms. \label{RCLpE}

\begin{remark} \label{Remark-p<2<p}
\emph{Note that for $1 \le p < 2$ we have $0 < p/2 < 1$ so that we
are forced to use Haagerup $L_p$ spaces in the definition
(\ref{Equation-Conditional-Dual}). On the other hand, let us note
that the assumption $x \in \mathcal{N} \mathrm{D}^{1/p}
\mathcal{M}$ (resp. $x \in \mathcal{M} \mathrm{D}^{1/p}
\mathcal{N}$) is not needed to define the norm of $x$ in
$L_p^r(\mathcal{M}, \mathsf{E})$ (resp. $L_p^c(\mathcal{M},
\mathsf{E})$) for $2 \le p \le \infty$. Indeed, given $x \in
L_p(\mathcal{M})$ we can define $$\|x\|_{L_p^r(\mathcal{M},
\mathsf{E})} = \big\| \mathsf{E}(xx^*)
\big\|_{L_{p/2}(\mathcal{N})}^{1/2} \quad \mbox{and} \quad
\|x\|_{L_p^c(\mathcal{M}, \mathsf{E})} = \big\| \mathsf{E}(x^*x)
\big\|_{L_{p/2}(\mathcal{N})}^{1/2}.$$ In that case
$\mathsf{E}(xx^*)$ and $\mathsf{E}(x^*x)$ are well-defined and
$$\max \Big\{ \big\| \mathsf{E}(xx^*)^{1/2}
\big\|_{L_p(\mathcal{N})}, \big\| \mathsf{E}(x^*x)^{1/2}
\big\|_{L_p(\mathcal{N})} \Big\} \le \|x\|_{L_p(\mathcal{M})}.$$
Hence, by the density of $\mathcal{N} \mathrm{D}^{1/p}
\mathcal{M}$ and $\mathcal{M} \mathrm{D}^{1/p} \mathcal{N}$ in
$L_p(\mathcal{M})$, we obtain the same closure as in
(\ref{Equation-Conditional-Dual}). However, in the case $1 \le p <
2$ the conditional expectation $\mathsf{E}$ is no longer
continuous on $L_{p/2}(\mathcal{M})$ so that we need this
alternative definition.}
\end{remark}

The duality of $L_p^r(\mathcal{M}, \mathsf{E})$ and
$L_p^c(\mathcal{M}, \mathsf{E})$ was studied in \cite{J1}. For the
moment we just need to know that, given $1 < p < \infty$, the
following isometries hold via the anti-linear duality bracket
$\langle x,y \rangle = \mbox{tr}(x^*y)$
\begin{equation} \label{Equation-Dual-Row-Column}
\begin{array}{rcl}
L_p^r(\mathcal{M}, \mathsf{E})^* & = & L_{p'}^r(\mathcal{M},
\mathsf{E}), \\ L_p^c(\mathcal{M}, \mathsf{E})^* & = &
L_{p'}^c(\mathcal{M}, \mathsf{E}).
\end{array}
\end{equation}

\begin{lemma} \label{Lemma-Norm-Conditional-Sup}
If $2 \le p \le \infty$ and $1/2 = 1/p + 1/s$, we have
\begin{eqnarray*}
\|x\|_{L_p^r(\mathcal{M}, \mathsf{E})} & = & \sup \Big\{ \|\alpha
x\|_{L_2(\mathcal{M})} \, \big| \ \|\alpha\|_{L_s(\mathcal{N})}
\le 1 \Big\}, \\ \|x\|_{L_p^c(\mathcal{M}, \mathsf{E})} & = & \sup
\Big\{ \|x \beta\|_{L_2(\mathcal{M})} \, \big| \
\|\beta\|_{L_s(\mathcal{N})} \le 1 \Big\}.
\end{eqnarray*}
\end{lemma}

\begin{proof} Given $x \in L_p^r(\mathcal{M}, \mathsf{E})$, the
operator $\mathsf{E}(xx^*)$ is positive. Hence
\begin{eqnarray*}
\|x\|_{L_p^r(\mathcal{M}, \mathsf{E})}^2 & = & \sup \Big\{
\mbox{tr} \big( a \mathsf{E}(xx^*) \big) \, \big| \ a \ge 0, \
\|a\|_{L_{(p/2)'}(\mathcal{N})} \le 1 \Big\} \\ & = & \sup \Big\{
\mbox{tr} \big( \alpha \mathsf{E}(xx^*) \alpha^* \big) \, \big| \
\|\alpha^* \alpha\|_{L_{(p/2)'}(\mathcal{N})} \le 1 \Big\}
\\ & = & \sup \Big\{ \mbox{tr} \big( \alpha xx^* \alpha^*
\big) \, \big| \ \|\alpha^* \alpha\|_{L_{(p/2)'}(\mathcal{N})} \le
1 \Big\} \\ & = & \sup \Big\{ \|\alpha x\|_{L_2(\mathcal{M})}^2 \,
\big| \ \|\alpha\|_{L_{2(p/2)'}(\mathcal{N})} \le 1 \Big\}.
\end{eqnarray*}
Finally we recall that $s = 2 (p/2)'$. The proof for
$L_p^c(\mathcal{M}, \mathsf{E})$ is entirely similar.
\end{proof}

As usual, we shall write $L_2^r(\mathcal{M}) =
\mathcal{B}(L_2(\mathcal{M}), \C)$ and $L_2^c(\mathcal{M}) =
\mathcal{B}(\C, L_2(\mathcal{M}))$ to denote the row and column
Hilbert spaces over $L_2(\mathcal{M})$. Both spaces embed
isometrically in $\mathcal{B}(L_2(\mathcal{M}))$. Hence, they
admit a natural \emph{o.s.s.} Given any index $2 \le p \le
\infty$, we generalize these spaces as follows. According to the
definition of Haagerup $L_p$ spaces, we may consider the
contractive inclusions $$\begin{array}{rcl} j_r: x \in \mathcal{M}
& \mapsto & \mathrm{D}^{\frac12} x \in L_2(\mathcal{M}), \\ j_c: x
\in \mathcal{M} & \mapsto & x \mathrm{D}^{\frac12} \in
L_2(\mathcal{M}). \end{array}$$ Then we define for  $2 \le p \le
\infty$
\begin{equation} \label{Et-LpAsimet}
\begin{array}{rcl}
L_p^r(\mathcal{M}) & = & \big[ j_r(\mathcal{M}),
L_2^r(\mathcal{M}) \big]_{\frac2p} \quad \mbox{with} \quad \|x\|_0
= \big\| \mathrm{D}^{-1/2} x \big\|_{\mathcal{M}}, \\
[5pt] L_p^c(\mathcal{M}) & = & \big[ j_c(\mathcal{M}),
L_2^c(\mathcal{M}) \big]_{\frac2p} \quad \mbox{with} \quad \|x\|_0
= \big\| x \mathrm{D}^{-1/2} \big\|_{\mathcal{M}}.
\end{array}
\end{equation}
Note that $L_p^r(\mathcal{M}) = L_p(\mathcal{M}) =
L_p^c(\mathcal{M})$ as Banach spaces.

\vskip5pt

Let $\mathcal{N}$ be a von Neumann subalgebra of $\mathcal{M}$.
Then, given $2 \le u,q,v \le \infty$, we consider the closed ideal
$\mathcal{I}_1$ in the Haagerup tensor product $L_u^r(\mathcal{N})
\otimes_h L_q^c(\mathcal{M})$ generated by the differences $x
\gamma \otimes y - x \otimes \gamma y$, with $x \in
L_u^r(\mathcal{N})$, $y \in L_q^c(\mathcal{M})$ and $\gamma \in
\mathcal{N}$. Similarly, we consider the closed ideal
$\mathcal{I}_2$ in $L_q^r(\mathcal{M}) \otimes_h
L_v^c(\mathcal{N})$ generated by the differences $x \gamma \otimes
y - x \otimes \gamma y$, with $x \in L_q^r(\mathcal{M})$, $y \in
L_v^c(\mathcal{N})$ and $\gamma \in \mathcal{N}$. Then, we define
the amalgamated Haagerup tensor products as the quotients
\begin{eqnarray*}
L_u^r(\mathcal{N}) \otimes_{\mathcal{N},h} L_q^c(\mathcal{M}) & =
& L_u^r(\mathcal{N}) \otimes_h L_q^c(\mathcal{M}) \big/
\mathcal{I}_1, \\ L_q^r(\mathcal{M}) \otimes_{\mathcal{N},h}
L_v^c(\mathcal{N}) & = & L_q^r(\mathcal{M}) \otimes_h
L_v^c(\mathcal{N}) \big/ \mathcal{I}_2.
\end{eqnarray*}
Let $a_1, a_2, \ldots, a_m \in L_u(\mathcal{N})$ and $b_1, b_2,
\ldots, b_m \in L_q(\mathcal{M})$. Using that the Haagerup tensor
product commutes with complex interpolation, it is not difficult
to check that the following identities hold
\begin{eqnarray*}
\Big\| \sum_{k=1}^m a_k \otimes e_{1k} \Big\|_{L_u^r(\mathcal{N})
\otimes_h R_m} & = & \Big\| \big( \sum_{k=1}^m a_k a_k^*
\big)^{1/2} \Big\|_{L_u(\mathcal{N})}, \\ \Big\| \sum_{k=1}^m
e_{k1} \otimes b_k \Big\|_{C_m \otimes_h L_q^c(\mathcal{M})} & = &
\Big\| \big( \sum_{k=1}^m b_k^* b_k \big)^{1/2}
\Big\|_{L_q(\mathcal{M})}.
\end{eqnarray*}
Let us write $\mathcal{S}_r(a) = (\summ_k a_k a_k^*)^{1/2}$ and
$\mathcal{S}_c(b) = (\summ_k b_k^* b_k )^{1/2}$. Then, following
the definition of the Haagerup tensor product, the norm of $x$ in
$L_u^r(\mathcal{N}) \otimes_{\mathcal{N}, h} L_q^c(\mathcal{M})$
can be written as $$\|x\|_{u \cdot q} = \inf \Big\{ \big\|
\mathcal{S}_r(a) \big\|_{L_u(\mathcal{N})} \big\| \mathcal{S}_c(b)
\big\|_{L_q(\mathcal{M})} \, \big| \ x \simeq \sum_{k=1}^m a_k
\otimes b_k \Big\},$$ where $\simeq$ means that the difference
belongs to $\mathcal{I}_1$. Similarly, given $e_1, e_2, \ldots,
e_m$ in $L_q^r(\mathcal{M})$ and $f_1, f_2, \ldots, f_m$ in
$L_v^c(\mathcal{N})$, we can write the norm of an element $x$ in
the space $L_q^r(\mathcal{M}) \otimes_{\mathcal{N}, h}
L_v^c(\mathcal{N})$ as follows $$\|x\|_{q \cdot v} = \inf \Big\{
\big\| \mathcal{S}_r(e) \big\|_{L_q(\mathcal{M})} \big\|
\mathcal{S}_c(f) \big\|_{L_v(\mathcal{N})} \, \big| \ x \simeq
\sum_{k=1}^m e_k \otimes f_k \Big\}.$$

\begin{lemma} \label{Lemma-Norm-Amalgamated-OS}
The following identities hold $$\|x\|_{u \cdot q} = \inf_{x =
\alpha y} \|\alpha\|_{L_u(\mathcal{N})}
\|y\|_{L_q(\mathcal{M})},$$ $$\|x\|_{q \cdot v} = \inf_{x = y
\beta} \|y\|_{L_q(\mathcal{M})} \|\beta\|_{L_v(\mathcal{N})}.$$
\end{lemma}

\begin{proof} The upper estimates $$\|x\|_{u \cdot q} \le \inf_{x = \alpha
y} \|\alpha\|_{L_u(\mathcal{N})} \|y\|_{L_q(\mathcal{M})},$$
$$\|x\|_{q \cdot v} \le \inf_{x = y \beta}
\|y\|_{L_q(\mathcal{M})} \|\beta\|_{L_v(\mathcal{N})},$$ are
trivial. We only prove the converse for the first identity since
the second identity can be derived in the same way. Let us assume
that $\|x\|_{u \cdot q} < 1$. Then we can find $a_1, a_2, \ldots,
a_m \in L_u(\mathcal{N})$ and $b_1, b_2, \ldots, b_m \in
L_q(\mathcal{M})$ such that $$x \simeq \sum_{k=1}^m a_k \otimes
b_k \quad \mbox{and} \quad \big\| \mathcal{S}_r(a)
\big\|_{L_u(\mathcal{N})} \big\| \mathcal{S}_c(b)
\big\|_{L_q(\mathcal{M})} < 1.$$ If $\mathrm{D}$ denotes the
density associated to $\varphi$, we define
$$\widetilde{\mathcal{S}}_r(a) = \Big( \sum_{k=1}^m a_ka_k^* +
\delta \mathrm{D}^{2/u} \Big)^{1/2} \quad \mbox{with} \quad \delta
> 0.$$ Since $\mbox{supp} \, \mathrm{D} = 1$,
$\widetilde{\mathcal{S}}_r(a)$ is invertible and we can define
$\alpha_k$ by $a_k = \widetilde{\mathcal{S}}_r(a) \alpha_k$ so
that $$\big\| \mathcal{S}_r(\alpha) \big\|_{L_\infty(\mathcal{N})}
= \Big\| \big( \sum_{i=1}^m \alpha_i \alpha_i^* \big)^{1/2}
\Big\|_{L_{\infty}(\mathcal{N})} \le 1.$$ The amalgamation over
$\mathcal{N}$ allows us to write
$$\sum_{k=1}^m a_k \otimes b_k \simeq \widetilde{\mathcal{S}}_r(a)
\otimes \big( \sum_{k=1}^m \alpha_k b_k \big)$$ for any possible
decomposition of $x$. Then, we have $$\Big\| \sum_{k=1}^m \alpha_k
b_k \Big\|_{L_q(\mathcal{M})} \le \big\| \mathcal{S}_r(\alpha)
\big\|_{L_{\infty}(\mathcal{N})} \big\| \mathcal{S}_c(b)
\big\|_{L_q(\mathcal{M})} \le \big\| \mathcal{S}_c(b)
\big\|_{L_q(\mathcal{M})}.$$ Moreover, applying the triangle
inequality we obtain $$\big\| \widetilde{\mathcal{S}}_r(a)
\big\|_{L_u(\mathcal{N})}^2 = \Big\| \sum_{k=1}^m a_ka_k^* +
\delta \mathrm{D}^{2/u} \Big\|_{L_{u/2}(\mathcal{N})} \le \big\|
\mathcal{S}_r(a) \big\|_{L_u(\mathcal{N})}^2 + \delta$$ In
summary, we have $$\inf_{x = \alpha y}
\|\alpha\|_{L_u(\mathcal{N})} \|y\|_{L_q(\mathcal{M})} \le
\sqrt{\big\| \mathcal{S}_r(a) \big\|_u^2 + \delta} \ \big\|
\mathcal{S}_c(b) \big\|_{L_q(\mathcal{M})}.$$ This holds for any
decomposition $x \simeq \summ_k a_k \otimes b_k$. We conclude
letting $\delta \to 0$. \end{proof}

\begin{proposition} \label{Proposition-Norm-Conditional-Inf}
Let $2 \le p < \infty$ and $2 < s \le \infty$ related by $1/2 =
1/p + 1/s$. Then, we have the following isometries $$L_{p'}^r
(\mathcal{M}, \mathsf{E}) = L_p^r(\mathcal{M}, \mathsf{E})^* =
L_s^r(\mathcal{N}) \otimes_{\mathcal{N}, h} L_2^c(\mathcal{M}),$$
$$L_{p'}^c (\mathcal{M}, \mathsf{E}) = L_p^c(\mathcal{M},
\mathsf{E})^* = L_2^r(\mathcal{M}) \otimes_{\mathcal{N},h}
L_s^c(\mathcal{N}).$$
\end{proposition}

\begin{proof} As we have already said, the isometries $$L_{p'}^r
(\mathcal{M}, \mathsf{E}) = L_p^r(\mathcal{M}, \mathsf{E})^*,$$
$$L_{p'}^c (\mathcal{M}, \mathsf{E}) = L_p^c(\mathcal{M},
\mathsf{E})^*,$$ were proved in \cite{J1}. Now let us consider an
element $x \in L_s^r(\mathcal{N}) \otimes_{\mathcal{N}, h}
L_2^c(\mathcal{M})$ and let $x = \alpha y$ be a decomposition with
$\alpha \in L_s(\mathcal{N})$ and $y \in L_2(\mathcal{M})$. Then,
(\ref{Equation-Dual-Row-Column}) and Lemma
\ref{Lemma-Norm-Conditional-Sup} provide the following inequality
$$\|x\|_{L_{p'}^r(\mathcal{M}, \mathsf{E})} = \sup \Big\{
\mbox{tr} \big( \alpha^* z y^* \big) \, \big| \
\sup_{\|\gamma\|_{L_s(\mathcal{N})} \le 1} \|\gamma
z\|_{L_2(\mathcal{M})} \le 1 \Big\} \le
\|\alpha\|_{L_s(\mathcal{N})} \|y\|_{L_2(\mathcal{M})}.$$
According to Lemma \ref{Lemma-Norm-Amalgamated-OS}, this proves
that $$id: x \in L_s^r(\mathcal{N}) \otimes_{\mathcal{N}, h}
L_2^c(\mathcal{M}) \mapsto x \in L_{p'}^r(\mathcal{M},
\mathsf{E})$$ is a contraction. Reciprocally, let us consider an
element of the form $x = \xi \mathrm{D}^{1/p'} a$ with $\xi \in
\mathcal{N}$ and $a \in \mathcal{M}$. Then, given any $\delta
> 0$, we consider the decomposition $x = \alpha_{\delta}
y_{\delta}$ with $\alpha_{\delta}$ and $y_{\delta}$ given by
$$\alpha_{\delta} = \Big( \xi \mathrm{D}^{\frac{1}{p'}}
\mathsf{E}(aa^* + \delta 1) \mathrm{D}^{\frac{1}{p'}} \xi^*
\Big)^{\gamma/2} \qquad \mbox{and} \qquad y_{\delta} =
\alpha_{\delta}^{-1} x$$ with $$1 - \gamma = \frac{p'}{2} =
\frac{s \gamma}{2}.$$ Then, we have
$$\|\alpha_{\delta}\|_{L_s(\mathcal{N})} = \big\| \xi
\mathrm{D}^{\frac{1}{p'}} \mathsf{E}(aa^* + \delta 1)
\mathrm{D}^{\frac{1}{p'}} \xi^* \big\|_{L_{s
\gamma/2}(\mathcal{N})}^{\gamma/2} = \big\| \xi
\mathrm{D}^{\frac{1}{p'}} \mathsf{E}(aa^* + \delta 1)
\mathrm{D}^{\frac{1}{p'}} \xi^*
\big\|_{L_{p'/2}(\mathcal{N})}^{\gamma/2}$$ and
\begin{eqnarray*}
\|y_{\delta}\|_{L_2(\mathcal{M})} & = & \mbox{tr} \Big( xx^* \big[
\xi \mathrm{D}^{\frac{1}{p'}} \mathsf{E}(aa^* + \delta 1)
\mathrm{D}^{\frac{1}{p'}} \xi^* \big]^{-\gamma} \Big)^{1/2} \\ & =
& \mbox{tr} \Big( \xi \mathrm{D}^{\frac{1}{p'}} \mathsf{E}(aa^*)
\mathrm{D}^{\frac{1}{p'}} \xi^* \big[ \xi
\mathrm{D}^{\frac{1}{p'}} \mathsf{E}(aa^* + \delta 1)
\mathrm{D}^{\frac{1}{p'}} \xi^* \big]^{-\gamma} \Big)^{1/2} \\ &
\le & \mbox{tr} \Big( \xi \mathrm{D}^{\frac{1}{p'}}
\mathsf{E}(aa^* + \delta 1) \mathrm{D}^{\frac{1}{p'}} \xi^* \big[
\xi \mathrm{D}^{\frac{1}{p'}} \mathsf{E}(aa^* + \delta 1)
\mathrm{D}^{\frac{1}{p'}} \xi^* \big]^{-\gamma} \Big)^{1/2}.
\end{eqnarray*}
In particular, $$\|y_{\delta}\|_{L_2(\mathcal{M})} \le \big\| \xi
\mathrm{D}^{\frac{1}{p'}} \mathsf{E}(aa^* + \delta 1)
\mathrm{D}^{\frac{1}{p'}} \xi^*
\big\|_{L_{p'/2}(\mathcal{N})}^{1/2-\gamma/2}.$$ We finally
conclude $$\inf_{x = \alpha y} \|\alpha\|_{L_s(\mathcal{N})}
\|y\|_{L_2(\mathcal{M})} \le \lim_{\delta \rightarrow 0}
\|\alpha_{\delta}\|_{L_s(\mathcal{N})}
\|y_{\delta}\|_{L_2(\mathcal{M})} \le \|x\|_{L_{p'}^r(\mathcal{M},
\mathsf{E})}.$$ Note that the case $p=2$ degenerates since we
obtain $\gamma = 0$ and $s = \infty$. However, this is a trivial
case since it suffices to consider the decomposition given by
$\alpha = 1$ and $y=x$. In summary, we have seen that the norms of
$L_s^r(\mathcal{N}) \otimes_{\mathcal{N}, h} L_2^c(\mathcal{M})$
and of $L_{p'}^r(\mathcal{M}, \mathsf{E})$ coincide on
$$\mathcal{A}_{p'} = \mathcal{N} \, \mathrm{D}^{1/p'}
\mathcal{M}.$$ Moreover, $\mathcal{A}_{p'} = \mathcal{N}
\mathrm{D}^{1/s} \mathrm{D}^{1/2} \mathcal{M}$ is clearly dense in
$L_s^r(\mathcal{N}) \otimes_{\mathcal{N}, h} L_2^c(\mathcal{M})$.
Therefore, since the density of $\mathcal{A}_{p'}$ in
$L_{p'}^r(\mathcal{M}, \mathsf{E})$ follows by definition, we
obtain the desired isometry. The arguments for the column case are
exactly the same. \end{proof}

\begin{remark} \label{Remark-OSS-Row-Column}
\emph{Proposition \ref{Proposition-Norm-Conditional-Inf} provides
an \emph{o.s.s.} for $L_p^r(\mathcal{M}, \mathsf{E})$ and
$L_p^c(\mathcal{M}, \mathsf{E})$ when $1 < p \le 2$. Moreover, by
anti-linear duality we also obtain a natural \emph{o.s.s.} for
$L_p^r(\mathcal{M}, \mathsf{E})$ and $L_p^c(\mathcal{M},
\mathsf{E})$ when $2 \le p < \infty$. However, we shall be
interested only on the Banach space structure of these spaces
rather than the operator space one. Therefore, we shall use the
simpler notation
$$L_u(\mathcal{N}) L_q(\mathcal{M}) \quad \mbox{and} \quad
L_q(\mathcal{M}) L_v(\mathcal{N})$$ to denote the underlying
Banach space of $L_u^r(\mathcal{N}) \otimes_{\mathcal{N},h}
L_q^c(\mathcal{M})$ or $L_q^r(\mathcal{M}) \otimes_{\mathcal{N},h}
L_v^c(\mathcal{N})$.}
\end{remark}

We conclude by giving one more characterization of the norm of
$L_p^r(\mathcal{M}, \mathsf{E})$ and $L_p^c(\mathcal{M},
\mathsf{E})$ for $2 \le p \le \infty$. To that aim, we shall
denote by $\mathsf{Lm}_{\mathcal{N}}(\mathrm{X}_1, \mathrm{X}_2)$
and $\mathsf{Rm}_{\mathcal{N}}(\mathrm{X}_1, \mathrm{X}_2)$ the
spaces of left and right $\mathcal{N}$-module maps between
$\mathrm{X}_1$ and $\mathrm{X}_2$.

\begin{proposition} \label{Proposition-Dual-Module}
Let $\mathcal{N}$ be a von Neumann subalgebra of $\mathcal{M}$.
Then, given any three indices $2 \le u,q,v < \infty$, we have the
following isometric isomorphisms $$\begin{array}{rclcl}
\displaystyle \big( L_u(\mathcal{N}) L_q(\mathcal{M}) \big)^* & =
& \mathsf{Rm}_{\mathcal{N}} (L_u(\mathcal{N}),
L_{q'}(\mathcal{M})) & = & \mathsf{Lm}_{\mathcal{N}}
(L_q(\mathcal{M}), L_{u'}(\mathcal{N})),
\\ \displaystyle \big( L_q(\mathcal{M})
L_v(\mathcal{N}) \big)^* & = & \mathsf{Lm}_{\mathcal{N}}
(L_v(\mathcal{N}), L_{q'}(\mathcal{M})) & = &
\mathsf{Rm}_{\mathcal{N}} (L_q(\mathcal{M}), L_{v'}(\mathcal{N})).
\end{array}$$
\end{proposition}

\begin{proof} The arguments we shall be using hold for both isometries.
Hence, we only prove the first one. Let us consider a linear
functional $\Phi: L_u(\mathcal{N}) L_q(\mathcal{M}) \rightarrow
\C$ and let $\Psi$ denote the associated bilinear map $$\Psi:
L_u(\mathcal{N}) \times L_q(\mathcal{M}) \rightarrow \C,$$ defined
by $\Psi(\alpha,x) = \Phi(\alpha \otimes x)$. Clearly, we have
$$\|\Psi\| = \sup \Big\{ \big| \Phi(\alpha \otimes x) \big| \,
\big| \ \|\alpha\|_{L_u(\mathcal{N})}, \|x\|_{L_q(\mathcal{M})}
\le 1 \Big\} = \|\Phi\|.$$ On the other hand, we use the isometry
$$\mathcal{B}(L_u(\mathcal{N}) \times L_q(\mathcal{M}), \C) =
\mathcal{B}(L_u(\mathcal{N}), L_{q'}(\mathcal{M}))$$ defined by
$\Psi(\alpha,x) = \mbox{tr} \big( \mathrm{T}(\alpha)x \big)$. By
the amalgamation over $\mathcal{N}$, any such linear map
$\mathrm{T}$ is a right $\mathcal{N}$-module map (i.e.
$\mathrm{T}(\alpha \beta) = \mathrm{T}(\alpha) \beta$ for all
$\beta \in \mathcal{N}$). Indeed, given $\beta \in \mathcal{N}$,
we have $$\mbox{tr} \big( \mathrm{T}(\alpha \beta) x \big) =
\Phi(\alpha \beta \otimes x) = \Phi(\alpha \otimes \beta x) =
\mbox{tr} \big( \mathrm{T}(\alpha) \beta x \big)$$ for all $x \in
L_q(\mathcal{M})$, which implies our claim. Reciprocally, if
$\mathrm{T}: L_u(\mathcal{N}) \rightarrow L_{q'}(\mathcal{M})$ is
a right $\mathcal{N}$-module map, we know by the same argument
that $\mathrm{T}$ arises from a bilinear map $\Psi$ satisfying
$\Psi(\alpha \beta, x) = \Psi(\alpha, \beta x)$, which corresponds
to a linear functional $\Phi$ on $L_u(\mathcal{N})
L_q(\mathcal{M})$ with the same norm. This proves the first
identity. For the second we proceed in a similar way by using the
adjoint map $\mathrm{T}^*$ instead of $\mathrm{T}$. \end{proof}

\begin{remark}
\emph{Note that the proof given above still works when we only
require one of the two indices to be finite. In particular, we
cover all combinations appearing for $L_p^r(\mathcal{M},
\mathsf{E})$ and $L_p^c(\mathcal{M}, \mathsf{E})$ since in that
case the $\mathcal{M}$-index is always $2$.}
\end{remark}

\chapter{Amalgamated $L_p$ spaces}
\label{Section2}

Let $\mathcal{M}$ be a von Neumann algebra equipped with a
\emph{n.f.} state $\varphi$ and let us consider a von Neumann
subalgebra $\mathcal{N}$ of $\mathcal{M}$. In what follows we
shall work with indices $(u,q,v)$ satisfying the following
property
\begin{equation} \label{Equation-Indices-1}
1 \le q \le \infty \quad \mbox{and} \quad 2 \le u, v \le \infty
\quad \mbox{and} \quad \frac{1}{u} + \frac{1}{q} + \frac{1}{v} =
\frac{1}{p} \le 1.
\end{equation}
Let us define $\mathcal{N}_u L_q(\mathcal{M}) \mathcal{N}_v$
\label{NLqN} to be the space $L_q(\mathcal{M})$ equipped with
$$|||x|||_{u \cdot q \cdot v} = \inf_{x = \alpha y \beta} \Big\{
\big\| \mathrm{D}^{\frac{1}{u}} \alpha \big\|_{L_u(\mathcal{N})}
\|y\|_{L_q(\mathcal{M})} \big\| \beta \mathrm{D}^{\frac{1}{v}}
\big\|_{L_v(\mathcal{N})} \, \big| \ \alpha, \beta \in
\mathcal{N}, y \in L_q(\mathcal{M}) \Big\}.$$

In the sequel it will be quite useful to have a geometric
representation of the indices $(u,q,v)$ satisfying property
(\ref{Equation-Indices-1}). To that aim, we consider the variables
$(1/u,1/v,1/q)$ in the Euclidean space $\R^3$. Then, we note that
the points satisfying (\ref{Equation-Indices-1}) are given by the
intersection of the simplex $0 \le 1/u + 1/v + 1/q \le 1$ with the
prism $0 \le \min(1/u,1/v) \le \max(1/u,1/v) \le 1/2$. This gives
rise to a solid $\mathsf{K}$ sketched below.

\begin{figure}[ht!]

\begin{center}

\includegraphics[width=10cm]{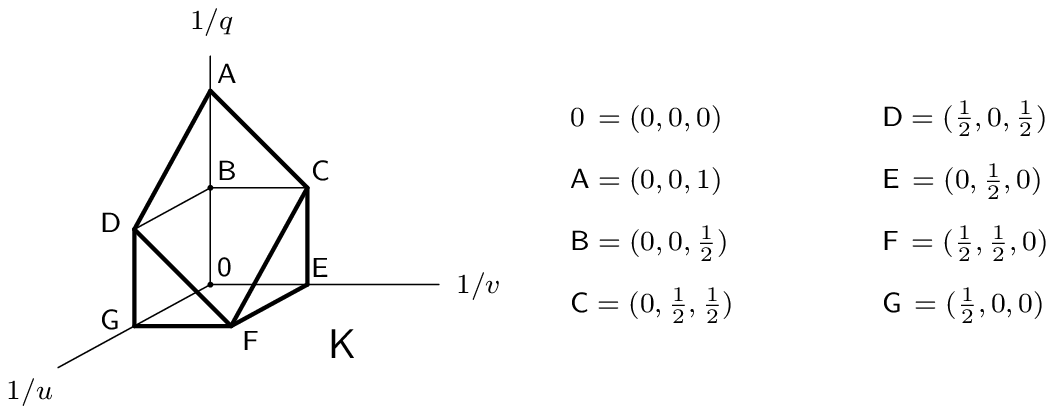}

\end{center}

\caption{\textsc{The solid $\mathsf{K}$.}}

\label{SolidK}

\end{figure}

%\vskip5pt

Now, let $(u,q,v)$ be any three indices satisfying property
(\ref{Equation-Indices-1}). Then we define the \emph{amalgamated
$L_p$ space} $L_u(\mathcal{N}) L_q(\mathcal{M}) L_v(\mathcal{N})$
\label{AmalLp} as the set of elements $x \in L_p(\mathcal{M})$,
with $1/p = 1/u + 1/q + 1/v$, such that there exists a
factorization $x = ayb$ with $a \in L_u(\mathcal{N})$, $y \in
L_q(\mathcal{M})$ and $b \in L_v(\mathcal{N})$. For any such
element $x$, we define
$$\|x\|_{u \cdot q \cdot v} = \inf \Big\{ \|a\|_{L_u(\mathcal{N})}
\|y\|_{L_q(\mathcal{M})} \|b\|_{L_v(\mathcal{N})} \, \big| \ x =
ayb \Big\}.$$ Chapters \ref{Section2}, \ref{Section3} and
\ref{Section4} are devoted to prove that $L_u(\mathcal{N})
L_q(\mathcal{M}) L_v(\mathcal{N})$ is a Banach space for any
$(u,q,v)$ satisfying (\ref{Equation-Indices-1}) and to study the
complex interpolation and duality of such spaces. In the process,
it will be quite relevant to know that the spaces $\mathcal{N}_u
L_q(\mathcal{M}) \mathcal{N}_v$ embed isometrically in
$L_u(\mathcal{N}) L_q(\mathcal{M}) L_v(\mathcal{N})$ as dense
subspaces. This will be essentially the aim of this chapter.

\begin{remark}
\emph{In order to justify the restriction $2 \le u,v \le \infty$
on (\ref{Equation-Indices-1}), let us show how the triangle
inequality might fail when this restriction is not satisfied. Let
us take $\mathcal{M} = \ell_{\infty}(\mathcal{N})$, then we claim
that $$L_w(\mathcal{N}) L_{\infty}(\mathcal{M})
L_{\infty}(\mathcal{N}) \quad \mbox{and} \quad
L_{\infty}(\mathcal{N}) L_{\infty}(\mathcal{M}) L_w(\mathcal{N})$$
are not normed for $1 \le w < 2$. Indeed, given a sequence $x =
(x_n)_{n \ge 1}$ in $L_{\infty}(\mathcal{M})$ and elements $a,b
\in L_w(\mathcal{N})$, we have by definition the following
inequalities $$\|ax\|_{w \cdot \infty \cdot \infty} \le
\|a\|_{L_w(\mathcal{N})} \|x\|_{L_{\infty}(\mathcal{M})} \quad
\mbox{and} \quad \|xb\|_{\infty \cdot \infty \cdot w} \le
\|x\|_{L_{\infty}(\mathcal{M})} \|b\|_{L_w(\mathcal{N})}.$$ In
particular, if these spaces were normed they should contain
contractively the projective tensor product $L_w(\mathcal{N})
\otimes_{\pi} L_{\infty}(\mathcal{M})$. However, by the
noncommutative Menchoff-Rademacher inequality in \cite{DJ}, we may
find $(z_n) \in L_w(\mathcal{N}) \otimes_{\pi}
L_{\infty}(\mathcal{M})$ of the form $$z_n = \sum_{k=1}^n
\varepsilon_k \frac{x_n}{1 + \log k} \quad \mbox{with} \quad x_k
\in L_w(\mathcal{N})$$ and such that $(z_n) \notin
L_w(\mathcal{N}) L_{\infty}(\mathcal{M}) L_{\infty}(\mathcal{N})$.
Moreover, taking adjoints we also may find a sequence $(z_n)$ such
that $(z_n) \in L_w(\mathcal{N}) \otimes_{\pi}
L_{\infty}(\mathcal{M}) \setminus L_{\infty}(\mathcal{N})
L_{\infty}(\mathcal{M}) L_w(\mathcal{N})$.}
\end{remark}

\begin{example}
\emph{As noted in the Introduction, several noncommutative
function spaces arise as particular cases of our notion of
amalgamated $L_p$ space. Let us mention four particularly relevant
examples:}
\begin{itemize}
\item[(a)] \emph{The noncommutative $L_p$ spaces arise as
$$L_q(\mathcal{M}) = L_{\infty}(\mathcal{N}) L_q(\mathcal{M})
L_{\infty}(\mathcal{N}).$$ Note that these spaces are represented
by the segment $0\mathsf{A}$ in Figure I.}

\vskip3pt

\item[(b)] \emph{If $p \le q$ and $1/r = 1/p - 1/q$, the spaces
$L_p(\mathcal{N}_1; L_q(\mathcal{N}_2))$ arise as
$$L_{2r}(\mathcal{N}_1) L_q(\mathcal{N}_1 \bar\otimes
\mathcal{N}_2) L_{2r}(\mathcal{N}_1).$$ When the index $p$ is
fixed, these spaces are represented in Figure I as segments
parallel to the upper face $\mathsf{ACDF}$ whose projection into
the plane $xy$ goes in the direction of the diagonal.}

\vskip3pt

\item[(c)] \emph{By Proposition
\ref{Proposition-Norm-Conditional-Inf}, $L_p^r(\mathcal{M},
\mathsf{E})$ and $L_p^c(\mathcal{M}, \mathsf{E})$ $(1 \le p \le
2)$ arise as
\begin{eqnarray*}
L_p^r(\mathcal{M}, \mathsf{E}) & = & L_s(\mathcal{N})
L_2(\mathcal{M}) L_{\infty}(\mathcal{N}) \quad \mbox{with} \quad
1/p = 1/2 + 1/s, \\ L_p^c(\mathcal{M}, \mathsf{E}) & = &
L_{\infty}(\mathcal{N}) L_2(\mathcal{M}) L_s(\mathcal{N}) \quad
\mbox{with} \quad 1/p = 1/2 + 1/s.
\end{eqnarray*}
These spaces are represented by the segments $\mathsf{BD}$ and
$\mathsf{BC}$ in Figure I. As we shall see in Chapter
\ref{Section4}, the spaces $L_p(\mathcal{M}; R_p^n)$ and
$L_p(\mathcal{M}; C_p^n)$ arise as particular cases of the latter
ones.}

\vskip3pt

\item[(d)] \emph{If $\mathcal{M} = \mathcal{N}$, the so-called
asymmetric noncommutative $L_p$ spaces arise as
$$L_{(u,v)}(\mathcal{M}) = L_u(\mathcal{M})
L_{\infty}(\mathcal{M}) L_v(\mathcal{M}).$$ These spaces are
represented by the square $0\mathsf{EFG}$ in Figure I. The reader
is refereed to \cite{JP} and Chapter \ref{Section7} below for a
detailed exposition of the main properties of these spaces and
their vector-valued analogues.}
\end{itemize}
\end{example}

\section{Haagerup's construction}
\label{Subsection-Preliminaries}

\label{Uffe2}

We now briefly sketch a well-known unpublished result due to
Haagerup \cite{H2} which will be essential for our further
purposes. Let us consider a $\sigma$-finite von Neumann algebra
$\mathcal{M}$ equipped with a distinguished \emph{n.f.} state
$\varphi$ and let us define the discrete multiplicative group
$$\mathrm{G} = \bigcup_{n \in \N} 2^{-n} \Z.$$ Then we can
construct the crossed product $\mathcal{R} = \mathcal{M}
\rtimes_{\sigma} \mathrm{G}$ in the same way we constructed the
crossed product $\mathcal{M} \rtimes_{\sigma} \R$. That is, if
$\mathcal{H}$ is the Hilbert space provided by the GNS
construction applied to $\varphi$, $\mathcal{R}$ is generated by
$\pi: \mathcal{M} \to \mathcal{B}(L_2(\mathrm{G}; \mathcal{H}))$
and $\lambda: \mathrm{G} \to \mathcal{B}(L_2(\mathrm{G};
\mathcal{H}))$ where
$$\big( \pi(x) \xi \big) (g) = \sigma_{-g}(x) \xi(g) \qquad
\mbox{and} \qquad \big( \lambda(h) \xi \big) (g) = \xi(g-h).$$ By
the faithfulness of $\pi$ we are allowed to identify $\mathcal{M}$
with $\pi(\mathcal{M})$. Then, a generic element in $\mathcal{R}$
has the form $\sum_g x_g \lambda(g)$ with $x_g \in \mathcal{M}$
and we have the conditional expectation
$$\mathsf{E}_{\mathcal{M}}: \sum_{g \in \mathrm{G}} x_g \lambda(g)
\in \mathcal{R} \mapsto x_0 \in \mathcal{M}.$$ This gives rise to
the state $$\widehat{\varphi}: \sum_{g \in \mathrm{G}} x_g
\lambda(g) \in \mathcal{R} \mapsto \varphi \circ
\mathsf{E}_{\mathcal{M}} \Big( \sum_{g \in \mathrm{G}} x_g
\lambda(g) \Big) = \varphi(x_0) \in \C.$$ According to \cite{H2},
$\mathcal{R}$ is the closure of a union of finite von Neumann
algebras $$\bigcup_{k \ge 1} \mathcal{M}_k$$ where
$(\mathcal{M}_k)_{k \ge 1}$ is directed by inclusion
$\mathcal{M}_1 \subset \mathcal{M}_2 \subset \ldots$ and each
$\mathcal{M}_k$ satisfies
\begin{equation} \label{Equation-Boundedness-RN-Derivative-2}
c_1(k) 1_{\mathcal{M}_k} \le \mathrm{D}_{\varphi_k} \le c_2(k)
1_{\mathcal{M}_k}
\end{equation}
for some constants $0 < c_1(k) \le c_2(k) < \infty$. Here
$\varphi_k$ denotes the restriction of $\widehat{\varphi}$ to
$\mathcal{M}_k$ and $\mathrm{D}_{\varphi_k}$ stands for the
corresponding density. Moreover, it also follows from \cite{H2}
that we can find for any integer $k \ge 1$ conditional
expectations $$\mathcal{E}_k: \mathcal{R} \to \mathcal{M}_k$$ such
that the following limit holds for every $\hat{x} \in
\mathcal{R}$, any $1 \le p < \infty$ and all $0 \le \eta \le 1$
\begin{equation} \label{Equation-Limit-Expectation}
\Big\| \mathrm{D}_{\widehat{\varphi}}^{(1-\eta)/p} \big(\hat{x} -
\mathcal{E}_k(\hat{x}) \big) \,
\mathrm{D}_{\widehat{\varphi}}^{\eta/p} \Big\|_{L_p(\mathcal{R})}
\longrightarrow 0 \qquad \mbox{as} \qquad k \rightarrow \infty.
\end{equation}

\begin{remark}
\emph{The $\sigma$-finiteness assumption might be dropped if we
replaced sequences of finite von Neumann algebras by nets.
However, it suffices for our aims to consider the $\sigma$-finite
case.}
\end{remark}

\vskip3pt

In the following result we provide the first application of this
construction. As usual, we write $\mathcal{S}$ for the strip of
complex numbers $z \in \C$ with $0 < \mathrm{Re}(z) < 1$ where we
consider the decomposition $\partial \mathcal{S} =
\partial_0 \cup \partial_1$ of its boundary into the sets
\label{Band} $$\partial_0 = \Big\{ z \in \C \, \big| \
\mbox{Re}(z) = 0 \Big\} \quad \mbox{and} \quad \partial_1 = \Big\{
z \in \C \, \big| \ \mbox{Re}(z) = 1 \Big\}.$$ Let $\mathrm{X}$ be
a Banach space and let $x$ be an element of $\mathrm{X}$. Then,
given $0 < \theta < 1$, we write
$\mathcal{A}(\mathrm{X},\theta,x)$ for set of bounded analytic
functions $f: \mathcal{S} \rightarrow \mathrm{X}$ (i.e. bounded
and continuous on $\mathcal{S} \cup
\partial \mathcal{S}$ and analytic on $\mathcal{S}$) which
satisfy $f(\theta) = x$.

\begin{lemma} \label{Lemma-Kosaki}
Let $x$ be an element in a von Neumann algebra $\mathcal{M}$.
Then, given $0 \le \eta \le 1$ and $1/p_{\theta} = (1-\theta)/p$
for some $0 < \theta < 1$ and some $1 \le p < \infty$, we have
$$\big\| \mathrm{D}^{\frac{1-\eta}{p_{\theta}}} x
\mathrm{D}^{\frac{\eta}{p_{\theta}}}
\big\|_{L_{p_{\theta}}(\mathcal{M})} \!\! = \inf_{} \left\{ \max
\Big\{ \! \sup_{z \in \partial_0} \big\|
\mathrm{D}^{\frac{1-\eta}{p}} f(z) \mathrm{D}^{\frac{\eta}{p}}
\big\|_{L_p(\mathcal{M})}, \sup_{z \in \partial_1} \big\| f(z)
\big\|_{L_{\infty}(\mathcal{M})} \Big\} \right\},$$ where the
infimum runs over all analytic functions $f: \mathcal{S}
\rightarrow \mathcal{M}$ in the set
$\mathcal{A}(\mathcal{M},\theta,x)$.
\end{lemma}

\begin{proof} The upper estimate follows by Kosaki's interpolation. To
prove the lower estimate, we assume by homogeneity that the norm
on the left is $1$. Then we begin with the case where
$\mathcal{M}$ is a finite von Neumann algebra equipped with a
\emph{n.f.} state $\varphi_0$ such that the density
$\mathrm{D}_{\varphi_0}$ satisfies
(\ref{Equation-Boundedness-RN-Derivative-2}). That is, $$c_1
1_{\mathcal{M}} \le \mathrm{D}_{\varphi_0} \le c_2
1_{\mathcal{M}}$$ for some positive numbers $0 < c_1 \le c_2 <
\infty$. Note that our assumption implies $\mathrm{D}_{\varphi_0}
\in \mathcal{M}$. In particular, the function $z \to
\mathrm{D}_{\varphi_0}^{\lambda z}$ is analytic for any scalar
$\lambda \in \C$. Then it is easy to find an optimal function.
Indeed, since $$x(\eta,\theta) =
\mathrm{D}_{\varphi_0}^{(1-\eta)/p_{\theta}} x \,
\mathrm{D}_{\varphi_0}^{\eta/p_{\theta}} \in \mathcal{M}$$ we find
by polar decomposition a partial isometry $\omega$ such that
$x(\eta,\theta) = \omega |x(\eta, \theta)|$. Our optimal function
is defined as $$f(z) = \mathrm{D}_{\varphi_0}^{-
\frac{(1-\eta)(1-z)}{p}} \omega \, |x(\eta,
\theta)|^{\frac{p_{\theta}(1-z)}{p}} \, \mathrm{D}_{\varphi_0}^{-
\frac{\eta (1-z)}{p}}.$$ Our assumption on
$\mathrm{D}_{\varphi_0}$ implies the boundedness and analyticity
of $f$. On the other hand, it is easy to check that $f(\theta) =
x$ so that $f \in \mathcal{A}(\mathcal{M}, \theta,x)$. Then,
recalling that both $|x(\eta, \theta)|^{i \lambda}$ and
$\mathrm{D}_{\varphi_0}^{i \lambda}$ are unitaries for any
$\lambda \in \R$, we have
\begin{eqnarray*}
\sup_{z \in \partial_0} \Big\|
\mathrm{D}_{\varphi_0}^{\frac{1-\eta}{p}} f(z)
\mathrm{D}_{\varphi_0}^{\frac{\eta}{p}} \Big\|_{L_p(\mathcal{M})}
& = & \sup_{z \in \partial_0} \Big\|
\mathrm{D}_{\varphi_0}^{\frac{(1-\eta) z}{p}} \omega |x(\eta,
\theta)|^{\frac{p_{\theta}}{p}} |x(\eta,
\theta)|^{-\frac{p_{\theta}z}{p}}
\mathrm{D}_{\varphi_0}^{\frac{\eta z}{p}}
\Big\|_{L_p(\mathcal{M})} \\ & \le & \Big\| |x(\eta,
\theta)|^{\frac{p_{\theta}}{p}} \Big\|_{L_p(\mathcal{M})} = \Big\|
\mathrm{D}_{\varphi_0}^{\frac{1-\eta}{p_{\theta}}} x
\mathrm{D}_{\varphi_0}^{\frac{\eta}{p_{\theta}}}
\Big\|_{L_{p_{\theta}}(\mathcal{M})}^{\frac{1}{1-\theta}} = 1.
\end{eqnarray*}
Similarly, we have on $\partial_1$ $$\sup_{z \in \partial_1}
\big\| f(z) \big\|_{L_{\infty}(\mathcal{M})} = \sup_{z \in
\partial_1} \Big\| \mathrm{D}_{\varphi_0}^{\frac{i (1-\eta)
\mathrm{Im} z}{p}} \omega |x(\eta, \theta)|^{- \frac{i p_{\theta}
\mathrm{Im} z}{p}} \mathrm{D}_{\varphi_0}^{\frac{i \eta
\mathrm{Im} z}{p}} \Big\|_{L_{\infty}(\mathcal{M})} \le 1.$$ This
completes the proof for finite von Neumann algebras satisfying
(\ref{Equation-Boundedness-RN-Derivative-2}). In the case of a
general von Neumann algebra $\mathcal{M}$ we use the Haagerup
construction sketched above. Let us introduce the shorter notation
\begin{eqnarray*}
\mathsf{N}_1 (\mathcal{M}, \varphi, x) & = & \big\|
\mathrm{D}^{\frac{1-\eta}{p_{\theta}}} x \,
\mathrm{D}^{\frac{\eta}{p_{\theta}}}
\big\|_{L_{p_{\theta}}(\mathcal{M})}, \\ \mathsf{N}_2
(\mathcal{M}, \varphi, x) & = & \inf_{} \left\{ \max \Big\{ \!
\sup_{z \in \partial_0} \big\| \mathrm{D}^{\frac{1-\eta}{p}} f(z)
\mathrm{D}^{\frac{\eta}{p}} \big\|_{L_p(\mathcal{M})}, \sup_{z \in
\partial_1} \big\| f(z) \big\|_{L_{\infty}(\mathcal{M})} \Big\}
\right\}.
\end{eqnarray*}
We are interested in proving $\mathsf{N}_1 (\mathcal{M}, \varphi,
x) \ge \mathsf{N}_2 (\mathcal{M}, \varphi, x)$. Combining property
(\ref{Equation-Boundedness-RN-Derivative-2}) with the first part
of this proof, we deduce that any $x \in \mathcal{M}$ satisfies
the following identities
\begin{equation} \label{Equation-N1N2}
\mathsf{N}_1 (\mathcal{M}_k, \varphi_k, \mathcal{E}_k(x)) =
\mathsf{N}_2 (\mathcal{M}_k, \varphi_k, \mathcal{E}_k(x)) \quad
\mbox{for all} \quad k \ge 1.
\end{equation}
By the triangle inequality $$\mathsf{N}_2 (\mathcal{M}, \varphi,
x) \le \limsup_{k \to \infty} \mathsf{N}_2 (\mathcal{M}, \varphi,
\mathsf{E}_{\mathcal{M}} (\mathcal{E}_k(x))) + \mathsf{N}_2
(\mathcal{M}, \varphi, x - \mathsf{E}_{\mathcal{M}}
(\mathcal{E}_k(x))).$$ Applying the contractivity of
$\mathsf{E}_{\mathcal{M}}$, (\ref{Equation-Limit-Expectation}) and
(\ref{Equation-N1N2})
\begin{eqnarray*}
\limsup_{k \to \infty} \mathsf{N}_2 (\mathcal{M}, \varphi,
\mathsf{E}_{\mathcal{M}} (\mathcal{E}_k(x))) & \le & \limsup_{k
\to \infty} \mathsf{N}_2 (\mathcal{M}_k, \varphi_k,
\mathcal{E}_k(x))
\\ & = & \limsup_{k \to \infty} \mathsf{N}_1 (\mathcal{M}_k,
\varphi_k, \mathcal{E}_k(x))
\\ & = & \mathsf{N}_1 (\mathcal{M}, \varphi, x)
\end{eqnarray*}
It remains to see that $\mathsf{N}_2 (\mathcal{M}, \varphi, x -
\mathsf{E}_{\mathcal{M}} (\mathcal{E}_k(x)))$ is arbitrary small.
First we note that $$\mathsf{N}_2 (\mathcal{M}, \varphi, x -
\mathsf{E}_{\mathcal{M}} (\mathcal{E}_k(x))) = \mathsf{N}_2
(\mathcal{M}, \varphi, \mathsf{E}_{\mathcal{M}} (x -
\mathcal{E}_k(x))) \le \mathsf{N}_2 (\mathcal{R},
\widehat{\varphi}, x - \mathcal{E}_k(x)).$$ Then we consider the
bounded analytic functions $$f_k: z \in \mathcal{S} \mapsto
\mathsf{m}_{0k}^{1-z} \mathsf{m}_{1k}^z \big( x - \mathcal{E}_k(x)
\big) \in \mathcal{R}$$ with the constants $\mathsf{m}_{0k},
\mathsf{m}_{1k}$ defined by
\begin{eqnarray*}
\mathsf{m}_{0k} & = & \Big\|
\mathrm{D}_{\widehat{\varphi}}^{(1-\eta)/p} \big(x -
\mathcal{E}_k(x) \big) \, \mathrm{D}_{\widehat{\varphi}}^{\eta/p}
\Big\|_{L_p(\mathcal{R})}^{-1/2}, \\ \mathsf{m}_{1k} & = & \Big\|
\mathrm{D}_{\widehat{\varphi}}^{(1-\eta)/p} \big(x -
\mathcal{E}_k(x) \big) \, \mathrm{D}_{\widehat{\varphi}}^{\eta/p}
\Big\|_{L_p(\mathcal{R})}^{(1-\theta)/2\theta}.
\end{eqnarray*}
Note that $f_k \in \mathcal{A} \big(\mathcal{R}, \theta, x -
\mathcal{E}_k(x) \big)$ since $\mathsf{m}_{0k}^{1-\theta}
\mathsf{m}_{1k}^{\theta} = 1$. On the other hand, $$\sup_{z \in
\partial_0} \Big\| \mathrm{D}_{\widehat{\varphi}}^{(1-\eta)/p}
f_k(z) \mathrm{D}_{\widehat{\varphi}}^{\eta/p}
\Big\|_{L_p(\mathcal{R})} = \mathsf{m}_{0k} \Big\|
\mathrm{D}_{\widehat{\varphi}}^{(1-\eta)/p} \big(x -
\mathcal{E}_k(x) \big) \, \mathrm{D}_{\widehat{\varphi}}^{\eta/p}
\Big\|_{L_p(\mathcal{R})}.$$ Similarly, we have on $\partial_1$
$$\sup_{z \in \partial_1} \big\| f_k(z)
\big\|_{L_{\infty}(\mathcal{R})} = \mathsf{m}_{1k} \big\| x -
\mathcal{E}_k(x) \big\|_{L_{\infty}(\mathcal{R})} \le 2
\mathsf{m}_{1k} \|x\|_{\mathcal{M}}.$$ Therefore, the proof is
completed since from (\ref{Equation-Limit-Expectation}) both terms
tend to $0$ with $k$. \end{proof}

\section{Triangle inequality on $\partial_{\infty} \mathsf{K}$}

\label{BoundaryK}

Let $\partial_{\infty} \mathsf{K}$ be the subset of the boundary
of $\mathsf{K}$ given by the intersection of $\mathsf{K}$ with the
coordinate planes (i.e. the union of the plane regions
$0\mathsf{ACE}$, $0\mathsf{ADG}$ and $0\mathsf{EFG}$). This set
will appear repeatedly in what follows. Note that the indices
$(u,q,v)$ which are represented in Figure I by the points of
$\partial_{\infty} \mathsf{K}$ are those satisfying
(\ref{Equation-Indices-1}) and
\begin{equation} \label{Equation-Indices-2}
\min \Big\{ 1/u, 1/q, 1/v \Big\} = 0.
\end{equation}

\vskip3pt

In the following result, we apply a complex interpolation trick
based on the operator-valued version of Szeg\"{o}'s classical
factorization theorem. This result is due to Devinatz \cite{D}. A
precise statement of Devinatz's theorem adapted to our aims can be
found in Pisier's paper \label{Devth} \cite{P0}, where he also
applies it in the context of complex interpolation. We note in
passing that a more general result (combining previous results due
to Devinatz, Helson \& Lowdenslager, Sarason and Wiener \& Masani)
can be found in the survey \cite{PX2}, see Corollary 8.2. This
interpolation technique will be used frequently in the sequel.

\begin{lemma} \label{Lemma-Triangle-Inequality}
If $(1/u,1/v,1/q) \in \partial_{\infty} \mathsf{K}$,
$\mathcal{N}_u L_q(\mathcal{M}) \mathcal{N}_v$ is a normed space.
\end{lemma}

\begin{proof} Assuming (\ref{Equation-Indices-1}), H\"{o}lder
inequality gives for $1/p = 1/u + 1/q + 1/v$ $$\big\|
\mathrm{D}^{\frac{1}{u}} x \mathrm{D}^{\frac{1}{v}}
\big\|_{L_p(\mathcal{M})} \le |||x|||_{u \cdot q \cdot v} \qquad
\mbox{for all} \qquad x \in \mathcal{N}_u L_q(\mathcal{M})
\mathcal{N}_v.$$ Therefore, $|||x|||_{u \cdot q \cdot v} = 0$
implies $x = 0$. Since the homogeneity over $\R_+$ is clear, it
remains to show that the triangle inequality holds in this case.
Our assumption $(1/u,1/v,1/q) \in
\partial_{\infty} \mathsf{K}$ reduces the possibilities to those
in which at least one of the indices is infinite. When $1/q = 0$,
Pisier's factorization argument (\emph{c.f.} Lemma 3.5 in
\cite{P2}) suffices to obtain the triangle inequality. Indeed, it
can be checked that this factorization provides the estimate
$$|||x_1 + x_2|||_{u \cdot q \cdot v} \le 2^{1/q} \big(
|||x_1|||_{u \cdot q \cdot v} + |||x_2|||_{u \cdot q \cdot v}
\big).$$ Hence, it remains to consider the cases $1/u = 0$ and
$1/v = 0$. Both can be treated with the same arguments so that we
only show the triangle inequality for $1/u = 0$. In that case, the
left term $\mathcal{N}_{\infty}$ is irrelevant in
$\mathcal{N}_{\infty} L_q(\mathcal{M}) \mathcal{N}_v$ so that we
shall ignore it in what follows. We need to consider two different
cases.

\vskip3pt

\noindent \textsc{Case I.} We first assume that $1 \le q \le 2$.
Let $x_1, x_2, \ldots, x_m$ be a finite sequence of vectors in
$L_q(\mathcal{M}) \mathcal{N}_v$ satisfying $|||x_k|||_{q \cdot v}
< 1$ and let us also consider a finite sequence of positive
numbers $\lambda_1, \lambda_2, \ldots, \lambda_m$ with $\summ_k
\lambda_k = 1$. Then, it clearly suffices to see that
$$\Big|\Big|\Big| \sum_{k=1}^m \lambda_k x_k \Big|\Big|\Big|_{q
\cdot v} \le 1.$$ By hypothesis, we may assume that $x_k = y_k
\beta_k,$ with $$\max \Big\{ \|y_k\|_{L_q(\mathcal{M})}, \big\|
\beta_k \mathrm{D}^{\frac{1}{v}} \big\|_{L_v(\mathcal{N})} \Big\}
< 1.$$ On the other hand, since $1 \le q \le 2$ we can consider $2
\le q_1 \le \infty$ defined by $$\frac{1}{q_1} + \frac{1}{2} =
\frac{1}{p} = \frac{1}{q} + \frac{1}{v}.$$ It is not difficult to
check that
\begin{eqnarray*}
L_q(\mathcal{M}) & = & \big[ L_p(\mathcal{M}),
L_{q_1}(\mathcal{M}) \big]_{2/v}, \\ L_v(\mathcal{N}) & = & \big[
L_{\infty}(\mathcal{N}) \, , \, L_2(\mathcal{N}) \big]_{2/v}.
\end{eqnarray*}
Then, by the complex interpolation method we can find bounded
analytic functions $$f_{1k}: \mathcal{S} \to L_p(\mathcal{M}) +
L_{q_1}(\mathcal{M})$$ satisfying $f_{1k}(2/v) = y_k$ for all
$k=1,2,\ldots,m$ and
\begin{equation} \label{Equation-Boundary-0}
\begin{array}{c}
\displaystyle \sup_{z \in \partial_0} \big\| f_{1k}(z)
\big\|_{L_p(\mathcal{M})} \ < 1, \\ \displaystyle \sup_{z \in
\partial_1} \big\| f_{1k}(z) \big\|_{L_{q_1}(\mathcal{M})} < 1.
\end{array}
\end{equation}
Similarly, Lemma \ref{Lemma-Kosaki} provides us with bounded
analytic functions $$f_{2k}: \mathcal{S} \rightarrow \mathcal{N}$$
satisfying $f_{2k}(2/v) = \beta_k$ for all $k=1,2,\ldots,m$ and
$$\begin{array}{r} \displaystyle \sup_{z \in \partial_0} \big\|
f_{2k}(z) \big\|_{L_{\infty}(\mathcal{N})} < 1, \\ \displaystyle
\sup_{z \in
\partial_1} \big\| f_{2k}(z) \mathrm{D}^{\frac{1}{2}}
\big\|_{L_2(\mathcal{N})} \ < 1.
\end{array}$$ Given $\delta > 0$, we define the following function
on $\partial \mathcal{S}$ $$\mathrm{W}(z) = \left\{
\begin{array}{ll} 1 & \mbox{if} \ z \in \partial_0, \\ \delta 1 +
\summ_k \lambda_k f_{2k}(z)^* f_{2k}(z) & \mbox{if} \ z \in
\partial_1.
\end{array} \right.$$ According to Devinatz's factorization
theorem \cite{P0}, there exists a bounded analytic function
$\mathrm{w}: \mathcal{S} \rightarrow \mathcal{N}$ with bounded
analytic inverse and satisfying the following identity on
$\partial \mathcal{S}$ $$\mathrm{w}(z)^* \mathrm{w}(z) =
\mathrm{W}(z).$$ Let us consider the bounded analytic function
$$g(z) = \Big( \sum_{k=1}^m \lambda_k \, f_{1k}(z) f_{2k}(z) \Big)
\mathrm{w}^{-1}(z).$$ Then, since $\mathrm{w}(z)^* \mathrm{w}(z)
\equiv 1$ on $\partial_0$, we have $$\sup_{z \in \partial_0}
\big\| g(z) \big\|_{L_p(\mathcal{M})} \le \sup_{z \in
\partial_0} \sum_{k=1}^m \lambda_k \big\| f_{1k}(z)
\big\|_{L_p(\mathcal{M})} \big\| f_{2k}(z) \mathrm{w}^{-1}(z)
\big\|_{L_{\infty}(\mathcal{N})} < 1.$$ On the other hand, we can
write $$g(z) = \sum_{k=1}^m \sqrt{\lambda_k} \, f_{1k}(z)
\gamma_k(z) \quad \mbox{with} \quad \sqrt{\lambda_k} f_{2k}(z) =
\gamma_k(z) \mathrm{w}(z).$$ Note that, according to the
definition of $\mathrm{w}$, we have for any $z \in
\partial_1$
\begin{equation} \label{Equation-Indentity-Bound}
\sum_{k=1}^m \gamma_k(z)^* \gamma_k(z) \le 1.
\end{equation}
Therefore, the following estimate follows from
(\ref{Equation-Boundary-0}) and (\ref{Equation-Indentity-Bound})
\begin{eqnarray*}
\sup_{z \in \partial_1} \big\| g(z) \big\|_{L_{q_1}(\mathcal{M})}
& \le & \sup_{z \in \partial_1} \Big\| \Big( \sum_{k=1}^m
\lambda_k f_{1k}(z) f_{1k}(z)^* \Big)^{1/2}
\Big\|_{L_{q_1}(\mathcal{M})} \\ & \times & \sup_{z \in
\partial_1} \Big\| \Big( \sum_{k=1}^m \gamma_k(z)^* \gamma_k(z)
\Big)^{1/2} \Big\|_{L_\infty(\mathcal{N})} \\ & \le & \sup_{z \in
\partial_1} \Big( \sum_{k=1}^m \lambda_k \big\| f_{1k}(z)
\big\|_{L_{q_1}(\mathcal{M})}^2 \Big)^{1/2} < 1.
\end{eqnarray*}
Note that the last estimate uses the triangle inequality on
$L_{q_1/2}(\mathcal{M})$ and that we are allowed to do so since
$q_1 \ge 2$. Combining the estimates obtained so far and applying
Kosaki's interpolation we obtain $$\big\| g(2/v)
\big\|_{L_q(\mathcal{M})} < 1.$$ On the other hand, we have
\begin{eqnarray*}
\sup_{z \in \partial_0} \big\| \mathrm{w}(z)
\big\|_{L_{\infty}(\mathcal{N})} & = & \sup_{z \in
\partial_0} \big\| \mathrm{w}(z)^*
\mathrm{w}(z) \big\|_{L_{\infty}(\mathcal{N})}^{1/2} = 1, \\
\sup_{z \in \partial_1} \big\| \mathrm{w}(z)
\mathrm{D}^{\frac{1}{2}} \big\|_{L_2(\mathcal{N})}^2 & = & \sup_{z
\in \partial_1} \big\| \mathrm{D}^{\frac{1}{2}} \mathrm{w}(z)^*
\mathrm{w}(z) \mathrm{D}^{\frac{1}{2}} \big\|_{L_1(\mathcal{N})}
\\ & \le & \sup_{z \in \partial_1} \sum_{k=1}^m \lambda_k \big\|
f_{2k}(z) \mathrm{D}^{\frac{1}{2}} \big\|_{L_2(\mathcal{N})}^2 +
\delta < 1 + \delta.
\end{eqnarray*}
Again, Kosaki's interpolation provides the estimate $$\big\|
\mathrm{w}(2/v) \mathrm{D}^{\frac{1}{v}} \big\|_{L_v(\mathcal{N})}
\le 1 + \delta.$$ In summary, recalling that $$\sum_{k=1}^m
\lambda_k x_k = \sum_{k=1}^m \lambda_k \, f_{1k}(2/v) f_{2k}(2/v)
= g(2/v) \mathrm{w}(2/v),$$ we obtain from our previous estimates
that $$\Big|\Big|\Big| \sum_{k=1}^m \lambda_k x_k
\Big|\Big|\Big|_{q \cdot v} < 1 + \delta.$$ Thus, the triangle
inequality follows by letting $\delta \rightarrow 0$ in the
expression above.

\vskip3pt

\noindent \textsc{Case II.} It remains to consider the case $2 \le
q \le \infty$. This case is simpler. Indeed, given any family of
vectors $x_1, x_2, \ldots, x_m$ and scalars $\lambda_1, \lambda_2,
\ldots, \lambda_m$ as above (with $x_k = y_k \beta_k$ and $y_k,
\beta_k$ satisfying the same inequalities), we define for $\delta
> 0$ $$\mathcal{S}_{\beta} = \Big(
\sum_{k=1}^m \lambda_k \beta_k^* \beta_k + \delta 1 \Big)^{1/2}
\qquad \mbox{so that} \qquad \beta_k = b_k \mathcal{S}_{\beta},$$
for some $b_1, b_2, \ldots, b_m \in \mathcal{N}$. Then we have a
factorization $$\sum_{k=1}^m \lambda_k x_k = \Big( \sum_{k=1}^m
\lambda_k y_k b_k \Big) \mathcal{S}_{\beta}.$$ Now, since $2 \le q
\le \infty$, we have triangle inequality in $L_{q/2}(\mathcal{M})$
and
\begin{eqnarray*}
\Big\| \sum_{k=1}^m \lambda_k y_k b_k \Big\|_{L_q(\mathcal{M})} &
\le & \Big\| \Big( \sum_{k=1}^m \lambda_k y_k y_k^* \Big)^{1/2}
\Big\|_{L_q(\mathcal{M})} \, \Big\| \Big( \sum_{k=1}^m \lambda_k
b_k^* b_k \Big)^{1/2} \Big\|_{L_{\infty}(\mathcal{N})} \\ & \le &
\Big( \sum_{k=1}^m \lambda_k \|y_k\|_{L_q(\mathcal{M})}^2
\Big)^{1/2} \, \Big\| \mathcal{S}_{\beta}^{-1} \Big( \sum_{k=1}^m
\lambda_k \beta_k^* \beta_k \Big) \mathcal{S}_{\beta}^{-1}
\Big\|_{L_{\infty}(N)}^{1/2},
\end{eqnarray*}
which is clearly bounded by $1$. On the other hand,
\begin{eqnarray*}
\big\| \mathcal{S}_{\beta} \mathrm{D}^{\frac{1}{v}}
\big\|_{L_v(\mathcal{N})}^2 & = & \Big\| \mathrm{D}^{\frac{1}{v}}
\Big( \sum_{k=1}^m \lambda_k \beta_k^* \beta_k + \delta 1 \Big)
\mathrm{D}^{\frac{1}{v}} \Big\|_{L_{v/2}(\mathcal{N})} \\ & \le &
\sum_{k=1}^m \lambda_k \big\| \beta_k \mathrm{D}^{\frac{1}{v}}
\big\|_{L_v(\mathcal{N})}^2 + \delta < 1 + \delta.
\end{eqnarray*}
Therefore, the triangle inequality follows one more time by
letting $\delta \to 0$. \end{proof}

The following will be a key point in the proof of Theorem
\ref{TheoremA1}, the main result in Chapter \ref{Section3}. Note
that we state it under the assumption that $\mathcal{N}_u
L_q(\mathcal{M}) \mathcal{N}_v$ is a normed space and that for the
moment we only know it (from Lemma
\ref{Lemma-Triangle-Inequality}) whenever $(1/u,1/v,1/q) \in
\partial_{\infty} \mathsf{K}$. However, since eventually we shall need to
apply this result for any point in $\mathsf{K}$, we state it in
full generality.

\begin{proposition} \label{Proposition-Isometry-Triangle-Complete}
If $\mathcal{N}_u L_q(\mathcal{M}) \mathcal{N}_v$ is a normed
space, we have
\begin{itemize}
\item[i)] The following map is an isometry $$j_{u,v}: x \in
\mathcal{N}_u L_q(\mathcal{M}) \mathcal{N}_v \mapsto
\mathrm{D}^{\frac{1}{u}} x \mathrm{D}^{\frac{1}{v}} \in
L_u(\mathcal{N}) L_q(\mathcal{M}) L_v(\mathcal{N}).$$ \item[ii)]
$L_u(\mathcal{N}) L_q(\mathcal{M}) L_v(\mathcal{N})$ is a Banach
space which completes $\mathcal{N}_u L_q(\mathcal{M})
\mathcal{N}_v$.
\end{itemize}
\end{proposition}

\begin{proof} It is clear that $$\|j_{u,v}(x)\|_{u \cdot q \cdot v} \le
|||x|||_{u \cdot q \cdot v} \qquad \mbox{for all} \qquad x \in
\mathcal{N}_u L_q(\mathcal{M}) \mathcal{N}_v.$$ Let us see that
the reverse inequality holds. Assume it does not hold, then we can
find $x_0 \in \mathcal{N}_u L_q(\mathcal{M}) \mathcal{N}_v$ such
that $|||x_0|||_{u \cdot q \cdot v} > 1$ and $\|j_{u,v}(x_0)\|_{u
\cdot q \cdot v} < 1$. By the Hahn-Banach theorem (here we use the
assumption that the space $\mathcal{N}_u L_q(\mathcal{M})
\mathcal{N}_v$ is normed) there exists a norm one functional
$$\varphi: \mathcal{N}_u L_q(\mathcal{M}) \mathcal{N}_v
\rightarrow \C$$ satisfying
\begin{itemize}
\item[(a)] $\, \varphi(x_0) = |||x_0|||_{u \cdot q \cdot v} > 1$.
\vskip4pt \item[(b)] $\big| \varphi \big( \alpha y \beta \big)
\big| \le \big\| \mathrm{D}^{\frac{1}{u}} \alpha
\big\|_{L_u(\mathcal{N})} \|y\|_{L_q(\mathcal{M})} \big\| \beta
\mathrm{D}^{\frac{1}{v}} \big\|_{L_v(\mathcal{N})}.$
\end{itemize}

\vskip3pt

\noindent Note that any $y \in L_q(\mathcal{M})$ provides a
densely defined bilinear map $$\Phi_y: \big(
\mathrm{D}^{\frac{1}{u}} \alpha, \beta \mathrm{D}^{\frac{1}{v}}
\big) \in L_u(\mathcal{N}) \times L_v(\mathcal{N}) \mapsto \varphi
\big( \alpha y \beta \big) \in \C,$$ which satisfies $\|\Phi_y\|
\le \|y\|_q$ and $$\Phi_{n_1yn_2} \big( \mathrm{D}^{\frac{1}{u}}
\alpha, \beta \mathrm{D}^{\frac{1}{v}} \big) = \Phi_y \big(
\mathrm{D}^{\frac{1}{u}} \alpha n_1, n_2 \beta
\mathrm{D}^{\frac{1}{v}} \big) \quad \mbox{for all} \quad n_1, n_2
\in \mathcal{N}.$$ On the other hand, since $\|j_{u,v}(x_0)\|_{u
\cdot q \cdot v} < 1$ we must have $j_{u,v}(x_0) = a_0 y_0 b_0$
with $$\max \Big\{ \|a_0\|_{L_u(\mathcal{N})},
\|y_0\|_{L_q(\mathcal{M})}, \|b_0\|_{L_v(\mathcal{N})} \Big\} <
1.$$ If we consider the invertible elements $$a = \big( a_0a_0^* +
\delta \mathrm{D}^{2/u} \big)^{1/2} \qquad \mbox{and} \qquad b =
\big( b_0^*b_0 + \delta \mathrm{D}^{2/v} \big)^{1/2},$$ we can
write $j_{u,v}(x_0) = a a^{-1} a_0 y_0 b_0 b^{-1} b = a y b$.
Moreover, for $\delta > 0$ small enough $$\max \Big\{
\|a\|_{L_u(\mathcal{N})}, \|y\|_{L_q(\mathcal{M})},
\|b\|_{L_v(\mathcal{N})} \Big\} < 1.$$ Since $a^2 \ge \delta
\mathrm{D}^{2/u}$ and $b^2 \ge \delta \mathrm{D}^{2/v}$, there
exist bounded elements $\alpha, \beta \in \mathcal{N}$ with
$$\mathrm{D}^{1/u} = a \alpha \qquad \mbox{and} \qquad
\mathrm{D}^{1/v} = \beta b.$$ Let us denote by $e$ the left
support of $\alpha$ and by $f$ the right support of $\beta$. We
note that the right support of $\alpha$ and the left support of
$\beta$ is 1. Then, we use polar decomposition to find strictly
positive densities $d_1 \in e \mathcal{N} e$ and $d_2 \in f
\mathcal{N} f$ and partial isometries $w_1$ and $w_2$ such that
$$\alpha = d_1 w_1 \qquad \mbox{and} \qquad \beta = w_2 d_2.$$
Note that $w_1^*w_1 = 1$, $w_1 w_1^* = e$, $w_2^* w_2 = f$ and
$w_2 w_2^* = 1$. Then we observe that $$ayb = j_{u,v}(x_0) =
\mathrm{D}^{1/u} x_0 \mathrm{D}^{1/v} = a \alpha x_0 \beta b \
\Rightarrow \ y = \alpha x_0 \beta.$$ This yields $$eyf = \alpha
x_0 \beta = d_1 w_1 x_0 w_2 d_2.$$ In particular, taking spectral
projections $e_n = 1_{[\frac{1}{n},n]}(d_1)$ and $f_n =
1_{[\frac{1}{n},n]}(d_2)$ $$w_1^* e_n d_1^{-1} eyf d_2^{-1} f_n
w_2^* = w_1^* e_n w_1 x_0 w_2 f_n w_2^*.$$ This implies that
\begin{eqnarray*}
\big| \varphi( w_1^* e_n w_1 x_0 w_2 f_n w_2^*) \big| & = & \big|
\Phi_{w_1^* e_n w_1 x_0 w_2 f_n w_2^*} \big(
\mathrm{D}^{\frac{1}{u}}, \mathrm{D}^{\frac{1}{v}} \big) \big| \\
& = & \big| \Phi_{w_1^* e_n d_1^{-1} eyf d_2^{-1} f_n w_2^*} \big(
\mathrm{D}^{\frac{1}{u}}, \mathrm{D}^{\frac{1}{v}} \big) \big| \\
& = & \big| \Phi_{w_1^* e_n d_1^{-1} eyf d_2^{-1} f_n w_2^*} \big(
a d_1 w_1, w_2 d_2 b \big) \big| \\ & = & \big| \Phi_{eyf} \big( a
e_n, f_n b \big) \big| \\ & \le & \|a e_n\|_u \|eyf\|_q \|f_n
b\|_v \\ & \le & \|a\|_u \|y\|_q \|b\|_v < 1.
\end{eqnarray*}
On the other hand, $(w_1^* e_n w_1)$ (resp. $(w_2 f_n w_2^*)$)
converges to $w_1^* w_1 = 1$ (resp. $w_2 w_2^* = 1$) strongly. By
Lemma 2.3 in \cite{J1}, this implies that $(\mathrm{D}^{1/u} w_1^*
e_n w_1)$ (resp. $(w_2 f_n w_2^* \mathrm{D}^{1/v})$) converges to
$\mathrm{D}^{1/u}$ (resp. $\mathrm{D}^{1/v}$) in the norm of
$L_u(\mathcal{N})$ (resp. $L_v(\mathcal{N})$). This combined with
the continuity of $\Phi_{x_0}$ gives
\begin{eqnarray*}
\big| \varphi(x_0) \big| & = & \lim_{n \to \infty} \big|
\Phi_{x_0} \big( \mathrm{D}^{\frac{1}{u}} w_1^* e_n w_1, w_2 f_n
w_2^* \mathrm{D}^{\frac{1}{v}} \big) \big| \\ & = & \lim_{n \to
\infty} \big| \varphi \big( w_1^* e_n w_1 x_0 w_2 f_n w_2^* \big)
\big| \le 1.
\end{eqnarray*}
This contradicts condition (a). Therefore, the map $j_{u,v}$
defines an isometry and the proof of i) is completed. Next, we see
that $\| \ \|_{u \cdot q \cdot v}$ is a norm \lq outside\rq${}$
the space $\mathcal{N}_u L_q(\mathcal{M}) \mathcal{N}_v$. The
homogeneity over $\R_+$ is clear and the positive definiteness
follows as in Lemma \ref{Lemma-Triangle-Inequality}. Thus, it
suffices to show that the triangle inequality holds. As we did
above, we begin by taking a family $x_1, x_2, \ldots , x_m$ of
elements in $L_u(\mathcal{N}) L_q(\mathcal{M}) L_v(\mathcal{N})$
satisfying $\|x_k\|_{u \cdot q \cdot v} < 1$ and a collection of
positive numbers $\lambda_1, \lambda_2, \ldots , \lambda_m$ with
$\summ_k \lambda_k = 1$. By hypothesis, we may assume that $x_k =
a_k y_k b_k$ with $$\max \Big\{ \|a_k\|_{L_u(\mathcal{N})},
\|y_k\|_{L_q(\mathcal{M})}, \|b_k\|_{L_v(\mathcal{N})} \Big\} <
1.$$ Given any $\xi > 1$ and by the density of
$\mathrm{D}^{\frac{1}{u}} \mathcal{N}$ (resp. $\mathcal{N}
\mathrm{D}^{\frac{1}{v}}$) in $L_u(\mathcal{N})$ (resp.
$L_v(\mathcal{N})$), it is not difficult to check that both $a_k$
and $b_k$ can be written in the following way
\begin{eqnarray*}
a_k & = & \sum_{i=0}^{\infty} \mathrm{D}^{\frac{1}{u}} a_{ik}
\qquad \mbox{with} \qquad \big\| \mathrm{D}^{\frac{1}{u}} a_{ik}
\big\|_{L_u(\mathcal{N})} \le \xi^{-i}, \\ b_k & = &
\sum_{j=0}^{\infty} b_{jk} \mathrm{D}^{\frac{1}{v}} \qquad
\mbox{with} \qquad \big\| b_{jk} \mathrm{D}^{\frac{1}{v}}
\big\|_{L_v(\mathcal{N})} \le \xi^{-j}.
\end{eqnarray*}
This gives rise to $$\sum_{k=1}^m \lambda_k x_k =
\sum_{i,j=0}^{\infty} \mathrm{D}^{\frac{1}{u}} \Big( \sum_{k=1}^m
\lambda_k a_{ik} y_k b_{jk} \Big) \mathrm{D}^{\frac{1}{v}}.$$
Assuming the triangle inequality on $\mathcal{N}_u
L_q(\mathcal{M}) \mathcal{N}_v$, we can write
$$\mathrm{D}^{\frac{1}{u}} \Big( \sum_{k=1}^m \lambda_k a_{ik} y_k
b_{jk} \Big) \mathrm{D}^{\frac{1}{v}} = \mathsf{A}_{ij}
\mathsf{Y}_{ij} \mathsf{B}_{ij}$$ where
$$\|\mathsf{A}_{ij}\|_{L_u(\mathcal{N})} \le \xi^{-\frac{i+j}{4}},
\quad \|\mathsf{Y}_{ij}\|_{L_q(\mathcal{M})} \le
\xi^{-\frac{i+j}{2}}, \quad \|\mathsf{B}_{ij}\|_{L_v(\mathcal{N})}
\le \xi^{-\frac{i+j}{4}}.$$ Then, we define for $\delta
> 0$ $$\mathcal{S}_r(\mathsf{A}) = \Big(
\sum_{i,j=0}^{\infty} \mathsf{A}_{ij} \mathsf{A}_{ij}^* + \delta
\mathrm{D}^{\frac{2}{u}} \Big)^{1/2} \quad \mbox{and} \quad
\mathcal{S}_c(\mathsf{B}) = \Big( \sum_{i,j=0}^{\infty}
\mathsf{B}_{ij}^* \mathsf{B}_{ij} + \delta
\mathrm{D}^{\frac{2}{v}} \Big)^{1/2},$$ so that there exist
$\alpha_{ij}, \beta_{ij} \in \mathcal{N}$ given by
$\mathcal{S}_r(\mathsf{A}) \alpha_{ij} = \mathsf{A}_{ij}$ and
$\beta_{ij} \mathcal{S}_c(\mathsf{B}) = \mathsf{B}_{ij}$. Hence,
$$\sum_{k=1}^m \lambda_k x_k  = \mathcal{S}_r(\mathsf{A}) \Big(
\sum_{i,j=0}^{\infty} \alpha_{ij} \mathsf{Y}_{ij} \beta_{ij} \Big)
\mathcal{S}_c(\mathsf{B}).$$ Moreover, we have $$\big\|
\mathcal{S}_r(\mathsf{A}) \big\|_{L_u(\mathcal{N})} \le \Big(
\delta + \sum_{i,j=0}^{\infty}
\|\mathsf{A}_{ij}\|_{L_u(\mathcal{N})}^2 \Big)^{1/2} \le
\sqrt{\delta} + \frac{1}{1 - 1/\sqrt{\xi}}.$$ The same estimate
applies to $\mathcal{S}_c(\mathsf{B})$. The middle term satisfies
\begin{eqnarray*}
\Big\| \sum_{i,j=0}^{\infty} \alpha_{ij} \mathsf{Y}_{ij}
\beta_{ij} \Big\|_q & \le & \Big\| \Big( \sum_{i,j=0}^{\infty}
\alpha_{ij} \alpha_{ij}^* \Big)^{\frac{1}{2}} \Big\|_{\infty}
\Big( \sum_{i,j=0}^{\infty} \|\mathsf{Y}_{ij}\|_q^q \
\Big)^{\frac{1}{q}} \, \Big\| \Big( \sum_{i,j=0}^{\infty}
\beta_{ij}^* \beta_{ij} \Big)^{\frac{1}{2}} \Big\|_{\infty} \\ &
\le & \Big( \sum_{i,j=0}^{\infty} \|\mathsf{Y}_{ij}\|_q^q \
\Big)^{\frac{1}{q}} \le \Big( \frac{1}{1 - 1/\xi^{q/2}}
\Big)^{2/q}
\end{eqnarray*}
In summary, $$\Big\| \sum_{k=1}^m \lambda_k x_k \Big\|_{u \cdot q
\cdot v} \le \big\| \mathcal{S}_r(\mathsf{A})
\big\|_{L_u(\mathcal{N})} \Big\| \sum_{i,j=0}^{\infty} \alpha_{ij}
\mathsf{Y}_{ij} \beta_{ij} \Big\|_{L_q(\mathcal{M})} \big\|
\mathcal{S}_c(\mathsf{B}) \big\|_{L_v(\mathcal{N})}
\longrightarrow 1$$ as $\delta \to 0$ and $\xi \to \infty$. This
proves the triangle inequality. To prove completeness we use again
a geometric series argument. Let $x_1, x_2, \ldots $ be a
countable family of elements in $L_u(\mathcal{N}) L_q(\mathcal{M})
L_v(\mathcal{N})$ with $\|x_k\|_{u \cdot q \cdot v} < 4^{-k}$ for
any integer $k \ge 1$. That is, we have $x_k = 4^{-k} a_k y_k b_k$
with $$\max \Big\{ \|a_k\|_{L_u(\mathcal{N})},
\|y_k\|_{L_q(\mathcal{M})}, \|b_k\|_{L_v(\mathcal{N})} \Big\} <
1.$$ Then, by a well known characterization of completeness, it
suffices to see that the sum $\summ_k x_k$ belongs to
$L_u(\mathcal{N}) L_q(\mathcal{M}) L_v(\mathcal{N})$. We use again
the same factorization trick $$\sum_{k=1}^{\infty} x_k =
\mathcal{S}_r(a) \Big( \sum_{k=1}^{\infty} 2^{-k} \alpha_k y_k
\beta_k \Big) \mathcal{S}_c(b)$$ where $$\begin{array}{lclclcl}
\mathcal{S}_r(a) & = & \displaystyle \Big( \sum_{k=1}^{\infty}
2^{-k} a_ka_k^* + \delta \mathrm{D}^{\frac{2}{u}} \Big)^{1/2} &
\quad \mbox{and} \quad & \mathcal{S}_r(a) \alpha_k & = & 2^{-k/2}
a_k, \\ \mathcal{S}_c(b) & = & \displaystyle \Big(
\sum_{k=1}^{\infty} 2^{-k} b_k^*b_k \, + \delta
\mathrm{D}^{\frac{2}{v}} \Big)^{1/2} & \quad \mbox{and} \quad &
\beta_k \mathcal{S}_c(b) & = & 2^{-k/2} b_k.
\end{array}$$ In particular, we just need to show that
$\mathcal{S}_r(a) \in L_u(\mathcal{N})$, $\mathcal{S}_c(b) \in
L_v(\mathcal{N})$ and the middle term belongs to
$L_q(\mathcal{M})$. However, this follows again by applying the
same estimates as above, details are left to the reader. To
conclude, we just need to show that $\mathcal{N}_u
L_q(\mathcal{M}) \mathcal{N}_v$ is dense in $L_u(\mathcal{N})
L_q(\mathcal{M}) L_v(\mathcal{N})$. However, this follows easily
from the density of $\mathrm{D}^{1/u} \mathcal{N}$ (resp.
$\mathcal{N} \mathrm{D}^{1/v}$) in $L_u(\mathcal{N})$ (resp.
$L_v(\mathcal{N})$) and the triangle inequality proved above.
\end{proof}

\section{A metric structure on the solid $\mathsf{K}$}

At this point, we are not able to prove that $\mathcal{N}_u
L_q(\mathcal{M}) \mathcal{N}_v$ is a normed space for any point
$(1/u,1/v,1/q)$ in the solid $\mathsf{K}$. However, we need at
least to know that we have a metric space structure. In fact, we
shall prove that $||| \ |||_{u \cdot q \cdot v}$ is always a
$\gamma$-norm for some $0 < \gamma \le 1$.

\begin{lemma}
If $(1/u,1/v,1/q) \in \mathsf{K}$, there exists $0 < \gamma \le 1$
such that $$|||x_1 + x_2|||_{u \cdot q \cdot v}^{\gamma} \le
|||x_1|||_{u \cdot q \cdot v}^{\gamma} + |||x_2|||_{u \cdot q
\cdot v}^{\gamma} \quad \mbox{for all} \quad x_1,x_2 \in
\mathcal{N}_u L_q(\mathcal{M}) \mathcal{N}_v.$$
\end{lemma}

\begin{proof} According to Lemma \ref{Lemma-Triangle-Inequality}, we can
assume in what follows that $$(1/u,1/v,1/q) \in \mathsf{K}
\setminus \partial_{\infty} \mathsf{K}.$$ As we did to prove the
triangle inequality, let $x_1, x_2, \ldots, x_m \in \mathcal{N}_u
L_q(\mathcal{M}) \mathcal{N}_v$ be a sequence of vectors
satisfying $|||x_k|||_{u \cdot q \cdot v} < 1$ and let $\lambda_1,
\lambda_2, \ldots, \lambda_m \in \R_+$ with sum $\summ_k \lambda_k
= 1$. Then it suffices to show that $$\Big|\Big|\Big| \sum_{k=1}^m
\lambda_k x_k \Big|\Big|\Big|_{u \cdot q \cdot v}^{\gamma} \le
\sum_{k=1}^m \lambda_k^{\gamma} \quad \mbox{for some} \quad 0 <
\gamma \le 1.$$ By hypothesis, we may assume that $x_k = \alpha_k
y_k \beta_k,$ with $$\max \Big\{ \big\| \mathrm{D}^{\frac{1}{u}}
\alpha_k \big\|_{L_u(\mathcal{N})}, \|y_k\|_{L_q(\mathcal{M})},
\big\| \beta_k \mathrm{D}^{\frac{1}{v}} \big\|_{L_v(\mathcal{N})}
\Big\} < 1.$$ By Figure I we can always find $1 \le u_1,q_0,v_1
\le \infty$ and $0 < \theta < 1$ such that
\begin{itemize}
\item[(a)] We have $1/q_0 = 1/u + 1/q + 1/v = 1/u_1 + 1/v_1$.
\item[(b)] We have $\big( 1/u, 1/v, 1/q \big) = \big( 0, 0,
(1-\theta)/q_0 \big) + \big( \theta/u_1, \theta/v_1, 0 \big)$.
\end{itemize}
Indeed, we first consider the plane $\mathsf{P}$ parallel to
$\mathsf{ACDF}$ containing $(1/u,1/v,1/q)$. The point
$(0,0,1/q_0)$ is the intersection of $\mathsf{P}$ with the segment
$0\mathsf{A}$. Then, we consider the line $\mathsf{L}$ passing
through $(0,0,1/q_0)$ and $(1/u,1/v,1/q)$. Then, the point
$(1/u_1,1/v_1,0)$ is the intersection of $\mathsf{L}$ with the
coordinate plane $z=0$. Note that the point $(1/u_1,1/v_1,0)$ is
not necessarily in $\mathsf{K}$. However, according to (a) it
always satisfies $1/u_1+1/v_1 \le 1$. Note also that, since we
have excluded the points in $\partial_{\infty} \mathsf{K}$ at the
beginning of this proof, we can always assume that $0 < \theta <
1$. Then it follows from (b) that
\begin{eqnarray*}
L_u(\mathcal{N}) & = & \big[ L_{\infty}(\mathcal{N}) \, , \,
L_{u_1}(\mathcal{N}) \big]_{\theta}, \\ L_q(\mathcal{M}) & = &
\big[ L_{q_0}(\mathcal{M}), L_{\infty}(\mathcal{M})
\big]_{\theta}, \\ L_v(\mathcal{N}) & = & \big[
L_{\infty}(\mathcal{N}) \, , \, L_{v_1}(\mathcal{N})
\big]_{\theta}.
\end{eqnarray*}
By the complex interpolation method, we can find bounded analytic
functions $$f_{2k}: \mathcal{S} \rightarrow L_{q_0}(\mathcal{M}) +
L_{\infty}(\mathcal{M})$$ satisfying $f_{2k}(\theta) = y_k$ for
$k=1,2,\ldots,m$ and such that
\begin{equation} \label{Equation-Boundary-1}
\max \Big\{ \sup_{z \in \partial_0} \big\| f_{2k}(z)
\big\|_{L_{q_0}(\mathcal{M})}, \sup_{z \in \partial_1} \big\|
f_{2k}(z) \big\|_{L_{\infty}(\mathcal{M})} \Big\} < 1.
\end{equation}
Similarly, by Lemma \ref{Lemma-Kosaki} we also have bounded
analytic functions $$f_{1k}, f_{3k}: \mathcal{S} \to \mathcal{N}$$
for any $k=1,2, \ldots,m$ satisfying $(f_{1k}(\theta),
f_{3k}(\theta)) = (\alpha_k, \beta_k)$ and
\begin{equation} \label{Equation-Boundary-2}
\begin{array}{rcr}
\displaystyle \sup_{z \in \partial_0} \ \max \Big\{ \big\|
f_{1k}(z) \big\|_{L_{\infty}(\mathcal{N})}, \big\| f_{3k}(z)
\big\|_{L_{\infty}(\mathcal{N})} \Big\} < 1, \\ \displaystyle
\sup_{z \in \partial_1} \ \max \Big\{ \big\|
\mathrm{D}^{\frac{1}{u_1}} f_{1k}(z)
\big\|_{L_{u_1}(\mathcal{N})}, \big\| f_{3k}(z)
\mathrm{D}^{\frac{1}{v_1}} \big\|_{L_{v_1}(\mathcal{N})} \Big\} <
1.
\end{array}
\end{equation}
Given $\delta > 0$, we consider the following functions on the
boundary $$\mathrm{W}_1(z) = \left\{
\begin{array}{ll} 1 & \mbox{if} \ z \in \partial_0,
\\ \delta 1 + \summ_k \lambda_k f_{1k}(z) f_{1k}(z)^* &
\mbox{if} \ z \in \partial_1,
\end{array} \right.$$ $$\mathrm{W}_3(z) = \left\{
\begin{array}{ll} 1 & \mbox{if} \ z \in \partial_0,
\\ \delta 1 + \summ_k \lambda_k f_{3k}(z)^* f_{3k}(z) &
\mbox{if} \ z \in \partial_1.
\end{array} \right.$$ According to Devinatz's factorization
theorem \cite{P0}, there exist bounded analytic functions
$\mathrm{w}_1, \mathrm{w}_3: \mathcal{S} \rightarrow \mathcal{N}$
with bounded analytic inverse and satisfying the following
identities on $\partial \mathcal{S}$ $$\mathrm{w}_1(z)
\mathrm{w}_1(z)^* = \mathrm{W}_1(z),$$ $$\mathrm{w}_3(z)^*
\mathrm{w}_3(z) = \mathrm{W}_3(z).$$ Then we can write
$$\sum_{k=1}^m \lambda_k x_k = \mathrm{w}_1(\theta) \Big[
\mathrm{w}_1^{-1}(\theta) \Big( \sum_{k=1}^m \lambda_k
f_{1k}(\theta) f_{2k}(\theta) f_{3k}(\theta) \Big)
\mathrm{w}_3^{-1}(\theta) \Big] \mathrm{w}_3(\theta).$$ Let us
estimate the norms of the three factors above. First we clearly
have
\begin{eqnarray*}
\sup_{z \in \partial_0} \big\| \mathrm{w}_1(z)
\big\|_{L_{\infty}(\mathcal{N})} & = & \sup_{z \in \partial_0}
\big\| \mathrm{w}_1(z) \mathrm{w}_1(z)^*
\big\|_{L_{\infty}(\mathcal{N})}^{1/2} = 1, \\ \sup_{z \in
\partial_0} \big\| \mathrm{w}_3(z)
\big\|_{L_{\infty}(\mathcal{N})} & = & \sup_{z \in \partial_0}
\big\| \mathrm{w}_3(z)^* \mathrm{w}_3(z)
\big\|_{L_{\infty}(\mathcal{N})}^{1/2} = 1.
\end{eqnarray*}
On the other hand, since $L_p(\mathcal{N})$ is always a
$\min(p,1)$-normed space
\begin{eqnarray*}
\lefteqn{\sup_{z \in \partial_1} \big\| \mathrm{D}^{\frac{1}{u_1}}
\mathrm{w}_1(z) \big\|_{L_{u_1}(\mathcal{N})}^{u_1} = \sup_{z \in
\partial_1} \big\| \mathrm{D}^{\frac{1}{u_1}} \mathrm{w}_1(z)
\mathrm{w}_1(z)^* \mathrm{D}^{\frac{1}{u_1}}
\big\|_{L_{u_1/2}(\mathcal{N})}^{u_1/2}} \\ & \le & \sup_{z \in
\partial_1} \left\{ \begin{array}{ll} \Big(
\sum_{k=1}^m \lambda_k \big\| \mathrm{D}^{\frac{1}{u_1}} f_{1k}(z)
\big\|_{L_{u_1}(\mathcal{N})}^2 + \delta \Big)^{u_1/2}, &
\mbox{if} \ u_1 \ge 2, \\ \sum_{k=1}^m \lambda_k^{u_1/2} \big\|
\mathrm{D}^{\frac{1}{u_1}} f_{1k}(z)
\big\|_{L_{u_1}(\mathcal{N})}^{u_1} + \delta^{u_1/2}, & \mbox{if}
\ u_1 < 2. \end{array} \right. \\ & < & \max \left\{
(1+\delta)^{u_1/2}, \delta^{u_1/2} + \sum_{k=1}^m
\lambda_k^{u_1/2} \right\}.
\end{eqnarray*}
Similarly, we have $$\sup_{z \in \partial_1} \big\|
\mathrm{w}_3(z) \mathrm{D}^{\frac{1}{v_1}}
\big\|_{L_{v_1}(\mathcal{N})}^{v_1} < \max \left\{
(1+\delta)^{v_1/2}, \delta^{v_1/2} + \sum_{k=1}^m
\lambda_k^{v_1/2} \right\}.$$ In summary, we obtain by Kosaki's
interpolation
\begin{eqnarray*}
\big\| \mathrm{D}^{\frac{1}{u}} \mathrm{w}_1(\theta)
\big\|_{L_u(\mathcal{N})} & < & \max \left\{ \sqrt{1+\delta},
\Big( \delta^{u_1/2} + \sum_{k=1}^m \lambda_k^{u_1/2}
\Big)^{1/u_1} \right\},
\\ \big\| \mathrm{w}_3(\theta) \mathrm{D}^{\frac{1}{v}}
\big\|_{L_v(\mathcal{N})} & < & \max \left\{ \sqrt{1+\delta},
\Big( \delta^{v_1/2} + \sum_{k=1}^m \lambda_k^{v_1/2}
\Big)^{1/v_1} \right\}.
\end{eqnarray*}
As we have pointed out above, we have $1/u_1 + 1/v_1 \le 1$. In
particular, $u_1/2$ and $v_1/2$ can not be simultaneously less
than $1$. Thus, at least one of the two sums above must be less or
equal than $1$. This means that taking $$\gamma = \min \Big\{ 1,
u_1/2, v_1/2 \Big\},$$ we deduce $$\lim_{\delta \to 0} \Big(
\big\| \mathrm{D}^{\frac{1}{u}} \mathrm{w}_1(\theta)
\big\|_{L_u(\mathcal{N})} \, \big\| \mathrm{w}_3(\theta)
\mathrm{D}^{\frac{1}{v}} \big\|_{L_v(\mathcal{N})} \Big) < \Big(
\sum_{k=1}^m \lambda_k^{\gamma} \Big)^{1/2\gamma} \le \Big(
\sum_{k=1}^m \lambda_k^{\gamma} \Big)^{1/\gamma}.$$ In particular,
it suffices to see that $$\Big\| \mathrm{w}_1^{-1}(\theta) \Big(
\sum_{k=1}^m \lambda_k f_{1k}(\theta) f_{2k}(\theta)
f_{3k}(\theta) \Big) \mathrm{w}_3^{-1}(\theta)
\Big\|_{L_q(\mathcal{M})} \le 1.$$ Let us consider the bounded
analytic function $$g(z) = \mathrm{w}_1^{-1}(z) \Big( \sum_{k=1}^m
\lambda_k f_{1k}(z) f_{2k}(z) f_{3k}(z) \Big)
\mathrm{w}_3^{-1}(z).$$ Since $\mathrm{w}_1(z) \mathrm{w}_1(z)^* =
1 = \mathrm{w}_3(z)^* \mathrm{w}_3(z)$ for all $z \in \partial_0$,
we clearly have
\begin{eqnarray*}
\sup_{z \in \partial_0} \big\| g(z) \big\|_{L_{q_0}(\mathcal{M})}
& = & \sup_{z \in \partial_0} \Big\| \sum_{k=1}^m \lambda_k
f_{1k}(z) f_{2k}(z) f_{3k}(z) \Big\|_{L_{q_0}(\mathcal{M})} \\ &
\le & \sup_{z \in \partial_0} \sum_{k=1}^m \lambda_k \big\|
f_{1k}(z) \big\|_{L_{\infty}(\mathcal{N})} \big\| f_{2k}(z)
\big\|_{L_{q_0}(\mathcal{M})} \big\| f_{3k}(z)
\big\|_{L_{\infty}(\mathcal{N})}.
\end{eqnarray*}
Then it follows from (\ref{Equation-Boundary-1}) and
(\ref{Equation-Boundary-2}) that the expression above is bounded
above by $1$. On the other hand, since $\mathrm{w}_1$ and
$\mathrm{w}_3$ are invertible, we can define the functions
$h_{1k}, h_{3k}: \mathcal{S} \to \mathcal{N}$ by the relations
$$\sqrt{\lambda_k} f_{1k}(z) = \mathrm{w}_1(z) h_{1k}(z) \quad
\mbox{and} \quad \sqrt{\lambda_k} f_{3k}(z) = h_{3k}(z)
\mathrm{w}_3(z).$$ Then it is easy to check that for any $z \in
\partial_1$
\begin{equation} \label{Equation-h1k-h3k}
\sum_{k=1}^m h_{1k}(z) h_{1k}(z)^* \le 1 \quad \mbox{and} \quad
\sum_{k=1}^m h_{3k}(z)^* h_{3k}(z) \le 1.
\end{equation}
Therefore, we obtain from (\ref{Equation-Boundary-1}) and
(\ref{Equation-h1k-h3k}) that
\begin{eqnarray*}
\sup_{z \in \partial_1} \big\| g(z)
\big\|_{L_{\infty}(\mathcal{M})} & = & \sup_{z \in
\partial_1} \Big\| \sum_{k=1}^m h_{1k}(z) f_{2k}(z) h_{3k}(z)
\Big\|_{L_{\infty}(\mathcal{M})} \\ & \le & \sup_{z \in
\partial_1} \Big\| \Big( \sum_{k=1}^m h_{1k}(z) h_{1k}(z)^*
\Big)^{1/2} \Big\|_{L_\infty(\mathcal{N})} \\ & \times & \sup_{z
\in \partial_1} \max_{1 \le k \le m} \big\| f_{2k}(z)
\big\|_{L_{\infty}(\mathcal{M})} \\ & \times & \sup_{z \in
\partial_1} \Big\| \Big( \sum_{k=1}^m h_{3k}(z)^* h_{3k}(z)
\Big)^{1/2} \Big\|_{L_\infty(\mathcal{N})} < 1.
\end{eqnarray*}
By Kosaki's interpolation we have $\big\|g(\theta)
\big\|_{L_q(\mathcal{M})} \le 1$. This completes the proof.
\end{proof}

Arguing as in the proof of Proposition
\ref{Proposition-Isometry-Triangle-Complete} ii), we can see that
any amalgamated space $L_u(\mathcal{N}) L_q(\mathcal{M})
L_v(\mathcal{N})$ is a $\gamma$-Banach space and that $j_{u,v}
\big( \mathcal{N}_u L_q(\mathcal{M}) \mathcal{N}_v \big)$ (with
$j_{u,v}$ as defined in Proposition
\ref{Proposition-Isometry-Triangle-Complete}) is a dense subspace.
Indeed, we just need to repeat the two given arguments (using
geometric series) conveniently rewritten with exponents $\gamma$
everywhere. We leave the details to the reader. Let us state this
result for future reference.

\begin{proposition} \label{Proposition-gamma-Norm}
If $(1/u,1/v,1/q) \in \mathsf{K}$, there exists $0 < \gamma \le 1$
such that
\begin{itemize}
\item[i)] $L_u(\mathcal{N}) L_q(\mathcal{M}) L_v(\mathcal{N})$ is
a $\gamma$-Banach space. \item[ii)] $j_{u,v} \big( \mathcal{N}_u
L_q(\mathcal{M}) \mathcal{N}_v \big)$ is a dense subspace of
$L_u(\mathcal{N}) L_q(\mathcal{M}) L_v(\mathcal{N})$.
\end{itemize}
\end{proposition}

\begin{observation} \label{Observation-Lack-gamma-Norm}
\emph{We do not claim at this moment that $j_{u,v}$ is an
isometry.}
\end{observation}

\chapter{An interpolation theorem}
\label{Section3}

In this chapter we prove that the solid $\mathsf{K}$ is an
interpolation family on the indices $(u,q,v)$. Of course, we need
to know a priori that $L_u(\mathcal{N}) L_q(\mathcal{M})
L_v(\mathcal{N})$ is a Banach space for any $(1/u,1/v,1/q)$ in the
solid $\mathsf{K}$. In fact, as we shall see the proofs of both
results depend on the proof of the other. Indeed, let us consider
a parameter $0 \le \tau \le 1$. According to the terminology
employed in Figure I, let us define $\mathsf{K}_{\tau}$ to be the
intersection of $\mathsf{K}$ with the plane $\mathsf{P}_{\tau}$
which contains the point $(0,0,\tau)$ and is parallel to the upper
face $\mathsf{ACDF}$ of $\mathsf{K}$. Roughly speaking, we first
prove that $$\mathcal{K}_{\tau} = \Big\{ L_u(\mathcal{N})
L_q(\mathcal{M}) L_v(\mathcal{N}) \, \big| \ (1/u,1/v,1/q) \in
\mathsf{K}_{\tau} \Big\}$$ is an interpolation family on the
indices $(u,q,v)$ for any $0 \le \tau \le 1$ and with ending
points lying on $\partial_{\infty} \mathsf{K}$. Note that we
already proved in Chapter \ref{Section2} that the spaces
associated to the points in $\partial_{\infty} \mathsf{K}$ are
Banach spaces. Then, we use this result to prove that every point
in $\mathsf{K}$ corresponds to a Banach space and derive our main
interpolation theorem. More concretely, we first prove the
following result.

\begin{lemma} \label{Lemma0}
Let us assume that
\begin{itemize}
\item[i)] $(1/u_j,1/v_j,1/q_j) \in \partial_{\infty} \mathsf{K}$
for $j=0,1$. \vskip4pt \item[ii)] $1/u_0 + 1/q_0 + 1/v_0 = 1/u_1 +
1/q_1 + 1/v_1$.
\end{itemize} \vskip3pt
Then $L_{u_{\theta}}(\mathcal{N}) L_{q_{\theta}}(\mathcal{M})
L_{v_{\theta}}(\mathcal{N})$ is a Banach space isometrically
isomorphic to $$\mathrm{X}_{\theta}(\mathcal{M}) = \Big[
L_{u_0}(\mathcal{N}) L_{q_0}(\mathcal{M}) L_{v_0}(\mathcal{N}),
L_{u_1}(\mathcal{N}) L_{q_1}(\mathcal{M}) L_{v_1}(\mathcal{N})
\Big]_{\theta}^{\null}.$$
\end{lemma}

We know from Lemma \ref{Lemma-Triangle-Inequality} and Proposition
\ref{Proposition-Isometry-Triangle-Complete} that the
interpolation pairs considered in Lemma \ref{Lemma0} are made of
Banach spaces. After the proof of Lemma \ref{Lemma0} we shall show
that $L_u(\mathcal{N}) L_q(\mathcal{M}) L_v(\mathcal{N})$ is a
Banach space for any $(1/u,1/v,1/q)$ in the solid $\mathsf{K}$ and
we shall deduce our main result in this chapter.

\begin{theorem} \label{TheoremA1}
The amalgamated space $L_u(\mathcal{N}) L_q(\mathcal{M})
L_v(\mathcal{N})$ is a Banach space for any $(1/u,1/v,1/q) \in
\mathsf{K}$. Moreover, if $(1/u_j,1/v_j,1/q_j) \in \mathsf{K}$ for
$j=0,1$, the space $L_{u_{\theta}}(\mathcal{N})
L_{q_{\theta}}(\mathcal{M}) L_{v_{\theta}}(\mathcal{N})$ is
isometrically isomorphic to $$\mathrm{X}_{\theta}(\mathcal{M}) =
\Big[L_{u_0}(\mathcal{N}) L_{q_0}(\mathcal{M})
L_{v_0}(\mathcal{N}), L_{u_1}(\mathcal{N}) L_{q_1}(\mathcal{M})
L_{v_1}(\mathcal{N}) \Big]_{\theta}^{\null}.$$
\end{theorem}

Note that the pairs of Banach spaces considered in Lemma
\ref{Lemma0} and Theorem \ref{TheoremA1} are compatible for
complex interpolation. Indeed, since the amalgamated space
$L_{u_j}(\mathcal{N}) L_{q_j}(\mathcal{M}) L_{v_j}(\mathcal{N})$
is continuously injected (by means of H\"{o}lder inequality) in
$L_{p_j}(\mathcal{M})$ for $j=0,1$ and $1/p_j = 1/u_j + 1/q_j +
1/v_j$, any such pair lives in the sum $L_{p_0}(\mathcal{M}) +
L_{p_1}(\mathcal{M})$. The main difficulty which appears to prove
Lemma \ref{Lemma0} lies in the fact that the intersection of
amalgamated spaces is quite difficult to describe in the general
case. Therefore, any attempt to work on a dense subspace meets
this obstacle. We shall avoid this difficulty by proving this
interpolation result in the particular case of finite von Neumann
algebras whose density is bounded above and below. Under this
assumption, we will be able to find a nice dense subspace to work
with. Then we shall use the Haagerup construction \cite{H2}
sketched in Section \ref{Subsection-Preliminaries} (suitably
modified to work in the present context) to get the result in the
general case. After that, the proof of Theorem \ref{TheoremA1}
follows easily by using similar techniques.

\section{Finite von Neumann algebras}

We begin by proving Lemma \ref{Lemma0} for a finite von Neumann
algebra $\mathcal{M}$ equipped with a $\emph{n.f.}$ state
$\varphi$ with respect to which, the corresponding density
$\mathrm{D}$ satisfies the following property for some positive
constants $0 < c_1 \le c_2 < \infty$
\begin{equation} \label{Equation-Boundedness-Density}
c_1 1 \le \mathrm{D} \le c_2 1.
\end{equation}

\begin{lemma} \label{Lemma-Interpolation-Density}
Assume that $(1/u_j,1/v_j,1/q_j) \in \partial_{\infty} \mathsf{K}$
for $j=0,1$. Then, given any $0 < \theta < 1$, the space
$L_{q_0}(\mathcal{M}) \cap L_{q_1}(\mathcal{M})$ is dense in
$$\mathrm{X}_{\theta}(\mathcal{M}) = \Big[ L_{u_0}(\mathcal{N})
L_{q_0}(\mathcal{M}) L_{v_0}(\mathcal{N}), L_{u_1}(\mathcal{N})
L_{q_1}(\mathcal{M}) L_{v_1}(\mathcal{N})
\Big]_{\theta}^{\null}.$$
\end{lemma}

\begin{proof} Let us write $\Delta$ for the intersection
$$L_{u_0}(\mathcal{N}) L_{q_0}(\mathcal{M}) L_{v_0}(\mathcal{N})
\cap L_{u_1}(\mathcal{N}) L_{q_1}(\mathcal{M})
L_{v_1}(\mathcal{N}).$$ According to the complex interpolation
method, we know that $\Delta$ is dense in
$\mathrm{X}_{\theta}(\mathcal{M})$. In particular, it suffices to
show that we have density in $\Delta$ with respect to the norm of
$\mathrm{X}_{\theta}(\mathcal{M})$. Let $x$ be an element in
$\Delta$ so that we have decompositions $$\begin{array}{rcl} x & =
& a_0 \tilde{y}_0 b_0 \\ x & = & a_1 \tilde{y}_1 b_1 \end{array}
\qquad \mbox{where} \qquad a_j \in L_{u_j}(\mathcal{N}), \
\tilde{y}_j \in L_{q_j}(\mathcal{M}), \ b_j \in
L_{v_j}(\mathcal{N}).$$

\noindent \textsc{Case I.} Let us first assume that $\max \big\{
u_0,u_1,v_0,v_1 \big\} < \infty$. Then, we can define
\begin{eqnarray*}
a & = & \sum_{j=0,1} \big( a_ja_j^* \big)^{u_j/2} + \delta
\mathrm{D},
\\ b & = & \sum_{j=0,1} \big( b_j^* \, b_j \big)^{v_j/2} \hskip1pt
+ \delta \mathrm{D}.
\end{eqnarray*}
Note that, since $2/u_j$ and $2/v_j$ are in $(0,1]$, we have
$$\left. \begin{array}{lcl} \big( a_ja_j^* \big)^{u_j/2} & \le & a
\\ \big( b_j^* \, b_j \big)^{v_j/2} & \le & b \end{array} \right\}
\Rightarrow \left\{ \begin{array}{rcl} a_ja_j^* & \le & a^{2/u_j},
\\ b_j^* \, b_j & \le & b^{2/v_j}.
\end{array} \right.$$ Therefore, if we define $\alpha_j,
\beta_j$ by $a_j = a^{1/u_j} \alpha_j$ and $b_j = \beta_j
b^{1/v_j}$, we have $\alpha_j \alpha_j^* \le 1$ and $\beta_j^*
\beta_j \le 1$. Thus, $\alpha_j, \beta_j \in \mathcal{N}$ for
$j=0,1$ and we can set $y_j = \alpha_j \tilde{y}_j \beta_j$ in
$L_{q_j}(\mathcal{M})$ so that
\begin{equation} \label{Equation-Factorization}
x = a^{1/u_j} y_j b^{1/v_j}
\end{equation}
for $j=0,1$. Now we consider the spectral projections $$e_n =
1_{[\frac{1}{n},n]}(a) \quad \mbox{and} \quad f_n =
1_{[\frac{1}{n},n]}(b).$$ Then, since $e_n x f_n = (e_n a^{1/u_j})
y_j (b^{1/v_j} f_n)$ we find $e_n x f_n \in L_{q_0}(\mathcal{M})
\cap L_{q_1}(\mathcal{M})$. Thus, we have to show that $e_n x f_n$
tends to $x$ in the norm of $\mathrm{X}_{\theta}(\mathcal{M})$.
Let us write $$x - e_n x f_n = (1-e_n)x + e_n x (1-f_n).$$ Then it
is clear that
\begin{eqnarray*}
\big\| (1-e_n)x \big\|_{\Delta} & = & \max \Big\{ \big\| (1-e_n)
a^{\frac{1}{u_j}} y_j b^{\frac{1}{v_j}} \big\|_{u_j \cdot q_j
\cdot v_j} \, \big| \ j=0,1 \Big\} \\ & \le & \max \Big\{ \big\|
(1-e_n) a^{\frac{1}{u_j}} \big\|_{u_j} \|y_j\|_{q_j} \big\|
b^{\frac{1}{v_j}} \big\|_{v_j} \, \big| \ j=0,1 \Big\}.
\end{eqnarray*}
The term on the right tends to $0$ as $n \to \infty$ since
$\max(u_0,u_1) < \infty$. Therefore, since the norm of
$\mathrm{X}_{\theta}(\mathcal{M})$ is controlled by the norm of
$\Delta$ when $0 < \theta < 1$, we obtain that $(1-e_n)x$ tends to
$0$ in the norm of $\mathrm{X}_{\theta}(\mathcal{M})$. The second
term $e_n x (1-f_n)$ can be estimated in the same way. This
concludes the proof of Case I.

\vskip3pt

\noindent \textsc{Case II.} Now let us assume that there exists
some indices among $u_0,u_1,v_0,v_1$ which are infinite. Then we
can assume without lost of generality that $a_j$ (resp. $b_j$) is
$1$ whenever $u_j$ (resp. $v_j$) is infinite. According to this
and assuming the convention $1^{\infty} = 1$, the previous
definition of $a$ and $b$ still makes sense. Moreover, the
relations obtained in (\ref{Equation-Factorization}) also hold in
this case. Therefore, we need to prove again that $e_n x f_n$
tends to $x$ with respect to the norm of
$\mathrm{X}_{\theta}(\mathcal{M})$. We only prove the convergence
for the term $(1 - e_n) x$ since again the second one can be
estimated in the same way. By the three lines lemma, we have
$$\big\|(1-e_n)x \big\|_{\mathrm{X}_{\theta}(\mathcal{M})} \le
\big\|(1-e_n)x \big\|_{u_0 \cdot q_0 \cdot v_0}^{1-\theta}
\big\|(1-e_n)x \big\|_{u_1 \cdot q_1 \cdot v_1}^{\theta}.$$ Now,
let us recall that both norms on the right hand side are uniformly
bounded on $n \ge 1$. In particular, since $0 < \theta < 1$, it
suffices to prove that $e_n x$ tends to $x$ with respect to the
norm of $L_{u_j}(\mathcal{N}) L_{q_j}(\mathcal{M})
L_{v_j}(\mathcal{N})$ for either $j=0$ or $j=1$. There are only
three possible situations:
\begin{itemize}
\item[(a)] Assume $\min(u_0,u_1) < \infty$. Let us suppose
(w.l.o.g.) that $\min(u_0,u_1) = u_0$. Then, we apply the estimate
$$\qquad \big\| (1-e_n) x \big\|_{u_0 \cdot q_0 \cdot v_0} \le \big\|
(1-e_n) a^{\frac{1}{u_0}} \big\|_{L_{u_0}(\mathcal{N})}
\|y_0\|_{L_{q_0}(\mathcal{M})} \big\| b^{\frac{1}{v_0}}
\big\|_{L_{v_0}(\mathcal{N})}.$$ \item[(b)] Assume $\min(q_0,q_1)
< \infty$. Let us suppose (w.l.o.g.) that $\min(q_0,q_1) = q_0$.
Then, we apply the estimate
$$\qquad \big\| (1-e_n) x \big\|_{u_0 \cdot q_0 \cdot v_0} \le \big\|
a^{\frac{1}{u_0}} \big\|_{L_{u_0}(\mathcal{N})} \|(1-e_n)
y_0\|_{L_{q_0}(\mathcal{M})} \big\| b^{\frac{1}{v_0}}
\big\|_{L_{v_0}(\mathcal{N})}.$$ \item[(c)] Finally, assume that
neither (a) nor (b) hold. Let us suppose (w.l.o.g.) that
$\min(v_0,v_1) = v_0$. In this case we have $u_0 = q_0 = u_1 = q_1
= \infty$ and $\Delta = L_{\infty}(\mathcal{N})
L_{\infty}(\mathcal{M}) L_{v_1}(\mathcal{N})$. Therefore,
$\mathcal{M}$ ($= L_{q_0}(\mathcal{M}) \cap L_{q_1}(\mathcal{M})$)
is dense in $\Delta$ (by Proposition
\ref{Proposition-Isometry-Triangle-Complete}) and thereby in
$\mathrm{X}_{\theta}(\mathcal{M})$.
\end{itemize}
\end{proof}

\begin{remark}
\emph{Note that the finiteness of $\mathcal{M}$ is used in the
proof of Lemma \ref{Lemma-Interpolation-Density} to ensure that
$e_n a^{1/u_j}$ and $b^{1/v_j} f_n$ are in $\mathcal{N}$. In
particular, in the case of general von Neumann algebras we shall
need to use a different approach.}
\end{remark}

\textsc{Proof of Lemma \ref{Lemma0}.} The lower estimate is easy.
Indeed, let $$\mathrm{T}_0: L_{u_0}(\ENE) \times
L_{q_0}(\mathcal{M}) \times L_{v_0}(\ENE) \rightarrow
L_{u_0}(\mathcal{N}) L_{q_0}(\mathcal{M}) L_{v_0}(\mathcal{N}),$$
$$\mathrm{T}_1: L_{u_1}(\ENE) \times L_{q_1}(\mathcal{M}) \times
L_{v_1}(\ENE) \rightarrow L_{u_1}(\mathcal{N})
L_{q_1}(\mathcal{M}) L_{v_1}(\mathcal{N}),$$ be the multilinear
maps given by $\mathrm{T}_0(a,y,b) = ayb = \mathrm{T}_1(a,y,b)$.
It is clear that both $\mathrm{T}_0$ and $\mathrm{T}_1$ are
contractive. Hence, the lower estimate follows by multilinear
interpolation. For the converse, we proceed in two steps.

\vskip5pt

\noindent \textsc{Step 1.} Let $x$ be of the form $$x =
\mathrm{D}^{1/u_{\theta}} \alpha y \beta \mathrm{D}^{1/v_{\theta}}
\quad \mbox{with} \quad \alpha, \beta \in \mathcal{N}, \ y \in
L_{q_{\theta}}(\mathcal{M}).$$ Using the isometry $j$ defined in
Proposition \ref{Proposition-Isometry-Triangle-Complete}, we have
$$x \in j_{u_{\theta},v_{\theta}} \Big( \mathcal{N}_{u_{\theta}}
L_{q_{\theta}} (\mathcal{M}) \mathcal{N}_{v_{\theta}} \Big).$$ We
claim that it suffices to see that
\begin{equation} \label{Equation-Basic-Estimate}
|||j_{u_{\theta},v_{\theta}}^{-1}(x)|||_{u_{\theta} \cdot
q_{\theta} \cdot v_{\theta}} \le
\|x\|_{\mathrm{X}_{\theta}(\mathcal{M})}
\end{equation}
for all $x$ of the form considered above. Indeed, let us assume
that inequality (\ref{Equation-Basic-Estimate}) holds. Then,
combining this inequality with the lower estimate proved above, we
obtain
\begin{equation} \label{Equation-Inequality-Equality}
|||j_{u_{\theta},v_{\theta}}^{-1}(x)|||_{u_{\theta} \cdot
q_{\theta} \cdot v_{\theta}} \le
\|x\|_{\mathrm{X}_{\theta}(\mathcal{M})} \le \|x\|_{u_{\theta}
\cdot q_{\theta} \cdot v_{\theta}} \le
|||j_{u_{\theta},v_{\theta}}^{-1}(x)|||_{u_{\theta} \cdot
q_{\theta} \cdot v_{\theta}}.
\end{equation}
Note that, since $\mathrm{X}_{\theta}(\mathcal{M})$ is a Banach
space, we deduce that $\mathcal{N}_{u_{\theta}}
L_{q_{\theta}}(\mathcal{M}) \mathcal{N}_{v_{\theta}}$ satisfies
the triangle inequality. Then, applying again Proposition
\ref{Proposition-Isometry-Triangle-Complete} we deduce that the
space $L_{u_{\theta}}(\mathcal{N}) L_{q_{\theta}}(\mathcal{M})
L_{v_{\theta}}(\mathcal{N})$ is the completion of
$\mathcal{N}_{u_{\theta}} L_{q_{\theta}}(\mathcal{M})
\mathcal{N}_{v_{\theta}}$. Therefore, it follows that for any $x
\in L_{u_{\theta}}(\mathcal{N}) L_{q_{\theta}}(\mathcal{M})
L_{v_{\theta}}(\mathcal{N})$, we have $$\|x\|_{u_{\theta} \cdot
q_{\theta} \cdot v_{\theta}} =
\|x\|_{\mathrm{X}_{\theta}(\mathcal{M})}.$$ In particular,
$L_{u_{\theta}}(\mathcal{N}) L_{q_{\theta}}(\mathcal{M})
L_{v_{\theta}}(\mathcal{N})$ embeds isometrically in
$\mathrm{X}_{\theta}(\mathcal{M})$. To conclude, it remains to see
that both spaces are the same. However, according to Lemma
\ref{Lemma-Interpolation-Density}, we conclude that
$L_{u_{\theta}}(\mathcal{N}) L_{q_{\theta}}(\mathcal{M})
L_{v_{\theta}}(\mathcal{N})$ is norm dense in
$\mathrm{X}_{\theta}(\mathcal{M})$.

\vskip5pt

\noindent \textsc{Step 2.} Now we prove inequality
(\ref{Equation-Basic-Estimate}). We use one more time the
interpolation trick based on Devinatz's factorization theorem
\cite{D}. Again, we refer the reader to Pisier's paper \cite{P0}
for a precise statement of Devinatz's theorem adapted to our aims.
Let $x$ be an element of the form $$x = \mathrm{D}^{1/u_{\theta}}
\alpha y \beta \mathrm{D}^{1/v_{\theta}} \quad \mbox{with} \quad
\alpha, \beta \in \mathcal{N}, \ y \in
L_{q_{\theta}}(\mathcal{M}).$$ Assume the norm of $x$ in
$\mathrm{X}_{\theta}(\mathcal{M})$ is less than one. Then, the
complex interpolation method provides a bounded analytic function
$$f: \mathcal{S} \longrightarrow L_{u_0}(\mathcal{N})
L_{q_0}(\mathcal{M}) L_{v_0}(\mathcal{N}) + L_{u_1}(\mathcal{N})
L_{q_1}(\mathcal{M}) L_{v_1}(\mathcal{N})$$ satisfying $f(\theta)
= x$ and
\begin{equation} \label{Equation-Boundary}
\begin{array}{rcl} \displaystyle \sup_{z \in
\partial_0} \|f(z)\|_{u_0 \cdot q_0 \cdot v_0} & < & 1, \\
\displaystyle \sup_{z \in \partial_1} \|f(z)\|_{u_1 \cdot q_1
\cdot v_1} & < & 1.
\end{array}
\end{equation}
On the other hand, we are assuming that the density $\mathrm{D}$
satisfies the boundedness condition
(\ref{Equation-Boundedness-Density}). In particular, $z \in
\mathcal{S} \mapsto \mathrm{D}^{\lambda z} \in \mathcal{M}$ is a
bounded analytic function for any $\lambda \in \C$. Therefore
(multiplying if necessary on the left and on the right by certain
powers of $\mathrm{D}^z$ and its inverses), we may assume that $f$
has the form $$f(z) = \mathrm{D}^{\frac{1-z}{u_0}}
\mathrm{D}^{\frac{z}{u_1}} f_1(z) \mathrm{D}^{\frac{z}{v_1}}
\mathrm{D}^{\frac{1-z}{v_0}},$$ where $f_1: \mathcal{S} \to
L_{q_0}(\mathcal{M}) + L_{q_1}(\mathcal{M})$ is bounded analytic.
Hence, we deduce from the boundary conditions
(\ref{Equation-Boundary}) and Proposition
\ref{Proposition-Isometry-Triangle-Complete} that $f$ can be
written on the boundary $\partial \mathcal{S}$ as follows
$$f(z) = \mathrm{D}^{\frac{1-z}{u_0}} \mathrm{D}^{\frac{z}{u_1}}
g_1(z) g_2(z) g_3(z) \mathrm{D}^{\frac{z}{v_1}}
\mathrm{D}^{\frac{1-z}{v_0}},$$ where $g_1,g_3:
\partial \mathcal{S} \to \mathcal{N}$ and $g_2: \partial_j \to
L_{q_j}(\mathcal{M})$ satisfy the following estimates
\begin{equation} \label{Equation-Boundary-Bounds}
\begin{array}{rcl}
\displaystyle \sup_{z \in \partial_0} \max \Big\{ \big\|
\mathrm{D}^{\frac{1}{u_0}} g_1(z) \big\|_{u_0}, \|g_2(z)\|_{q_0},
\big\| g_3(z) \mathrm{D}^{\frac{1}{v_0}} \big\|_{v_0} \Big\} & < &
1, \\ \displaystyle \sup_{z \in \partial_1} \max \Big\{ \big\|
\mathrm{D}^{\frac{1}{u_1}} g_1(z) \big\|_{u_1}, \|g_2(z)\|_{q_1},
\big\| g_3(z) \mathrm{D}^{\frac{1}{v_1}} \big\|_{v_1} \Big\} & < &
1.
\end{array}
\end{equation}
Given any $\delta > 0$, we define the following functions on the
boundary
\begin{eqnarray*}
\mathrm{W}_1: z \in \partial \mathcal{S} & \mapsto & g_1(z)
g_1(z)^* + \delta 1 \in \mathcal{N}, \\ \mathrm{W}_3: z \in
\partial \mathcal{S} & \mapsto & g_3(z)^* g_3(z) + \delta 1 \in
\mathcal{N}.
\end{eqnarray*}
According to Devinatz's factorization theorem \cite{P0}, we can
find invertible bounded analytic functions $\mathrm{w}_1,
\mathrm{w}_3: \mathcal{S} \rightarrow \mathcal{N}$ with bounded
analytic inverse and satisfying the following identities on
$\partial \mathcal{S}$
\begin{equation} \label{Equation-Devinatz}
\begin{array}{rcl}
\mathrm{w}_1(z) \mathrm{w}_1(z)^* & = & \mathrm{W}_1(z), \\
\mathrm{w}_3(z)^* \mathrm{w}_3(z) & = & \mathrm{W}_3(z).
\end{array}
\end{equation}
Then, we consider the factorization $$f(z) = h_1(z) h_2(z)
h_3(z)$$ with $h_2(z) = h_1^{-1}(z) f(z) h_3^{-1}(z)$ and $h_1,
h_3$ given by
\begin{eqnarray*}
h_1(z) & = & \mathrm{D}^{\frac{1-z}{u_0}}
\mathrm{D}^{\frac{z}{u_1}} \mathrm{w}_1(z)
\mathrm{D}^{-\frac{z}{u_1}} \mathrm{D}^{\frac{z}{u_0}}, \\ h_3(z)
& = & \mathrm{D}^{\frac{z}{v_0}} \mathrm{D}^{-\frac{z}{v_1}}
\mathrm{w}_3(z) \mathrm{D}^{\frac{z}{v_1}}
\mathrm{D}^{\frac{1-z}{v_0}}.
\end{eqnarray*}
Note than our original hypothesis
(\ref{Equation-Boundedness-Density}) implies the boundedness and
analyticity of $h_1, h_2, h_3$. Then, recalling that
$\mathrm{D}^{\omega}$ is a unitary for any $\omega \in \C$ such
that $\mbox{Re} \, \omega = 0$ and that $2 \le u_j, v_j \le
\infty$ (so that we have triangle inequality on
$L_{u_j/2}(\mathcal{N})$ and $L_{v_j/2}(\mathcal{N})$), we obtain
from (\ref{Equation-Boundary-Bounds}) the following estimates for
$h_1$ and $h_3$
\begin{eqnarray*}
\sup_{z \in \partial_0} \big\| h_1(z)
\big\|_{L_{u_0}(\mathcal{N})}^2 & = & \sup_{z \in \partial_0}
\big\| \mathrm{D}^{\frac{1}{u_0}} \mathrm{w}_1(z)
\mathrm{w}_1(z)^* \mathrm{D}^{\frac{1}{u_0}}
\big\|_{L_{u_0/2}(\mathcal{N})} \\ & \le & \sup_{z \in \partial_0}
\big\| \mathrm{D}^{\frac{1}{u_0}} g_1(z)
\big\|_{L_{u_0}(\mathcal{N})}^2 + \delta < 1 + \delta,
\\ \sup_{z \in \partial_0} \big\| h_3(z)
\big\|_{L_{v_0}(\mathcal{N})}^2 & = & \sup_{z \in \partial_0}
\big\| \mathrm{D}^{\frac{1}{v_0}} \mathrm{w}_3(z)^*
\mathrm{w}_3(z) \mathrm{D}^{\frac{1}{v_0}}
\big\|_{L_{v_0/2}(\mathcal{N})} \\ & \le & \sup_{z \in
\partial_0} \big\| g_3(z) \mathrm{D}^{\frac{1}{v_0}}
\big\|_{L_{v_0}(\mathcal{N})}^2 + \delta < 1 + \delta,
\\ \sup_{z \in \partial_1} \big\| h_1(z) \mathrm{D}^{\frac{1}{u_1} -
\frac{1}{u_0}} \big\|_{L_{u_1}(\mathcal{N})}^2 & = & \sup_{z \in
\partial_1} \big\| \mathrm{D}^{\frac{1}{u_1}} \mathrm{w}_1(z)
\mathrm{w}_1(z)^* \mathrm{D}^{\frac{1}{u_1}}
\big\|_{L_{u_1/2}(\mathcal{N})} \\ & \le & \sup_{z \in
\partial_1} \big\| \mathrm{D}^{\frac{1}{u_1}} g_1(z)
\big\|_{L_{u_1}(\mathcal{N})}^2 + \delta < 1 + \delta,
\\ \sup_{z \in \partial_1} \big\| \mathrm{D}^{\frac{1}{v_1} -
\frac{1}{v_0}} h_3(z) \big\|_{L_{v_1}(\mathcal{N})}^2 & = &
\sup_{z \in \partial_1} \big\| \mathrm{D}^{\frac{1}{v_1}}
\mathrm{w}_3(z)^* \mathrm{w}_3(z) \mathrm{D}^{\frac{1}{v_1}}
\big\|_{L_{v_1/2}(\mathcal{N})} \\ & \le & \sup_{z \in
\partial_1} \big\| g_3(z) \mathrm{D}^{\frac{1}{v_1}}
\big\|_{L_{v_1}(\mathcal{N})}^2 + \delta < 1 + \delta.
\end{eqnarray*}
By Kosaki's interpolation, we obtain
\begin{eqnarray*}
\big\| h_1(\theta) \mathrm{D}^{\frac{\theta}{u_1} -
\frac{\theta}{u_0}} \big\|_{L_{u_{\theta}}(\mathcal{N})} & < &
\sqrt{1 + \delta}, \\ \big\| \mathrm{D}^{\frac{\theta}{v_1} -
\frac{\theta}{v_0}} h_3(\theta)
\big\|_{L_{v_{\theta}}(\mathcal{N})} \hskip1pt & < & \sqrt{1 +
\delta}.
\end{eqnarray*}
On the other hand, using again the unitarity of
$\mathrm{D}^{\omega}$ when $\mbox{Re} \, \omega = 0$, we have
\begin{eqnarray*}
\lefteqn{\sup_{z \in \partial_0} \big\| h_2(z)
\big\|_{L_{q_0}(\mathcal{M})}^2} \\ & = & \sup_{z \in
\partial_0} \big\| \mathrm{w}_1^{-1}(z) g_1(z)
g_2(z) g_3(z) \mathrm{w}_3^{-1}(z) \big\|_{L_{q_0}(\mathcal{M})}^2
\\ & \le & \sup_{z \in \partial_0} \big\| \mathrm{w}_1^{-1}(z)
g_1(z) \big\|_{L_{\infty}(\mathcal{N})}^2 \big\| g_2(z)
\big\|_{L_{q_0}(\mathcal{M})}^2 \big\| g_3(z) \mathrm{w}_3^{-1}(z)
\big\|_{L_{\infty}(\mathcal{N})}^2
\\ & < & \sup_{z \in \partial_0} \big\| \mathrm{w}_1^{-1}(z)
g_1(z) g_1(z)^* \mathrm{w}_1^{-1}(z)^*
\big\|_{L_{\infty}(\mathcal{N})} \big\| \mathrm{w}_3^{-1}(z)^*
g_3(z)^* g_3(z) \mathrm{w}_3^{-1}(z)
\big\|_{L_{\infty}(\mathcal{N})},
\end{eqnarray*}
\begin{eqnarray*}
\lefteqn{\sup_{z \in
\partial_1} \big\| \mathrm{D}^{\frac{1}{u_0} - \frac{1}{u_1}}
h_2(z) \mathrm{D}^{\frac{1}{v_0} - \frac{1}{v_1}}
\big\|_{L_{q_1}(\mathcal{M})}^2} \\ & = & \sup_{z \in
\partial_1} \big\| \mathrm{w}_1^{-1}(z) g_1(z)
g_2(z) g_3(z) \mathrm{w}_3^{-1}(z) \big\|_{L_{q_1}(\mathcal{M})}^2
\\ & \le & \sup_{z \in \partial_1} \big\| \mathrm{w}_1^{-1}(z)
g_1(z) \big\|_{L_{\infty}(\mathcal{N})}^2 \big\| g_2(z)
\big\|_{L_{q_1}(\mathcal{M})}^2 \big\| g_3(z) \mathrm{w}_3^{-1}(z)
\big\|_{L_{\infty}(\mathcal{N})}^2
\\ & < & \sup_{z \in \partial_1} \big\| \mathrm{w}_1^{-1}(z)
g_1(z) g_1(z)^* \mathrm{w}_1^{-1}(z)^*
\big\|_{L_{\infty}(\mathcal{N})} \big\| \mathrm{w}_3^{-1}(z)^*
g_3(z)^* g_3(z) \mathrm{w}_3^{-1}(z)
\big\|_{L_{\infty}(\mathcal{N})}.
\end{eqnarray*}
Then we combine (\ref{Equation-Devinatz}) with the inequalities
$$g_1(z) g_1(z)^* \le g_1(z) g_1(z)^* + \delta 1,$$ $$g_3(z)^*
g_3(z) \le g_3(z)^* g_3(z) + \delta 1,$$ to obtain by Kosaki's
interpolation $$\big\| \mathrm{D}^{\frac{\theta}{u_0} -
\frac{\theta}{u_1}} h_2(\theta) \mathrm{D}^{\frac{\theta}{v_0} -
\frac{\theta}{v_1}} \big\|_{L_{q_{\theta}}(\mathcal{M})} < 1.$$ In
summary, we have obtained the following factorization $$f(\theta)
= \big( \mathrm{D}^{\frac{1}{u_{\theta}}} \mathrm{w}_1(\theta)
\big) \big( \mathrm{D}^{\frac{\theta}{u_0} - \frac{\theta}{u_1}}
h_2(\theta) \mathrm{D}^{\frac{\theta}{v_0} - \frac{\theta}{v_1}}
\big) \big( \mathrm{w}_3(\theta) \mathrm{D}^{\frac{1}{v_{\theta}}}
\big) = j_1(\theta) j_2(\theta) j_3(\theta) ,$$ with
$$\|j_1(\theta)\|_{L_{u_{\theta}}(\mathcal{N})} < \sqrt{1 +
\delta}, \quad \|j_2(\theta)\|_{L_{q_{\theta}}(\mathcal{M})} < 1,
\quad \|j_3(\theta)\|_{L_{v_{\theta}}(\mathcal{N})} < \sqrt{1 +
\delta}.$$ Therefore, (\ref{Equation-Basic-Estimate}) follows by
letting $\delta \rightarrow 0$ in $$\big| \big| \big|
j_{u_{\theta}, v_{\theta}}^{-1} \big( f(\theta) \big) \big| \big|
\big|_{u_{\theta} \cdot q_{\theta} \cdot v_{\theta}} < 1 +
\delta.$$ This concludes the proof of Lemma \ref{Lemma0} for
finite von Neumann algebras. \fin

\begin{remark} \label{Remark-Slice-NO}
\emph{Note that condition ii) of Lemma \ref{Lemma0} is not needed
at any point in the proof given above for finite von Neumann
algebras. This will be crucial in the proof of Theorem
\ref{TheoremA1} below.}
\end{remark}

\section{Conditional expectations on $\partial_{\infty} \mathsf{K}$}
\label{Section4.2}

Before the proof of Lemma \ref{Lemma0} for general von Neumann
algebras, we need some preliminary results in order to adapt
Haagerup's construction to the present context. We begin with a
technical lemma and some auxiliary interpolation results.

\begin{lemma} \label{Lemma-Analyticity}
Let $\mathcal{M}$ be a von Neumann algebra equipped with a n.f.
state $\varphi$ and let $\mathrm{D}$ be the associated density.
Let us consider a bounded analytic function $f: \mathcal{S}
\rightarrow \mathcal{M}$. Then, given $1 \le p < \infty$, the
following functions are also bounded analytic
\begin{eqnarray*}
h_1: z \in \mathcal{S} & \mapsto & \mathrm{D}^{(1-z)/p} f(z)
\mathrm{D}^{z/p} \in L_p(\mathcal{M}),
\\ h_2: z \in \mathcal{S} & \mapsto &
\mathrm{D}^{z/p} f(z) \mathrm{D}^{(1-z)/p} \in L_p(\mathcal{M}).
\end{eqnarray*}
\end{lemma}

\begin{proof} The arguments to be used hold for both $h_1$ and $h_2$.
Hence, we only prove the assertion for $h_1$. The continuity and
boundedness of $h_1$ on the closure of $\mathcal{S}$ is trivial.
To prove the analyticity of $h_1$ we may clearly assume that the
function $f$ is a finite power series $$f(z) = \sum_{k=1}^n x_k
z^k \quad \mbox{with} \quad x_k \in \mathcal{M}.$$ In particular,
it suffices to see that for a fixed element $x_0 \in \mathcal{M}$,
the function $$h_0: z \in \mathcal{S} \mapsto \mathrm{D}^{(1-z)/p}
x_0 \mathrm{D}^{z/p} \in L_p(\mathcal{M})$$ is analytic. If $x_0
\in \mathcal{M}_a$ is an analytic element this is clear. Assume
that $x_0$ is not an analytic element. According to
Pedersen-Takesaki \cite{PT}, the net $(x_{\gamma})_{\gamma > 0}
\subset \mathcal{M}_a$ given by $$x_{\gamma} =
\sqrt{\frac{\gamma}{\pi}} \int_{\R} \sigma_t(x_0) \exp(- \gamma
t^2) \, dt$$ converges strongly to $x_0$ as $\gamma \rightarrow
\infty$. Then Lemma 2.3 in \cite{J1} gives that $$h_{\gamma}(z) =
\mathrm{D}^{(1-z)/p} x_{\gamma} \mathrm{D}^{z/p}$$ converges
pointwise to $h_0(z)$ in the norm of $L_p(\mathcal{M})$. Now, let
us consider a linear functional $\varphi: L_p(\mathcal{M})
\rightarrow \C$ and a cycle $\Gamma$ in $\mathcal{S}$ homologous
to zero with respect to $\mathcal{S}$. Then, since
$\|x_{\gamma}\|_{\mathcal{M}} \le \|x_0\|_{\mathcal{M}}$ for all
$\gamma > 0$, the dominated convergence theorem gives
$$\int_{\Gamma} \varphi(h_0(z)) dz = \lim_{\gamma \rightarrow
\infty} \int_{\Gamma} \varphi(h_{\gamma}(z)) dz = 0.$$ Thus,
$\varphi(h_0)$ is analytic for any linear functional $\varphi:
L_p(\mathcal{M}) \rightarrow \C$ and so is $h_0$. \end{proof}

\begin{lemma} \label{Lemma-Interpolation-Partial0K}
If $2 \le u,v \le \infty$, we have the following isometries
\begin{eqnarray*}
\big[ L_u(\mathcal{N}) L_{\infty}(\mathcal{M}), L_u(\mathcal{N})
L_2(\mathcal{M}) \big]_{\theta} & = & L_u(\mathcal{N})
L_{2/\theta}(\mathcal{M}),
\\ \big[ L_{\infty}(\mathcal{M})
L_v(\mathcal{N}), L_2(\mathcal{M}) L_v(\mathcal{N}) \big]_{\theta}
& = & L_{2/\theta}(\mathcal{M}) L_v(\mathcal{N}).
\end{eqnarray*}
\end{lemma}

\begin{proof} We shall only prove that
$$\mathrm{X}_{\theta}(\mathcal{M}) = \big[ L_u(\mathcal{N})
L_{\infty}(\mathcal{M}), L_u(\mathcal{N}) L_2(\mathcal{M})
\big]_{\theta} = L_u(\mathcal{N}) L_{2/\theta}(\mathcal{M}).$$ The
other isometry can be proved in the same way. The lower estimate
follows by multilinear interpolation just as in the proof of Lemma
\ref{Lemma0} for finite von Neumann algebras. To prove the upper
estimate we first note that, according to Proposition
\ref{Proposition-Isometry-Triangle-Complete}, we have a dense
subset $$\mathrm{D}^{1/u} \mathcal{N} \big( L_2(\mathcal{M}) \cap
L_{\infty}(\mathcal{M}) \big) \subset L_u(\mathcal{N})
L_{2/\theta}(\mathcal{M}).$$ Since this subset is also dense in
$\mathrm{X}_{\theta}(\mathcal{M})$ (note that the intersection in
this case is $L_u(\mathcal{N}) L_{\infty}(\mathcal{M})$) , it
suffices to show that for any $x$ of the form $\mathrm{D}^{1/u} a
y$ with $a \in \mathcal{N}$ and $y \in L_2(\mathcal{M}) \cap
L_{\infty}(\mathcal{M})$, we have the following inequality
\begin{equation} \label{Equation-Enough}
\|x\|_{u \cdot \frac{2}{\theta}} \le
\|x\|_{\mathrm{X}_{\theta}(\mathcal{M})}.
\end{equation}
Assume that the norm of $x$ in $\mathrm{X}_{\theta}(\mathcal{M})$
is less than $1$. Then we can find a bounded analytic function $f:
\mathcal{S} \to L_u(\mathcal{N}) L_{\infty}(\mathcal{M}) +
L_u(\mathcal{N}) L_2(\mathcal{M})$ satisfying $f(\theta) = x$ and
the inequalities
\begin{equation} \label{Equation-Boundary3}
\begin{array}{lcl} \displaystyle \sup_{z \in
\partial_0} \big\| f(z) \big\|_{u \cdot \infty} & < & 1, \\
\displaystyle \sup_{z \in \partial_1} \big\| f(z) \big\|_{u \cdot
2} \ & < & 1.
\end{array}
\end{equation}
Moreover, by our previous considerations we can assume that $f$
has the form $$f(z) = \mathrm{D}^{\frac{1}{u}} f_1(z),$$ where
$f_1: \mathcal{S} \to L_2(\mathcal{M}) \cap
L_{\infty}(\mathcal{M})$ is a bounded analytic function. Then we
deduce from the boundary conditions (\ref{Equation-Boundary3}) and
Proposition \ref{Proposition-Isometry-Triangle-Complete} that $f$
can be written on $\partial \mathcal{S}$ as follows $$f(z) =
\mathrm{D}^{\frac{1}{u}} g_1(z) g_2(z),$$ with $g_1: \partial
\mathcal{S} \to \mathcal{N}$ and $g_2: \partial \mathcal{S} \to
L_2(\mathcal{M}) \cap L_{\infty}(\mathcal{M})$ satisfying $$\max
\Big\{ \sup_{z \in \partial \mathcal{S}} \big\|
\mathrm{D}^{\frac{1}{u}} g_1(z) \big\|_{L_u(\mathcal{N})}, \sup_{z
\in \partial_0} \big\| g_2(z) \big\|_{L_{\infty}(\mathcal{M})},
\sup_{z \in \partial_1} \big\| g_2(z) \big\|_{L_2(\mathcal{M})}
\Big\} < 1.$$ Now we apply Devinatz's theorem \cite{P0} to
$\mathrm{W}: z \in \partial \mathcal{S} \to g_1(z) g_1(z)^* +
\delta 1$, so that we find an invertible bounded analytic function
$\mathrm{w}: \mathcal{S} \to \mathcal{N}$ satisfying $\mathrm{w}
\mathrm{w}^* = \mathrm{W}$ on the boundary. Then we consider the
factorization $$f(z) = \big( \mathrm{D}^{\frac{1}{u}}
\mathrm{w}(z) \big) \big( \mathrm{w}^{-1}(z) f_1(z) \big).$$
Clearly both factors are bounded analytic and we have
\begin{eqnarray*}
\sup_{z \in \partial \mathcal{S}} \big\| \mathrm{D}^{\frac{1}{u}}
\mathrm{w}(z) \big\|_{L_u(\mathcal{N})}^2 \ & = & \sup_{z \in
\partial \mathcal{S}} \big\| \mathrm{D}^{\frac{1}{u}}
\mathrm{w}(z) \mathrm{w}(z)^* \mathrm{D}^{\frac{1}{u}}
\big\|_{L_{u/2}(\mathcal{N})} < 1 + \delta, \\ \sup_{z \in
\partial_0} \big\| \mathrm{w}^{-1}(z) f_1(z)
\big\|_{L_{\infty}(\mathcal{M})} & \le & \sup_{z \in
\partial_0} \big\| \mathrm{w}^{-1}(z) g_1(z)
\big\|_{L_{\infty}(\mathcal{N})} \big\| g_2(z)
\big\|_{L_{\infty}(\mathcal{M})} < 1, \\ \sup_{z \in
\partial_1} \big\| \mathrm{w}^{-1}(z) f_1(z)
\big\|_{L_2(\mathcal{M})} \ & \le & \sup_{z \in
\partial_1} \big\| \mathrm{w}^{-1}(z) g_1(z)
\big\|_{L_{\infty}(\mathcal{N})} \big\| g_2(z)
\big\|_{L_2(\mathcal{M})} \ < 1.
\end{eqnarray*}
Finally, by Kosaki's interpolation and letting $\delta \to 0$, we
obtain inequality (\ref{Equation-Enough}). \end{proof}

\begin{lemma} \label{Lemma-Interpolation-Partial0K-2}
Assume $1/u + 1/2 = 1/p = 1/2 + 1/v$ and $1/q_{\theta} =
(1-\theta)/p + \theta/2$. Then, we have the following isometries
\begin{eqnarray*}
\big[ L_p(\mathcal{M}), L_u(\mathcal{N}) L_2(\mathcal{M})
\big]_{\theta} & = & L_{u/\theta}(\mathcal{N})
L_{q_{\theta}}(\mathcal{M}), \\ \big[ L_p(\mathcal{M}),
L_2(\mathcal{M}) L_v(\mathcal{N}) \big]_{\theta} & = &
L_{q_{\theta}}(\mathcal{M}) L_{v/\theta}(\mathcal{N}).
\end{eqnarray*}
\end{lemma}

\begin{proof} One more time, the lower estimate follows by multilinear
interpolation and we shall only prove the first isometry
$$\mathrm{X}_{\theta}(\mathcal{M}) = \big[ L_p(\mathcal{M}),
L_u(\mathcal{N}) L_2(\mathcal{M}) \big]_{\theta} =
L_{u/\theta}(\mathcal{N}) L_{q_{\theta}}(\mathcal{M}).$$ Let us
point out that H\"{o}lder inequality gives $$\Delta =
L_p(\mathcal{M}) \cap L_u(\mathcal{N}) L_2(\mathcal{M}) =
L_u(\mathcal{N}) L_2(\mathcal{M}).$$ On the other hand, it follows
from Lemma \ref{Lemma-Density-Analytic} that $$\mathrm{D}^{1/p}
\mathcal{M}_a = \mathrm{D}^{1/u} \mathcal{N}_a \mathcal{M}_a
\mathrm{D}^{1/2} = \mathrm{D}^{\theta/u} \mathcal{N}_a
\mathcal{M}_a \mathrm{D}^{1/q_{\theta}}.$$ Hence $\mathrm{D}^{1/p}
\mathcal{M}_a$ is dense in $\Delta$ and in
$L_{u/\theta}(\mathcal{N}) L_{q_{\theta}}(\mathcal{M})$ so that it
suffices to see $$\|x\|_{\frac{u}{\theta} \cdot q_{\theta}} \le
\|x\|_{\mathrm{X}_{\theta}(\mathcal{M})}$$ for any element $x$ of
the form $\mathrm{D}^{1/p} y$ for some $y \in \mathcal{M}_a$.
Assume that the norm of $x$ in $\mathrm{X}_{\theta}(\mathcal{M})$
is less than $1$. Then, according to the considerations above, we
can find a bounded analytic function $f: \mathcal{S} \to
L_p(\mathcal{M}) + L_u(\mathcal{N}) L_2(\mathcal{M})$ of the form
$$f(z) = \mathrm{D}^{1/p} f_1(z) \quad \mbox{with} \quad f_1:
\mathcal{S} \to \mathcal{M}_a \quad \mbox{bounded analytic},$$
satisfying $f(\theta) = x$ and such that
\begin{equation} \label{Equation-Boundary4}
\max \Big\{ \sup_{z \in
\partial_0} \big\| f(z) \big\|_p \, , \sup_{z \in
\partial_1} \big\| f(z) \big\|_{u \cdot 2} \Big\} < 1.
\end{equation}
Moreover, since $f_1$ takes values in $\mathcal{M}_a$ we can
rewrite $f$ as $$f(z) = \mathrm{D}^{\frac{z}{u} + \frac{1-z}{p}}
f_2(z) \mathrm{D}^{\frac{z}{2}} \quad \mbox{with} \quad f_2:
\mathcal{S} \to \mathcal{M}_a \quad \mbox{bounded analytic}.$$ In
particular, $f_{\mid_{\partial_1}}$ has the form $$f(1+it) =
\mathrm{D}^{\frac{1}{u}} \sigma_{-t/2} \big( f_2(1+it) \big)
\mathrm{D}^{\frac{1}{2}}.$$ According to Proposition
\ref{Proposition-Isometry-Triangle-Complete} and the boundary
estimate for $f$ on $\partial_1$, we can write $$f(1+it) = \,
\mathrm{D}^{\frac{1}{u}} g_1(1+it) g_2(1+it),$$ with $g_1:
\partial_1 \to \mathcal{N}$ and $g_2: \partial_1 \to
L_2(\mathcal{M})$ satisfying
\begin{equation} \label{Equation-Boundary5}
\max \Big\{ \sup_{z \in
\partial_1} \big\| \mathrm{D}^{\frac{1}{u}} g_1(z)
\big\|_{L_u(\mathcal{N})}, \sup_{z \in \partial_1} \big\| g_2(z)
\big\|_{L_2(\mathcal{M})} \Big\} < 1.
\end{equation}
By Devinatz's theorem, we can consider an invertible bounded
analytic function $\mathrm{w}: \mathcal{S} \to \mathcal{N}$
satisfying $\mathrm{w}(z) \mathrm{w}(z)^* = \mathrm{W}(z)$ for all
$z \in
\partial \mathcal{S}$ where this time the function $\mathrm{W}:
\partial \mathcal{S} \to \mathcal{N}$ is given by $$\mathrm{W}(z)
= \left\{ \begin{array}{ll} 1, & \mbox{if} \ z \in \partial_0, \\
\sigma_{-\mathrm{Im}z/u}\big( g_1(z)g_1(z)^* \big) + \delta 1, &
\mbox{if} \ z \in \partial_1, \end{array} \right.$$ with $z =
\mathrm{Re}z + i \mathrm{Im}z$. Then we factorize $f$ as $f(z) =
h_1(z) \mathrm{D}^{-\frac{1}{u}} h_2(z)$ where
\begin{eqnarray*}
h_1(z) & = & \mathrm{D}^{\frac{z}{u}} \mathrm{w}(z)
\mathrm{D}^{\frac{1-z}{u}}, \\ h_2(z) & = &
\mathrm{D}^{\frac{z}{u}} \mathrm{w}^{-1}(z)
\mathrm{D}^{\frac{1-z}{u}} \mathrm{D}^{\frac{1-z}{2}} f_2(z)
\mathrm{D}^{\frac{z}{2}}.
\end{eqnarray*}
Note that Lemma \ref{Lemma-Analyticity} provides the boundedness
and analyticity of $h_1$ and $h_2$. As usual, we need to estimate
the norms of $h_1$ and $h_2$ on the boundary. We begin with the
estimates for $h_1$. On $\partial_0$ we have
\begin{eqnarray*}
\sup_{z \in \partial_0} \big\| h_1(z) \mathrm{D}^{-\frac{1}{u}}
\big\|_{L_{\infty}(\mathcal{N})} & = & \sup_{t \in \R} \big\|
\sigma_{t/u} \big( \mathrm{w}(it) \big)
\big\|_{L_{\infty}(\mathcal{N})} \\ & = & \sup_{z \in \partial_0}
\big\| \mathrm{w}(z) \big\|_{L_{\infty}(\mathcal{N})} \\ & = &
\sup_{z \in \partial_0} \big\| \mathrm{w}(z) \mathrm{w}(z)^*
\big\|_{L_{\infty}(\mathcal{N})}^{1/2} = 1.
\end{eqnarray*}
On $\partial_1$ we apply the boundary condition
(\ref{Equation-Boundary5})
\begin{eqnarray*}
\sup_{z \in \partial_1} \big\| h_1(z) \big\|_{L_u(\mathcal{N})}^2
& = & \sup_{t \in \R} \, \big\| \sigma_{t/u} \big(
\mathrm{D}^{\frac{1}{u}} \mathrm{w}(1+it) \big)
\big\|_{L_u(\mathcal{N})}^2 \\ & = & \sup_{z \in
\partial_1} \big\| \mathrm{D}^{\frac{1}{u}} \mathrm{w}(z)
\big\|_{L_u(\mathcal{N})}^2 \\ & = & \sup_{z \in \partial_1}
\big\| \mathrm{D}^{\frac{1}{u}} \mathrm{w}(z) \mathrm{w}(z)^*
\mathrm{D}^{\frac{1}{u}} \big\|_{L_{u/2}(\mathcal{N})} \\ & \le &
\sup_ {z \in \partial_1} \big\| \sigma_{-\mathrm{Im}z/u} \big(
\mathrm{D}^{\frac{1}{u}} g_1(z) g_1(z)^* \mathrm{D}^{\frac{1}{u}}
\big) \big\|_{L_{u/2}(\mathcal{N})} + \delta \\ & = & \sup_ {z \in
\partial_1} \big\| \mathrm{D}^{\frac{1}{u}} g_1(z)
\big\|_{L_u(\mathcal{N})}^2 + \delta < 1 + \delta.
\end{eqnarray*}
For the estimate of $h_2$ on $\partial_0$ we use
(\ref{Equation-Boundary4})
\begin{eqnarray*}
\sup_{z \in \partial_0} \big\| h_2(z) \big\|_{L_p(\mathcal{M})} &
= & \sup_{t \in \R} \, \big\| \sigma_{t/u} \big(
\mathrm{w}^{-1}(it) \mathrm{D}^{\frac{1}{u}} \big) \sigma_{-t/2}
\big( \mathrm{D}^{\frac{1}{2}} f_2(it) \big)
\big\|_{L_p(\mathcal{M})} \\ & = & \sup_{t \in \R} \, \big\|
\sigma_{t/u} \big( \mathrm{w}^{-1}(it) \big) \sigma_{t/u} \big(
\sigma_{-t/p} \big( \mathrm{D}^{\frac{1}{p}} f_2(it) \big) \big)
\big\|_{L_p(\mathcal{M})} \\ & = & \sup_{t \in \R} \, \big\|
\mathrm{w}^{-1}(it) \sigma_{-t/p} \big( \mathrm{D}^{\frac{1}{p}}
f_2(it) \big) \big\|_{L_p(\mathcal{M})} \\ & \le & \sup_{t \in \R}
\, \big\| \big( \mathrm{w}(it) \mathrm{w}(it)^* \big)^{-1}
\big\|_{L_{\infty}(\mathcal{N})}^{1/2} \big\|
\mathrm{D}^{\frac{1}{p}} f_2(it) \big\|_{L_p(\mathcal{M})} \\ & =
& \sup_{t \in \R} \, \big\| \sigma_{-t/2} \big(
\mathrm{D}^{\frac{1}{p}} f_2(it) \big) \big\|_{L_p(\mathcal{M})}
\\ & = & \sup_ {z \in \partial_0} \big\| f(z)
\big\|_{L_p(\mathcal{M})} < 1.
\end{eqnarray*}
Finally, we use again (\ref{Equation-Boundary5}) to estimate $h_2$
on $\partial_1$
\begin{eqnarray*}
\lefteqn{\sup_{z \in \partial_1} \big\| \mathrm{D}^{-\frac{1}{u}}
h_2(z) \big\|_{L_2(\mathcal{M})}} \\ & = & \sup_{t \in \R} \,
\big\| \sigma_{t/u} \big( \mathrm{w}^{-1}(1+it) \big)
\sigma_{-t/2} \big( f_2(1+it) \mathrm{D}^{\frac{1}{2}} \big)
\big\|_{L_2(\mathcal{M})}
\\ & = & \sup_{t \in \R} \, \big\|
\sigma_{t/u} \big( \mathrm{w}^{-1}(1+it) \big) g_1(1+it) g_2(1+it)
\big\|_{L_2(\mathcal{M})} \\ & \le & \sup_{t \in \R} \, \big\|
\sigma_{t/u} \big( \mathrm{w}^{-1}(1+it) \big) g_1(1+it)
\big\|_{L_{\infty}(\mathcal{N})} \big\| g_2(1+it)
\big\|_{L_2(\mathcal{M})} \\ & < & \sup_{t \in \R} \, \big\|
\sigma_{t/u} \big( \mathrm{w}^{-1}(1+it) \big) g_1(1+it)
g_1(1+it)^* \sigma_{t/u} \big( \mathrm{w}^{-1}(1+it) \big)^*
\big\|_{L_{\infty}(\mathcal{N})}^{1/2} \\ & = & \sup_{t \in \R} \,
\big\| \mathrm{w}^{-1}(1+it) \sigma_{-t/u} \big( g_1(1+it)
g_1(1+it)^* \big) \mathrm{w}^{-1}(1+it)^*
\big\|_{L_{\infty}(\mathcal{N})}^{1/2}
\\ & \le & \sup_{z \in \partial_1}
\big\| \mathrm{w}^{-1}(z) \mathrm{w}(z) \mathrm{w}(z)^*
\mathrm{w}^{-1}(z)^* \big\|_{L_{\infty}(\mathcal{N})}^{1/2} = 1.
\end{eqnarray*}
In summary, by Kosaki's interpolation we find $$\big\|h_1(\theta)
\mathrm{D}^{-\frac{1-\theta}{u}}
\big\|_{L_{u/\theta}(\mathcal{N})} < \sqrt{1+\delta} \quad
\mbox{and} \quad \big\| \mathrm{D}^{-\frac{\theta}{u}} h_2(\theta)
\big\|_{L_{q_{\theta}}} < 1.$$ Therefore, since we have
$$f(\theta) = \big( h_1(\theta) \mathrm{D}^{-\frac{1-\theta}{u}}
\big) \big( \mathrm{D}^{-\frac{\theta}{u}} h_2(\theta) \big),$$
the result follows by letting $\delta \to 0$. This completes the
proof. \end{proof}

\begin{remark}
\emph{For the proof of Lemmas \ref{Lemma-Interpolation-Partial0K}
and \ref{Lemma-Interpolation-Partial0K-2} it has been essential to
have an explicit description of the intersection $\Delta$ of the
interpolation pair. The lack of this description in the general
case is what forces us to use Haagerup's construction to extend
the result from finite to general von Neumann algebras.}
\end{remark}

Let us consider a so-called commutative square of conditional
expectations. That is, on one hand we have a von Neumann
subalgebra $\mathcal{N}_1$ of $\mathcal{M}_1$. On the other, we
consider a von Neumann subalgebra $\mathcal{N}_2$ of
$\mathcal{M}_2$ such that $\mathcal{M}_2$ (resp. $\mathcal{N}_2$)
is a von Neumann subalgebra of $\mathcal{M}_1$ (resp.
$\mathcal{N}_1$) so that the diagram in Figure II below commutes.
The interpolation results in Lemmas
\ref{Lemma-Interpolation-Partial0K} and
\ref{Lemma-Interpolation-Partial0K-2}  will allow us to show that
the conditional expectation $\mathsf{E}_{\mathcal{M}}$ extends to
a contractive projection on the amalgamated spaces which
correspond to points $(1/u,1/v,1/q)$ in the set $\partial_{\infty}
\mathsf{K}$.

\begin{figure}[ht!]

\begin{center}

\includegraphics[width=2.5cm]{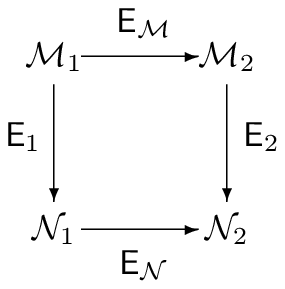}

\end{center}

\vskip-15pt

\caption{\textsc{A commutative square}.}

\end{figure}

\vskip5pt

\begin{observation} \label{Observation-Restriction-Square}
\emph{Let us assume that the commutative square above satisfies
that $\mathsf{E}_1(\mathcal{M}_2) \subset \mathcal{N}_2$. Then we
claim that the restriction ${\mathsf{E}_1}_{\mid_{\mathcal{M}_2}}$
of $\mathsf{E}_1$ to $\mathcal{M}_2$ coincides with
$\mathsf{E}_2$. Indeed, let us consider an element $x \in
\mathcal{M}_2$. Then, since $\mathsf{E}_1(x) \in \mathcal{N}_2$,
we have $$\mathsf{E}_1(x) = \mathsf{E}_{\mathcal{N}} \circ
\mathsf{E}_1(x) = \mathsf{E}_2 \circ \mathsf{E}_{\mathcal{M}}(x) =
\mathsf{E}_2(x).$$ Similarly, we have
${\mathsf{E}_{\mathcal{M}}}_{\mid_{\mathcal{N}_1}} =
\mathsf{E}_{\mathcal{N}}$ whenever
$\mathsf{E}_{\mathcal{M}}(\mathcal{N}_1) \subset \mathcal{N}_2$.
In what follows we shall assume these conditions on the
commutative squares we are using. In fact, the commutative squares
obtained from the Haagerup construction defined below will satisfy
these assumptions.}
\end{observation}

\begin{proposition} \label{Proposition-Expectation-1}
If $2 \le u,v \le \infty$, $\mathsf{E}_{\mathcal{M}}$ extends to
contractions
\begin{eqnarray*}
\mathsf{E}_{\mathcal{M}}: L_u(\mathcal{N}_1) L_2(\mathcal{M}_1) &
\to & L_u(\mathcal{N}_2) L_2(\mathcal{M}_2), \\
\mathsf{E}_{\mathcal{M}}: L_2(\mathcal{M}_1) L_v(\mathcal{N}_1) &
\to & L_2(\mathcal{M}_2) L_v(\mathcal{N}_2).
\end{eqnarray*}
\end{proposition}

\begin{proof} Let us note that for any index $2 \le p \le \infty$,
the natural inclusion $$j: L_p^r(\mathcal{M}_2, \mathsf{E}_2) \to
L_p^r(\mathcal{M}_1, \mathsf{E}_1)$$ is an isometry. Indeed, since
the space $L_p(\mathcal{N}_2)$ embeds isometrically in
$L_p(\mathcal{N}_1)$ and the conditional expectation
$\mathsf{E}_2$ is the restriction
${\mathsf{E}_1}_{\mid_{\mathcal{M}_2}}$ of $\mathsf{E}_1$ to
$\mathcal{M}_2$ (see Observation
\ref{Observation-Restriction-Square}), the following identity
holds for any $x \in L_p^r(\mathcal{M}_2, \mathsf{E}_2)$
$$\|j(x)\|_{L_p^r(\mathcal{M}_1, \mathsf{E}_1)} = \big\|
\mathsf{E}_1(xx^*)^{1/2} \big\|_{L_p(\mathcal{N}_1)} = \big\|
\mathsf{E}_2(xx^*)^{1/2} \big\|_{L_p(\mathcal{N}_2)} =
\|x\|_{L_p^r(\mathcal{M}_2, \mathsf{E}_2)}.$$ Now let us consider
the index $2 \le u \le \infty$ defined by $1/2 = 1/u + 1/p$. When
$u > 2$ we have $p < \infty$ and Proposition
\ref{Proposition-Norm-Conditional-Inf} gives $$L_u(\mathcal{N}_j)
L_2(\mathcal{M}_j) = L_p^r(\mathcal{M}_j, \mathsf{E}_j)^* \quad
\mbox{for} \quad j=1,2.$$ Therefore, our map
$\mathsf{E}_{\mathcal{M}}: L_u(\mathcal{N}_1) L_2(\mathcal{M}_1)
\to L_u(\mathcal{N}_2) L_2(\mathcal{M}_2)$ coincides with the
adjoint of $j$ so that we obtain a contraction. Finally, for $u=2$
we just need to note that $L_u(\mathcal{N}_j) L_2(\mathcal{M}_j)$
embeds isometrically in $L_p^r(\mathcal{M}_j, \mathsf{E}_j)^*$ for
$j=0,1$. Indeed, the first part of the proof of Proposition
\ref{Proposition-Norm-Conditional-Inf} holds even for $p=\infty$.
Therefore, in this case our map is a restriction of $j^*$ to a
closed subspace of $L_p^r(\mathcal{M}_1, \mathsf{E}_1)^*$ and
hence contractive. The proof of the contractivity of the second
mapping follows in the same way after replacing
$L_p^r(\mathcal{M}, \mathsf{E})$ by $L_p^c(\mathcal{M},
\mathsf{E})$.
\end{proof}

\begin{proposition} \label{Proposition-Expectation-2}
If $2 \le u,v \le \infty$, $\mathsf{E}_{\mathcal{M}}$ extends to a
contraction $$\mathsf{E}_{\mathcal{M}}: L_u(\mathcal{N}_1)
L_{\infty}(\mathcal{M}_1) L_v(\mathcal{N}_1) \to
L_u(\mathcal{N}_2) L_{\infty}(\mathcal{M}_2) L_v(\mathcal{N}_2).$$
\end{proposition}

\begin{proof} According to Proposition
\ref{Proposition-Isometry-Triangle-Complete}, it suffices to prove
the assertion on the dense subspace $j_{uv} \big( \mathcal{N}_{1u}
L_{\infty}(\mathcal{M}_1) \mathcal{N}_{1v} \big)$. Thus, let $x$
be an element in such subspace so that it decomposes as $x=
\mathrm{D}^{1/u} \alpha y_0 \beta \mathrm{D}^{1/v}$ with $\alpha,
\beta \in \mathcal{N}_1$ and $y_0 \in \mathcal{M}_1$. Let us
define the operators $$a = \big( \alpha \alpha^* + \delta 1
\big)^{1/2}, \quad b = \big( \beta^* \beta + \delta 1 \big)^{1/2},
\quad y = a^{-1} \alpha y_0 \beta b^{-1},$$ so that $x =
\mathrm{D}^{1/u} a y b \mathrm{D}^{1/v}$. Then, we proceed as in
\cite{JR} by writing $\mathsf{E}_{\mathcal{M}}(x)$ as follows
$$\mathsf{E}_{\mathcal{M}} (x) = \mathrm{\mathrm{D}}^{1/u}
\mathsf{E}_{\mathcal{M}}(a^2)^{1/2} \Big[
\mathsf{E}_{\mathcal{M}}(a^2)^{-1/2} \mathsf{E}_{\mathcal{M}}
(ayb) \mathsf{E}_{\mathcal{M}}(b^2)^{-1/2} \Big]
\mathsf{E}_{\mathcal{M}}(b^2)^{1/2} \mathrm{D}^{1/v}.$$ We clearly
have \begin{equation} \label{Equation-Right-Left}
\begin{array}{lcl} \big\| \mathrm{D}^{\frac{1}{u}}
\mathsf{E}_{\mathcal{M}}(a^2)^{1/2} \big\|_{L_u(\mathcal{N}_2)} &
\le & \Big( \big\| \mathrm{D}^{\frac{1}{u}} \alpha
\big\|_{L_u(\mathcal{N}_1)}^2 + \delta \Big)^{1/2}, \\ \big\|
\mathsf{E}_{\mathcal{M}}(b^2)^{1/2} \, \mathrm{D}^{\frac{1}{v}}
\big\|_{L_v(\mathcal{N}_2)} & \le & \Big( \big\| \beta
\mathrm{D}^{\frac{1}{v}} \big\|_{L_v(\mathcal{N}_1)}^2  + \delta
\Big)^{1/2}.
\end{array}
\end{equation}
On the other hand, let us consider for a moment a von Neumann
subalgebra $\mathcal{N}$ of a given von Neumann algebra
$\mathcal{M}$ and let $\mathsf{E}: \mathcal{M} \to \mathcal{N}$ be
the corresponding conditional expectation. Let us consider a
positive element $\gamma \in \mathcal{M}_+$ satisfying $\gamma \ge
\delta 1$. Then we define the map $$\Lambda_{\gamma}: w \in
\mathcal{M} \to \mathsf{E}(\gamma^2)^{-1/2} \mathsf{E}(\gamma w
\gamma) \mathsf{E}(\gamma^2)^{-1/2} \in \mathcal{N}.$$ According
to \cite{JR}, this is a completely positive map so that
$\|\Lambda_{\gamma}\| = \|\Lambda_{\gamma}(1)\| = 1$. Thus,
$\Lambda_{\gamma}$ is a contraction. Then we apply this result to
$$\gamma = \Big( \begin{array}{cc} a & 0 \\ 0 & b \end{array}
\Big).$$ More concretely, $\Lambda_{\gamma}: \mathrm{M}_2 \otimes
\mathcal{M}_1 \to \mathrm{M}_2 \otimes \mathcal{M}_2$ is given by
$$\Lambda_{\gamma} (w) = \Big[
\widetilde{\mathsf{E}}_{\mathcal{M}} \Big( \begin{array}{cl} a^2 &
0 \\ 0 & b^2 \end{array} \Big) \Big]^{-\frac{1}{2}}
\widetilde{\mathsf{E}}_{\mathcal{M}} (\gamma w \gamma) \ \Big[
\widetilde{\mathsf{E}}_{\mathcal{M}} \Big(
\begin{array}{cl} a^2 & 0 \\ 0 & b^2 \end{array}
\Big) \Big]^{-\frac{1}{2}} \ \ \mbox{with} \ \
\widetilde{\mathsf{E}}_{\mathcal{M}} = id \otimes
\mathsf{E}_{\mathcal{M}}.$$ Then we observe that
$$\Lambda_{\gamma} \Big(
\begin{array}{cc} 0 & y \\ 0 & 0 \end{array}
\Big) = \left( \begin{array}{cc} 0 &
\mathsf{E}_{\mathcal{M}}(a^2)^{-1/2} \mathsf{E}_{\mathcal{M}} (a y
b) \mathsf{E}_{\mathcal{M}}(b^2)^{-1/2}  \\ 0 & 0
\end{array} \right).$$ In particular, since
$\|y\|_{\mathcal{M}_1} \le \big\| a^{-1} \alpha
\big\|_{\mathcal{N}_1} \|y_0\|_{\mathcal{M}_1} \big\| \beta b^{-1}
\big\|_{\mathcal{N}_1} \le \|y_0\|_{\mathcal{M}_1}$, we deduce
\begin{equation} \label{Equation-Middle}
\big\| \mathsf{E}_{\mathcal{M}}(a^2)^{-1/2}
\mathsf{E}_{\mathcal{M}} (a y b)
\mathsf{E}_{\mathcal{M}}(b^2)^{-1/2}
\big\|_{L_{\infty}(\mathcal{M}_2)} \le
\|y_0\|_{L_{\infty}(\mathcal{M}_1)}.
\end{equation}
In summary, according to (\ref{Equation-Right-Left}) and
(\ref{Equation-Middle}) $$\big\| \mathrm{D}^{\frac{1}{u}}
\mathsf{E}_{\mathcal{M}}(\alpha y_0 \beta)
\mathrm{D}^{\frac{1}{v}} \big\|_{u \cdot \infty \cdot v} \le \Big(
\big\| \mathrm{D}^{\frac{1}{u}} \alpha
\big\|_{L_u(\mathcal{N}_1)}^2 + \delta \Big)^{\frac{1}{2}}
\|y_0\|_{L_{\infty}(\mathcal{M}_1)} \Big( \big\| \beta
\mathrm{D}^{\frac{1}{v}} \big\|_{L_v(\mathcal{N}_1)}^2 + \delta
\Big)^{\frac{1}{2}}.$$ The proof is completed by letting $\delta
\to 0$ and taking the infimum on the right. \end{proof}

\begin{corollary} \label{Corollary-Expectation}
If $(1/u,1/v,1/q) \in \partial_{\infty} \mathsf{K}$,
$\mathsf{E}_{\mathcal{M}}$ extends to a contraction
$$\mathsf{E}_{\mathcal{M}}: L_u(\mathcal{N}_1) L_q(\mathcal{M}_1)
L_v(\mathcal{N}_1) \to L_u(\mathcal{N}_2) L_q(\mathcal{M}_2)
L_v(\mathcal{N}_2).$$
\end{corollary}

\begin{proof} The case $1/q = 0$ has already been considered in Proposition
\ref{Proposition-Expectation-2}. Thus, assume (without lost of
generality) that $1/v=0$. Then, we know from Propositions
\ref{Proposition-Expectation-1} and
\ref{Proposition-Expectation-2} that $$\begin{array}{llcl}
\mathsf{E}_{\mathcal{M}}: & L_u(\mathcal{N}_1) L_2(\mathcal{M}_1)
& \to & L_u(\mathcal{N}_2) L_2(\mathcal{M}_2), \\
\mathsf{E}_{\mathcal{M}}: & L_u(\mathcal{N}_1)
L_{\infty}(\mathcal{M}_1) & \to & L_u(\mathcal{N}_2)
L_{\infty}(\mathcal{M}_2), \end{array}$$ are contractions.
Therefore, it follows from Lemma
\ref{Lemma-Interpolation-Partial0K} that the same holds for
$$\mathsf{E}_{\mathcal{M}}: L_u(\mathcal{N}_1) L_q(\mathcal{M}_1)
\to L_u(\mathcal{N}_2) L_q(\mathcal{M}_2)$$ when $2 \le q \le
\infty$. Finally, it remains to consider the case $1 \le q \le 2$.
Given $1 \le q \le 2$ and $2 \le u \le \infty$ such that $1/u +
1/q \le 1$, we consider the index defined by $1/p = 1/u + 1/q$.
Then, according to Lemma \ref{Lemma-Interpolation-Partial0K-2} we
have $$L_u(\mathcal{N}_1) L_q(\mathcal{M}_1) = \big[
L_p(\mathcal{M}_1), L_{u_1}(\mathcal{N}_1) L_2(\mathcal{M}_1)
\big]_{\theta},$$ for some $0 \le \theta \le 1$ and some $2 \le
u_1 \le \infty$. Note that the fact that $u_1$ can be chosen so
that $2 \le u_1 \le \infty$ follows by a quick look at Figure I.
By Proposition \ref{Proposition-Expectation-1},
$\mathsf{E}_{\mathcal{M}}$ is contractive on
$L_{u_1}(\mathcal{N}_1) L_2(\mathcal{M}_1)$ and of course on
$L_p(\mathcal{M}_1)$. Therefore, the result follows by complex
interpolation. \end{proof}

\section{General von Neumann algebras I}
\label{Section4.3}

In this section we prove Lemma \ref{Lemma0} under the assumption
$\min(q_0,q_1) < \infty$. However, most of the arguments we are
giving here will be useful for the remaining case. Before starting
the proof we fix some notation. As in Lemma \ref{Lemma0}, we are
given a von Neumann subalgebra $\mathcal{N}$ of $\mathcal{M}$ and
we denote by $\mathsf{E}: \mathcal{M} \to \mathcal{N}$ the
corresponding conditional expectation. According to the Haagerup
construction sketched in Section \ref{Subsection-Preliminaries},
we have
$$\mathcal{R}_{\mathcal{M}} = \mathcal{M} \rtimes_{\sigma}
\mathrm{G} = \overline{\bigcup_{k \ge 1} \mathcal{M}_k}.$$ As
above, we consider the conditional expectation
$$\mathsf{E}_{\mathcal{M}}: \sum_{g \in \mathrm{G}} x_g \lambda(g)
\in \mathcal{R}_{\mathcal{M}} \mapsto x_0 \in \mathcal{M}$$ and
the state $\widehat{\varphi} = \varphi \circ
\mathsf{E}_{\mathcal{M}}$ on $\mathcal{R}_{\mathcal{M}}$.
Moreover, we have conditional expectations
$\mathcal{E}_{\mathcal{M}_k}: \mathcal{R_M} \to \mathcal{M}_k$ for
each $k \ge 1$ so that (if $\varphi_k$ denotes the restriction of
$\widehat{\varphi}$ to $\mathcal{M}_k$ and
$\mathrm{D}_{\varphi_k}$ stands for the corresponding density) the
family of finite von Neumann subalgebras $(\mathcal{M}_k)_{k \ge
1}$ satisfy
\begin{equation} \label{Equation-Finiteness}
c_1(k) 1_{\mathcal{M}_k} \le \mathrm{D}_{\varphi_k} \le c_2(k)
1_{\mathcal{M}_k}.
\end{equation}
Moreover, given $0 \le \eta \le 1$ and $1 \le p < \infty$, we have
\begin{equation} \label{Equation-Limit}
\lim_{k \to \infty} \Big\|
\mathrm{D}_{\widehat{\varphi}}^{(1-\eta)/p} \big( \hat{x} -
\mathcal{E}_k(\hat{x}) \big)
\mathrm{D}_{\widehat{\varphi}}^{\eta/p}
\Big\|_{L_p(\mathcal{R}_{\mathcal{M}})} = 0
\end{equation}
for any $\hat{x} \in \mathcal{R}_{\mathcal{M}}$. Clearly, we can
consider another Haagerup construction for the von Neumann
subalgebra $\mathcal{N}$ so that the analogous properties hold. In
summary, we sketch the situation in Figure III below.

\vskip1cm

\begin{figure}[ht!]

\begin{center}

\includegraphics[width=5cm]{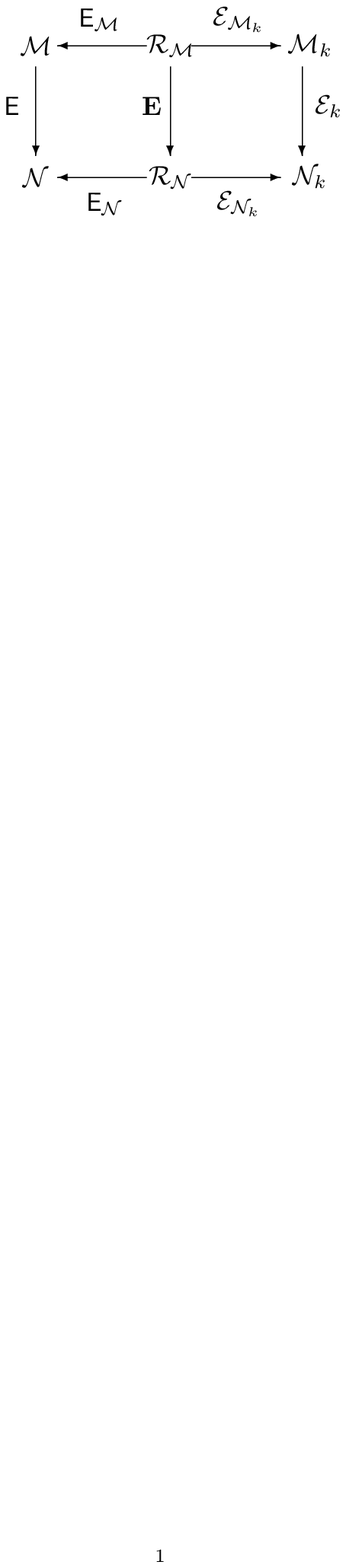}

\end{center}

\vskip-20cm

\caption{\textsc{Haagerup's construction}.}

\end{figure}

\vskip5pt

Now let us consider two points $(1/u_0, 1/v_0, 1/q_0)$ and
$(1/u_1, 1/v_1, 1/q_1)$ in $\partial_{\infty} \mathsf{K}$. Then it
is clear that the natural inclusion mappings $$id:
L_{u_0}(\mathcal{N}) L_{q_0}(\mathcal{M}) L_{v_0}(\mathcal{N}) \to
L_{u_0}(\mathcal{R_N}) L_{q_0}(\mathcal{R_M})
L_{v_0}(\mathcal{R_N}),$$ $$id: L_{u_1}(\mathcal{N})
L_{q_1}(\mathcal{M}) L_{v_1}(\mathcal{N}) \to
L_{u_1}(\mathcal{R_N}) L_{q_1}(\mathcal{R_M})
L_{v_1}(\mathcal{R_N}),$$ are contractions. In particular,
defining for $0 \le \theta \le 1$
\begin{eqnarray*}
\mathrm{X}_{\theta}(\mathcal{M}) \!\! & = & \!\! \Big[
L_{u_0}(\mathcal{N}) L_{q_0}(\mathcal{M}) L_{v_0}(\mathcal{N}),
L_{u_1}(\mathcal{N}) L_{q_1}(\mathcal{M}) L_{v_1}(\mathcal{N})
\Big]_{\theta}, \\ \mathrm{X}_{\theta}(\mathcal{R_M}) \!\! & = &
\!\! \Big[ L_{u_0}(\mathcal{R_N}) L_{q_0}(\mathcal{R_M})
L_{v_0}(\mathcal{R_N}), L_{u_1}(\mathcal{R_N})
L_{q_1}(\mathcal{R_M}) L_{v_1}(\mathcal{R_N}) \Big]_{\theta},
\end{eqnarray*}
we obtain by interpolation that $id:
\mathrm{X}_{\theta}(\mathcal{M}) \to
\mathrm{X}_{\theta}(\mathcal{R_M})$ is also a contraction. Now,
considering the commutative square on the right of Figure III, it
follows from Corollary \ref{Corollary-Expectation} that for each
$k \ge 1$ we have contractions $$\mathcal{E}_{\mathcal{M}_k}:
L_{u_0}(\mathcal{R_N}) L_{q_0}(\mathcal{R_M})
L_{v_0}(\mathcal{R_N}) \to L_{u_0}(\mathcal{N}_k)
L_{q_0}(\mathcal{M}_k) L_{v_0}(\mathcal{N}_k),$$
$$\mathcal{E}_{\mathcal{M}_k}: L_{u_1}(\mathcal{R_N})
L_{q_1}(\mathcal{R_M}) L_{v_1}(\mathcal{R_N}) \to
L_{u_1}(\mathcal{N}_k) L_{q_1}(\mathcal{M}_k)
L_{v_1}(\mathcal{N}_k).$$ Moreover, $\mathcal{M}_k$ satisfies
(\ref{Equation-Finiteness}) for each $k \ge 1$. In particular,
since we have already proved that Lemma \ref{Lemma0} holds in this
case, we obtain by complex interpolation a contraction
$\mathcal{E}_{\mathcal{M}_k}: \mathrm{X}_{\theta}(\mathcal{R_M})
\to L_{u_{\theta}}(\mathcal{N}_k) L_{q_{\theta}}(\mathcal{M}_k)
L_{v_{\theta}}(\mathcal{N}_k)$. In summary, writing
$$\mathcal{E}_{\mathcal{M}_k} = \mathcal{E}_{\mathcal{M}_k} \circ
id,$$ we have found a contraction
\begin{equation} \label{Equation-Contraction-1}
\mathcal{E}_{\mathcal{M}_k}: \mathrm{X}_{\theta}(\mathcal{M}) \to
L_{u_{\theta}}(\mathcal{N}_k) L_{q_{\theta}}(\mathcal{M}_k)
L_{v_{\theta}}(\mathcal{N}_k).
\end{equation}

On the other hand, regarding the space
$L_{u_{\theta}}(\mathcal{N}_k) L_{q_{\theta}}(\mathcal{M}_k)
L_{v_{\theta}}(\mathcal{N}_k)$ as a subspace of
$L_{u_{\theta}}(\mathcal{R_N}) L_{q_{\theta}}(\mathcal{R_M})
L_{v_{\theta}}(\mathcal{R_N})$, we can consider the restriction of
the conditional expectation $\mathsf{E}_{\mathcal{M}}$ to
$\mathcal{M}_k$ and define the following map
$$\mathsf{E}_{\mathcal{M}}: L_{u_{\theta}}(\mathcal{N}_k)
L_{q_{\theta}}(\mathcal{M}_k) L_{v_{\theta}}(\mathcal{N}_k) \to
L_{u_{\theta}}(\mathcal{N}) L_{q_{\theta}}(\mathcal{M})
L_{v_{\theta}}(\mathcal{N}).$$ The following technical result will
be the key to prove Lemma \ref{Lemma0}.

\begin{proposition} \label{Proposition-Contraction-Expectation}
If $k \ge 1$ and $0 \le \theta \le 1$, $\mathsf{E}_{\mathcal{M}}$
extends to a contraction $$\mathsf{E}_{\mathcal{M}}:
L_{u_{\theta}}(\mathcal{N}_k) L_{q_{\theta}}(\mathcal{M}_k)
L_{v_{\theta}}(\mathcal{N}_k) \to L_{u_{\theta}}(\mathcal{N})
L_{q_{\theta}}(\mathcal{M}) L_{v_{\theta}}(\mathcal{N}).$$
\end{proposition}

\begin{proof} It clearly suffices to see our assertion on the dense
subspace $$\mathcal{A}_k = \mathrm{D}_{\varphi_k}^{1/u_{\theta}}
\mathcal{N}_{k,a} \mathrm{D}_{\varphi_k}^{1/2q_{\theta}}
\mathcal{M}_{k,a} \mathrm{D}_{\varphi_k}^{1/2q_{\theta}}
\mathcal{N}_{k,a} \mathrm{D}^{1/v_{\theta}},$$ where
$\mathcal{M}_{k,a}$ and $\mathcal{N}_{k,a}$ denote the
$*$-algebras of analytic elements in $\mathcal{M}_k$ and
$\mathcal{N}_k$ respectively. Thus, let us consider an element $x$
in $\mathcal{A}_k$ and assume that the norm of $x$ in
$L_{u_{\theta}}(\mathcal{N}_k) L_{q_{\theta}}(\mathcal{M}_k)
L_{v_{\theta}}(\mathcal{N}_k)$ is less that $1$. Then, since we
know that Lemma \ref{Lemma0} holds for finite von Neumann algebras
satisfying (\ref{Equation-Finiteness}), we can find a bounded
analytic function $f: \mathcal{S} \to \mathcal{A}_k$ satisfying
$f(\theta) = x$ and the boundary estimate $$\max \Big\{ \sup_{z
\in \partial_0} \big\| f(z) \big\|_{u_0 \cdot q_0 \cdot v_0},
\sup_{z \in \partial_1} \big\| f(z) \big\|_{u_1 \cdot q_1 \cdot
v_1} \Big\} < 1.$$ Then we use the hypothesis $1/u_{\theta} +
1/q_{\theta} + 1/v_{\theta} = 1/u_0+ 1/v_0 + 1/q_0$ and Lemma
\ref{Lemma-Density-Analytic} to rewrite $\mathcal{A}_k$ as follows
$$\mathcal{A}_k = \mathrm{D}_{\varphi_k}^{1/u_0 + 1/2q_0}
\mathcal{M}_{k,a} \mathrm{D}_{\varphi_k}^{1/2q_0 + 1/v_0}.$$
Therefore, multiplying if necessary on the left and on the right
by certain powers of $\mathrm{D}_{\varphi_k}^z$ and its inverses,
we may assume that $f$ has the following form (recall that we are
using here the property (\ref{Equation-Finiteness}))
\begin{equation} \label{Equation-Form-of-f}
f(z) = \mathrm{D}_{\varphi_k}^{\frac{1-z}{u_0}}
\mathrm{D}_{\varphi_k}^{\frac{z}{u_1}}
\mathrm{D}_{\varphi_k}^{\frac{1-z}{2q_0}}
\mathrm{D}_{\varphi_k}^{\frac{z}{2q_1}} \, f_1(z) \,
\mathrm{D}_{\varphi_k}^{\frac{z}{2q_1}}
\mathrm{D}_{\varphi_k}^{\frac{1-z}{2q_0}}
\mathrm{D}_{\varphi_k}^{\frac{z}{v_1}}
\mathrm{D}_{\varphi_k}^{\frac{1-z}{v_0}},
\end{equation}
with $f_1: \mathcal{S} \to \mathcal{M}_{k,a}$ bounded analytic. In
particular, we have
\begin{eqnarray*}
\sup_{z \in \partial_0} \Big\|
\mathrm{D}_{\varphi_k}^{\frac{1}{u_0} + \frac{1}{2q_0}} f_1(z)
\mathrm{D}_{\varphi_k}^{\frac{1}{2q_0} + \frac{1}{v_0}}
\Big\|_{u_0 \cdot q_0 \cdot v_0} & < & 1, \\ \sup_{z \in
\partial_1} \Big\| \mathrm{D}_{\varphi_k}^{\frac{1}{u_1} +
\frac{1}{2q_1}} f_1(z) \mathrm{D}_{\varphi_k}^{\frac{1}{2q_1} +
\frac{1}{v_1}} \Big\|_{u_1 \cdot q_1 \cdot v_1} & < & 1.
\end{eqnarray*}
Note that here we have used the fact that
$\mathrm{D}_{\varphi_k}^w$ is a unitary for any $w \in \C$ such
that $\mbox{Re} \, w = 0$. Now we apply Corollary
\ref{Corollary-Expectation} to the commutative square on the left
of Figure III. In other words, we have contractions
$$\mathsf{E}_{\mathcal{M}}: L_{u_j}(\mathcal{N}_k)
L_{q_j}(\mathcal{M}_k) L_{v_j}(\mathcal{N}_k) \to
L_{u_j}(\mathcal{N}) L_{q_j}(\mathcal{M}) L_{v_j}(\mathcal{N})$$
since $(1/u_j,1/v_j,1/q_j) \in \partial_{\infty} \mathsf{K}$ for
$j=0,1$. Therefore,
\begin{eqnarray*}
\sup_{z \in \partial_0} \Big\| \mathrm{D}_{\varphi}^{\frac{1}{u_0}
+ \frac{1}{2q_0}} \mathsf{E}_{\mathcal{M}}(f_1(z))
\mathrm{D}_{\varphi}^{\frac{1}{2q_0} + \frac{1}{v_0}} \Big\|_{u_0
\cdot q_0 \cdot v_0} & < & 1, \\ \sup_{z \in
\partial_1} \Big\| \mathrm{D}_{\varphi}^{\frac{1}{u_1} +
\frac{1}{2q_1}} \mathsf{E}_{\mathcal{M}}(f_1(z))
\mathrm{D}_{\varphi}^{\frac{1}{2q_1} + \frac{1}{v_1}} \Big\|_{u_1
\cdot q_1 \cdot v_1} & < & 1.
\end{eqnarray*}
Then, according to Proposition
\ref{Proposition-Isometry-Triangle-Complete}, we can find
functions $g_1,g_3: \partial \mathcal{S} \to \mathcal{N}$ and
$g_2: \partial_j \to L_{q_j}(\mathcal{M})$ such that
$$\mathrm{D}_{\varphi}^{1/2q_j} \mathsf{E}_{\mathcal{M}}(f_1(z))
\mathrm{D}_{\varphi}^{1/2q_j} = g_1(z) g_2(z) g_3(z) \quad
\mbox{for all} \quad z \in \partial_j$$ and satisfying the
following boundary estimates
\begin{eqnarray*}
\sup_{z \in \partial_0} \Big\{ \big\|
\mathrm{D}_{\varphi}^{\frac{1}{u_0}} g_1(z)
\big\|_{L_{u_0}(\mathcal{N})}, \big\| g_2(z)
\big\|_{L_{q_0}(\mathcal{M})}, \big\| g_3(z)
\mathrm{D}_{\varphi}^{\frac{1}{v_0}} \big\|_{L_{v_0}(\mathcal{N})}
\Big\} & < & 1, \\ \sup_{z \in \partial_1} \Big\{ \big\|
\mathrm{D}_{\varphi}^{\frac{1}{u_1}} g_1(z)
\big\|_{L_{u_1}(\mathcal{N})}, \big\| g_2(z)
\big\|_{L_{q_1}(\mathcal{M})}, \big\| g_3(z)
\mathrm{D}_{\varphi}^{\frac{1}{v_1}} \big\|_{L_{v_1}(\mathcal{N})}
\Big\} & < & 1.
\end{eqnarray*}
Then we can proceed as in the proof of Lemma \ref{Lemma0} for
finite von Neumann algebras. Namely, we apply Devinatz's
factorization theorem to the functions
\begin{eqnarray*}
\mathrm{W}_1: z \in \partial \mathcal{S} & \mapsto &
\sigma_{\mathrm{Im}z/2q_1 - \mathrm{Im}z/2q_0} \big( g_1(z)
g_1(z)^* \big) + \delta 1 \in \mathcal{N},
\\ \mathrm{W}_3: z \in \partial \mathcal{S} & \mapsto &
\sigma_{\mathrm{Im}z/2q_0 - \mathrm{Im}z/2q_1} \big( g_3(z)^*
g_3(z) \big) + \delta 1 \in \mathcal{N},
\end{eqnarray*}
so that we can find invertible bounded analytic functions
$\mathrm{w}_1, \mathrm{w}_3: \mathcal{S} \to \mathcal{N}$
satisfying
\begin{eqnarray*}
\mathrm{w}_1(z) \mathrm{w}_1(z)^* & = & \mathrm{W}_1(z), \\
\mathrm{w}_3(z)^* \mathrm{w}_3(z) & = & \mathrm{W}_3(z),
\end{eqnarray*}
for all $z \in \partial \mathcal{S}$. Then we consider the
factorization $$\mathsf{E}_{\mathcal{M}} (f(z)) = h_1(z) h_2(z)
h_3(z)$$ with $h_2(z) = h_1^{-1}(z) \mathsf{E}_{\mathcal{M}}
(f(z)) h_3^{-1}(z)$ and $h_1,h_3$ given by
\begin{eqnarray*}
h_1(z) & = & \mathrm{D}_{\varphi}^{\frac{1-z}{u_0}}
\mathrm{D}_{\varphi}^{\frac{z}{u_1}} \mathrm{w}_1(z)
\mathrm{D}_{\varphi}^{-\frac{z}{u_1}}
\mathrm{D}_{\varphi}^{\frac{z}{u_0}}, \\ h_3(z) & = &
\mathrm{D}_{\varphi}^{\frac{z}{v_0}}
\mathrm{D}_{\varphi}^{-\frac{z}{v_1}} \mathrm{w}_3(z)
\mathrm{D}_{\varphi}^{\frac{z}{v_1}}
\mathrm{D}_{\varphi}^{\frac{1-z}{v_0}}.
\end{eqnarray*}
To estimate the norm of $h_j(\theta)$ for $j=1,2,3$ in the
corresponding $L_p$ space we proceed as in the proof of Lemma
\ref{Lemma0} for finite von Neumann algebras. That is, we first
show the boundedness and analyticity of $h_1, h_2, h_3$. This
enables us to estimate the norm on the boundary $\partial
\mathcal{S}$ and apply Kosaki's interpolation. Let us start with
the function $h_1$. To that aim we note that
\begin{itemize}
\item[(a)] If $u_0 \le u_1$, we can write $$h_1(z) =
\mathrm{D}_{\varphi}^{\frac{1}{u_1}}
\mathrm{D}_{\varphi}^{(\frac{1}{u_0}-\frac{1}{u_1})(1-z)}
\mathrm{w}_1(z)
\mathrm{D}_{\varphi}^{(\frac{1}{u_0}-\frac{1}{u_1})z}.$$
\item[(b)] If $u_0 \ge u_1$, we can write $$h_1(z) = \Big(
\mathrm{D}_{\varphi}^{\frac{1}{u_0}}
\mathrm{D}_{\varphi}^{(\frac{1}{u_1}-\frac{1}{u_0})z}
\mathrm{w}_1(z)
\mathrm{D}_{\varphi}^{(\frac{1}{u_1}-\frac{1}{u_0})(1-z)} \Big)
\mathrm{D}_{\varphi}^{(\frac{1}{u_0}-\frac{1}{u_1})}.$$
\end{itemize}
By Lemma \ref{Lemma-Analyticity}, $h_1$ is either bounded analytic
or can be written as $$h_1(z) = j_1(z)
\mathrm{D}_{\varphi}^{(\frac{1}{u_0}-\frac{1}{u_1})},$$ with $j_1$
bounded analytic. In any case, Kosaki's interpolation gives
$$\big\| h_1(\theta) \mathrm{D}_{\varphi}^{\frac{\theta}{u_1} -
\frac{\theta}{u_0}} \big\|_{L_{u_{\theta}}(\mathcal{N})} \le \max
\Big\{ \sup_{z \in
\partial_0} \big\| h_1(z) \big\|_{L_{u_0}(\mathcal{N})},
\sup_{z \in \partial_1} \big\| h_1(z)
\mathrm{D}_{\varphi}^{\frac{1}{u_1} - \frac{1}{u_0}}
\big\|_{L_{u_1}(\mathcal{N})} \Big\}.$$ Moreover, arguing as in
the proof of Lemma \ref{Lemma0} for finite von Neumann algebras
\begin{eqnarray*}
\sup_{z \in \partial_0} \big\| h_1(z)
\big\|_{L_{u_0}(\mathcal{N})} & < & \sqrt{1 + \delta},
\\ \sup_{z \in \partial_1} \big\| h_1(z)
\mathrm{D}_{\varphi}^{\frac{1}{u_1} - \frac{1}{u_0}}
\big\|_{L_{u_1}(\mathcal{N})} & < & \sqrt{1 + \delta}.
\end{eqnarray*}
Indeed, the only significant difference (with respect to the proof
for finite von Neumann algebras) is that we can not assume any
longer that $\mathrm{D}_{\varphi}^w$ is a unitary for a purely
imaginary complex number $w$. However, this difficulty is easily
avoided by using the norm-invariance property of the one-parameter
modular automorphism group. In fact, that is the reason why we
used the modular automorphism group in the definition of
$\mathrm{W}_1$ and $\mathrm{W}_3$, see also the proof of Lemma
\ref{Lemma-Interpolation-Partial0K-2}. Our previous estimates give
rise to
\begin{equation} \label{Equation-Estimacion-Uno}
\big\| h_1(\theta) \mathrm{D}_{\varphi}^{\frac{\theta}{u_1} -
\frac{\theta}{u_0}} \big\|_{L_{u_{\theta}}(\mathcal{N})} < \sqrt{1
+ \delta}.
\end{equation}
Similarly, we have
\begin{equation} \label{Equation-Estimacion-Dos}
\big\| \mathrm{D}_{\varphi}^{\frac{\theta}{v_1} -
\frac{\theta}{v_0}} h_3(\theta)
\big\|_{L_{v_{\theta}}(\mathcal{N})} \hskip1pt < \sqrt{1 +
\delta}.
\end{equation}
Finally, we consider the function $h_2$. To study the boundedness
and analyticity of $h_2$, we recall the expression for $f$
obtained in (\ref{Equation-Form-of-f}). Then, we can rewrite $h_2$
in the following way $$h_2(z) =
\mathrm{D}_{\varphi}^{\frac{z}{u_1}}
\mathrm{D}_{\varphi}^{-\frac{z}{u_0}} \mathrm{w}_1^{-1}(z)
\mathrm{D}_{\varphi}^{\frac{1-z}{2q_0}}
\mathrm{D}_{\varphi}^{\frac{z}{2q_1}}
\mathsf{E}_{\mathcal{M}}(f_1(z))
\mathrm{D}_{\varphi}^{\frac{z}{2q_1}}
\mathrm{D}_{\varphi}^{\frac{1-z}{2q_0}} \mathrm{w}_3^{-1}(z)
\mathrm{D}_{\varphi}^{-\frac{z}{v_0}}
\mathrm{D}_{\varphi}^{\frac{z}{v_1}}.$$ Using one more time that
$1/u_0+1/q_0+1/v_0 = 1/u_1+1/q_1+1/v_1$, we can write

\vskip3pt

\begin{itemize}
\item[(a)] If $\max(u_0,u_1) = u_0$ and $\max(v_0,v_1) = v_0$, we
have
\begin{eqnarray*}
h_2(z) & = & \mathrm{D}_{\varphi}^{\frac{z}{u_1} - \frac{z}{u_0}}
\mathrm{w}_1^{-1}(z) \mathrm{D}_{\varphi}^{\frac{1-z}{u_1} -
\frac{1-z}{u_0}} \\ & \times & \mathsf{E}_{\mathcal{M}} \Big(
\mathrm{D}_{\varphi_k}^{\frac{1-z}{u_0} - \frac{1-z}{u_1} +
\frac{1-z}{2 q_0} + \frac{z}{2 q_1}} f_1(z)
\mathrm{D}_{\varphi_k}^{\frac{z}{2 q_1} + \frac{1-z}{2 q_0} -
\frac{1-z}{v_1} + \frac{1-z}{v_0}} \Big) \\ & \times &
\mathrm{D}_{\varphi}^{\frac{1-z}{v_1} - \frac{1-z}{v_0}}
\mathrm{w}_3^{-1}(z) \mathrm{D}_{\varphi}^{\frac{z}{v_1}
-\frac{z}{v_0}}.
\end{eqnarray*}
The first and third terms on the right hand side are bounded
analytic by Lemma \ref{Lemma-Analyticity}. The middle term is
clearly bounded analytic in $L_{q_1}(\mathcal{M})$. Indeed, it
follows easily from the fact that $\mathrm{D}_{\varphi_k}$
satisfies (\ref{Equation-Finiteness}).

\vskip3pt

\item[(b)] If $\max(u_0,u_1) = u_1$ and $\max(v_0,v_1) = v_1$, we
proceed as above with $$\mathrm{D}_{\varphi}^{\frac{1}{u_0} -
\frac{1}{u_1}} h_2(z) \mathrm{D}_{\varphi}^{\frac{1}{v_0} -
\frac{1}{v_1}}.$$

\vskip3pt

\item[(c)] If $\max(u_0,u_1) = u_0$ and $\max(v_0,v_1) = v_1$, we
proceed as above with $$h_2(z) \mathrm{D}_{\varphi}^{\frac{1}{v_0}
-\frac{1}{v_1}}.$$

\vskip3pt

\item[(d)] If $\max(u_0,u_1) = u_1$ and $\max(v_0,v_1) = v_0$, we
proceed as above with $$\mathrm{D}_{\varphi}^{\frac{1}{u_0}
-\frac{1}{u_1}} h_2(z).$$
\end{itemize}

In any of the possible situations considered, Kosaki's
interpolation provides the following estimate
\begin{eqnarray*}
\lefteqn{\big\| \mathrm{D}_{\varphi}^{\frac{\theta}{u_0} -
\frac{\theta}{u_1}} h_2(\theta)
\mathrm{D}_{\varphi}^{\frac{\theta}{v_0} - \frac{\theta}{v_1}}
\big\|_{L_{q_{\theta}}(\mathcal{M})}} \\ & \le & \max \Big\{
\sup_{z \in \partial_0} \big\| h_2(z)
\big\|_{L_{q_0}(\mathcal{M})}, \sup_{z \in \partial_1} \big\|
\mathrm{D}_{\varphi}^{\frac{1}{u_0} - \frac{1}{u_1}} h_2(z)
\mathrm{D}_{\varphi}^{\frac{1}{v_0} - \frac{1}{v_1}}
\big\|_{L_{q_1}(\mathcal{M})} \Big\}.
\end{eqnarray*}
On the other hand, arguing one more time as in the proof of Lemma
\ref{Lemma0},
\begin{eqnarray*}
\sup_{z \in \partial_0} \big\| h_2(z)
\big\|_{L_{q_0}(\mathcal{M})} & < & 1, \\ \sup_{z \in
\partial_1} \big\| \mathrm{D}_{\varphi}^{\frac{1}{u_0} -
\frac{1}{u_1}} h_2(z) \mathrm{D}_{\varphi}^{\frac{1}{v_0} -
\frac{1}{v_1}} \big\|_{L_{q_1}(\mathcal{M})} & < & 1.
\end{eqnarray*}
In particular,
\begin{equation} \label{Equation-Estimacion-Tres}
\big\| \mathrm{D}_{\varphi}^{\frac{\theta}{u_0} -
\frac{\theta}{u_1}} h_2(\theta)
\mathrm{D}_{\varphi}^{\frac{\theta}{v_0} - \frac{\theta}{v_1}}
\big\|_{L_{q_{\theta}}(\mathcal{M})} < 1.
\end{equation}
In summary, we have found a factorization of
$\mathsf{E}_{\mathcal{M}} (x)$ $$\mathsf{E}_{\mathcal{M}}
(f(\theta)) = \big( h_1(\theta)
\mathrm{D}_{\varphi}^{\frac{\theta}{u_1} - \frac{\theta}{u_0}}
\big) \big( \mathrm{D}_{\varphi}^{\frac{\theta}{u_0} -
\frac{\theta}{u_1}} h_2(\theta)
\mathrm{D}_{\varphi}^{\frac{\theta}{v_0} - \frac{\theta}{v_1}}
\big) \big( \mathrm{D}_{\varphi}^{\frac{\theta}{v_1} -
\frac{\theta}{v_0}} h_3(\theta) \big)$$ which, according to
(\ref{Equation-Estimacion-Uno}, \ref{Equation-Estimacion-Dos},
\ref{Equation-Estimacion-Tres}), provides the estimate $$\big\|
\mathsf{E}_{\mathcal{M}}(x) \big\|_{u_{\theta} \cdot q_{\theta}
\cdot v_{\theta}} < 1+ \delta.$$ Therefore, the proof is concluded
by letting $\delta \to 0$ in the inequality above. \end{proof}

\begin{corollary} \label{Corollary-Contraction-Expectation}
If $k \ge 1$ and $0 \le \theta \le 1$, we have $$\big\|
\mathsf{E}_{\mathcal{M}} \mathcal{E}_{\mathcal{M}_k} (x)
\big\|_{u_{\theta} \cdot q_{\theta} \cdot v_{\theta}} \le
\|x\|_{\mathrm{X}_{\theta}(\mathcal{M})} \quad \mbox{for all}
\quad x \in \mathrm{X}_{\theta}(\mathcal{M}).$$
\end{corollary}

\begin{proof} The result follows automatically from
(\ref{Equation-Contraction-1}) and Proposition
\ref{Proposition-Contraction-Expectation}. \end{proof}

\textsc{Proof of Lemma \ref{Lemma0}.} As in Lemmas
\ref{Lemma-Interpolation-Partial0K} and
\ref{Lemma-Interpolation-Partial0K-2} and also as in the proof for
finite von Neumann algebras, the lower estimate follows by using
multilinear interpolation. This means that we have a contraction
\begin{equation} \label{Equation-Multilinear}
id: L_{u_{\theta}}(\mathcal{N}) L_{q_{\theta}}(\mathcal{M})
L_{v_{\theta}}(\mathcal{N}) \to \mathrm{X}_{\theta}(\mathcal{M}).
\end{equation}
To prove the converse, let us consider the space $$\mathcal{A} =
\mathrm{D}_{\varphi}^{1/u_{\theta}} \mathcal{N}
\mathrm{D}_{\varphi}^{1/2q_{\theta}} \mathcal{M}
\mathrm{D}_{\varphi}^{1/2q_{\theta}} \mathcal{N}
\mathrm{D}_{\varphi}^{1/v_{\theta}}.$$ According to Proposition
\ref{Proposition-gamma-Norm} and Corollary
\ref{Corollary-Contraction-Expectation}, given $x \in \mathcal{A}$
we have $$\|x\|_{u_{\theta} \cdot q_{\theta} \cdot
v_{\theta}}^{\gamma} \le
\|x\|_{\mathrm{X}_{\theta}(\mathcal{M})}^{\gamma} + \lim_{k \to
\infty} \big\| x - \mathsf{E}_{\mathcal{M}}
\mathcal{E}_{\mathcal{M}_k}(x) \big\|_{u_{\theta} \cdot q_{\theta}
\cdot v_{\theta}}^{\gamma}.$$ Therefore, we need to see that the
limit above is $0$. To that aim, we use that $x \in \mathcal{A}$
so that we can write $$x = \mathrm{D}_{\varphi}^{1/u_{\theta}} y
\mathrm{D}_{\varphi}^{1/v_{\theta}} \quad \mbox{with} \quad y \in
L_{q_{\theta}}(\mathcal{M}).$$ Then we use the contractivity of
$\mathsf{E}_{\mathcal{M}}$ on $L_{q_{\theta}}(\mathcal{R_M})$ to
obtain
\begin{eqnarray*}
\lim_{k \to \infty} \big\| x - \mathsf{E}_{\mathcal{M}}
\mathcal{E}_{\mathcal{M}_k}(x) \big\|_{u_{\theta} \cdot q_{\theta}
\cdot v_{\theta}} & = & \lim_{k \to \infty} \big\|
\mathrm{D}_{\varphi}^{1/u_{\theta}} \big( y -
\mathsf{E}_{\mathcal{M}} \mathcal{E}_{\mathcal{M}_k}(y) \big)
\mathrm{D}_{\varphi}^{1/v_{\theta}} \big\|_{u_{\theta} \cdot
q_{\theta} \cdot v_{\theta}} \\ & \le & \lim_{k \to \infty} \big\|
y - \mathsf{E}_{\mathcal{M}} \mathcal{E}_{\mathcal{M}_k}(y)
\big\|_{L_{q_{\theta}}(\mathcal{M})} \\ & = & \lim_{k \to \infty}
\big\| \mathsf{E}_{\mathcal{M}} \big( y -
\mathcal{E}_{\mathcal{M}_k}(y) \big)
\big\|_{L_{q_{\theta}}(\mathcal{M})} \\ & \le & \lim_{k \to
\infty} \big\| y - \mathcal{E}_{\mathcal{M}_k}(y)
\big\|_{L_{q_{\theta}}(\mathcal{R_M})}.
\end{eqnarray*}
Then, recalling our hypothesis $\min(q_0,q_1) < \infty$ assumed at
the beginning of this paragraph, we deduce that $q_{\theta} <
\infty$ for any $0 < \theta < 1$. Therefore, according to
(\ref{Equation-Limit}) we conclude that the limit above is $0$. In
particular, we have $$\|x\|_{u_{\theta} \cdot q_{\theta} \cdot
v_{\theta}} = \|x\|_{\mathrm{X}_{\theta}(\mathcal{M})} \quad
\mbox{for all} \quad x \in \mathcal{A}.$$ Now, using the density
of $\mathcal{A}$ in $L_{u_{\theta}}(\mathcal{N})
L_{q_{\theta}}(\mathcal{M}) L_{v_{\theta}}(\mathcal{N})$ and
(\ref{Equation-Multilinear}) we deduce that the same holds for any
$x \in L_{u_{\theta}}(\mathcal{N}) L_{q_{\theta}}(\mathcal{M})
L_{v_{\theta}}(\mathcal{N})$. Hence, it remains to see that
$L_{u_{\theta}}(\mathcal{N}) L_{q_{\theta}}(\mathcal{M})
L_{v_{\theta}}(\mathcal{N})$ is dense in
$\mathrm{X}_{\theta}(\mathcal{M})$. To that aim, it suffices to
prove that the subspace $$\mathsf{E}_{\mathcal{M}}
\mathcal{E}_{\mathcal{M}_k} (\mathrm{X}_{\theta}(\mathcal{M}))
\subset L_{u_{\theta}}(\mathcal{N}) L_{q_{\theta}}(\mathcal{M})
L_{v_{\theta}}(\mathcal{N})$$ is dense in
$\mathrm{X}_{\theta}(\mathcal{M})$. Moreover, since the
intersection space $$\Delta = L_{u_0}(\mathcal{N})
L_{q_0}(\mathcal{M}) L_{v_0}(\mathcal{N}) \cap
L_{u_1}(\mathcal{N}) L_{q_1}(\mathcal{M}) L_{v_1}(\mathcal{N})$$
is dense in $\mathrm{X}_{\theta}(\mathcal{M})$ for any $0 < \theta
< 1$, we just need to approximate any element $x$ in $\Delta$ by
an element in $\mathsf{E}_{\mathcal{M}}
\mathcal{E}_{\mathcal{M}_k} (\mathrm{X}_{\theta}(\mathcal{M}))$
with respect to the norm of $\mathrm{X}_{\theta}(\mathcal{M})$.
Using one more time that $\min(q_0,q_1) < \infty$ we assume
(without lost of generality) that $ q_0 < \infty$. Then, applying
the three lines lemma $$\big\| x - \mathsf{E}_{\mathcal{M}}
\mathcal{E}_{\mathcal{M}_k} (x)
\big\|_{\mathrm{X}_{\theta}(\mathcal{M})} \le \big\| x -
\mathsf{E}_{\mathcal{M}} \mathcal{E}_{\mathcal{M}_k} (x)
\big\|_{u_0 \cdot q_0 \cdot v_0}^{1-\theta} \big\| x -
\mathsf{E}_{\mathcal{M}} \mathcal{E}_{\mathcal{M}_k} (x)
\big\|_{u_1 \cdot q_1 \cdot v_1}^{\theta},$$ we just need to show
that the first term on the right tends to $0$ as $k \to \infty$
while the second term is uniformly bounded on $k$. The uniform
boundedness follows from Corollary
\ref{Corollary-Contraction-Expectation} since $$\big\| x -
\mathsf{E}_{\mathcal{M}} \mathcal{E}_{\mathcal{M}_k} (x)
\big\|_{u_1 \cdot q_1 \cdot v_1} \le \|x\|_{u_1 \cdot q_1 \cdot
v_1} + \big\| \mathsf{E}_{\mathcal{M}} \mathcal{E}_{\mathcal{M}_k}
(x) \big\|_{u_1 \cdot q_1 \cdot v_1} \le 2 \|x\|_{u_1 \cdot q_1
\cdot v_1}.$$ For the first term, we pick $y \in
L_{q_0}(\mathcal{M})$ so that $$\big\| x -
\mathrm{D}_{\varphi}^{1/u_0} y \mathrm{D}_{\varphi}^{1/v_0}
\big\|_{u_0 \cdot q_0 \cdot v_0} \le \delta.$$ Then we have
\begin{eqnarray*}
\big\| x - \mathsf{E}_{\mathcal{M}} \mathcal{E}_{\mathcal{M}_k}
(x) \big\|_{u_0 \cdot q_0 \cdot v_0} & \le & \big\| x -
\mathrm{D}_{\varphi}^{1/u_0} y \mathrm{D}_{\varphi}^{1/v_0}
\big\|_{u_0 \cdot q_0 \cdot v_0} \\ & + & \big\|
\mathrm{D}_{\varphi}^{1/u_0} \big( y - \mathsf{E}_{\mathcal{M}}
\mathcal{E}_{\mathcal{M}_k}(y) \big) \mathrm{D}_{\varphi}^{1/v_0}
\big\|_{u_0 \cdot q_0 \cdot v_0} \\ & + & \big\|
\mathrm{D}_{\varphi}^{1/u_0} \mathsf{E}_{\mathcal{M}}
\mathcal{E}_{\mathcal{M}_k} (y) \mathrm{D}_{\varphi}^{1/v_0} -
\mathsf{E}_{\mathcal{M}} \mathcal{E}_{\mathcal{M}_k} (x)
\big\|_{u_0 \cdot q_0 \cdot v_0}.
\end{eqnarray*}
In particular, according to Corollary
\ref{Corollary-Contraction-Expectation}
\begin{eqnarray*}
\big\| x - \mathsf{E}_{\mathcal{M}} \mathcal{E}_{\mathcal{M}_k}
(x) \big\|_{u_0 \cdot q_0 \cdot v_0} & \le & 2 \delta + \big\|
\mathrm{D}_{\varphi}^{1/u_0} \big( y - \mathsf{E}_{\mathcal{M}}
\mathcal{E}_{\mathcal{M}_k}(y) \big) \mathrm{D}_{\varphi}^{1/v_0}
\big\|_{u_0 \cdot q_0 \cdot v_0} \\ & \le & 2 \delta + \big\| y -
\mathsf{E}_{\mathcal{M}} \mathcal{E}_{\mathcal{M}_k}(y)
\big\|_{L_{q_0}(\mathcal{M})} \\ & = & 2 \delta + \big\|
\mathsf{E}_{\mathcal{M}} \big( y - \mathcal{E}_{\mathcal{M}_k}(y)
\big) \big\|_{L_{q_0}(\mathcal{M})} \\ & \le & 2 \delta + \big\| y
- \mathcal{E}_{\mathcal{M}_k}(y) \big\|_{L_{q_0}(\mathcal{R_M})}.
\end{eqnarray*}
Finally, since $q_0 < \infty$ by hypothesis, we know that the
second term on the right tends to $0$ as $k \to \infty$. Then, we
let $\delta \to 0$. This completes the proof of Lemma \ref{Lemma0}
for general von Neumann algebras with $\min(q_0,q_1) < \infty$.
\fin

\section{General von Neumann algebras II}

To complete the proof of Lemma \ref{Lemma0} we have to study the
case $\min(q_0,q_1) = \infty$. Note that the proof above fails in
this case since (\ref{Equation-Limit}) does not hold for
$p=\infty$. However, Corollary
\ref{Corollary-Contraction-Expectation} is still valid in this
case so that it suffices to see that
\begin{equation} \label{Equation-Enough-Infty}
\|x\|_{u_{\theta} \cdot \infty \cdot v_{\theta}} \le \sup_{k \ge
1} \big\| \mathsf{E}_{\mathcal{M}} \mathcal{E}_{\mathcal{M}_k} (x)
\big\|_{u_{\theta} \cdot \infty \cdot v_{\theta}} \quad \mbox{for
all} \quad x \in \mathrm{X}_{\theta}(\mathcal{M}).
\end{equation}
Indeed, since the lower estimate follows one more time by
multilinear interpolation, inequality
(\ref{Equation-Enough-Infty}) and Corollary
\ref{Corollary-Contraction-Expectation} are enough to conclude the
proof of Lemma \ref{Lemma0}. In order to prove
(\ref{Equation-Enough-Infty}) we shall need to consider the spaces
\label{RCLp(amal)qE}
$$L_p^r(\mathcal{M}, \mathsf{E}) \oM L_q^c(\mathcal{M},
\mathsf{E}) = \Big\{ \summ_k w_{1k} w_{2k} \, \big| \ w_{1k} \in
L_p^r(\mathcal{M}, \mathsf{E}), \ w_{2k} \in L_q^c(\mathcal{M},
\mathsf{E}) \Big\}$$ for $2 \le p,q \le \infty$ and equipped with
$$\|y\|_{r_p \cdot c_q} = \inf \left\{ \Big\| \Big( \summ_k \mathsf{E}(w_{1k} w_{1k}^*) \Big)^{1/2} \Big\|_{L_p(\mathcal{N})}
\Big\| \Big( \summ_k \mathsf{E}(w_{2k}^* w_{2k}) \Big)^{1/2}
\Big\|_{L_q(\mathcal{N})} \right\}$$ where the infimum runs over
all possible decompositions $$y = \summ_k w_{1k} w_{2k}.$$ It is
not hard to check that $L_p^r(\mathcal{M}, \mathsf{E}) \oM
L_q^c(\mathcal{M}, \mathsf{E})$ is a normed space, see e.g. Lemma
3.5 in \cite{J1} for a similar result. The notation $\oM$ is
motivated by the fact that the norm given above comes from an
amalgamated Haagerup tensor product, we refer the reader to
Chapter \ref{Section6} below for a more detailed explanation. The
following result is the key to conclude the proof of Lemma
\ref{Lemma0}. We use the well known Grothendieck-Pietsch version
of the Hahn-Banach theorem, see \cite{P} for more on this topic.

\label{GroPi}

\begin{theorem} \label{Theorem-Grothendieck}
Let $2 < p,q,u,v < \infty$ related by $$1/u + 1/p = 1/2 = 1/v +
1/q.$$ Then, we have the following isometry via the anti-linear
bracket $\langle x,y \rangle = \mathrm{tr} \big(x^*y \big)$
$$L_u(\mathcal{N}) L_{\infty}(\mathcal{M}) L_v(\mathcal{N}) =
\big( L_p^r(\mathcal{M}, \mathsf{E}) \oM L_q^c(\mathcal{M},
\mathsf{E}) \big)^*.$$
\end{theorem}

\begin{proof} Given $x= ayb$ in $L_u(\mathcal{N})
L_{\infty}(\mathcal{M}) L_v(\mathcal{N})$, H\"{o}lder inequality
gives
\begin{eqnarray*}
\lefteqn{\Big| \mathrm{tr} \Big( x^* \summ_k w_{1k}w_{2k} \Big)
\Big|} \\ & = & \Big| \mathrm{tr} \Big( y^* \Big[ \summ_k a^*
w_{1k} \otimes e_{1k} \Big] \Big[ \summ_k w_{2k} b^* \otimes
e_{k1} \Big] \Big) \Big| \\ & \le & \Big\|a^* \Big( \summ_k w_{1k}
\otimes e_{1k} \Big) \Big\|_2 \|y\|_{\infty} \Big\| \Big( \summ_k
w_{2k} \otimes e_{k1} \Big) b^* \Big\|_2
\\ & = & \mbox{tr} \Big( a a^* \Big[ \summ_k \mathsf{E} (w_{1k} w_{1k}^*)
\Big] \Big)^{1/2} \|y\|_{\infty} \mbox{tr} \Big( \Big[ \summ_k
\mathsf{E}(w_{2k}^* w_{2k}) \Big] b^* b \Big)^{1/2} \\ & \le &
\|a\|_u \|y\|_{\infty} \|b\|_v \Big\| \Big( \summ_k \mathsf{E}
(w_{1k}w_{1k}^*) \Big)^{1/2} \Big\|_p \Big\| \Big( \summ_k
\mathsf{E} (w_{2k}^*w_{2k}) \Big)^{1/2} \Big\|_q.
\end{eqnarray*}
Thus, taking the infimum on the right, we have a contraction $$x
\in L_u(\mathcal{N}) L_{\infty}(\mathcal{M}) L_v(\mathcal{N})
\mapsto \mbox{tr}( x^* \, \cdot) \in \big( L_p^r(\mathcal{M},
\mathsf{E}) \oM L_q^c(\mathcal{M}, \mathsf{E}) \big)^*.$$ To prove
the converse, we take a norm one functional $\varphi$ on
$L_p^r(\mathcal{M}, \mathsf{E}) \oM L_q^c(\mathcal{M},
\mathsf{E})$. If $1/s = 1/p + 1/q$ it is clear that (see Remark
\ref{Remark-p<2<p}) $$L_s(\mathcal{M}) = L_p(\mathcal{M})
L_q(\mathcal{M}) \to L_p^r(\mathcal{M}, \mathsf{E}) \oM
L_q^c(\mathcal{M}, \mathsf{E})$$ is a dense contractive inclusion.
In particular, we can assume that there exists $$x \in
L_{s'}(\mathcal{M}) = L_u(\mathcal{M}) L_v(\mathcal{M})$$
satisfying
\begin{equation} \label{Equation-Trace-Form}
\varphi(y) = \varphi_x(y) = \mbox{tr}(x^*y).
\end{equation}
To conclude, it suffices to see that $\|x\|_{u \cdot \infty \cdot
v} \le 1$. Let us consider a finite family $y_1, y_2, \ldots, y_m$
in the dense subspace $L_p(\mathcal{M}) L_q(\mathcal{M})$ with
decompositions $$y_k = w_{1k} w_{2k}.$$ Since $\varphi_x$ has norm
one $$\Big| \summ_k \varphi_x \big( w_{1k} w_{2k} \big) \Big| \le
\Big\| \summ_k \mathsf{E}(w_{1k} w_{1k}^*)
\Big\|_{L_{p/2}(\mathcal{N})}^{1/2} \Big\| \summ_k \mathsf{E} (
w_{2k}^* w_{2k}) \Big\|_{L_{q/2}(\mathcal{N})}^{1/2}.$$ Moreover,
since the right hand side remains unchanged under multiplication
with unimodular complex numbers $z_k \in \mathbb{T}$, we have the
following inequality $$\summ_k \big| \varphi_x \big( w_{1k} w_{2k}
\big) \big| \le \Big\| \summ_k \mathsf{E}(w_{1k} w_{1k}^*)
\Big\|_{L_{p/2}(\mathcal{N})}^{1/2} \Big\| \summ_k \mathsf{E} (
w_{2k}^* w_{2k}) \Big\|_{L_{q/2}(\mathcal{N})}^{1/2}.$$ Now we
consider the unit balls in $L_u(\mathcal{N})$ and
$L_v(\mathcal{N})$
\begin{eqnarray*}
\mathsf{B}_1 & = & \Big\{ \alpha \in L_u(\mathcal{N}) \, \big| \
\|\alpha\|_{L_u(\mathcal{N})} \le 1 \Big\}, \\ \mathsf{B}_2 & = &
\Big\{ \beta \in L_v(\mathcal{N}) \ \big| \
\|\beta\|_{L_v(\mathcal{N})} \le 1 \Big\}.
\end{eqnarray*}
By the arithmetic-geometric mean inequality
\begin{eqnarray} \label{Equation-Arith-Geo}
\lefteqn{\quad \summ_k \big| \varphi_x \big( w_{1k} w_{2k} \big)
\big|} \\ \nonumber & \le & \frac{1}{2} \Big( \sup_{\alpha \in
\mathsf{B}_1} \summ_k \mbox{tr} \big( \alpha \mathsf{E}(w_{1k}
w_{1k}^*) \alpha^* \big) + \sup_{\beta \in \mathsf{B}_2} \summ_k
\mbox{tr} \big( \beta^* \mathsf{E} (w_{2k}^* w_{2k}) \beta \big)
\Big) \\ \nonumber & = & \frac{1}{2} \Big( \sup_{\alpha \in
\mathsf{B}_1} \summ_k \big\| \alpha w_{1k}
\big\|_{L_2(\mathcal{M})}^2 + \sup_{\beta \in \mathsf{B}_2}
\summ_k \big\| w_{2k} \beta \big\|_{L_2(\mathcal{M})}^2 \Big).
\end{eqnarray}
Note that $\mathsf{B}_1$ and $\mathsf{B}_2$ are compact when
equipped with the $\sigma(L_u(\mathcal{N}), L_{u'}(\mathcal{N}))$
and the $\sigma(L_v(\mathcal{N}), L_{v'}(\mathcal{N}))$ topologies
respectively. Now, labelling $(w_{1k})_{k \ge 1}$ and $(w_{2k})_{k
\ge 1}$ by $\mathsf{w}_1$ and $\mathsf{w}_2$, we consider
$f_{\mathsf{w_1w_2}}: \mathsf{B}_1 \times \mathsf{B}_2 \to \R$
defined by $$f_{\mathsf{w_1w_2}}(\alpha,\beta) = \summ_k \big\|
\alpha w_{1k} \big\|_{L_2(\mathcal{M})}^2 + \summ_k \big\| w_{2k}
\beta \big\|_{L_2(\mathcal{M})}^2 - 2 \summ_k \big| \varphi_x
\big( w_{1k} w_{2k} \big) \big|.$$ This gives rise to the cone
$$\mathsf{C}_+ = \Big\{ f_{\mathsf{w_1w_2}} \in
\mathcal{C}(\mathsf{B}_1 \times \mathsf{B}_2) \, \big| \
w_{1k}w_{2k} \in L_p(\mathcal{M}) L_q(\mathcal{M}) \Big\}.$$ Then
we consider the open cone $$\mathsf{C}_- = \Big\{ f \in
\mathcal{C}(\mathsf{B}_1 \times \mathsf{B}_2) \, \big| \ \sup f <
0 \Big\}.$$ According to (\ref{Equation-Arith-Geo}), the cones
$\mathsf{C}_+$ and $\mathsf{C}_-$ are disjoint. Therefore, the
geometric Hahn-Banach theorem provides a norm one functional $\xi:
\mathcal{C}(\mathsf{B}_1 \times \mathsf{B}_2) \to \R$ satisfying
$$\xi(f_-) < \rho \le \xi(f_+)$$ for some $\rho \in \R$ and all
$(f_+,f_-) \in \mathsf{C}_+ \times \mathsf{C}_-$. Moreover, since
we are dealing with cones, it turns out that $\rho = 0$ and $\xi$
is a positive functional. Then, according to Riesz representation
theorem, there exists a unique (positive) Radon measure
$\mu_{\xi}$ on $\mathsf{B}_1 \times \mathsf{B}_2$ satisfying
\begin{equation} \label{Equation-Integral-Riesz}
\xi(f) = \int_{\mathsf{B}_1 \times \mathsf{B}_2} f \, d \mu_{\xi}
\quad \mbox{for all} \quad f \in \mathcal{C}(\mathsf{B}_1 \times
\mathsf{B}_2).
\end{equation}
In fact, since $\xi$ is a norm one positive functional,
$\mu_{\xi}$ is a probability measure. Now we use that
$\xi_{\mid_{\mathsf{C}_+}}$ takes values in $\R_+$ and
(\ref{Equation-Integral-Riesz}) to obtain the following inequality
\begin{eqnarray*}
2 \summ_k \big| \varphi_x(w_{1k}w_{2k}) \big| & \le & \summ_k
\int_{\mathsf{B}_1 \times \mathsf{B}_2} \mbox{tr} \big( w_{1k}
w_{1k}^* \alpha^* \alpha \big) \, d \mu_{\xi}(\alpha,\beta) \\ & +
& \summ_k \int_{\mathsf{B}_1 \times \mathsf{B}_2} \mbox{tr} \big(
w_{2k}^* w_{2k} \beta \beta^* \big) \, d \mu_{\xi}(\alpha,\beta)
\\ & = & \summ_k \big\| \alpha_0 w_{1k}
\big\|_{L_2(\mathcal{M})}^2 + \summ_k \big\| w_{2k} \beta_0
\big\|_{L_2(\mathcal{M})}^2,
\end{eqnarray*}
where $(\alpha_0,\beta_0) \in \mathsf{B}_1 \times \mathsf{B}_2$
are given by
\begin{eqnarray*}
\alpha_0 & = & \Big( \int_{\mathsf{B}_1 \times \mathsf{B}_2}
\alpha^* \alpha \, d \mu_{\xi}(\alpha,\beta) \Big)^{1/2} \in
\mathsf{B}_1, \\ \beta_0 & = & \Big( \int_{\mathsf{B}_1 \times
\mathsf{B}_2} \beta \beta^* \, d \mu_{\xi}(\alpha,\beta)
\Big)^{1/2} \in \mathsf{B}_2.
\end{eqnarray*}
Then, using the identity $2 rs = \inf_{\gamma > 0} \, (\gamma r)^2
+ (s/\gamma)^2$, we conclude $$\summ_k \big|
\varphi_x(w_{1k}w_{2k}) \big| \le \Big( \summ_k \big\| \alpha_0
w_{1k} \big\|_{L_2(\mathcal{M})}^2 \Big)^{1/2} \Big( \summ_k
\big\| w_{2k} \beta_0 \big\|_{L_2(\mathcal{M})}^2 \Big)^{1/2}.$$
In particular, given any pair $(w_1,w_2) \in L_p(\mathcal{M})
\times L_q(\mathcal{M})$, we have
\begin{equation} \label{Equation-Bound-Bilinear}
\big| \varphi_x(w_1w_2) \big| \le \big\| \alpha_0 w_1
\big\|_{L_2(\mathcal{M})} \big\| w_2 \beta_0
\big\|_{L_2(\mathcal{M})}.
\end{equation}
Let us write $q_{\alpha_0}$ and $q_{\beta_0}$ for the support
projections of $\alpha_0$ and $\beta_0$. Then we define
\begin{eqnarray*}
d_{\alpha_0} & = & \alpha_0^u + (1-q_{\alpha_0}) \mathrm{D}
(1-q_{\alpha_0}), \\ d_{\beta_0} & = & \beta_0^v \hskip0.5pt +
(1-q_{\beta_0}\hskip0.5pt ) \mathrm{D} (1-q_{\beta_0}).
\end{eqnarray*}
Note that $\phi_{\alpha_0} = \mbox{tr}(d_{\alpha_0} \cdot)$ and
$\phi_{\beta_0} = \mbox{tr}(d_{\beta_0} \cdot)$ are \emph{n.f.}
finite weights on $\mathcal{M}$. In particular, by Theorem
\ref{Theorem-Haagerup-Kosaki} we know that $$d_{\alpha_0}^{1/2}
\mathcal{M} \to L_2(\mathcal{M}) \quad \mbox{and} \quad
\mathcal{M} d_{\beta_0}^{1/2} \to L_2(\mathcal{M})$$ are dense
inclusions. Therefore, since $$\begin{array}{rclclcl} q_{\alpha_0}
d_{\alpha_0}^{1/2} \mathcal{M} & = & \alpha_0^{u/2} \mathcal{M} &
= & \alpha_0 \alpha_0^{u/p} \mathcal{M} & \subset & \alpha_0
L_p(\mathcal{M}),
\\ \mathcal{M} d_{\beta_0}^{1/2} q_{\beta_0} & = &
\mathcal{M} \beta_0^{v/2} & = & \mathcal{M} \beta_0^{v/q} \beta_0
& \subset & L_q(\mathcal{M}) \beta_0,
\end{array}$$ it follows that $\alpha_0 L_p(\mathcal{M})$ (resp.
$L_q(\mathcal{M}) \beta_0$) is dense in $q_{\alpha_0}
L_2(\mathcal{M})$ (resp. $L_2(\mathcal{M}) q_{\beta_0}$). Hence we
can consider the linear map $$T_x: q_{\alpha_0} L_2(\mathcal{M})
\to q_{\beta_0} L_2(\mathcal{M})$$ determined by the relation
$$\big\langle \beta_0 w_2^*, T_x(\alpha_0 w_1) \big\rangle =
\mbox{tr} \big( w_2 \beta_0 T_x(\alpha_0 w_1) \big) =
\varphi_x(w_1w_2).$$ According to (\ref{Equation-Bound-Bilinear}),
$T_x$ is contractive. Moreover, $T_x$ is clearly a right
$\mathcal{M}$ module map so that it commutes with the right action
on $\mathcal{M}$. This means that there exists a contraction $m
\in \mathcal{M}$ satisfying $T_x(\alpha_0 w_1) = m \alpha_0 w_1$.
Finally, applying (\ref{Equation-Trace-Form}) we deduce the
following identity $$\mbox{tr} \big( x^* w_1w_2) =
\varphi_x(w_1w_2) = \mbox{tr} \big( T_x(\alpha_0 w_1) w_2 \beta_0
\big) = \mbox{tr} \big( \beta_0 m \alpha_0 w_1 w_2 \big),$$ which
holds for any pair $(w_1,w_2) \in L_p(\mathcal{M}) \times
L_q(\mathcal{M})$. Therefore, by the density of $L_p(\mathcal{M})
L_q(\mathcal{M})$ in $L_p^r(\mathcal{M}, \mathsf{E}) \oM
L_q^c(\mathcal{M}, \mathsf{E})$ we have $$x = \alpha_0^* m^*
\beta_0^*.$$ Then, since $(\alpha_0, \beta_0) \in \mathsf{B}_1
\times \mathsf{B}_2$ and $m$ is contractive, we have $\|x\|_ {u
\cdot \infty \cdot v} \le 1$. \end{proof}

\begin{observation} \label{Observation-max(u,v)}
\emph{With a slight change in the arguments used, we can see that
Theorem \ref{Theorem-Grothendieck} holds for any $(u,v) \in
[2,\infty] \times [2,\infty]$ such that $\max(u,v) > 2$. Indeed,
by symmetry it suffices to see that}
\begin{itemize}
\item[(a)] $L_u(\mathcal{N}) L_\infty(\mathcal{M})
L_2(\mathcal{N}) = \big( L_p^r(\mathcal{M},\mathsf{E})
\ten_{\mathcal{M}} L_\infty^c(\mathcal{M},\mathsf{E}) \big)^*$
\emph{for any} $2 < u \le \infty$. \item[(b)] $L_u(\mathcal{N})
L_\infty(\mathcal{M}) L_\infty(\mathcal{N}) = \big(
L_p^r(\mathcal{M},\mathsf{E}) \ten_{\mathcal{M}}
L_2^c(\mathcal{M},\mathsf{E}) \big)^*$ \emph{for any} $2 < u \le
\infty$.
\end{itemize}
\emph{Since the proofs are similar, we only prove (a). Recalling
that $L_p(\mathcal{M}) L_{\infty}(\mathcal{M})$ is norm dense in
$L_p^r(\mathcal{M},\mathsf{E}) \ten_{\mathcal{M}}
L_\infty^c(\mathcal{M},\mathsf{E})$, we deduce that every norm one
functional $\varphi: L_p^r(\mathcal{M}, \mathsf{E})
\ten_{\mathcal{M}} L_\infty^c(\mathcal{M}, \mathsf{E}) \to \C$ is
given by $\varphi(y) = \varphi_x(y) = \mbox{tr}(x^*y)$ for some $x
\in L_{p'}(\mathcal{M})$. Using one more time the
Grothendieck-Pietsch separation trick we get
\[ \big| \mbox{tr}(x^*w_1w_2) \big| \le \big\| \alpha_0 w_1 \big\|_2
\psi \big( \mathsf{E}(w_2^*w_2) \big)^{1/2} \] for some $\alpha_0$
in the unit ball of $L_u(\mathcal{N})$ and $\psi \in
\mathcal{N}^*$. Let $(e_{\alpha})$ be a net such that $e_{\al} \to
1$ strongly so that $\lim_{\alpha} \varphi (e_{\al}y) =
\varphi_{n}(y)$ gives the normal part. Now replace $y$ by
$ye_{\al}$ and we get in the limit \[ \big| \mbox{tr}(x^*w_1w_2)
\big| \le \big\| \alpha_0 w_1 \big\|_2 \varphi_n\big(
\mathsf{E}(w_2^*w_2) \big)^{1/2}.\] Thus $|\mbox{tr}(x^* w_1w_2)|
\le \|\alpha_0w_1\|_2 \|w_2\beta_0\|_2$ and we may continue as in
Theorem \ref{Theorem-Grothendieck}.}
\end{observation}

\textsc{Proof of Lemma \ref{Lemma0}.} As we already pointed out at
the beginning of this section, we have to prove inequality
(\ref{Equation-Enough-Infty}). Before doing it we recall that the
indices $(u_{\theta}, v_{\theta})$ satisfy
\begin{equation} \label{Equation-uv-Theta}
2 < \max(u_{\theta}, v_{\theta}) \quad \mbox{and} \quad
\min(u_{\theta}, v_{\theta}) < \infty.
\end{equation}
for any $0 < \theta < 1$. Indeed, otherwise we would have
$u_{\theta} = v_{\theta} = 2$ or $u_{\theta} = v_{\theta} =
\infty$. However, $(1/2,1/2,0)$ and $(0,0,0)$ are extreme points
of $\mathsf{K} \cap \{z=0\}$ and this is not possible for $0 <
\theta < 1$. The first inequality in (\ref{Equation-uv-Theta})
allows us to apply Theorem \ref{Theorem-Grothendieck} after
Observation \ref{Observation-max(u,v)}, while the second
inequality will be used below. Now, let $x$ be an element of
$L_{u_{\theta}}(\mathcal{N}) L_{\infty}(\mathcal{M})
L_{v_{\theta}}(\mathcal{N})$. According to Theorem
\ref{Theorem-Grothendieck}, we know that we can find a finite
family $(w_{1k},w_{2k}) \in L_{p_{\theta}}(\mathcal{M}) \times
L_{q_{\theta}}(\mathcal{M})$ with $$1/u_{\theta} + 1/p_{\theta} =
1/2 = 1/v_{\theta} + 1/q_{\theta}$$ so that
\begin{equation} \label{Equation-Norm-One}
\max \left\{ \Big\| \Big( \summ_k \mathsf{E}(w_{1k}w_{1k}^*)
\Big)^{1/2} \Big\|_{L_{p_{\theta}} (\mathcal{N})}, \Big\|\Big(
\summ_k \mathsf{E}(w_{2k}^*w_{2k}) \Big)^{1/2}
\Big\|_{L_{q_{\theta}}(\mathcal{N})} \right\} \le 1,
\end{equation} and $$\|x\|_{u_{\theta} \cdot \infty \cdot
v_{\theta}} \le \Big| \mbox{tr} \Big( x^* \summ_k w_{1k}w_{2k}
\Big) \Big| + \delta.$$ Moreover, given $1/s = 1/u_{\theta} +
1/v_{\theta}$, we have
\begin{eqnarray*}
\Big| \mbox{tr} \Big( \big[ x - \mathsf{E}_{\mathcal{M}}
\mathcal{E}_{\mathcal{M}_k}(x) \big]^* \summ_j w_{1j}w_{2j} \Big)
\Big| & \le & \big\|x - \mathsf{E}_{\mathcal{M}}
\mathcal{E}_{\mathcal{M}_k}(x) \big\|_s \Big\| \summ_j w_{1j}
w_{2j} \Big\|_{s'}.
\end{eqnarray*}
According to (\ref{Equation-uv-Theta}) we have $1 < s < \infty$.
Then, it follows from (\ref{Equation-Limit}) that the first factor
on the right hand side tends to $0$ as $k \to \infty$ while the
second factor belongs to $L_{s'}(\mathcal{M})$. In conclusion, we
obtain the following estimate
\begin{eqnarray*}
\|x\|_{u_{\theta} \cdot \infty \cdot v_{\theta}} & \le & \lim_{k
\to \infty} \Big| \mbox{tr} \Big( \mathsf{E}_{\mathcal{M}}
\mathcal{E}_{\mathcal{M}_k}(x)^* \summ_k w_{1k}w_{2k} \Big) \Big|
+ \delta.
\end{eqnarray*}
Using (\ref{Equation-Norm-One}) and applying Theorem
\ref{Theorem-Grothendieck} one more time $$\|x\|_{u_{\theta} \cdot
\infty \cdot v_{\theta}} \le \sup_{k \ge 1} \big\|
\mathsf{E}_{\mathcal{M}} \mathcal{E}_{\mathcal{M}_k}(x)
\big\|_{u_{\theta} \cdot \infty \cdot v_{\theta}} + \delta.$$
Thus, (\ref{Equation-Enough-Infty}) follows for $x \in
L_{u_{\theta}}(\mathcal{N}) L_{\infty}(\mathcal{M})
L_{v_{\theta}}(\mathcal{N})$ by letting $\delta \to 0$. Finally,
as in the case $\min(q_0,q_1) < \infty$, it remains to see that
$L_{u_{\theta}}(\mathcal{N}) L_{\infty}(\mathcal{M})
L_{v_{\theta}}(\mathcal{N})$ is dense in
$\mathrm{X}_{\theta}(\mathcal{M})$. Here we also need a different
argument. Let us keep the notation $1/u_0 + 1/v_0 = 1/s = 1/u_1 +
1/v_1$. Then, we may assume that $(u_0,v_0) = (s,\infty)$ and
$(u_1,v_1) = (\infty,s)$. Indeed, if we conclude the proof in this
particular case, the general case follows from the reiteration
theorem for complex interpolation, see e.g. \cite{BL}. This can be
justified by means of Figure I, since the segment joining the
points $(1/u_0,1/v_0,0)$ and $(1/u_1,1/v_1,0)$ is always contained
in the segment joining $(1/s,0,0)$ and $(0,1/s,0)$. Thus, we
assume in what follows that $$\mathrm{X}_{\theta}(\mathcal{M}) =
\Big[ L_s(\mathcal{N}) L_{\infty}(\mathcal{M})
L_{\infty}(\mathcal{N}), L_{\infty}(\mathcal{N})
L_{\infty}(\mathcal{M}) L_s(\mathcal{N}) \Big]_{\theta}$$ so that
$1/u_{\theta} = (1-\theta)/s$ and $1/v_{\theta} = \theta/s$. By
the density of $$\Delta = L_s(\mathcal{N}) L_{\infty}(\mathcal{M})
L_{\infty}(\mathcal{N}) \cap L_{\infty}(\mathcal{N})
L_{\infty}(\mathcal{M}) L_s(\mathcal{N})$$ in
$\mathrm{X}_{\theta}(\mathcal{M})$, it suffices to approximate any
element $x \in \Delta$. In particular, we can write $x=a_0b_0$ and
$x=b_1a_1$ where $a_0,a_1 \in L_s(\mathcal{N})$ and $b_0, b_1 \in
\mathcal{M}$. Moreover, we can assume that $a_0 = a_1$. Indeed,
taking $$a = \big( a_0a_0^* + a_1^*a_1 + \delta
\mathrm{D}_{\varphi}^{2/s} \big)^{1/2}$$ we have $x = a a^{-1} a_0
b_0 = a c_0$ and $x = b_1 a_1 a^{-1} a = c_1 a$ with
$$\|c_j\|_{\mathcal{M}} \le \|b_j\|_{\mathcal{M}} \quad \mbox{for}
\quad j=0,1.$$ Then $c_1= a c_0 a^{-1}$ and we claim that
$a^{\theta} c_0 a^{-\theta}$ is in $\mathcal{M}$ for all
$0<\theta<1$. Indeed, let $\phi$ be the \emph{n.f.} finite weight
$\phi(\cdot) = \mbox{tr}(a^s \cdot)$ and let $\mathcal{M}
\rtimes_{\sigma^{\phi}} \mathbb{R}$ be the crossed product with
respect to the modular automorphism group associated to $\phi$.
Let us consider the spectral projection $p_n = 1_{[1/n,n]}(a)$.
Then, for a fixed integer $n$ the function $$f_n(z)= p_n a^z c_0
a^{-z} p_n = p_n a^z p_n c_0 p_n a^{-z} p_n$$ is analytic so that
$$\|f_n(\theta)\| \le \sup_{t \in \R} \big\|
\sigma_{t/s}^{\phi}(c_0) \big\|^{1-\theta} \big\|
\sigma_{t/s}^{\phi}(ac_0a^{-1}) \big\|^{\theta} =
\|c_0\|^{1-\theta} \|a c_0 a^{-1}\|^{\theta} = \|c_0\|^{1-\theta}
\|c_1\|^{\theta}.$$ Sending $n$ to infinity, we deduce that
$a^{\theta} c_0 a^{-\theta}$ is a bounded element of $\mathcal{M}
\rtimes_{\sigma^{\phi}} \mathbb{R}$. Moreover, using the dual
action with respect to $\phi$, we find that $a^{\theta} c_0
a^{-\theta}$ belongs to $\mathcal{M}$. Therefore, we obtain $$x=
a^{1-\theta} a^{\theta} c_0 a^{-\theta} a^{\theta} \in
L_{u_{\theta}}(\mathcal{N}) L_{\infty}(\mathcal{M})
L_{v_{\theta}}(\mathcal{N}).$$ We have seen that the intersection
space $\Delta$ is included in $L_{u_{\theta}}(\mathcal{N})
L_{\infty}(\mathcal{M}) L_{v_{\theta}}(\mathcal{N})$. Hence, the
result follows since $\Delta$ is dense in
$\mathrm{X}_{\theta}(\mathcal{M})$. The proof of Lemma
\ref{Lemma0} (for any von Neumann algebra) is therefore completed.
\fin

\section{Proof of the main interpolation theorem}

To prove Theorem \ref{TheoremA1}, we need to know a priori that
$L_u(\mathcal{N}) L_q(\mathcal{M}) L_v(\mathcal{N})$ is a Banach
space for any indices $(u,q,v)$ associated to a point
$(1/u,1/v,1/q) \in \mathsf{K}$. This is a simple consequence of
Lemma \ref{Lemma0}. Indeed, according to Lemma
\ref{Lemma-Triangle-Inequality} and Proposition
\ref{Proposition-Isometry-Triangle-Complete}, we know that our
assertion is true for any $(1/u,1/v,1/q) \in \partial_{\infty}
\mathsf{K}$. In particular, it follows from Lemma \ref{Lemma0}
that $$L_{u_{\theta}}(\mathcal{N}) L_{q_{\theta}}(\mathcal{M})
L_{v_{\theta}}(\mathcal{N})$$ is a Banach space for any $0 \le
\theta \le 1$ whenever $(u_j,q_j,v_j) \in \partial_{\infty}
\mathsf{K}$ for $j=0,1$ and $1/u_0 + 1/q_0 + 1/v_0 = 1/u_1 + 1/q_1
+ 1/v_1$. In other words, according to the notation introduced at
the beginning of this chapter, this condition holds whenever
$(1/u_j,1/v_j,1/q_j) \in \mathsf{K}_{\tau} \cap
\partial_{\infty} \mathsf{K}$ for $j=0,1$ and some $0 \le \tau \le 1$.
Therefore, it suffices to see that $\mathsf{K}_{\tau}$ is the
convex hull of $\mathsf{K}_{\tau} \cap
\partial_{\infty} \mathsf{K}$ for any $0 \le \tau \le 1$. However,
this follows easily from Figure I. Note that $\mathsf{K}_{\tau}$
is either a point ($\tau=0$), a triangle ($0 < \tau \le 1/2$), a
pentagon ($1/2 < \tau < 1$) or a parallelogram ($\tau=1$).

\begin{observation} \label{Observation-Subspace-Normed}
\emph{In fact, in the case of finite von Neumann algebras, Lemma
\ref{Lemma0} provides more information. Namely, according to
(\ref{Equation-Inequality-Equality}) and the fact that
$\mathsf{K}_{\tau}$ is the convex hull of $\mathsf{K}_{\tau} \cap
\partial_{\infty} \mathsf{K}$, we deduce that $\mathcal{N}_u
L_q(\mathcal{M}) \mathcal{N}_v$ is a normed space when equipped
with $||| \ |||_{u \cdot q \cdot v}$ for any $(1/u,1/v,1/q) \in
\mathsf{K}$. Moreover, $\mathcal{N}_u L_q(\mathcal{M})
\mathcal{N}_v$ embeds isometrically in $L_u(\mathcal{N})
L_q(\mathcal{M}) L_v(\mathcal{N})$ as a dense subspace, something
we did not know up to now (see Observation
\ref{Observation-Lack-gamma-Norm}). This means that Lemma
\ref{Lemma-Triangle-Inequality} and Proposition
\ref{Proposition-Isometry-Triangle-Complete} hold for any point in
$\mathsf{K}$ in the case of finite von Neumann algebras.}
\end{observation}

\vskip3pt

\textsc{Proof of Theorem \ref{TheoremA1}.} Now we are ready to
prove Theorem \ref{TheoremA1}. The arguments to be used follow the
same strategy used for the proof of Lemma \ref{Lemma0}. In
particular, we only need to point out how to proceed. In first
place, as usual, the lower estimate follows by multilinear
interpolation.

\vskip3pt

\noindent \textsc{Step 1.} Let us show the validity of Theorem
\ref{TheoremA1} for finite von Neumann algebras satisfying the
boundedness condition (\ref{Equation-Boundedness-Density}). First
we note that, once we know that $L_u(\mathcal{N}) L_q(\mathcal{M})
L_v(\mathcal{N})$ is always a Banach space, the proof of Lemma
\ref{Lemma-Interpolation-Density} is still valid for ending points
$(1/u_j,1/v_j,1/q_j)$ lying on $\mathsf{K} \setminus
\partial_{\infty} \mathsf{K}$. Then we follow the proof of Lemma
\ref{Lemma0} for finite von Neumann algebras verbatim to deduce
Theorem \ref{TheoremA1} in this case. Here is essential to observe
(as we did in Remark \ref{Remark-Slice-NO}) that the proof of
Lemma \ref{Lemma0} for finite von Neumann algebras does not use at
any point the restriction $1/u_0 + 1/q_0 + 1/v_0 = 1/u_1 + 1/q_1 +
1/v_1$. Note also that, in order to obtain the boundary estimates
(\ref{Equation-Boundary-Bounds}), Proposition
\ref{Proposition-Isometry-Triangle-Complete} is needed. Here is
where we apply Observation \ref{Observation-Subspace-Normed}. This
proves Theorem \ref{TheoremA1} for finite von Neumann algebras.

\vskip3pt

\noindent \textsc{Step 2.} The next goal is to show  that the
corresponding conditional expectations are contractive. First we
observe that, according to Lemma \ref{Lemma0}, we can extend the
validity of Corollary \ref{Corollary-Expectation} to any point
$(1/u,1/v,1/q) \in \mathsf{K}$ by complex interpolation. Here we
use again that $\mathsf{K}_{\tau} = \mbox{conv}
\big(\mathsf{K}_{\tau} \cap \partial_{\infty} \mathsf{K} \big)$.
Then, it is straightforward to see that Corollary
\ref{Corollary-Contraction-Expectation} also holds in this case.
Indeed, first we apply complex interpolation to obtain a
contraction $$\mathcal{E}_{\mathcal{M}_k}:
\mathrm{X}_{\theta}(\mathcal{M}) \to
\mathrm{X}_{\theta}(\mathcal{R_M}) \to
L_{u_{\theta}}(\mathcal{N}_k) L_{q_{\theta}}(\mathcal{M}_k)
L_{v_{\theta}}(\mathcal{N}_k).$$ Second, the contractivity of
$$\mathsf{E}_{\mathcal{M}}: L_{u_{\theta}}(\mathcal{N}_k)
L_{q_{\theta}}(\mathcal{M}_k) L_{v_{\theta}}(\mathcal{N}_k) \to
L_{u_{\theta}}(\mathcal{N}) L_{q_{\theta}}(\mathcal{M})
L_{v_{\theta}}(\mathcal{N})$$ follows since, as we have seen,
Corollary \ref{Corollary-Expectation} holds for any point
$(1/u,1/v,1/q) \in \mathsf{K}$.

\vskip3pt

\noindent \textsc{Step 3.} We now prove Theorem \ref{TheoremA1} in
the case $\min(q_0,q_1) < \infty$. It follows easily from Step 2.
Indeed, recalling again that the restriction $1/u_0 + 1/q_0 +
1/v_0 = 1/u_1 + 1/q_1 + 1/v_1$ is not used in the proof of Lemma
\ref{Lemma0} (once we know the validity of Corollary
\ref{Corollary-Contraction-Expectation}), the proof follows
verbatim.

\vskip3pt

\noindent \textsc{Step 4.} Finally, we consider the case
$\min(q_0,q_1) = \infty$. First we observe that the first half of
the proof of Lemma \ref{Lemma0} for this case holds for any two
ending points $\mathsf{p}_j = (1/u_j, 1/v_j, 0)$ in the square
$\mathsf{K} \cap \{z=0\}$. Thus it only remains to check that
$L_{u_{\theta}}(\mathcal{N}) L_{\infty}(\mathcal{M})
L_{v_{\theta}}(\mathcal{N})$ is dense in
$\mathrm{X}_{\theta}(\mathcal{M})$. Applying the reiteration
theorem as we did in the proof of Lemma \ref{Lemma0}, we may
assume that $\mathsf{p}_0$ and $\mathsf{p}_1$ are in the boundary
of $\mathsf{K} \cap \{z=0\}$. We have three possible situations.
First we assume that $\mathsf{p}_0$ and $\mathsf{p}_1$ \emph{live
in the same edge of $\mathsf{K} \cap \{z=0\}$}. Let $\Delta$ be
the intersection of the interpolation pair. In this case, we have
$\Delta = L_{u_0}(\mathcal{N}) L_{\infty}(\mathcal{M})
L_{v_0}(\mathcal{N})$ or $\Delta = L_{u_1}(\mathcal{N})
L_{\infty}(\mathcal{M}) L_{v_1}(\mathcal{N})$ since the points of
any edge of $\mathsf{K} \cap \{z=0\}$ are directed by inclusion.
In particular, we deduce $\Delta \subset L_{u_\theta}(\mathcal{N})
L_{\infty}(\mathcal{M}) L_{v_\theta}(\mathcal{N})$ from which the
result follows. If $\mathsf{p}_0$ and $\mathsf{p}_1$ \emph{live in
consecutive edges of $\mathsf{K} \cap \{z=0\}$}, we have four
choices according to the common vertex $\mathsf{v}$ of the
corresponding (consecutive) edges. Following Figure I we may have
$\mathsf{v} = 0, \mathsf{E}, \mathsf{F}, \mathsf{G}$. When
$\mathsf{v} = \mathsf{E}, \mathsf{G}$ we are back to the situation
above (one endpoint is contained in the other) and there is
nothing to prove. When $\mathsf{v} = 0$, we may assume w.l.o.g.
that
\begin{eqnarray*}
L_{u_0}(\mathcal{N}) L_{\infty}(\mathcal{M}) L_{v_0}(\mathcal{N})
& = & L_{s_0}(\mathcal{N}) L_{\infty}(\mathcal{M}), \\
L_{u_1}(\mathcal{N}) L_{\infty}(\mathcal{M}) L_{v_1}(\mathcal{N})
& = & L_{\infty}(\mathcal{M}) L_{s_1}(\mathcal{N}),
\end{eqnarray*}
for some $2 \le s_0, s_1 \le \infty$. Moreover, we may also assume
w.l.o.g. that $s_0 \le s_1$. This allows us to write $1/s_0 =
1/s_1 + 1/r$ for some index $2 \le r \le \infty$. In particular,
$L_{s_0}(\mathcal{N}) L_{\infty}(\mathcal{M}) = (
L_{s_1}(\mathcal{N}) L_r(\mathcal{N}) ) L_{\infty}(\mathcal{M})$.
Let $x \in \Delta$ so that $$x = \alpha \gamma m_0 = m_1 \beta$$
with $\alpha, \beta \in L_{s_1}(\mathcal{N})$, $\gamma \in
L_r(\mathcal{N})$ and $m_0, m_1 \in \mathcal{M}$. Taking $$a =
\big( \alpha\alpha^* + \beta^*\beta + \delta \mathrm{D}^{2/s_1}
\big)^{1/2},$$ we may write $x=ac_0=c_1a$ (so that
$c_1=ac_0a^{-1}$) with $$c_0=(a^{-1}\alpha\gamma) m_0 \in
L_r(\mathcal{N}) L_{\infty}(\mathcal{M}) \quad \mbox{and} \quad
c_1=m_1 \beta a^{-1} \in L_{\infty}(\mathcal{M}).$$ On the other
hand, $(1/r,0,0)$ and $(0,0,0)$ are in the same edge of
$\mathsf{K} \cap \{z=0\}$. Thus $[L_r(\mathcal{N})
L_{\infty}(\mathcal{M}), L_{\infty}(\mathcal{M})]_{\theta} =
L_{r_{\theta}}(\mathcal{N}) L_{\infty}(\mathcal{M})$ with
$1/r_{\theta} = (1-\theta)/r$. Using this interpolation result and
arguing as in the proof of Lemma \ref{Lemma0} for $\min(q_0,q_1) =
\infty$, we easily obtain that $a^{\theta} c_0 a^{-\theta} \in
L_{r_{\theta}}(\mathcal{N}) L_{\infty}(\mathcal{M})$. Thus we
deduce $$x = a^{1-\theta} a^{\theta} c_0 a^{-\theta} a^{\theta}
\in (L_{s_1/(1-\theta)}(\mathcal{N}) L_{r_{\theta}}(\mathcal{N}))
L_{\infty}(\mathcal{M}) L_{s_1/\theta}(\mathcal{N}).$$ Then, since
the latter space is $L_{u_{\theta}}(\mathcal{N})
L_{\infty}(\mathcal{M}) L_{v_{\theta}}(\mathcal{N})$, we have seen
that $\Delta$ is included in this space. This completes the proof
for $\mathsf{v} =0$. When $\mathsf{v} = \mathsf{F}$, we may assume
w.l.o.g. that
\begin{eqnarray*}
L_{u_0}(\mathcal{N}) L_{\infty}(\mathcal{M}) L_{v_0}(\mathcal{N})
& = & L_{s_0}(\mathcal{N}) L_{\infty}(\mathcal{M})
L_2(\mathcal{N}), \\ L_{u_1}(\mathcal{N}) L_{\infty}(\mathcal{M})
L_{v_1}(\mathcal{N}) & = & L_2(\mathcal{N})
L_{\infty}(\mathcal{M}) L_{s_1}(\mathcal{N}).
\end{eqnarray*}
Writing $1/2 = 1/s_0 + 1/ r_0 = 1/s_1 + 1/ r_1$ for some $2 \le
r_0,r_1 \le \infty$ we have
\begin{eqnarray*}
L_{s_0}(\mathcal{N}) L_{\infty}(\mathcal{M}) L_2(\mathcal{N}) & =
& L_{s_0}(\mathcal{N}) L_{\infty}(\mathcal{M})
L_{r_1}(\mathcal{N}) L_{s_1}(\mathcal{N}),
\\ L_2(\mathcal{N}) L_{\infty}(\mathcal{M}) L_{s_1}(\mathcal{N}) &
= & L_{s_0}(\mathcal{N}) L_{r_0}(\mathcal{N})
L_{\infty}(\mathcal{M}) L_{s_1}(\mathcal{N}).
\end{eqnarray*}
Thus, using our result for $\mathsf{v}=0$ we find
\begin{eqnarray*}
\Delta & = & L_{s_0}(\mathcal{N}) \big( L_{\infty}(\mathcal{M})
L_{r_1}(\mathcal{N}) \cap L_{r_0}(\mathcal{N})
L_{\infty}(\mathcal{M}) \big) L_{s_1}(\mathcal{N}) \\ & \subset &
L_{s_0}(\mathcal{N}) \big( L_{r_0/\theta}(\mathcal{N})
L_{\infty}(\mathcal{M}) L_{r_1/(1-\theta)}(\mathcal{N}) \big)
L_{s_1}(\mathcal{N}) \\ & \subset & L_{s_0(\theta)}(\mathcal{N})
L_{\infty}(\mathcal{M}) L_{s_1(\theta)}(\mathcal{N}),
\end{eqnarray*}
with $1/s_0(\theta)= (1-\theta)/s_0 + \theta/2$ and $1/s_1(\theta)
= (1-\theta)/2 + \theta/s_1$. This completes the proof for
consecutive edges. Finally, we assume that $\mathsf{p}_0$ and
$\mathsf{p}_1$ \emph{live in opposite edges of $\mathsf{K} \cap
\{z=0\}$}. Since the two possible situations are symmetric, we
only consider the case
\begin{eqnarray*}
L_{u_0}(\mathcal{N}) L_{\infty}(\mathcal{M}) L_{v_0}(\mathcal{N})
& = & L_{s_0}(\mathcal{N}) L_{\infty}(\mathcal{M}),
\\ L_{u_1}(\mathcal{N}) L_{\infty}(\mathcal{M})
L_{v_1}(\mathcal{N}) & = & L_{s_1}(\mathcal{N})
L_{\infty}(\mathcal{M}) L_2(\mathcal{N}).
\end{eqnarray*}
If $s_0 \ge s_1$ we clearly have $$\Delta = L_{u_0}(\mathcal{N})
L_{\infty}(\mathcal{M}) L_{v_0}(\mathcal{N}) \subset
L_{u_\theta}(\mathcal{N}) L_{\infty}(\mathcal{M})
L_{v_\theta}(\mathcal{N})$$ and there is nothing to prove. If $s_0
< s_1$ we have $1/s_0 = 1/s_1 + 1/r$ so that
\begin{eqnarray*}
\Delta & = & L_{s_1}(\mathcal{N}) \big( L_r(\mathcal{N})
L_{\infty}(\mathcal{M}) \cap L_{\infty}(\mathcal{M})
L_2(\mathcal{N}) \big) \\ & \subset & L_{s_1}(\mathcal{N}) \big(
L_{r/(1-\theta)}(\mathcal{N}) L_{\infty}(\mathcal{M})
L_{2/\theta}(\mathcal{N}) \big) \\ & = & L_{u_\theta}(\mathcal{N})
L_{\infty}(\mathcal{M}) L_{v_\theta}(\mathcal{N}).
\end{eqnarray*}
This proves the assertion for opposite edges and so the space
$L_{u_\theta}(\mathcal{N}) L_{\infty}(\mathcal{M})
L_{v_\theta}(\mathcal{N})$ is always dense in
$\mathrm{X}_{\theta}(\mathcal{M})$. The proof of Theorem
\ref{TheoremA1} is completed. \fin

\begin{remark}
\emph{The key points to see that the proof of Lemma \ref{Lemma0}
applies whenever we start with any two ending points
$(1/u_j,1/v_j,1/q_j)$ lying on $\mathsf{K}$ are the following:}
\begin{itemize}
\item[(a)] \emph{$L_{u_j}(\mathcal{N}) L_{q_j}(\mathcal{M})
L_{v_j}(\mathcal{N})$ is a Banach space.} \item[(b)] \emph{Lemma
\ref{Lemma-Triangle-Inequality} holds on $\mathsf{K}$ for finite
von Neumann algebras.} \item[(c)] \emph{Corollary
\ref{Corollary-Contraction-Expectation} also holds with ending
points in $\mathsf{K} \setminus
\partial_{\infty} \mathsf{K}$.}
\end{itemize}
\emph{Thus, it suffices to see that Lemma \ref{Lemma0} gives (a),
(b), (c) by complex interpolation. On the other hand, it is worthy
to explain with some more details why restriction $1/u_0 + 1/q_0 +
1/v_0 = 1/u_1 + 1/q_1 + 1/v_1$ can be dropped. The only two points
where this restriction is needed (apart from the case
$\min(q_0,q_1)=\infty$ which has been discussed in Step 4 above)
are in the proofs of Lemma \ref{Lemma-Interpolation-Partial0K-2}
and Proposition \ref{Proposition-Contraction-Expectation}.
However, Lemma \ref{Lemma-Interpolation-Partial0K-2} is only
needed to obtain Corollary \ref{Corollary-Expectation} (which we
have \emph{auto-improved} in Step 2 above by using Lemma
\ref{Lemma0}). Moreover, as we also pointed out in Step 2,
Proposition \ref{Proposition-Contraction-Expectation} now follows
from our improvement of Corollary \ref{Corollary-Expectation}.
Therefore, restriction $1/u_0 + 1/q_0 + 1/v_0 = 1/u_1 + 1/q_1 +
1/v_1$ can be ignored.}
\end{remark}

\chapter{Conditional $L_p$ spaces}
\label{Section4}

We conclude the first half of this paper by studying the duals of
amalgamated $L_p$ spaces and the subsequent applications of
Theorem \ref{TheoremA1}. Let us consider a von Neumann algebra
$\mathcal{M}$ equipped with a \emph{n.f.} state $\varphi$ and a
von Neumann subalgebra $\mathcal{N}$ of $\mathcal{M}$. Let
$\mathsf{E}: \mathcal{M} \rightarrow \mathcal{N}$ denote the
corresponding conditional expectation. We consider any three
indices $(u,p,v)$ such that $(1/u,1/v,1/p)$ belongs to
$\mathsf{K}$ and we define $1 \le s \le \infty$ by $1/s = 1/u +
1/p + 1/v$. Then, the \emph{conditional $L_p$ space}
\label{CondLp} $$L_{(u,v)}^p(\mathcal{M}, \mathsf{E})$$ is defined
as the completion of $L_p(\mathcal{M})$ with respect to the norm
$$\|x\|_{L_{(u,v)}^p(\mathcal{M}, \mathsf{E})} = \sup \Big\{
\|\alpha x \beta\|_{L_s(\mathcal{M})} \, \big| \
\|\alpha\|_{L_u(\mathcal{N})}, \|\beta\|_{L_v(\mathcal{N})} \le 1
\Big\}.$$

We shall show below that amalgamated and conditional $L_p$ spaces
are related by duality. According to our main result in Chapter
\ref{Section3}, this immediately provides interpolation isometries
of the form
\begin{equation} \label{Ecuacion-Cond-Int}
\Big[ L_{(u_0,v_0)}^{p_0}(\mathcal{M},\mathsf{E}),
L_{(u_1,v_1)}^{p_1}(\mathcal{M},\mathsf{E}) \Big]_{\theta} =
L_{(u_{\theta},v_{\theta})}^{p_{\theta}}(\mathcal{M}, \mathsf{E}).
\end{equation}
Our aim now is to explore these identities, since they will be
useful in the sequel.

\begin{example} \label{Remark-Conditional-I-II}
\emph{As in Chapter \ref{Section2}, several noncommutative
function spaces arise as particular cases of our notion of
conditional $L_p$ space. Let us mention four particularly relevant
examples:}
\begin{itemize}
\item[(a)] \emph{The noncommutative $L_p$ spaces arise as}
$$L_p(\mathcal{M}) = L_{(\infty,\infty)}^p(\mathcal{M}, \mathsf{E}).$$

\item[(b)] \emph{If $p \ge q$ and $1/r = 1/q - 1/p$, the spaces
$L_p(\mathcal{N}_1; L_q(\mathcal{N}_2))$ arise as
$$L_{(2r,2r)}^p(\mathcal{N}_1 \bar\otimes \mathcal{N}_2, \mathsf{E}),$$
where the conditional expectation $\mathsf{E}: \mathcal{N}_1
\bar\otimes \mathcal{N}_2 \to \mathcal{N}_1$ is $\mathsf{E} =
1_{\mathcal{N}_1} \otimes \varphi_{\mathcal{N}_2}$.}

\vskip3pt

\item[(c)] \emph{If $2 \le p \le \infty$ and $1/p + 1/q = 1/2$,
Lemma \ref{Lemma-Norm-Conditional-Sup} gives
\begin{eqnarray*} L_p^r(\mathcal{M},
\mathsf{E}) & = & L_{(q,\infty)}^p(\mathcal{M}, \mathsf{E}), \\
L_p^c(\mathcal{M}, \mathsf{E}) & = & L_{(\infty,q)}^p(\mathcal{M},
\mathsf{E}).
\end{eqnarray*}
As we shall see below, $L_p(\mathcal{M}; R_p^n)$ and
$L_p(\mathcal{M}; C_p^n)$ are particular cases.}

\vskip3pt

\item[(d)] \emph{In Chapter \ref{Section7} we will also identity
certain asymmetric noncommutative $L_p$ spaces as particular cases
of conditional $L_p$ spaces. We prefer in this case to leave the
details for Chapter \ref{Section7}.}
\end{itemize}
\end{example}

\section{Duality}

Note that given $(1/u,1/v,1/q) \in \mathsf{K}$, we usually take
$1/p = 1/u + 1/q + 1/v$. In the following it will be more
convenient to replace $p$ by $p'$, the index conjugate to $p$. Let
us consider the following restriction of
(\ref{Equation-Indices-1})
\begin{equation} \label{Equation-Indices-3}
1 < q < \infty \quad \mbox{and} \quad 2 < u, v \le \infty \quad
\mbox{and} \quad 0 < \frac{1}{u} + \frac{1}{q} + \frac{1}{v} =
\frac{1}{p'} < 1.
\end{equation}

\begin{theorem} \label{Theorem-Duality-Conditional}
Let $\mathsf{E}: \mathcal{M} \rightarrow \mathcal{N}$ denote the
conditional expectation of $\mathcal{M}$ onto $\mathcal{N}$ and
let $1 < p < \infty$ given by $1/q' = 1/u + 1/p + 1/v$, where the
indices $(u,q,v)$ satisfy $(\ref{Equation-Indices-3})$ and $q'$ is
conjugate to $q$. Then, the following isometric isomorphisms hold
via the anti-linear duality bracket $\langle x,y \rangle =
\mathrm{tr}(x^* y)$ $$\big( L_u(\mathcal{N}) L_q(\mathcal{M})
L_v(\mathcal{N}) \big)^* = L_{(u,v)}^p(\mathcal{M}, \mathsf{E}),$$
$$\big( L_{(u,v)}^p(\mathcal{M}, \mathsf{E}) \big)^* =
L_u(\mathcal{N}) L_q(\mathcal{M}) L_v(\mathcal{N}).$$
\end{theorem}

\begin{proof} Let us consider the map $$\Lambda_p: x \in
L_{(u,v)}^p(\mathcal{M}, \mathsf{E}) \mapsto \mbox{tr} \big( x^*
\cdot \big) \in \big( L_u(\mathcal{N}) L_q(\mathcal{M})
L_v(\mathcal{N}) \big)^*.$$ We first show that $\Lambda_p$ is an
isometry
\begin{eqnarray*}
\big\| \Lambda_p(x) \big\|_{(u \cdot q \cdot v)^*} & = & \sup
\Big\{ \big| \mbox{tr} (x^* y) \big| \, \big| \ \inf_{y = \alpha z
\beta} \|\alpha\|_{L_u(\mathcal{N})} \|z\|_{L_q(\mathcal{M})}
\|\beta\|_{L_v(\mathcal{N})} \le 1 \Big\}
\\ & = & \sup_{\null} \Big\{ \big| \mbox{tr} (\beta x^* \alpha z)
\big| \, \big| \ \|\alpha\|_{L_u(\mathcal{N})},
\|z\|_{L_q(\mathcal{M})}, \|\beta\|_{L_v(\mathcal{N})} \le 1
\Big\} \\ & = & \sup_{\null} \Big\{ \big\| \beta x^* \alpha
\big\|_{L_{q'}(\mathcal{M})} \hskip6pt  \big| \
\|\alpha\|_{L_u(\mathcal{N})} \le 1, \|\beta\|_{L_v(\mathcal{N})}
\le 1 \Big\} \\ & = & \sup_{\null} \Big\{ \big\| \alpha^* x
\beta^* \big\|_{L_{q'}(\mathcal{M})} \, \big| \
\|\alpha\|_{L_u(\mathcal{N})} \le 1, \|\beta\|_{L_v(\mathcal{N})}
\le 1 \Big\} \\ & = & \|x\|_{L_{(u,v)}^p(\mathcal{M},
\mathsf{E})}^{\null}.
\end{eqnarray*}
It remains to see that $\Lambda_p$ is surjective. To that aim we
use again the solid $\mathsf{K}$ in Figure I. Note that, since the
case $u = v = \infty$ is clear, we may assume that $\min(u,v) <
\infty$. In that case any point $(1/u,1/v,1/q)$ with $(u,q,v)$
satisfying (\ref{Equation-Indices-3}) lies in the interior of a
segment $\mathsf{S}$ contained in $\mathsf{K}$ and satisfying
\begin{itemize}
\item[(a)] One end point of $\mathsf{S}$ lies in the open interval
$(0,\mathsf{A})$. \item[(b)] The segment $\mathsf{S}$ belongs to a
plane parallel to $\mathsf{ACDF}$.
\end{itemize}
According to Theorem \ref{TheoremA1}, this means that
\begin{equation} \label{Equation-Reflexivity}
L_u(\mathcal{N}) L_q(\mathcal{M}) L_v(\mathcal{N}) = \Big[
L_{p'}(\mathcal{M}), L_{u_1}(\mathcal{N}) L_{q_1}(\mathcal{M})
L_{v_1}(\mathcal{N}) \Big]_{\theta},
\end{equation}
for some $0 < \theta < 1$ and some $(1/u_1,1/v_1,1/q_1) \in
\mathsf{K}$. Recalling that $1 < p < \infty$, we know that
$L_{p'}(\mathcal{M})$ is reflexive. In particular, the same holds
for the interpolation space in (\ref{Equation-Reflexivity}) and we
obtain the following isometric isomorphism $$\big(
L_u(\mathcal{N}) L_q(\mathcal{M}) L_v(\mathcal{N}) \big)^* = \Big[
L_p(\mathcal{M}), \big( L_{u_1}(\mathcal{N}) L_{q_1}(\mathcal{M})
L_{v_1}(\mathcal{N}) \big)^* \Big]_{\theta}.$$ Moreover, since $0
< \theta < 1$, we know that the intersection $$L_p(\mathcal{M})
\cap \big( L_{u_1}(\mathcal{N}) L_{q_1}(\mathcal{M})
L_{v_1}(\mathcal{N}) \big)^*$$ is norm dense in the space $\big(
L_u(\mathcal{N}) L_q(\mathcal{M}) L_v(\mathcal{N}) \big)^*.$ On
the other hand, recalling that $1/u_1 + 1/q_1 + 1/v_1 = 1/p'$, we
know from the definition of amalgamated spaces that
$$L_{p'}(\mathcal{M}) = L_{p'}(\mathcal{M}) + L_{u_1}(\mathcal{N})
L_{q_1}(\mathcal{M}) L_{v_1}(\mathcal{N}).$$ Hence,
$L_p(\mathcal{M}) = L_p(\mathcal{M}) \, \cap \, \big(
L_{u_1}(\mathcal{N}) L_{q_1}(\mathcal{M}) L_{v_1}(\mathcal{N})
\big)^*$ is norm dense in $$\big( L_u(\mathcal{N})
L_q(\mathcal{M}) L_v(\mathcal{N}) \big)^*.$$ Therefore,
$\Lambda_p$ has dense range since $$L_p(\mathcal{M}) \subset
L_{(u,v)}^p(\mathcal{M}, \mathsf{E}).$$ For the second part, we
use from (\ref{Equation-Reflexivity}) that $L_u(\mathcal{N})
L_q(\mathcal{M}) L_v(\mathcal{N})$ is reflexive.
\end{proof}

\begin{remark} \label{Remark-Isometric-Duality}
\emph{The first part of the proof of Theorem
\ref{Theorem-Duality-Conditional} holds for any point
$(1/u,1/v,1/q)$ in the solid $\mathsf{K}$. In particular, we
always have an isometric embedding $$L_{(u,v)}^p(\mathcal{M},
\mathsf{E}) \longrightarrow \big( L_u(\mathcal{N})
L_q(\mathcal{M}) L_v(\mathcal{N}) \big)^*.$$}
\end{remark}

\begin{remark}
\emph{Note that the indices excluded in Theorem
\ref{Theorem-Duality-Conditional} by the restriction imposed by
property (\ref{Equation-Indices-3}) are the natural ones. For
instance, the last restriction $0 < 1/u + 1/q + 1/v < 1$ only
affects conditional/amalgamated $L_1$ and $L_{\infty}$ spaces,
which are not expected to be reflexive. Moreover, the spaces
$L_u(\mathcal{N}) L_{\infty}(\mathcal{M}) L_v(\mathcal{N})$ are
not reflexive in general. Indeed, let us consider the particular
case in which $\mathcal{N}$ is the complex field. These spaces
collapse into $L_{\infty}(\mathcal{M})$ which is not reflexive.}
\end{remark}

\section{Conditional $L_{\infty}$ spaces}

Among the non-reflexive conditional spaces, we concentrate on some
properties of conditional $L_{\infty}$ spaces that will be needed
in the second half of this paper. Note that, given indices
$(u,q,v)$ satisfying (\ref{Equation-Indices-1}) with $1/u + 1/q +
1/v = 1$, we have defined the space
$$L_{(u,v)}^{\infty}(\mathcal{M}, \mathsf{E})$$ as the completion
of $L_{\infty}(\mathcal{M})$ with respect to the norm
$$\|x\|_{L_{(u,v)}^{\infty}(\mathcal{M}, \mathsf{E})} = \sup
\Big\{ \|\alpha x \beta \|_{L_{q'}(\mathcal{M})} \, \big| \
\|\alpha\|_{L_u(\mathcal{N})}, \|\beta\|_{L_v(\mathcal{N})} \le 1
\Big\}.$$ According to Remark \ref{Remark-Isometric-Duality}, we
know that this space embeds isometrically in \label{GothicL}
$$\mathcal{L}_{(u,v)}^{\infty}(\mathcal{M}, \mathsf{E}) = \big(
L_u(\mathcal{N}) L_q(\mathcal{M}) L_v(\mathcal{N}) \big)^*.$$

\begin{proposition} \label{Proposition-Cond-Linfty}
The following properties hold:
\begin{itemize}
\item[i)] $\mathcal{L}_{(u,v)}^{\infty}(\mathcal{M}, \mathsf{E})$
is contractively included in $L_{q'}(\mathcal{M})$.

\vskip2pt

\item[ii)] $L_{\infty}(\mathcal{M})$ and
$L_{(u,v)}^{\infty}(\mathcal{M}, \mathsf{E})$ are weak$^*$ dense
subspaces of $\mathcal{L}_{(u,v)}^{\infty}(\mathcal{M},
\mathsf{E})$.
\end{itemize}
\end{proposition}

\begin{proof} Let us consider the map $$j: \varphi \in
\mathcal{L}_{(u,v)}^{\infty}(\mathcal{M}, \mathsf{E}) \to \varphi
\big( \mathrm{D}^{\frac{1}{u}} \cdot \mathrm{D}^{\frac{1}{v}}
\big) \in L_{q'}(\mathcal{M}).$$ By Proposition
\ref{Proposition-Isometry-Triangle-Complete}, the map $j$ is
clearly injective. On the other hand,
\begin{eqnarray*}
\lefteqn{\big\| \varphi \big( \mathrm{D}^{\frac{1}{u}} \cdot
\mathrm{D}^{\frac{1}{v}} \big) \big\|_{L_{q'}(\mathcal{M})}} \\ &
= & \sup \Big\{ \big| \varphi \big( \mathrm{D}^{\frac{1}{u}} y
\mathrm{D}^{\frac{1}{v}} \big) \big| \, \big| \ \|y\|_q \le 1
\Big\} \\ & \le & \sup \Big\{ \big| \varphi \big(
\mathrm{D}^{\frac{1}{u}} \alpha y \beta \mathrm{D}^{\frac{1}{v}}
\big) \big| \, \big| \ \big\| \mathrm{D}^{\frac{1}{u}} \alpha
\big\|_{L_u(\mathcal{N})}, \|y\|_q, \big\| \beta
\mathrm{D}^{\frac{1}{v}} \big\|_{L_v(\mathcal{N})} \le 1 \Big\}
\\ & = & \|\varphi\|_{\mathcal{L}_{(u,v)}^{\infty}(\mathcal{M},
\mathsf{E})}.
\end{eqnarray*}
The last identity follows from Proposition
\ref{Proposition-Isometry-Triangle-Complete}. This shows that $j$
is a contraction. For the second part, it suffices to see that
$L_{\infty}(\mathcal{M})$ is weak$^*$ dense. However, it is clear
from the definition of amalgamated spaces that the inclusion map
$$L_u(\mathcal{N}) L_q(\mathcal{M}) L_v(\mathcal{N})
\longrightarrow L_1(\mathcal{M})$$ is injective. In particular,
taking adjoints we get the announced weak$^*$ density. \end{proof}

\section{Interpolation results and applications}

In this last section we consider some interesting particular cases
of the dual version \eqref{Ecuacion-Cond-Int} of Theorem
\ref{TheoremA1}. One of the applications we shall consider
generalizes Pisier's interpolation result \cite{P0} and Xu's
recent extension \cite{X}.

\begin{theorem} \label{Theorem-Interpolation-Conditional1}
Let $\mathcal{N}$ be a von Neumann subalgebra of $\mathcal{M}$ and
let $\mathsf{E}: \mathcal{M} \rightarrow \mathcal{N}$ be the
corresponding conditional expectation. Assume that $(u_j,q_j,v_j)$
satisfy $(\ref{Equation-Indices-3})$ for $j=0,1$ and that $1/u_j +
1/q_j + 1/v_j = 1/p_j'$. Then, if $p_j$ denotes the conjugate
index to $p_j'$, the following isometric isomorphism holds $$\Big[
L_{(u_0,v_0)}^{p_0}(\mathcal{M},\mathsf{E}),
L_{(u_1,v_1)}^{p_1}(\mathcal{M},\mathsf{E}) \Big]_{\theta} =
L_{(u_{\theta},v_{\theta})}^{p_{\theta}}(\mathcal{M},
\mathsf{E}).$$ Moreover, if $2 \le u_j, v_j \le \infty$ for
$j=0,1$, we also have $$\Big[
\mathcal{L}_{(u_0,v_0)}^{\infty}(\mathcal{M},\mathsf{E}),
\mathcal{L}_{(u_1,v_1)}^{\infty}(\mathcal{M},\mathsf{E})
\Big]^{\theta} =
\mathcal{L}_{(u_{\theta},v_{\theta})}^{\infty}(\mathcal{M},
\mathsf{E}).$$
\end{theorem}

\begin{proof} The first part follows automatically from Theorem
\ref{TheoremA1} and Theorem \ref{Theorem-Duality-Conditional}. The
second part follows from Theorem \ref{TheoremA1} and the duality
properties which link the complex interpolation brackets $[ \ , \,
]_{\theta}$ and $[ \ , \, ]^{\theta}$, see e.g. \cite{BL}.
\end{proof}

Now we study some consequences of Theorem
\ref{Theorem-Interpolation-Conditional1}. We shall content
ourselves by exploring only the case $p_0=p_1$. This restriction
is motivated by the applications we are using in the successive
chapters. The last part of the following result requires to
introduce some notation. As usual, we shall write $R_p^n$ (resp.
$C_p^n$) to denote the interpolation space $[R_n, C_n]_{1/p}$
(resp. $[C_n,R_n]_{1/p}$), where $R_n$ and $C_n$ denote the
$n$-dimensional row and column Hilbert spaces. Alternatively, we
may define $R_p^n$ and $C_p^n$ as the first row and column
subspaces of the Schatten class $S_p^n$. On the other hand, given
an element $x_0$ in a von Neumann algebra $\mathcal{M}$, we shall
consider the mappings $L_{x_0}$ and $R_{x_0}$ on $\mathcal{M}$
defined respectively as follows $$L_{x_0}(x) = x_0 x \qquad
\mbox{and} \qquad R_{x_0}(x) = x x_0.$$

\begin{corollary} \label{Corollary-Interpolation-Conditional2}
Let $\mathcal{N}$ be a von Neumann subalgebra of $\mathcal{M}$ and
let $\mathsf{E}: \mathcal{M} \rightarrow \mathcal{N}$ be the
corresponding conditional expectation of $\mathcal{M}$ onto
$\mathcal{N}$. Then, we have the following isometric isomorphisms:
\begin{itemize}
\item[i)] If $2 \le p < \infty$ and $2 < q \le \infty$ are such
that $1/2 = 1/p + 1/q$, we have $$\begin{array}{rcl} \big[
L_p(\mathcal{M}), L_p^r(\mathcal{M}, \mathsf{E}) \big]_{\theta} &
= & L_{(s, \infty)}^p(\mathcal{M}, \mathsf{E}), \\ \big[
L_p(\mathcal{M}), L_p^c(\mathcal{M}, \mathsf{E}) \big]_{\theta} &
= & L_{(\infty,s)}^p(\mathcal{M}, \mathsf{E}), \end{array} \quad
\mbox{with} \quad \frac{1}{s} = \frac{\theta}{q}.$$ In the case
$p=\infty$, we obtain $$\begin{array}{rcl} \big[
L_{\infty}(\mathcal{M}), L_{\infty}^r(\mathcal{M}, \mathsf{E})
\big]_{\theta} & = & L_{(2/\theta, \infty)}^{\infty}(\mathcal{M},
\mathsf{E}), \\ \big[ L_{\infty}(\mathcal{M}),
L_{\infty}^c(\mathcal{M}, \mathsf{E}) \big]_{\theta} & = &
L_{(\infty,2/\theta)}^{\infty}(\mathcal{M}, \mathsf{E}).
\end{array}$$
\item[ii)] If $2 \le p < \infty$ and $2 < q \le \infty$ are such
that $1/2 = 1/p + 1/q$, we have $$\quad \big[
L_p^c(\mathcal{M},\mathsf{E}),L_p^r(\mathcal{M}, \mathsf{E})
\big]_{\theta} = L_{(u,v)}^p(\mathcal{M}, \mathsf{E})$$ with
$$(1/u,1/v) = (\theta/q, (1-\theta)/q).$$ \item[iii)] Let us define
for $2 \le p < \infty$
$$\mathrm{X}_{\theta}(\mathcal{M}) = \big[ L_p(\mathcal{M}; C_p^n),
L_p(\mathcal{M};R_p^n) \big]_{\theta}.$$ Then, if $1/w = 1/p +
1/v$ with $1/v = (1 - \theta)/q$, we deduce $$\Big\| \sum_{k=1}^n
x_k \otimes \delta_k \Big\|_{\mathrm{X}_{\theta}(\mathcal{M})}^2 =
\Big\| \sum_{k=1}^n L_{x_k} R_{x_k^*}: L_{v/2}(\mathcal{M})
\rightarrow L_{w/2}(\mathcal{M}) \Big\|.$$
\end{itemize}
\end{corollary}

\begin{proof} The assertions in the first part of i) and ii) follow from
Theorem \ref{Theorem-Interpolation-Conditional1} after the obvious
identifications, see Example \ref{Remark-Conditional-I-II}. For
the last part of i) we only prove the first identity since the
second one follows in the same way. According to Remark
\ref{Remark-Isometric-Duality}, $L_{\infty}^r(\mathcal{M},
\mathsf{E})$ is the closure of $L_{\infty}(\mathcal{M})$ in
$\mathcal{L}^{\infty}_{(2,\infty)} (\mathcal{M}, \mathsf{E})$.
Then, applying a well-known property of the complex method (see
e.g. \cite[Theorem 4.2.2]{BL}) we find
$$\big[ L_\infty(\mathcal{M}), L_{\infty}^r(\mathcal{M},
\mathsf{E}) \big]_{\theta} = \big[ L_\infty(\mathcal{M}),
\mathcal{L}^{\infty}_{(2,\infty)}(\mathcal{M}, \mathsf{E})
\big]_{\theta} \quad \mbox{for} \quad 0 < \theta < 1.$$ By Berg's
theorem, we have an isometric inclusion $$\big[
L_\infty(\mathcal{M}),
\mathcal{L}^{\infty}_{(2,\infty)}(\mathcal{M}, \mathsf{E})
\big]_{\theta} \subset \big[ L_\infty(\mathcal{M}),
\mathcal{L}^{\infty}_{(2,\infty)}(\mathcal{M}, \mathsf{E})
\big]^{\theta}.$$ Therefore given $x \in L_\infty(\mathcal{M})$,
Theorem \ref{Theorem-Interpolation-Conditional1} gives
\begin{eqnarray*}
\|x\|_{[L_\infty(\mathcal{M}),L_{\infty}^r(\mathcal{M},
\mathsf{E})]_{\theta}} & = & \|x\|_{[L_\infty(\mathcal{M}),
\mathcal{L}^{\infty}_{(2,\infty)}(\mathcal{M},
\mathsf{E})]_{\theta}} \\ & = & \|x\|_{[L_\infty(\mathcal{M}),
\mathcal{L}^{\infty}_{(2,\infty)}(\mathcal{M},
\mathsf{E})]^{\theta}} \\ & = &
\|x\|_{L_{(2/\theta,\infty)}^{\infty}(\mathcal{M}, \mathsf{E})}.
\end{eqnarray*}
The assertion then follows by a simple density argument $$\big[
L_\infty(\mathcal{M}),L_{\infty}^r(\mathcal{M}, \mathsf{E})
\big]_{\theta} = L_{(2/\theta,\infty)}^{\infty}(\mathcal{M},
\mathsf{E}).$$ Finally, for part iii) we consider the direct sum
$\mathcal{M}_{\oplus n} = \mathcal{M} \oplus \mathcal{M} \oplus
\cdots \oplus \mathcal{M}$ with $n$ terms and equipped with the
\emph{n.f.} state $\varphi_n (x_1, x_2, \cdots, x_n) = \frac{1}{n}
\sum_k \varphi(x_k)$. The natural conditional expectation is given
by $$\mathsf{E}_n: \sum_{k=1}^n x_k \otimes \delta_k \in
\mathcal{M}_{\oplus n} \mapsto \frac{1}{n} \sum_{k=1}^n x_k \in
\mathcal{M}.$$ It is clear that we have the isometries
\begin{eqnarray*}
L_p(\mathcal{M}; R_p^n) & = & \sqrt{n} \,
L_p^r(\mathcal{M}_{\oplus n}, \mathsf{E}_n), \\ L_p(\mathcal{M};
C_p^n) & = & \sqrt{n} \, L_p^c(\mathcal{M}_{\oplus n},
\mathsf{E}_n).
\end{eqnarray*}
According to ii) and the definition of the norm in
$L_{(u,v)}^p(\mathcal{M}_{\oplus n}, \mathsf{E}_n)$, we have
\begin{eqnarray*}
\lefteqn{\Big\| \sum_{k=1}^n x_k \otimes \delta_k
\Big\|_{\mathrm{X}_{\theta}(\mathcal{M})}^2}
\\ & = & \sup \left\{ n \Big\| \sum_{k=1}^n \alpha x_k
\beta \otimes \delta_k \Big\|_{L_2(\mathcal{M}_{\oplus n})}^2 \,
\big|
\ \|\alpha\|_{L_u(\EME)}, \|\beta\|_{L_v(\EME)} \le 1 \right\} \\
& = & \sup \left\{ \sum_{k=1}^n \mbox{tr} \big(\alpha x_k \beta
\beta^* x_k^* \alpha^* \big) \, \big| \ \|\alpha\|_{L_u(\EME)},
\|\beta\|_{L_v(\EME)} \le 1 \right\} \\ & = & \sup \left\{ \Big\|
\sum_{k=1}^n L_{x_k} R_{x_k^*} (\beta \beta^*)
\Big\|_{L_{(u/2)'}(\mathcal{M})} \, \big| \ \|\beta\|_{L_v(\EME)}
\le 1 \right\} \\ & = & \sup \left\{ \Big\| \sum_{k=1}^n L_{x_k}
R_{x_k^*} (\gamma) \Big\|_{L_{(u/2)'}(\mathcal{M})} \, \big| \
\gamma \ge 0, \ \|\gamma\|_{L_{v/2}(\EME)} \le 1 \right\}.
\end{eqnarray*}
Recalling that $\summ_k L_{x_k} R_{x_k^*}$ is positive and that
$(u/2)' = w/2$, we conclude. \end{proof}

\begin{observation} \label{Observation-Extension}
\emph{Arguing as in Corollary
\ref{Corollary-Interpolation-Conditional2}, we easily obtain}
$$\big[ L_{\infty}(\mathcal{M}), L_{(u,v)}^{\infty}(\mathcal{M},
\mathsf{E}) \big]_{\theta} = L_{(u/\theta,
v/\theta)}^{\infty}(\mathcal{M}, \mathsf{E}).$$ \emph{These
results shall be frequently used in the successive chapters.
Moreover, let us also mention that Corollary
\ref{Corollary-Interpolation-Conditional2} i) is needed in
\cite{JM} to study the noncommutative John-Nirenberg theorem.}
\end{observation}

At the time of this writing we do not know whether Corollary
\ref{Corollary-Interpolation-Conditional2} ii) extends to
$p=\infty$ in full generality. We will now show that the equality
holds when restricted to elements in $\mathcal{M}$. This result
will play a very important role in the sequel.

\begin{lemma} \label{Lemma-Last}
If $1 < p,q < \infty$ and $z \in \mathcal{M}$, we have $$\inf
\Big\{ \|x\|_{L^{\infty}_{(2p,\infty)}(\mathcal{M}, \mathsf{E})}
\, \|y\|_{L^{\infty}_{(\infty,2q)} (\mathcal{M}, \mathsf{E})} \,
\big| \ z = xy, \ x,y \in \mathcal{M} \Big\} \le
\|z\|_{L^{\infty}_{(2p,2q)}(\mathcal{M}, \mathsf{E})}.$$
\end{lemma}

\begin{proof} On $\mathcal{M}$ we define the norm $$\|z\|_h = \inf_{z=xy}
\|x\|_{L^{\infty}_{(2p,\infty)} (\mathcal{M}, \mathsf{E})} \,
\|y\|_{L^{\infty}_{(\infty,2q)} (\mathcal{M}, \mathsf{E})}$$ where
the infimum is taken over $x,y \in L_{\infty}(\mathcal{M})$. It is
easy to check that we do not need to consider sums here. Let us
assume that $\|z\|_h=1$. By the Hahn-Banach theorem there exists a
linear functional $\phi: \mathcal{M} \to \C$ such that $\phi(z)=1$
and
\begin{equation} \label{pp}
\Big| \summ_k \phi(x_k y_k) \Big| \le \sup_{\|a\|_{2p} \le 1}
\Big\| \summ_k a x_k x_k^* a^* \Big\|_p^{1/2} \sup_{\|b\|_{2q} \le
1} \Big\| \summ_k b^* y_k^*y_k b \Big\|_q^{1/2}.
\end{equation}
Note that $\|z\|_h \le \|z\|_{\infty}$ and thus $\phi$ is
continuous. It follows immediately that we may move the absolute
values in (\ref{pp}) inside. Thus, we get $$\summ_k |\phi(x_k
y_k)| \le \sup_{a,c} \, \mbox{tr} \Big( \summ_k a x_k x_k^* a^* c
\Big)^{1/2} \, \sup_{b,d} \, \mbox{tr} \Big(\summ_k b^* y_k^* y_k
b d \Big)^{1/2}.$$ Here we take the supremum over $(a,c,b,d)$ in
$$\mathsf{B}_{L_{2p}(\mathcal{N})} \times
\mathsf{B}_{L_{p'}(\mathcal{M})}^+ \times
\mathsf{B}_{L_{2q}(\mathcal{N})} \times
\mathsf{B}_{L_{q'}(\mathcal{M})}^+,$$ all equipped with the
weak$^*$ topology. Using the standard Grothendieck-Pietsch
separation argument as in Theorem \ref{Theorem-Grothendieck} we
obtain two probability measures $\mu_1$ and $\mu_2$ such that
$$|\phi(xy)| \le \Big( \int \mbox{tr}(axx^*a^*c) \, d\mu_1(a,c)
\Big)^{1/2} \Big( \int \mbox{tr}(b^*y^*ybd) \, d\mu_2(b,d)
\Big)^{1/2}.$$ Since $L_{2p}(\mathcal{N}) L_{p'}(\mathcal{M})
L_{2p}(\mathcal{N})$ and $L_{2q}(\mathcal{N}) L_{q'}(\mathcal{M})
L_{2q}(\mathcal{N})$ are Banach spaces, $$\alpha = \int a^*ca \,
d\mu_1(a,c) \quad \mbox{and} \quad \beta = \int bdb^* \,
d\mu_2(b,d)$$ are positive elements respectively in the unit balls
of $$L_{2p}(\mathcal{N}) L_{p'}(\mathcal{M}) L_{2p}(\mathcal{N})
\quad \mbox{and} \quad L_{2q}(\mathcal{N}) L_{q'}(\mathcal{M})
L_{2q}(\mathcal{N}).$$ Therefore we find $a_1,a_2 \in
\mathsf{B}_{L_{2p}(\mathcal{N})}$ and $c_1,c_2 \in
\mathsf{B}_{L_{2p'}(\mathcal{M})}$ such that $\alpha =
a_1c_1c_2a_2$. We deduce from the Cauchy-Schwartz inequality and
the arithmetic-geometric mean inequality that
\begin{eqnarray*}
\mbox{tr}(xx^*\alpha) & = & \big| \mbox{tr}(xx^*a_1c_1c_2a_2)
\big| \\ & = & \big| \mbox{tr}(c_2a_2xx^*a_1c_1) \big| \\ & \le &
\mbox{tr} \big(c_2a_2xx^*a_2^*c_2^* \big)^{1/2} \mbox{tr} \big(
c_1^*a_1^*xx^*a_1c_1 \big)^{1/2} \\ & \le & \mbox{tr} \Big( xx^*
\frac{a_1c_1c_1^*a_1^*+a_2^*c_2^*c_2a_2}{2} \Big).
\end{eqnarray*}
We could consider $a = (a_1^*a_1+a_2^*a_2)^{1/2}$ to deduce that
$$a^{-1} \frac{a_1c_1c_1^*a_1^*+a_2^*c_2^*c_2a_2}{2}a^{-1}$$ is a
positive element in $L_{p'}(\mathcal{M})$ of norm $\le 1$. This is
not enough for our purposes. However, following the proof of the
triangle inequality in Lemma \ref{Lemma-Triangle-Inequality}, we
may apply Devinatz's theorem one more time to find an operator $a
\in (1+\varepsilon) \mathsf{B}_{L_{2p}(\mathcal{N})}$ with full
support and $c \in \mathsf{B}_{L_{p'}(\mathcal{M})}$ such that
$$\frac{a_1c_1c_1^*a_1^*+a_2^*c_2^*c_2a_2}{2} = a^*ca.$$ We leave
the details to the interested reader. This implies that $c$ is
positive and $$\mbox{tr} (xx^*\alpha) \le \mbox{tr}(xx^*a^*ca) =
\|c^{1/2}ax\|_2^2.$$ The same argument for $\beta$ gives $b \in
(1+\varepsilon) \mathsf{B}_{L_{2q}(\mathcal{N})}$ and $d \in
\mathsf{B}_{L_{q'}(\mathcal{M})}^+$ such that
$$\mbox{tr}(y^*y\beta) \le \mbox{tr}(y^*ybdb^*) =
\|ybd^{1/2}\|_2^2.$$ This yields $$|\phi(xy)| \le \|c^{1/2}ax\|_2
\|ybd^{1/2}\|_2.$$ From this it is easy to find a contraction
$u\in \mathcal{M}$ such that $$\phi(xy) =
\mbox{tr}(uc^{\frac12}axybd^{\frac12}).$$ If $1/r = 1/2p + 1/2q$,
we deduce from H\"{o}lder's inequality that
\begin{eqnarray*}
\|z\|_h = |\phi(z)| = \big| \mbox{tr}(uc^{\frac12}azbd^{\frac12})
\big| \le \|azb\|_{L_r(\mathcal{M})} \le (1+\varepsilon)^2
\|z\|_{L_{(2p,2q)}^{\infty}(\mathcal{M}, \mathsf{E})}.
\end{eqnarray*} Finally, recalling that $\varepsilon>0$ is
arbitrary, the assertion follows by taking $\varepsilon \to 0$.
\end{proof}

\begin{corollary} \label{Corollary-Pisier-Conditional}
Assume that $1 \le p_0,p_1,q_0,q_1 \le \infty$ satisfy
$$\max(p_0,p_1), \max (q_0, q_1) > 1 \quad \mbox{and} \quad
\min(p_0,p_1), \min(q_0,q_1) < \infty.$$ Then, given $x \in
\mathcal{M}$ we have
$$\|x\|_{L^{\infty}_{(2p_\theta,2q_{\theta})}(\mathcal{M},
\mathsf{E})} = \|x\|_{[L^{\infty}_{(2p_0,2q_0)} (\mathcal{M},
\mathsf{E}), L^{\infty}_{(2p_1,2q_1)} (\mathcal{M},
\mathsf{E})]_{\theta}}.$$ In particular, for $\theta = 1/q$ we
obtain $$\|x\|_{[L_\infty^c(\mathcal{M},
\mathsf{E}),L_{\infty}^r(\mathcal{M}, \mathsf{E})]_{\theta}} =
\sup \Big\{ \|axb\|_{L_2(\mathcal{M})} \, \big| \
\|a\|_{L_{2q}(\mathcal{N})}, \, \|b\|_{L_{2q'}(\mathcal{N})} \le 1
\Big\}.$$
\end{corollary}

\begin{proof} The upper estimate is an easy application of trilinear
interpolation. For the converse we apply Lemma \ref{Lemma-Last} so
that for any $\varepsilon > 0$ we can always find a factorization
$x=x_1x_2$ satisfying $$\|x_1\|_{L^{\infty}_{(2p_{\theta},\infty)}
(\mathcal{M}, \mathsf{E})}
\|x_2\|_{L^{\infty}_{(\infty,2q_{\theta})} (\mathcal{M},
\mathsf{E})} \le (1+\varepsilon)
\|x\|_{L^{\infty}_{(2p_{\theta},2q_{\theta})}(\mathcal{M},
\mathsf{E})}.$$ According to Corollary
\ref{Corollary-Interpolation-Conditional2} i) we know that
$$L^{\infty}_{(2p_{\theta},\infty)} (\mathcal{M},
\mathsf{E}) = \big[ L_{\infty} (\mathcal{M}),
L_{\infty}^r(\mathcal{M}, \mathsf{E}) \big]_{1/p_{\theta}}.$$
Taking $\theta_0 = 1/p_0$ and $\theta_1 = 1/p_1$, the reiteration
theorem implies that
$$\big[ L_{\infty}(\mathcal{M}), L_{\infty}^r(\mathcal{M},
\mathsf{E}) \big]_{1/p_{\theta}} = \big[ [L_{\infty}(\mathcal{M}),
L_{\infty}^r(\mathcal{M}, \mathsf{E})]_{\theta_0},
[L_{\infty}(\mathcal{M}), L_{\infty}^r(\mathcal{M},
\mathsf{E})]_{\theta_1} \big]_{\theta}.$$ In particular, $$\big[
L_{\infty}(\mathcal{M}), L_{\infty}^r(\mathcal{M}, \mathsf{E})
\big]_{1/p_{\theta}} = \big[ L^{\infty}_{(2p_0,\infty)}
(\mathcal{M}, \mathsf{E}), L^{\infty}_{(2p_1,\infty)}
(\mathcal{M}, \mathsf{E}) \big]_{\theta}.$$ Therefore, we get
$$\|x_1\|_{[L^{\infty}_{(2p_0,\infty)} (\mathcal{M}, \mathsf{E}),
L^{\infty}_{(2p_1,\infty)} (\mathcal{M}, \mathsf{E})]_{\theta}} =
\|x_1\|_{L^{\infty}_{(2p_{\theta},\infty)} (\mathcal{M},
\mathsf{E})}.$$ Similarly, we have
$$\|x_2\|_{[L^{\infty}_{(\infty,2q_0)} (\mathcal{M}, \mathsf{E}),
L^{\infty}_{(\infty,2q_1)} (\mathcal{M}, \mathsf{E})]_{\theta}} =
\|x_2\|_{L^{\infty}_{(\infty,2q_{\theta})}(\mathcal{M},
\mathsf{E})} .$$ On the other hand the inequality
$$\|xy\|_{L^{\infty}_{(2p,2q)} (\mathcal{M}, \mathsf{E})} \le
\|x\|_{L^{\infty}_{(2p,\infty)} (\mathcal{M}, \mathsf{E})}
\|y\|_{L^{\infty}_{(\infty,2q)} (\mathcal{M}, \mathsf{E})}$$ holds
for all $1 \le p,q \le \infty$. Thus, by bilinear interpolation we
deduce
\begin{eqnarray*}
\lefteqn{\|x\|_{[L^{\infty}_{(2p_0,2q_0)} (\mathcal{M},
\mathsf{E}), L^{\infty}_{(2p_1,2q_1)} (\mathcal{M},
\mathsf{E})]_{\theta}}} \\ & \le &
\|x_1\|_{[L^{\infty}_{(2p_0,\infty)} (\mathcal{M}, \mathsf{E}),
L^{\infty}_{(2p_1,\infty)} (\mathcal{M}, \mathsf{E})]_{\theta}}
\|x_2\|_{[L^{\infty}_{(\infty,2q_0)} (\mathcal{M}, \mathsf{E}),
L^{\infty}_{(\infty,2q_1)} (\mathcal{M}, \mathsf{E})]_{\theta}}.
\end{eqnarray*}
Combining the previous estimates and taking $\varepsilon \to 0$ we
obtain the assertion. \end{proof}

\begin{remark} \label{Remark-Density-Discussion}
\emph{If $\mathcal{M}$ is dense in the intersection
$L_{\infty}^r(\mathcal{M}, \mathsf{E}) \cap
L_{\infty}^c(\mathcal{M}, \mathsf{E})$, then Corollary
\ref{Corollary-Interpolation-Conditional2} ii) extends to
$p=\infty$. This holds for instance when the inclusion
$\mathcal{N} \subset \mathcal{M}$ has finite index or when
$\mathcal{M} = \mathcal{A} \bar{\otimes} \mathcal{N}$ with
$\mathcal{A}$ finite-dimensional. However, at the time of this
writing it is not clear whether the density assumption is
satisfied for general conditional expectations. On the other hand,
we note that the validity of Corollary
\ref{Corollary-Interpolation-Conditional2} ii) for $p=\infty$ can
be understood as a \emph{conditional version} of Pisier's
interpolation result \cite{P0}. Indeed, taking $\theta = 1/q$ we
easily find that
\begin{eqnarray*}
\|x\|_{[L_{\infty}^c(\mathcal{M}, \mathsf{E}),
L_{\infty}^r(\mathcal{M}, \mathsf{E})]_{\theta}} & = & \sup \Big\{
\|x^*ax\|_{L_q(\mathcal{M})} \, \big| \ \|a\|_{L_q(\mathcal{N})}
\le 1 \Big\} \\ & = & \sup \Big\{ \|xbx^*\|_{L_{q'}(\mathcal{M})}
\, \big| \ \|b\|_{L_{q'}(\mathcal{N})} \le 1 \Big\}.
\end{eqnarray*}
Furthermore, according to Theorem
\ref{Theorem-Interpolation-Conditional1}, this also applies for $2
\le p \le \infty$ (we leave the details to the reader). In
particular, we also find a conditional version of Xu's
interpolation result \cite{X}. Moreover, when $1 \le p \le 2$ this
result follows from Theorem \ref{TheoremA1} instead of Theorem
\ref{Theorem-Interpolation-Conditional1}.}
\end{remark}

\chapter{Intersections of $L_p$ spaces}
\label{Section5}

The rest of this paper is devoted to study intersections of
certain generalized $L_p$ spaces. Our main goal is to prove a
noncommutative version of $(\Sigma_{pq})$, see the Introduction.
In this chapter we begin by proving certain interpolation results
for intersections. As usual we consider a von Neumann subalgebra
$\mathcal{N}$ of $\mathcal{M}$ with corresponding conditional
expectation $\mathsf{E}: \mathcal{M} \to \mathcal{N}$. Then given
a positive integer $n$, $1 \le q \le p \le \infty$ and $1/r =
1/q-1/p$, we define the following intersection spaces
\begin{eqnarray*}
\mathcal{R}_{2p,q}^n(\mathcal{M},\mathsf{E}) & = &
n^{\frac{1}{2p}} L_{2p}(\mathcal{M}) \cap n^{\frac{1}{2q}}
L_{(2r,\infty)}^{2p}(\mathcal{M}, \mathsf{E}), \\
\mathcal{C}_{2p,q}^n \, (\mathcal{M},\mathsf{E}) & = &
n^{\frac{1}{2p}} L_{2p}(\mathcal{M}) \cap n^{\frac{1}{2q}}
L_{(\infty,2r)}^{2p}(\mathcal{M}, \mathsf{E}).
\end{eqnarray*}
Our main result in this chapter shows that the two families of
intersection spaces considered above are interpolation scales in
the index $q$. That is, the intersections commute with the complex
interpolation functor. Indeed, we obtain the following
isomorphisms with relevant constants independent on $n$
\begin{equation} \label{Eq-Int-Int1}
\begin{array}{rcl}
\big[ \mathcal{R}_{2p,1}^n(\mathcal{M},\mathsf{E}),
\mathcal{R}_{2p,p}^n(\mathcal{M},\mathsf{E}) \big]_{\theta} &
\simeq & \mathcal{R}_{2p,q}^n(\mathcal{M},\mathsf{E}), \\ [5pt]
\big[ \hskip1.5pt \mathcal{C}_{2p,1}^n \,
(\mathcal{M},\mathsf{E}), \, \mathcal{C}_{2p,p}^n \,
(\mathcal{M},\mathsf{E}) \big]_{\theta} & \simeq & \,
\mathcal{C}_{2p,q}^n \, (\mathcal{M},\mathsf{E}),
\end{array}
\end{equation}
and with $1/q = 1-\theta + \theta/p$. Moreover, we shall also
prove that
\begin{equation} \label{Eq-Int-Int2}
\big[ \mathcal{R}_{2p,1}^n(\mathcal{M},\mathsf{E}),
\mathcal{C}_{2p,1}^n(\mathcal{M},\mathsf{E}) \big]_{\theta} \
\simeq \bigcap_{u,v \in \{ 2p', \infty \}} n^{\frac{1-\theta}{u} +
\frac{1}{2p} + \frac{\theta}{v}} \,
L_{(\frac{u}{1-\theta},\frac{v}{\theta})}^{2p} (\mathcal{M},
\mathsf{E}).
\end{equation}

\section{Free Rosenthal inequalities}

Our aim in this section is to present the free analogue given in
\cite{JPX} of Rosenthal inequalities \cite{Ro0}, where mean-zero
independent random variables are replaced by free random
variables. This will be one of the key tools needed for the proof
of the isomorphisms \eqref{Eq-Int-Int1} and \eqref{Eq-Int-Int2}.
For the sake of completeness we first recall the construction of
reduced amalgamated free products.

\subsection{Amalgamated free products}
\label{Subsubsection5.1.1}

The basics of free products without amalgamation can be found in
\cite{VDN}. Let $\mathcal{A}$ be a von Neumann algebra equipped
with a \emph{n.f.} state $\phi$ and let $\mathcal{N}$ be a von
Neumann subalgebra of $\mathcal{A}$. Let
$\mathsf{E}_{\mathcal{N}}: \mathcal{A} \rightarrow \mathcal{N}$ be
the corresponding conditional expectation onto $\mathcal{N}$. A
family $\mathsf{A}_1, \mathsf{A}_2, \ldots, \mathsf{A}_n$ of von
Neumann subalgebras of $\mathcal{A}$, having $\mathcal{N}$ as a
common subalgebra, is called \emph{freely independent}
\label{Libertad} over $\mathsf{E}_{\mathcal{N}}$ if
$$\mathsf{E}_{\mathcal{N}}(a_1 a_2 \cdots a_m) = 0$$ whenever
$\mathsf{E}_{\mathcal{N}}(a_k) = 0$ for all $1 \le k \le m$ and
$a_k \in \mathsf{A}_{j_k}$ with $j_1 \neq j_2 \neq \cdots \neq
j_m$. As in the scalar-valued case, operator-valued freeness
admits a Fock space representation.

\vskip5pt

We first assume that $\mathsf{A}_1, \mathsf{A}_2, \ldots,
\mathsf{A}_n$ are $\mathrm{C}^*$-algebras having $\mathcal{N}$ as
a common $\mathrm{C}^*$-subalgebra and that $\mathsf{E}_k =
{\mathsf{E}_{\mathcal{N}}}_{\mid_{\mathsf{A}_k}}$ are faithful
conditional expectations. Let us consider the mean-zero subspaces
$$\bubl_k = \Big\{ a_k \in \mathsf{A}_k \, \big| \
\mathsf{E}_k(a_k) = 0 \Big\}.$$ Following \cite{Dy,V}, we consider
the Hilbert $\mathcal{N}$-module $\bubl_{j_1} \otimes \bubl_{j_2}
\otimes \cdots \bubl_{j_m}$ (where the tensor products are
amalgamated over the von Neumann subalgebra $\mathcal{N}$)
equipped with the $\mathcal{N}$-valued inner product $$\big\langle
a_1 \otimes \cdots \otimes a_m, b_1 \otimes \cdots \otimes b_m
\big\rangle = \mathsf{E}_{j_m} \big( a_m^* \cdots
\mathsf{E}_{j_2}(a_2^* \, \mathsf{E}_{j_1}(a_1^* b_1) \, b_2)
\cdots b_m \big).$$ Then, the usual Fock space is replaced by the
Hilbert $\mathcal{N}$-module
$$\mathcal{H}_{\mathcal{N}} = \mathcal{N} \oplus \bigoplus_{m \ge
1} \bigoplus_{j_1 \neq j_2 \neq \cdots \neq j_m} \bubl_{j_1}
\otimes \bubl_{j_2} \otimes \cdots \otimes \bubl_{j_m}.$$ The
direct sums above are assumed to be $\mathcal{N}$-orthogonal. Let
$\mathcal{L(H_N)}$ stand for the algebra of adjointable maps on
$\mathcal{H_N}$. A linear right $\mathcal{N}$-module map
$\mathrm{T}: \mathcal{H_N} \rightarrow \mathcal{H_N}$ is called
\emph{adjointable} whenever there exists $\mathrm{S}:
\mathcal{H_N} \rightarrow \mathcal{H_N}$ such that
$$\langle x, \mathrm{T}y \rangle = \langle \mathrm{S} x,y \rangle
\qquad \mbox{for all} \qquad x,y \in \mathcal{H_N}.$$ Let us
recall how elements in $\mathsf{A}_k$ act on $\mathcal{H_N}$. We
decompose any $a_k \in \mathsf{A}_k$ as $$a_k = \bubla_k +
\mathsf{E}_k(a_k).$$ Any element in $\mathcal{N}$ acts on
$\mathcal{H_N}$ by left multiplication. Therefore, it suffices to
define the action of mean-zero elements. The $*$-homomorphism
$\pi_k: \mathsf{A}_k \rightarrow \mathcal{L(H_N)}$ is then defined
as follows $$\pi_k(\bubla_k) (b_{j_1} \otimes \cdots \otimes
b_{j_m}) = \left\{
\begin{array}{ll} \, \bubla_k \otimes b_{j_1} \otimes \cdots
\otimes b_{j_m}, & \mbox{if} \ k \neq j_1 \\ \\ \,
\mathsf{E}_k(\bubla_k b_{j_1}) \, b_{j_2} \otimes \cdots \otimes
b_{j_m} \ \oplus & \\ \big( \bubla_k b_{j_1} -
\mathsf{E}_k(\bubla_k b_{j_1}) \big) \otimes b_{j_2} \otimes
\cdots \otimes b_{j_m}, & \mbox{if} \ k = j_1.
\end{array} \right.$$ This definition also applies for the
empty word. Now, since the algebra $\mathcal{L(H_N)}$ is a
$\mathrm{C}^*$-algebra (\emph{c.f.} \cite{L}), we can define the
\emph{reduced $\mathcal{N}$-amalgamated free product}
$\mathrm{C}^*(*_{\mathcal{N}} \mathsf{A}_k)$ as the
$\mathrm{C}^*$-closure of linear combinations of operators of the
form $$\pi_{j_1}(a_1) \pi_{j_2}(a_2) \cdots \pi_{j_m}(a_m).$$

\vskip3pt

Now we consider the case in which $\mathcal{N}$ and $\mathsf{A}_1,
\mathsf{A}_2, \ldots, \mathsf{A}_n$ are von Neumann algebras. Let
$\varphi: \mathcal{N} \rightarrow \C$ be a \emph{n.f.} state on
$\mathcal{N}$. This provides us with the induced states
$\varphi_k: \mathsf{A}_k \rightarrow \C$ given by $\varphi_k =
\varphi \circ \mathsf{E}_k$. The Hilbert space $$L_2 \big(
\bubl_{j_1} \otimes \bubl_{j_2} \otimes \cdots \otimes
\bubl_{j_m}, \varphi \big)$$ is obtained from $\bubl_{j_1} \otimes
\bubl_{j_2} \otimes \cdots \otimes \bubl_{j_m}$ by considering the
inner product $$\big\langle a_1 \otimes \cdots \otimes a_m, b_1
\otimes \cdots \otimes b_m \big\rangle_{\varphi} = \varphi \Big(
\big\langle a_1 \otimes \cdots \otimes a_m, b_1 \otimes \cdots
\otimes b_m \big\rangle \Big).$$ Then we define the orthogonal
direct sum $$\mathcal{H}_{\varphi} = L_2(\mathcal{N}) \oplus
\bigoplus_{m \ge 1} \bigoplus_{j_1 \neq j_2 \neq \cdots \neq j_m}
L_2 \big( \bubl_{j_1} \otimes \bubl_{j_2} \otimes \cdots \otimes
\bubl_{j_m}, \varphi \big).$$ Let us consider the
$*$-representation $\lambda: \mathcal{L(H_N)} \rightarrow
\mathcal{B}(\mathcal{H}_{\varphi})$ defined by
$\lambda(\mathrm{T}) x = \mathrm{T} x$. The faithfulness of
$\lambda$ is implied by the fact that $\varphi$ is also faithful.
Let $\mathcal{M(H_N)}$ be the von Neumann algebra generated in
$\mathcal{B(H_{\varphi})}$ by $\mathcal{L(H_N)}$. Then, we define
the \emph{reduced $\mathcal{N}$-amalgamated free product}
$*_{\mathcal{N}} \mathsf{A}_k$ \label{Amalgamado} as the weak$^*$
closure of $\mathrm{C}^* (*_{\mathcal{N}} \mathsf{A}_k)$ in
$\mathcal{M(H_N)}$. After decomposing $$a_k = \bubla_k +
\mathsf{E}_k(a_k)$$ and identifying $\bubl_k$ with $\lambda \big(
\pi_k(\bubl_k) \big)$, we can think of $*_{\mathcal{N}}
\mathsf{A}_k$ as $$*_{\mathcal{N}} \mathsf{A}_k = \Big(
\mathcal{N} \oplus \bigoplus_{m \ge 1} \bigoplus_{j_1 \neq j_2
\neq \cdots \neq j_m} \bubl_{j_1} \bubl_{j_2} \cdots \bubl_{j_m}
\Big)''.$$ Let us consider the orthogonal projections
\begin{eqnarray*}
\mathcal{Q}_{\emptyset}: \mathcal{H}_{\varphi} & \rightarrow &
L_2(\mathcal{N}), \\ \mathcal{Q}_{j_1 \cdots j_m}:
\mathcal{H}_{\varphi} & \rightarrow & L_2 \big( \bubl_{j_1}
\otimes \bubl_{j_2} \otimes \cdots \otimes \bubl_{j_m}, \varphi
\big).
\end{eqnarray*}
If we also consider the projection $\mathcal{Q}_{\mathsf{A}_k} =
\mathcal{Q}_{\emptyset} + \mathcal{Q}_k$, the following mappings
\begin{eqnarray*}
\mathsf{E}_{\mathcal{N}}: x \in *_{\mathcal{N}} \mathsf{A}_k &
\mapsto & \mathcal{Q}_{\emptyset} x \mathcal{Q}_{\emptyset} \in
\mathcal{N},
\\ \mathcal{E}_{\mathsf{A}_k}: x \in *_{\mathcal{N}} \mathsf{A}_k
& \mapsto & \mathcal{Q}_{\mathsf{A}_k} x
\mathcal{Q}_{\mathsf{A}_k} \in \mathsf{A}_k,
\end{eqnarray*}
are faithful conditional expectations. Then, it turns out that
$\mathsf{A}_1, \mathsf{A}_2, \ldots, \mathsf{A}_n$ are von Neumann
subalgebras of $*_{\mathcal{N}} \mathsf{A}_k$ freely independent
over $\mathsf{E}_{\mathcal{N}}$. Reciprocally, if $\mathsf{A}_1,
\mathsf{A}_2, \ldots, \mathsf{A}_n$ is a collection of von Neumann
subalgebras of $\mathcal{A}$ freely independent over
$\mathsf{E}_{\mathcal{N}}: \mathcal{A} \rightarrow \mathcal{N}$
and generating $\mathcal{A}$, then $\mathcal{A}$ is isomorphic to
$*_{\mathcal{N}} \mathsf{A}_k$.

\vskip5pt

\begin{remark}
\emph{Let $\mathcal{N}_1$ and $\mathcal{N}_2$ be von Neumann
algebras and assume that $\mathcal{N}_2$ is equipped with a
\emph{n.f.} state $\varphi_2$. A relevant example of the
construction outlined above is the following. Let $\mathcal{A} =
\mathcal{N}_1 \otimes \mathcal{N}_2$ and let us consider the
conditional expectation $\mathsf{E}_{\mathcal{N}_1}: \mathcal{A}
\rightarrow \mathcal{N}_1$ defined by
$$\mathsf{E}_{\mathcal{N}_1}(x_1 \otimes x_2) = x_1 \otimes
\varphi_2(x_2) 1.$$ Assume that $\mathsf{A}_1, \mathsf{A}_2,
\ldots, \mathsf{A}_n$ are freely independent subalgebras of
$\mathcal{N}_2$ over $\varphi_2$. Then, it is well-known that
$\mathcal{N}_1 \otimes \mathsf{A}_1, \mathcal{N}_1 \otimes
\mathsf{A}_2, \ldots, \mathcal{N}_1 \otimes \mathsf{A}_n$ is a
family of freely independent subalgebras of $\mathcal{A}$ over
$\mathsf{E}_{\mathcal{N}_1}$, see e.g. Section 7 of \cite{J2}. In
particular, if $\mathsf{A}_1, \mathsf{A}_2, \ldots, \mathsf{A}_n$
generate $\mathcal{N}_2$, we obtain
\begin{equation} \label{Equation-CB-Amalgamated}
\mathcal{A} = \mathcal{N}_1 \otimes \big( \mathop{*}_{k=1}^n
\mathsf{A}_k \big) = *_{\mathcal{N}_1} ( \mathcal{N}_1 \otimes
\mathsf{A}_k).
\end{equation}}
\end{remark}

\subsection{A Rosenthal/Voiculescu type inequality}

In this paragraph we recall the free analogue \cite{JPX} of
Rosenthal inequalities \cite{Ro0} for mean-zero random variables
and prove a simple consequence of it. Let $\mathsf{A}_1,
\mathsf{A}_2, \ldots, \mathsf{A}_n$ be a family of von Neumann
algebras having $\mathcal{N}$ as a common von Neumann subalgebra
and let $*_\mathcal{N} \mathsf{A}_k$ be the corresponding
amalgamated free product. Given a family $a_1, a_2, \ldots, a_n$
in $*_\mathcal{N} \mathsf{A}_k$ we consider the row and column
conditional square functions
\begin{eqnarray*}
\mathcal{S}_{\mathrm{cond}}^r(a) & = & \Big( \sum_{k=1}^n
\mathsf{E}_{\mathcal{N}}(a_k a_k^*) \Big)^{1/2}, \\
\mathcal{S}_{\mathrm{cond}}^c(a) & = & \Big( \sum_{k=1}^n
\mathsf{E}_{\mathcal{N}}(a_k^* a_k) \Big)^{1/2}.
\end{eqnarray*}
\label{Dan}

\begin{theorem} \label{Theorem-Voiculescu}
If $2 \le p \le \infty$ and $a_1, a_2, \ldots, a_n \in
L_p(*_{\mathcal{N}} \mathsf{A}_k)$ with $a_k \in L_p(\bubl_k)$ for
$1 \le k \le n$, the following equivalence of norms holds with
relevant constants independent of $p$ or $n$
$$\Big\| \sum_{k=1}^n a_k \Big\|_p \sim \Big( \sum_{k=1}^n
\|a_k\|_p^p \Big)^{1/p} + \big\| \mathcal{S}_{\mathrm{cond}}^r(a)
\big\|_p + \big\| \mathcal{S}_{\mathrm{cond}}^c(a) \big\|_p.$$
\end{theorem}

This is the operator-valued/free analogue of Rosenthal's original
result. On the other hand, a noncommutative analogue was obtained
in \cite{JX4} for general algebras (non necessarily free
products), see also \cite{X2} for the notion of noncommutative
independence (called order independence) employed in it. Recalling
that freeness implies order independence, Theorem
\ref{Theorem-Voiculescu} follows from \cite{JX4} for $2 \le p <
\infty$. However, the constants there are not uniformly bounded as
$p \to \infty$, in sharp contrast with the situation in Theorem
\ref{Theorem-Voiculescu}. This is another example of an $L_p$
inequality involving independent random variables which only holds
in the limit case as $p \to \infty$ when considering their free
analogue. This constitutes a significant difference in this paper.
Theorem \ref{Theorem-Voiculescu} for $p=\infty$ was proved in
\cite{J2} and constitutes the operator valued extension of
Voiculescu's inequality \cite{V2}. Finally, we refer the reader to
\cite{JPX} for a generalization of Theorem
\ref{Theorem-Voiculescu}, where $a_1, a_2, \ldots, a_n$ are
replaced by certain words of a fixed degree $d \ge 1$.

\vskip5pt

The following result is an standard application for positive
random variables.

\begin{corollary} \label{Corollary-Voiculescu}
If $2 \le p \le \infty$ and $a_1, \ldots, a_n \in
L_p(*_{\mathcal{N}} \mathsf{A}_k)$ with $a_k \in
L_p(\mathsf{A}_k)$ for $1 \le k \le n$, the following equivalence
of norms holds with relevant constants independent of $p$ or $n$
\begin{eqnarray*}
\Big\| \big( \sum_{k=1}^n a_k a_k^* \big)^{1/2} \Big\|_p & \sim &
\Big( \sum_{k=1}^n \|a_k\|_p^p \Big)^{1/p} + \big\|
\mathcal{S}_{\mathrm{cond}}^r(a) \big\|_p, \\ \Big\| \big(
\sum_{k=1}^n a_k^* a_k \big)^{1/2} \Big\|_p & \sim & \Big(
\sum_{k=1}^n \|a_k\|_p^p \Big)^{1/p} + \big\|
\mathcal{S}_{\mathrm{cond}}^c(a) \big\|_p.
\end{eqnarray*}
\end{corollary}

\begin{proof} We we clearly have
$$\Big( \sum_{k=1}^n \|a_k\|_p^p \Big)^{1/p} \le \Big\| \big(
\sum_{k=1}^n a_k a_k^* \big)^{1/2} \Big\|_{p}.$$ Indeed, our claim
follows by complex interpolation from the trivial case $p =
\infty$ and the case $p = 2$, where the equality clearly holds.
This, together with the fact that $\mathsf{E}_{\mathcal{N}}$ is a
contraction on $L_{p/2} (*_{\mathcal{N}} \mathsf{A}_k)$, proves
the lower estimate with constant $2$. For the upper estimate, we
begin with the triangle inequality in $L_p(*_{\mathcal{N}}
\mathsf{A}_k; R_p^n)$ $$\Big\| \big( \sum_{k=1}^n a_k a_k^*
\big)^{1/2} \Big\|_p \le \Big\| \sum_{k=1}^n
\mathsf{E}_{\mathcal{N}}(a_k) \otimes e_{1k} \Big\|_p + \Big\|
\sum_{k=1}^n \bubla_k \otimes e_{1k} \Big\|_p = \mathrm{A} +
\mathrm{B}.$$ To estimate $\mathrm{A}$, we apply Kadison's
inequality (see Lemma \ref{Lemma-Properties-Expectation} i))
$$\mathrm{A} = \Big\| \Big( \sum_{k=1}^n
\mathsf{E}_{\mathcal{N}}(a_k) \mathsf{E}_{\mathcal{N}}(a_k)^*
\Big)^{1/2} \Big\|_p \le \Big\| \Big( \sum_{k=1}^n
\mathsf{E}_{\mathcal{N}}(a_k a_k^*) \Big)^{1/2} \Big\|_p = \big\|
\mathcal{S}_{\mathrm{cond}}^r(a) \big\|_p.$$ On the other hand,
according to \eqref{Equation-CB-Amalgamated} we can regard
$\bubla_k \otimes e_{1k}$ as an element of
$$S_p^n \big( L_p(*_{\mathcal{N}} \mathsf{A}_k) \big) = L_p \Big(
*_{S_{\infty}^n(\mathcal{N})} S_{\infty}^n(\mathsf{A}_k) \Big),$$
where $\mathsf{E}_{\mathcal{N}}$ is replaced by $1_{\mathrm{M}_n}
\otimes \mathsf{E}_{\mathcal{N}}$. Then, writing $x_k$ for
$\bubla_k \otimes e_{1k}$ and $\mathbf{E}_{\mathcal{N}}$ for
$1_{\mathrm{M}_n} \otimes \mathsf{E}_{\mathcal{N}}$, we estimate
$\mathrm{B}$ using the free Rosenthal inequalities in this bigger
space. Indeed, using $\|x_k\|_p = \|a_k\|_p$ we obtain
$$\mathrm{B} \sim \Big( \sum_{k=1}^n \|a_k\|_p^p \Big)^{1/p} + \Big\|
\Big( \sum_{k=1}^n \mathbf{E}_{\mathcal{N}} \big( x_k x_k^* \big)
\Big)^{1/2} \Big\|_p + \Big\| \Big( \sum_{k=1}^n
\mathbf{E}_{\mathcal{N}} \big( x_k^* x_k \big) \Big)^{1/2}
\Big\|_p.$$ Let us note that
$$\begin{array}{rclcl} \displaystyle \sum_{k=1}^n
\mathbf{E}_{\mathcal{N}} \big( x_k x_k^* \big) & = & \displaystyle
\sum_{k=1}^n e_{11}^{\null} \otimes \mathsf{E}_{\mathcal{N}} \big(
\bubla_k \bubla_k^* \big) & \le & \displaystyle \sum_{k=1}^n
e_{11}^{\null} \otimes \mathsf{E}_{\mathcal{N}} \big( a_k a_k^*
\big), \\ \displaystyle \sum_{k=1}^n \mathbf{E}_{\mathcal{N}}
\big( x_k^* x_k \big) & = & \displaystyle \sum_{k=1}^n \, e_{kk}
\otimes \mathsf{E}_{\mathcal{N}} \big( \bubla_k^* \bubla_k \big) &
\le & \displaystyle \sum_{k=1}^n \, e_{kk}  \otimes
\mathsf{E}_{\mathcal{N}} \big( a_k^* a_k \big). \end{array}$$ This
implies $$\Big\| \Big( \sum_{k=1}^n \mathbf{E}_{\mathcal{N}} \big(
x_k x_k^* \big) \Big)^{1/2} \Big\|_p \le \Big\| \Big( \sum_{k=1}^n
\mathsf{E}_{\mathcal{N}} (a_ka_k^*) \Big)^{1/2} \Big\|_p.$$ On the
other hand, the third term is controlled by
\begin{eqnarray*}
\Big\| \Big( \sum_{k=1}^n \mathbf{E}_{\mathcal{N}} \big( x_k^* x_k
\big) \Big)^{1/2} \Big\|_p & = & \Big\| \sum_{k=1}^n
\mathbf{E}_{\mathcal{N}} \big( x_k^* x_k \big) \Big\|_{p/2}^{1/2}
\\ & \le & \Big\| \sum_{k=1}^n e_{kk} \otimes \mathsf{E}_{\mathcal{N}}
\big( a_k^* a_k \big) \Big\|_{p/2}^{1/2} \\ & = & \Big(
\sum_{k=1}^n \big\| \mathsf{E}_\mathcal{N} \big( a_k^* a_k \big)
\big\|_{p/2}^{p/2} \Big)^{1/p} \le \Big( \sum_{k=1}^n \|a_k\|_p^p
\Big)^{1/p}.
\end{eqnarray*}
Combining the inequalities above we have the upper estimate with
constant $2$. \end{proof}

\section{Estimates for BMO type norms}

Apart from the free Rosenthal inequalities, our second key tool in
the proof of \eqref{Eq-Int-Int1} and \eqref{Eq-Int-Int2} will be
certain estimates that we develop in this section. Let us recall
the definition of several noncommutative Hardy spaces. Xu's survey
\cite{X2} contains a systematic exposition of these notions. Let
$\mathcal{M}$ be a von Neumann algebra equipped with a \emph{n.f.}
state $\varphi$. Let $\mathcal{M}_1, \mathcal{M}_2, \ldots$ be an
increasing filtration of von Neumann subalgebras of $\mathcal{M}$
which are invariant under the modular automorphism group on
$\mathcal{M}$ associated to $\varphi$. This allows us to consider
the corresponding conditional expectations $\mathcal{E}_n:
\mathcal{M} \to \mathcal{M}_n$ and noncommutative $L_p$
martingales $x = (x_n)_{n \ge 1}$ with martingale differences
$dx_1, dx_2, \ldots$ adapted to this filtration. The row and
column Hardy spaces $\Ha_p^r(\mathcal{M})$ and
$\Ha_p^c(\mathcal{M})$ are defined respectively as the closure of
the space of finite $L_p$ martingales with respect to the
following norms
\begin{eqnarray*}
\|x\|_{\Ha_p^r(\mathcal{M})} & = & \Big\| \big( \summ_k dx_k
dx_k^* \big)^{\frac12} \Big\|_p, \\ \|x\|_{\Ha_p^c(\mathcal{M})} &
= & \Big\| \big( \summ_k dx_k^* dx_k \big)^{\frac12} \Big\|_p.
\end{eqnarray*}
On the other hand, the space $\Ha_p^p(\mathcal{M})$ measures the
$p$-variation $$\|x\|_{\Ha_p^p(\mathcal{M})} = \Big( \summ_k
\|dx_k\|_p^p \Big)^{1/p}.$$ Finally, the conditional row and
column Hardy spaces for martingales $h_p^r(\mathcal{M})$ and
$h_p^c(\mathcal{M})$ are defined as the closure of the space of
finite $L_p$ martingales with respect to the following norms
\begin{eqnarray*}
\|x\|_{h_p^r(\mathcal{M})} & = & \Big\| \big( \summ_k
\mathcal{E}_{k-1} (dx_k dx_k^*) \big)^{\frac12} \Big\|_p, \\
\|x\|_{h_p^c(\mathcal{M})} & = & \Big\| \big( \summ_k
\mathcal{E}_{k-1} (dx_k^* dx_k) \big)^{\frac12} \Big\|_p.
\end{eqnarray*}

A further tool are maximal functions. The notion of maximal
function was introduced in \cite{J1} via the spaces
$L_p(\mathcal{M}; \ell_\infty)$. Here we are using variations of
these from Musat's paper \cite{M}. Given $2 \le p \le \infty$ we
define the spaces $L_p^r(\mathcal{M}; \ell_\infty)$ and
$L_p^c(\mathcal{M}; \ell_\infty)$ as the closure of the space of
bounded sequences in $L_p(\mathcal{M})$ with respect to the
following norms
\begin{eqnarray*}
\|(x_n)\|_{L_p^r(\EME; \ell_\infty)} & = & \big\| \sup_{n \ge 1}
\, x_n x_n^* \big\|_{p/2}^{1/2}, \\ \|(x_n)\|_{L_p^c(\EME;
\ell_\infty)} & = & \big\| \sup_{n \ge 1} \, x_n^* x_n
\big\|_{p/2}^{1/2}.
\end{eqnarray*}
This definition requires some explanation. Indeed, in the
noncommutative setting there is no obvious analogue for the
pointwise supremum of a family of positive operators. Therefore,
the above is to be understood in the sense of the suggestive
notation introduced in \cite{J1}. Among several characterizations
of the norm in $L_p(\mathcal{M}; \ell_\infty)$ of a sequence
$(z_n)_{n \ge 1}$ of positive operators, we outline the following
obtained by duality
\begin{equation} \label{Eq-Dual-Norm-Sup}
\big\| \sup_{n \ge 1} \, z_n \big\|_p = \sup \Big\{ \summ_n
\mbox{tr} (z_n w_n) \, \big| \ w_n \ge 0, \ \Big\| \summ_n w_n
\Big\|_{p'} \le 1 \Big\}.
\end{equation}
The reader is also referred to \cite{JX2} for a rather complete
exposition. One of the fundamental properties obtained in \cite{M}
of the spaces $L_p^r(\mathcal{M}; \ell_\infty)$ and
$L_p^c(\mathcal{M}; \ell_\infty)$ is that they form interpolation
families. Indeed, given $2 \le p_0, p_1 \le \infty$, we have the
following isometric isomorphisms for $1/p_\theta = (1-\theta)/p_0
+ \theta/p_1$
\begin{equation} \label{Eq-Interpol-Musat}
\begin{array}{rcl} \big[ L_{p_0}^r(\mathcal{M}; \ell_\infty),
L_{p_1}^r(\mathcal{M}; \ell_\infty) \big]_\theta & = &
L_{p_\theta}^r(\mathcal{M}; \ell_\infty), \\ [5pt] \big[
L_{p_0}^c(\mathcal{M}; \ell_\infty), L_{p_1}^c(\mathcal{M};
\ell_\infty) \big]_\theta & = & L_{p_\theta}^c(\mathcal{M};
\ell_\infty).
\end{array}
\end{equation}

\subsection{One-sided estimates}

We will simplify our arguments considerably by assuming that
$\mathcal{M}$ is finite and the density $\mathrm{D}$ associated to
the state $\varphi$ satisfies $c_1 1_\mathcal{M} \le \mathrm{D}
\le c_2 1_\mathcal{M}$. The general case will follow one more time
from Haagerup's approximation theorem.

\begin{lemma}
\label{Lemma-BMO1} Let $1 \le p \le 2 \le q \le \infty$ be such
that $1/p = 1/2 + 1/q$. Given $0 < \theta < 1$, let $x$ be a norm
one element in $[\Ha_p^p(\mathcal{M}),
\Ha_p^r(\mathcal{M})]_\theta$ and let us consider the indices $1
\le u \le 2 \le v \le \infty$ defined as follows
$$1/u = 1/p - \theta/q \quad \mbox{and} \quad 1/v = \theta/q.$$
Then we may find a positive element $a \in L_{v/2}(\mathcal{M})$
and $b_k \in L_u(\mathcal{M}_k)$ such that
$$dx_k = \mathcal{E}_k(a)^{\frac12} b_k \quad \mbox{and} \quad
\max \Big\{ \|a\|_{v/2}, \Big( \sum_{k \ge 1} \|b_k\|_u^u
\Big)^{1/u} \Big\} \le \sqrt{2}.$$ If $x$ belongs to
$[\Ha_p^p(\mathcal{M}), \Ha_p^c(\mathcal{M})]_\theta$, the same
conclusion holds with $dx_k = b_k^* \mathcal{E}_k(a)^{\frac12}$.
\end{lemma}

\begin{proof}
The last assertion follows from the first part of the statement by
taking adjoints. To prove the first assertion, we may assume by
approximation that $x$ is a finite martingale in
$L_p(\mathcal{M}_m)$ for some $m \ge 1$. Therefore, since $x$ is
of norm $1$ there exists an analytic function $f: \mathcal{S} \to
L_p(\mathcal{M}_m)$ of the form $$f(z) = \sum_{k=1}^m d_k(z),$$
which satisfies $f(\theta) = x$ and the estimate $$\max \left\{
\sup_{z \in \partial_0} \Big( \sum_{k=1}^m \|d_k(z)\|_p^p
\Big)^{\frac1p}, \, \sup_{z \in \partial_1} \Big\| \big(
\sum_{k=1}^m d_k(z) d_k(z)^* \big)^{\frac12} \Big\|_p \right\} \le
1.$$ Note that we are assuming that the $d_k$'s are also analytic
where $d_1(z), d_2(z), \ldots$ denote the martingale differences
of $f(z)$. Now we consider the following functions on the strip
for $1 \le k \le m$ $$g_k(z) = \left\{ \begin{array}{ll} 1, &
\mbox{if} \ z \in \partial_0, \\ \big( \sum_{j=1}^k d_j(z)
d_j(z)^* + \delta \mathrm{D}^{\frac2p} \big)^{\frac{p}{q}}, &
\mbox{if} \ z \in \partial_1.
\end{array} \right.$$ According to our original assumption on the
finiteness of $\mathcal{M}$ and the invertibility of $\mathrm{D}$,
we are in position to apply Devinatz's theorem. Indeed, here we
need Xu's modification, which can be found in Section 8 of Pisier
and Xu's survey \cite{PX2}. This provides us with analytic
function $h_k$ with analytic inverse and such that
\begin{equation} \label{Eq-Devinatz-Factor}
h_k(z) h_k(z)^* = g_k(z) \quad \mbox{for all} \quad z \in
\partial \mathcal{S}.
\end{equation}

\noindent \textsc{Step 1.} We claim that
\begin{equation} \label{Eq-Row-h}
\Big\| \sum_{k=1}^m \delta_k \otimes h_k(\theta)
\Big\|_{L_v^r(\mathcal{M}; \ell_\infty)} \le \big( 1 +
\delta^{\frac{p}{2}} \big)^{\frac{\theta}{q}},
\end{equation}
where $\delta$ appears in the definition of $g_k$. Indeed,
according to the interpolation isometries
\eqref{Eq-Interpol-Musat} and the three lines lemma, it suffices
to see that the following estimates hold
\begin{eqnarray}
\label{a1} \sup_{z \in \partial_0} \Big\| \sum_{k=1}^m \delta_k
\otimes h_k(z) \Big\|_{L_\infty^r(\mathcal{M}; \ell_\infty)} & \le
& 1, \\ \label{a2} \sup_{z \in \partial_1} \Big\| \sum_{k=1}^m
\delta_k \otimes h_k(z) \Big\|_{L_q^r(\mathcal{M}; \ell_\infty)} \
& \le & \big( 1 + \delta^{\frac{p}{2}} \big)^{\frac{1}{q}}.
\end{eqnarray}
To prove \eqref{a1} we first recall from
\eqref{Eq-Devinatz-Factor} and the definition of $g_k$ that
$h_k(z)$ is a unitary for any $z \in
\partial_0$ and any $1 \le k \le m$. Therefore, since Fubini's theorem
gives $L_\infty^r(\mathcal{M}; \ell_\infty) = L_\infty
(\ell_\infty(\mathcal{M}))$, we conclude
$$\sup_{z \in \partial_0} \Big\| \sum_{k=1}^m \delta_k \otimes
h_k(z) \Big\|_{L_\infty^r(\mathcal{M}; \ell_\infty)} = \sup_{z \in
\partial_0} \sup_{1 \le k \le m} \|h_k(z)\|_\infty = 1.$$ On the
other hand, if $z \in \partial_1$ we have the following estimate
from \eqref{Eq-Devinatz-Factor}
\begin{eqnarray*}
\Big\| \sum_{k=1}^m \delta_k \otimes h_k(z)
\Big\|_{L_q^r(\mathcal{M}; \ell_\infty)}^q & = & \Big\| \sup_{1
\le k \le m} \Big( \sum_{j=1}^k d_j(z) d_j(z)^* + \delta
\mathrm{D}^{\frac2p} \Big)^{p/q} \Big\|_{q/2}^{q/2} \\ & \le &
\Big\| \Big( \sum_{j=1}^m d_j(z) d_j(z)^* + \delta
\mathrm{D}^{\frac2p} \Big)^{p/q} \Big\|_{q/2}^{q/2} \\ & = &
\Big\| \sum_{j=1}^m d_j(z) d_j(z)^* + \delta \mathrm{D}^{\frac2p}
\Big\|_{p/2}^{p/2} \le 1 + \delta^{\frac{p}{2}},
\end{eqnarray*}
where the last inequality follows from the fact that
$L_{p/2}(\mathcal{M})$ is a $p/2$-normed space and also from the
right boundary estimate for the function $f$ given above. Thus,
inequality \eqref{a2} follows and the proof of \eqref{Eq-Row-h} is
completed.

\vskip5pt

\noindent \textsc{Step 2.} Let us consider the functions $w_k(z) =
h_k(z)^{-1} d_k(z)$. Now we claim that
\begin{equation} \label{Eq-lu-est}
\Big( \sum_{k=1}^m \|w_k(\theta)\|_u^u \Big)^{1/u} \le
\sqrt{2}^\theta \big( 1 + \delta^{\frac{p}{2}} \big)^{\theta/2}.
\end{equation}
According to \eqref{Eq-Devinatz-Factor}, we can write $h_k(z) =
g_k(z)^{\frac12} u_k(z)$ for some unitary $u_k(z)$. Thus, we
deduce that for $z \in \partial_0$ we have $w_k(z) = u_k(z)^*
d_k(z)$. This and the left boundary estimate for $f$ yield
$$\sup_{z \in
\partial_0} \Big( \sum_{k=1}^m \|w_k(z)\|_p^p \Big)^{1/p} \le 1.$$ The
interesting argument, based on the classical Fefferman-Stein
duality theorem, appears for $z \in \partial_1$. In that case we
have $h_k(z)^{-1} = u_k(z)^* g_k(z)^{-\frac12}$ and we find the
following estimate for any $z \in \partial_1$
\begin{eqnarray*}
\sum_{k=1}^m \|w_k(z)\|_2^2 & = & \sum_{k=1}^m \mbox{tr} \Big(
h_k(z)^{-1} d_k(z) d_k(z)^* h_k(z)^{-1*} \Big) \\ & = &
\sum_{k=1}^m \mbox{tr} \Big( u_k(z)^* g_k(z)^{-\frac12} d_k(z)
d_k(z)^* g_k(z)^{-\frac12} u_k(z) \Big) \\ & = & \sum_{k=1}^m
\mbox{tr} \Big( g_k(z)^{-\frac12} d_k(z) d_k(z)^*
g_k(z)^{-\frac12} \Big).
\end{eqnarray*}
Now we define the positive operators
\begin{eqnarray*}
\alpha_k & = & \Big( \sum_{j=1}^{k-1} d_j(z) d_j(z)^* + \delta
\mathrm{D}^{\frac2p} \Big)^{\frac{p}{2}}, \\ \beta_k & = & \Big(
\sum_{j=1}^{k} d_j(z) d_j(z)^* + \delta \mathrm{D}^{\frac2p}
\Big)^{\frac{p}{2}},
\end{eqnarray*}
and the indices $(s,t) = (2/q,2/p)$. Lemma 7.2 in \cite{JX} gives
$$\mbox{tr} \Big( \beta_k^{-s/2} (\beta_k^t - \alpha_k^t)
\beta_k^{-s/2} \Big) \le 2 \, \mbox{tr} \big( \beta_k - \alpha_k
\big).$$ According to our choice of $(s,t)$, we may rewrite this
inequality as follows $$\mbox{tr} \Big( g_k(z)^{-\frac12} d_k(z)
d_k(z)^* g_k(z)^{-\frac12} \Big) \le 2 \, \mbox{tr} \big( \beta_k
- \beta_{k-1} \big).$$ Summing up, we deduce that $$\sum_{k=1}^m
\|w_k(z)\|_2^2 \le 2 \, \Big\| \sum_{j=1}^{m} d_j(z) d_j(z)^* +
\delta \mathrm{D}^{\frac2p} \Big\|_{p/2}^{p/2} \le 2 \big( 1 +
\delta^{\frac{p}{2}} \big).$$ Finally, recalling that the
functions $w_k: \mathcal{S} \to L_p(\mathcal{M})$ are analytic
since are products of analytic functions, inequality
\eqref{Eq-lu-est} follows from the estimates given above and the
three lines lemma.

\vskip5pt

\noindent \textsc{Step 3.} For the moment, we have seen that
$$dx_k = h_k(\theta) w_k(\theta).$$ Now, recalling that
$$g_k: \partial \mathcal{S} \to L_{q/2}(\mathcal{M}_k)
\quad \mbox{for each} \quad 1 \le k \le m,$$ we conclude that
$h(\theta) = (h_1(\theta), h_2(\theta), \ldots, h_m(\theta))$ is
an \emph{adapted} sequence. That is, we have $h_k \in
L_v(\mathcal{M}_k)$ for all $1 \le k \le m$. On the other hand,
according to the definition of the space $L_v^r(\mathcal{M};
\ell_\infty)$, we may find for any $\delta > 0$ a positive
operator $a$ such that $h_k(\theta) h_k(\theta)^* \le a$ for $1
\le k \le m$ and
\begin{equation} \label{a3}
\|a\|_{v/2}^{1/2} \le (1 + \delta) \Big\| \sum_{k=1}^m \delta_k
\otimes h_k(\theta) \Big\|_{L_v^r(\mathcal{M}; \ell_\infty)} \le
\big( 1 + \delta^{\frac{p}{2}} \big)^{1 + \frac{\theta}{q}}.
\end{equation}
Moreover, since $h(\theta)$ is adapted, we have $h_k(\theta)
h_k(\theta)^* = \mathcal{E}_k (h_k(\theta) h_k(\theta)^*) \le
\mathcal{E}_k(a)$. This gives a contraction $\gamma_k \in
\mathcal{M}_k$ with $h_k(\theta) = \mathcal{E}_k(a)^{\frac12}
\gamma_k$. In particular, we deduce
\begin{equation} \label{a4}
dx_k = \mathcal{E}_k(a)^{\frac12} \gamma_k w_k(\theta) =
\mathcal{E}_k(a)^{\frac12} b_k.
\end{equation}
Finally, since $\gamma_k$ is a contraction
\begin{equation} \label{a5}
\Big( \sum_{k=1}^m \|b_k\|_u^u \Big)^{1/u} \le \Big( \sum_{k=1}^m
\|w_k(\theta)\|_u^u \Big)^{1/u} \le \sqrt{2}^\theta \big( 1 +
\delta^{\frac{p}{2}} \big)^{\theta/2}.
\end{equation}
Therefore, the assertion follows from \eqref{a3}, \eqref{a4} and
\eqref{a5} by letting $\delta \to 0^+$.
\end{proof}

Applying the anti-linear duality bracket $\langle x,y \rangle =
\mbox{tr} (x^*y)$, we consider in the following result an
immediate application of Lemma \ref{Lemma-BMO1} for the following
dual spaces
\begin{eqnarray*}
\mathcal{Z}_p^r(\mathcal{M},\theta) & = & \big[
\Ha_p^p(\mathcal{M}),
\Ha_p^r(\mathcal{M}) \big]_{\theta}^*, \\
\mathcal{Z}_p^c(\mathcal{M},\theta) & = & \big[
\Ha_p^p(\mathcal{M}), \Ha_p^c(\mathcal{M}) \big]_{\theta}^*.
\end{eqnarray*}

\begin{lemma} \label{Lemma-BMO2}
If $p,q,u,v$ are as above and $1/u+1/s=1$, we have
\begin{eqnarray*}
\|x\|_{\mathcal{Z}_p^r(\mathcal{M},\theta)}
& \le & 2 \, \sup_{\begin{subarray}{c} \|a\|_{v/2} \le 1 \\
a \ge 0 \end{subarray}} \Big( \summ_k \big\|
\mathcal{E}_k(a)^{\frac12} dx_k \big\|_s^s \Big)^{1/s}, \\
\|x\|_{\mathcal{Z}_p^c(\mathcal{M},\theta)}
& \le & 2 \, \sup_{\begin{subarray}{c} \|a\|_{v/2} \le 1 \\
a \ge 0
\end{subarray}} \Big( \summ_k \big\| dx_k
\mathcal{E}_k(a)^{\frac12} \big\|_s^s \Big)^{1/s}.
\end{eqnarray*}
\end{lemma}

\begin{proof}
There exists $y$ in the unit ball of $[\Ha_p^p(\mathcal{M}),
\Ha_p^r(\mathcal{M})]_\theta$ such that
$$\|x\|_{\mathcal{Z}_p^r(\mathcal{M},\theta)} = \big| \mbox{tr}
(xy^*) \big|.$$ On the other hand, according to Lemma
\ref{Lemma-BMO1} and using homogeneity, we may write $dy_k =
\mathcal{E}_k(a)^{1/2} b_k$ with $a$ being a positive operator in
the unit ball of $L_{v/2}(\mathcal{M})$ and $b_1, b_2, \ldots$
satisfying
\begin{equation} \label{Eq-bku}
\Big( \sum_{k \ge 1} \|b_k\|_u^u \Big)^{1/u} \le 2.
\end{equation}
With this decomposition we obtain the following estimate
\begin{eqnarray*}
\|x\|_{\mathcal{Z}_p^r(\mathcal{M},\theta)} & = & \Big| \summ_k
\mbox{tr} \big( dx_k dy_k^* \big) \Big| \\ & = & \Big| \summ_k
\mbox{tr} \big( dx_k  b_k^* \mathcal{E}_k(a)^{\frac12} \big) \Big|
\\ & \le & \Big( \summ_k \|b_k\|_u^u \Big)^{1/u} \Big( \summ_k \big\|
\mathcal{E}_k(a)^{\frac12} dx_k \big\|_s^s \Big)^{1/s}.
\end{eqnarray*}
The assertion follows from \eqref{Eq-bku}. The proof of the second
estimate is similar.
\end{proof}

In the following result and in the rest of this paper we shall
write $\mathrm{A} \lesssim \mathrm{B}$ to denote the existence of
an absolute constant $\mathrm{c}$ such that $\mathrm{A \le c \,
B}$ holds.

\begin{lemma} \label{Lemma-BMO6}
If $p,q,u,v$ are as above and $1/u+1/s=1$, we have
\begin{eqnarray*}
\sup_{\begin{subarray}{c} \|a\|_{v/2} \le 1 \\
a \ge 0 \end{subarray}} \Big( \summ_k \big\|
\mathcal{E}_k(a)^{\frac12} dx_k \big\|_s^s \Big)^{1/s} & \lesssim
&
\|x\|_{[\Ha_{p'}^{p'}(\mathcal{M}), \Ha_{p'}^r(\mathcal{M})]_\theta}, \\
\sup_{\begin{subarray}{c} \|a\|_{v/2} \le 1 \\
a \ge 0
\end{subarray}} \Big( \summ_k \big\| dx_k
\mathcal{E}_k(a)^{\frac12} \big\|_s^s \Big)^{1/s} & \lesssim &
\|x\|_{[\Ha_{p'}^{p'}(\mathcal{M}),
\Ha_{p'}^c(\mathcal{M})]_\theta}.
\end{eqnarray*}
\end{lemma}

\begin{proof}
We may find an analytic function $$f: \mathcal{S} \to
\Ha_{p'}^{p'}(\mathcal{M}) + \Ha_{p'}^r(\mathcal{M})$$ such that
$f(\theta) = x$ and $$\max \Big\{ \sup_{z \in \partial_0}
\|f(z)\|_{\Ha_{p'}^{p'}(\mathcal{M})}, \sup_{z \in \partial_1}
\|f(z)\|_{\Ha_{p'}^r(\mathcal{M})} \Big\} \lesssim
\|x\|_{[\Ha_{p'}^{p'}(\mathcal{M}),
\Ha_{p'}^r(\mathcal{M})]_\theta}.$$ By homogeneity, let us assume
that the right hand side above is $1$. Now we take positive
element $a$ in the unit ball of $L_{v/2}(\mathcal{M})$ so that we
may consider an analytic function $g: \mathcal{S} \to
L_\infty(\mathcal{M}) + L_{q/2}(\mathcal{M})$ satisfying
$g(\theta) = a$ and $$\max \Big\{ \sup_{z \in \partial_0}
\|g(z)\|_\infty, \sup_{z \in
\partial_1} \|g(z)\|_{q/2} \Big\} \le 1.$$ Then we construct the
analytic function $$h(z) = \summ_k \delta_k \otimes
d_k(f(\overline{z}))^* \mathcal{E}_k(g(z)) d_k(f(z)) \in
\ell_\infty \big( L_1(\mathcal{M}) \big).$$ For $z \in
\partial_0$ we have $$\|h(z)\|_{p'/2} \le \Big( \sum_{k \ge 1}
\|d_k(f(\overline{z}))\|_{p'}^{p'} \Big)^{1/p'} \sup_{k \ge 1}
\big\| \mathcal{E}_k(g(z)) \big\|_\infty \Big( \sum_{k \ge 1}
\|d_k(f(z))\|_{p'}^{p'} \Big)^{1/p'} \lesssim 1.$$ In the case $z
\in \partial_1$, we choose a factorization $g(z) = g_1(z) g_2(z)$
such that $$\|g_1(z)\|_q = \|g_2(z)\|_q = \sqrt{\|g(z)\|_{q/2}}
\le 1.$$ Then, H\"{o}lder inequality provides the following
estimate
\begin{eqnarray*}
\|h(z)\|_1 & = & \sum_{k \ge 1} \big\| \mathcal{E}_k \big(
d_k(f(\overline{z}))^* g_1(z) g_2(z) d_k(f(z)) \big) \big\|_1 \\ &
\le & \Big( \sum_{k \ge 1} \mathrm{tr} \big(
d_k(f(\overline{z}))^* g_1(z) g_1(z)^* d_k(f(\overline{z})) \big)
\Big)^{\frac12} \\ & \times & \Big( \sum_{k \ge 1} \mathrm{tr}
\big( d_k(f(z))^* g_2(z)^* g_2(z) d_k(f(z)) \big) \Big)^{\frac12}
\\ & \le & \Big\| \sum_{k \ge 1} d_k(f(\overline{z})) d_k(f(\overline{z}))^*
\Big\|_{p'/2}^{1/2} \Big\| \sum_{k \ge 1} d_k(f(z)) d_k(f(z))^*
\Big\|_{p'/2}^{1/2} \lesssim 1.
\end{eqnarray*}
Indeed, in order to apply H\"{o}lder inequality in the first
inequality above, we factorize the conditional expectation
$\mathcal{E}_{k-1}(a^*b)$ as the product $u_{k-1}(a)^* u_{k-1}(b)$
by a right $\mathcal{M}_{k-1}$-module map $u_{k-1}: \mathcal{M}
\to C_\infty(\mathcal{M}_{k-1})$, see e.g. \cite{J1,JPX}. By
complex interpolation ($2/s = 2(1-\theta)/p' + \theta$), we
conclude that $$\Big( \summ_k \big\| \mathcal{E}_k(a)^{\frac12}
dx_k \big\|_s^s \Big)^{1/s} = \|h(\theta)\|_{s/2}^{1/2} \lesssim
1.$$ The column version of this inequality follows by taking
adjoints.
\end{proof}

\begin{remark} \label{Remark-Equiv-Hardy}
\emph{When $1 < p \le 2$ we have}
\begin{eqnarray*}
\|x\|_{\mathcal{Z}_p^r(\mathcal{M},\theta)}
& \sim & \sup_{\begin{subarray}{c} \|a\|_{v/2} \le 1 \\
a \ge 0 \end{subarray}} \Big( \summ_k \big\|
\mathcal{E}_k(a)^{\frac12} dx_k \big\|_s^s \Big)^{1/s}, \\
\|x\|_{\mathcal{Z}_p^c(\mathcal{M},\theta)}
& \sim & \sup_{\begin{subarray}{c} \|a\|_{v/2} \le 1 \\
a \ge 0
\end{subarray}} \Big( \summ_k \big\| dx_k
\mathcal{E}_k(a)^{\frac12} \big\|_s^s \Big)^{1/s}.
\end{eqnarray*}
\emph{Indeed, since anti-linear duality is compatible with complex
interpolation via the analytic function
$\mbox{tr}(f(\overline{z})^*g(z))$, we find by reflexivity in the
case $1 < p \le 2$ the following isomorphisms}
\begin{equation} \label{Eq-DualH}
\begin{array}{rcl}
\mathcal{Z}_p^r(\mathcal{M},\theta) & \simeq &
[\Ha_{p'}^{p'}(\mathcal{M}), \Ha_{p'}^r(\mathcal{M})]_\theta, \\
[5pt] \mathcal{Z}_p^c(\mathcal{M},\theta) & \simeq &
[\Ha_{p'}^{p'}(\mathcal{M}), \Ha_{p'}^c(\mathcal{M})]_\theta.
\end{array}
\end{equation}
\emph{Therefore, the result follows from Lemma \ref{Lemma-BMO2}
and Lemma \ref{Lemma-BMO6}.}
\end{remark}

\begin{remark} \label{Remark-Results-CondHardy}
\emph{Lemmas \ref{Lemma-BMO1}, \ref{Lemma-BMO2}, \ref{Lemma-BMO6}
and Remark \ref{Remark-Equiv-Hardy} immediately generalize for the
row and column conditional Hardy spaces. Indeed, if we replace the
row and column Hardy spaces $\Ha_p^r(\mathcal{M})$ and
$\Ha_p^c(\mathcal{M})$ by their conditional analogues
$h_p^r(\mathcal{M})$ and $h_p^c(\mathcal{M})$ and the conditional
expectation $\mathcal{E}_k$ by $\mathcal{E}_{k-1}$, it can be
easily checked that the same arguments can be adapted to obtain
the conditioned results.}
\end{remark}

We shall also need to generalize the norms of $h_p^r(\mathcal{M})$
and $h_p^c(\mathcal{M})$ for arbitrary (non-necessarily adapted)
sequences $z_1, z_2, \ldots$ in $L_p(\mathcal{M})$ as follows
\begin{eqnarray*}
\Big\| \summ_k \delta_k \otimes z_k
\Big\|_{L_{\mathrm{cond}}^p(\mathcal{M}; \ell_2^r)} & = & \Big\|
\big( \summ_k
\mathcal{E}_{k-1} (z_k z_k^*) \big)^{\frac12} \Big\|_p, \\
\Big\| \summ_k \delta_k \otimes z_k
\Big\|_{L_{\mathrm{cond}}^p(\mathcal{M}; \ell_2^c)} & = & \Big\|
\big( \summ_k \mathcal{E}_{k-1} (z_k^* z_k) \big)^{\frac12}
\Big\|_p.
\end{eqnarray*}

\begin{lemma} \label{Lemma-BMO3}
Let $p,q$ be as above and let us consider the indices $(s,t)$
given by $1/s = (1-\eta)/p' + \eta/2$ and $1/t = \eta/q$ for some
parameter $0 < \eta < 1$. Then, the following estimates hold for
every martingale difference sequence $dx_1, dx_2, \ldots$ in
$L_{p'}(\mathcal{M})$
\begin{eqnarray*}
\Big\| \summ_k \delta_k \otimes dx_k a_k
\Big\|_{L_{\mathrm{cond}}^s(\mathcal{M}; \ell_2^r)} & \le & \Big(
\sum_{k \ge 1} \|a_k\|_t^t \Big)^{1/t}
\|x\|_{[h_{p'}^r(\mathcal{M}), \Ha_{p'}^{p'}(\mathcal{M})]_\eta},
\\ \Big\| \summ_k \delta_k \otimes a_k
dx_k \Big\|_{L_{\mathrm{cond}}^s(\mathcal{M}; \ell_2^c)} & \le &
\Big( \sum_{k \ge 1} \|a_k\|_t^t \Big)^{1/t}
\|x\|_{[h_{p'}^c(\mathcal{M}), \Ha_{p'}^{p'}(\mathcal{M})]_\eta}.
\end{eqnarray*}
\end{lemma}

\begin{proof}
We recall that
\begin{eqnarray*}
\big[ L_{\mathrm{cond}}^{p_0}(\mathcal{M}; \ell_2^r),
L_{\mathrm{cond}}^{p_1}(\mathcal{M}; \ell_2^r) \big]_{\theta} &
\subset & L_{\mathrm{cond}}^{p_\theta}(\mathcal{M}; \ell_2^r), \\
\big[ L_{\mathrm{cond}}^{p_0}(\mathcal{M}; \ell_2^c),
L_{\mathrm{cond}}^{p_1}(\mathcal{M}; \ell_2^c) \big]_{\theta} &
\subset & L_{\mathrm{cond}}^{p_\theta}(\mathcal{M}; \ell_2^c),
\end{eqnarray*}
hold isometrically. Indeed, we recall the factorization
$$\mathcal{E}_{k-1}(a^*b) = u_{k-1}(a)^*
u_{k-1}(b)$$ used in the proof of Lemma \ref{Lemma-BMO6} and the
resulting isometric embeddings
\begin{eqnarray*}
L_{\mathrm{cond}}^p(\mathcal{M}; \ell_2^r) & \subset & L_p \big(
\mathcal{M}; R_p(\N^2) \big), \\ L_{\mathrm{cond}}^p(\mathcal{M};
\ell_2^c) & \subset & L_p \big( \mathcal{M}; C_p(\N^2) \big).
\end{eqnarray*}
We refer the reader to \cite{J1} for a more detailed explanation.
According to our claim and by bilinear interpolation, it suffices
to prove the assertion for the extremal cases. When $\eta = 0$ we
have $(s,t) = (p', \infty)$ and the following estimate holds
$$\Big\| \big( \summ_k \mathcal{E}_{k-1} \big( dx_k a_k
a_k^* dx_k^* \big) \big)^{\frac12} \Big\|_{p'} \le \sup_{k \ge 1}
\|a_k\|_\infty \, \Big\| \big( \summ_k \mathcal{E}_{k-1} \big(
dx_k dx_k^* \big) \big)^{\frac12} \Big\|_{p'}.$$ On the other
hand, for $\eta = 1$ we have $(s,t) = (2,q)$ so that
\begin{eqnarray*}
\Big\| \big( \summ_k \mathcal{E}_{k-1} \big( dx_k a_k a_k^* dx_k^*
\big) \big)^{\frac12} \Big\|_2 & = & \Big( \summ_k \mbox{tr} \big(
a_k a_k^* dx_k^* dx_k \big) \Big)^{1/2} \\ & \le & \Big( \summ_k
\|a_k\|_q^2 \|dx_k\|_{p'}^2 \Big)^{1/2} \\ & \le & \Big( \summ_k
\|a_k\|_q^q \Big)^{1/q} \Big( \summ_k \|dx_k\|_{p'}^{p'}
\Big)^{1/p'}.
\end{eqnarray*}
This proves the first inequality. The second one follows by taking
adjoints.
\end{proof}

The following is the main result of this paragraph.

\begin{proposition} \label{Proposition-BMO}
If $p,q,u,v$ are as above and $1/u + 1/s = 1$, we have
\begin{eqnarray*}
\|x\|_{\mathcal{Z}_p^r(\mathcal{M},\theta)} & \le & c(p,\theta) \,
\max \Big\{ \|x\|_{\Ha_{p'}^{p'}(\mathcal{M})},
\|x\|_{[\Ha_{p'}^{p'}(\mathcal{M}), h_{p'}^r(\mathcal{M})]_\theta} \Big\}, \\
\|x\|_{\mathcal{Z}_p^c(\mathcal{M},\theta)} & \le & c(p,\theta) \,
\max \Big\{ \|x\|_{\Ha_{p'}^{p'}(\mathcal{M})},
\|x\|_{[\Ha_{p'}^{p'}(\mathcal{M}), h_{p'}^c(\mathcal{M})]_\theta}
\Big\}.
\end{eqnarray*}
Here the constant $c(p,\theta)$ satisfies $c(p,\theta) \sim 1$ as
$v \to \infty$ and
$$c(p,\theta) \sim \Big( \frac{p}{(2-p) \theta} - 1
\Big)^{-1/2} = \sqrt{\frac{2}{v-2}} \quad \mbox{as} \quad v \to
2.$$
\end{proposition}

\begin{proof}
According to Lemma \ref{Lemma-BMO2} we have
$$\|x\|_{\mathcal{Z}_p^r(\mathcal{M},\theta)}
\le 2 \, \sup_{\begin{subarray}{c} \|a\|_{v/2} \le 1 \\
a \ge 0 \end{subarray}} \Big( \summ_k \big\|
\mathcal{E}_k(a)^{\frac12} dx_k \big\|_s^s \Big)^{1/s}.$$ However,
since $2 \le s \le \infty$ we have
\begin{eqnarray}
\label{EqUse1} \lefteqn{\Big( \summ_k \big\|
\mathcal{E}_k(a)^{\frac12} dx_k \big\|_s^s \Big)^{1/s}} \\
\nonumber & = & \Big( \summ_k \big\| dx_k^* \big( da_k +
\mathcal{E}_{k-1}(a) \big) dx_k \big\|_{s/2}^{s/2} \Big)^{1/s}
\\ \nonumber & \le & \Big( \summ_k \big\| dx_k^* da_k dx_k
\big\|_{s/2}^{s/2} \Big)^{1/s} + \Big( \summ_k \big\| dx_k^*
\mathcal{E}_{k-1}(a) dx_k \big\|_{s/2}^{s/2} \Big)^{1/s}.
\end{eqnarray}
As we pointed above, Lemma \ref{Lemma-BMO6} and Remark
\ref{Remark-Equiv-Hardy} generalize to conditional Hardy spaces
after replacing $\mathcal{E}_k$ by $\mathcal{E}_{k-1}$. In other
words, the inequality below holds with absolute constants for $1
\le p \le 2$
\begin{equation} \label{Eq-Cond-Char-Hardy}
\sup_{\begin{subarray}{c} \|a\|_{v/2} \le 1 \\
a \ge 0 \end{subarray}} \Big( \summ_k \big\| dx_k^*
\mathcal{E}_{k-1}(a) dx_k \big\|_{s/2}^{s/2} \Big)^{1/s} \lesssim
\|x\|_{[\Ha_{p'}^{p'}(\mathcal{M}),
h_{p'}^r(\mathcal{M})]_\theta}.
\end{equation}
Therefore, it suffices to estimate the first term on the right of
\eqref{EqUse1}.

\vskip5pt

\noindent \textsc{Step 1}. We first assume $4 \le v \le \infty$.
Since $1/s = 1/v + 1/p'$
$$\Big( \summ_k \big\| dx_k da_k dx_k \big\|_{s/2}^{s/2} \Big)^{1/s}
\le \Big( \summ_k \|da_k\|_{v/2}^{v/2} \Big)^{1/v} \Big( \summ_k
\|dx_k\|_{p'}^{p'} \Big)^{1/p'}.$$ Then, complex interpolation
gives $$\Big( \summ_k \|da_k\|_{v/2}^{v/2} \Big)^{1/v} \le
\sqrt{2} \, \|a\|_{v/2}^{1/2} \le \sqrt{2} \quad \mbox{for} \quad
4 \le v \le \infty.$$ Indeed, our claim is trivial for the
extremal cases.

\vskip5pt

\noindent \textsc{Step 2}. The case $2 < v < 4$ is a little more
complicated. By the noncommutative Burkholder inequality
\cite{JX}, we may find a decomposition $da_k = d\alpha_k + d
\beta_k + d \gamma_k$ into three martingales satisfying the
following estimates
\begin{eqnarray}
\nonumber \Big( \summ_k \|d\alpha_k\|_{v/2}^{v/2} \Big)^{2/v} &
\le & c_{v} \|a\|_{v/2}, \\ \label{Eq-NCBurk-Est} \Big\| \big(
\summ_k \mathcal{E}_{k-1} (d\beta_k d\beta_k^*) \big)^{\frac12}
\Big\|_{v/2} & \le & c_{v} \|a\|_{v/2}, \\ \nonumber \Big\| \big(
\summ_k \mathcal{E}_{k-1} (d\gamma_k^* d\gamma_k) \big)^{\frac12}
\Big\|_{v/2} & \le & c_{v} \|a\|_{v/2},
\end{eqnarray}
where we know from \cite{Ran2} that $$c_v \lesssim \frac{1}{v-2}
\quad \mbox{for} \quad 2 \le v \le 4.$$ Since we need to estimate
the first term on the right of \eqref{EqUse1}, we decompose it
into three terms according to the martingale decomposition above.
The resulting term associated to $\alpha$ can be estimated as in
Step 1. For the second term (associated to $\beta$) we may find a
norm one element $(b_k)$ in $L_{(s/2)'}(\mathcal{M};
\ell_{(s/2)'})$ such that
\begin{equation} \label{Eq-RaizCuadrada}
\Big( \summ_k \big\| dx_k^* d\beta_k dx_k \big\|_{s/2}^{s/2}
\Big)^{1/s} = \sqrt{ \Big| \summ_k \mbox{tr} \big( b_k dx_k^*
d\beta_k dx_k\big) \Big|}.
\end{equation}
However, we have
\begin{eqnarray*}
\Big| \summ_k \mbox{tr} \big( b_k dx_k^* d\beta_k dx_k \big) \Big|
& = & \Big| \summ_k \mbox{tr} \big( d\beta_k dx_k b_k dx_k^* \big)
\Big| \\ & = & \Big| \summ_k \mbox{tr} \big(
\mathcal{E}_{k-1}(d\beta_k dx_k b_k dx_k^*) \big) \Big| \\ & \le &
\Big\| \big( \summ_k \mathcal{E}_{k-1} (d\beta_k d\beta_k^*)
\big)^{\frac12} \Big\|_{v/2} \\ & \times & \Big\| \big( \summ_k
\mathcal{E}_{k-1} (dx_k b_k^* dx_k^* dx_k b_k dx_k^* )
\big)^{\frac12} \Big\|_{(v/2)'} \\ & \le & c_v \|a\|_{v/2} \Big\|
\big( \summ_k \mathcal{E}_{k-1} (dx_k b_k^* dx_k^* dx_k b_k dx_k^*
) \big)^{\frac12} \Big\|_{(v/2)'}.
\end{eqnarray*}
Indeed, the first inequality above follows from H\"older
inequality after representing $\mathcal{E}_{k-1}(a^*b)$ as
$u_{k-1}(a)^* u_{k-1}(b)$ via the right module map $u_{k-1}$
considered in the proof of Lemma \ref{Lemma-BMO6}. The estimate
for the term associated to $\gamma$ is similar and yields the same
term with $b_k^*$ and $b_k$ exchanged. Now we have to estimate the
term
$$\Big\| \big( \summ_k \mathcal{E}_{k-1} (dx_k b_k^* dx_k^* dx_k
b_k dx_k^* ) \big)^{\frac12} \Big\|_{(v/2)'}.$$ Writing each $b_k$
as a linear combination $$b_k = \big( b_{k1} - b_{k2} \big) + i
\big( b_{k3} - b_{k4} \big)$$ of positive elements and allowing an
additional constant $2$, we may assume that the operators $b_1,
b_2, \ldots$ are positive. This consideration allows us to
construct the positive elements $z_k = (dx_k b_k
dx_k^*)^{\frac12}$. Then, recalling our assumption $2 < v < 4$, we
have $2 < (v/2)' < \infty$ and Lemma 5.2 of \cite{JX} gives for $t
= \frac12 (v/2)'$
$$\Big\| \summ_k \mathcal{E}_{k-1} (z_k^4) \Big\|_t \le \Big\|
\summ_k \mathcal{E}_{k-1} (z_k^2)
\Big\|_{2t}^{\frac{2(t-1)}{2t-1}} \Big( \summ_k \|z_k\|_{4t}^{4t}
\Big)^{\frac{1}{2t-1}}.$$ In our situation, this implies
\begin{eqnarray}
\label{Eq-ts} \lefteqn{\Big\| \big( \summ_k \mathcal{E}_{k-1}
(dx_k b_k^* dx_k^* dx_k b_k dx_k^* ) \big)^{\frac12}
\Big\|_{(v/2)'}} \\ \nonumber & \le & \Big\| \summ_k
\mathcal{E}_{k-1} (dx_k b_k dx_k^*)
\Big\|_{(v/2)'}^{\frac{t-1}{2t-1}} \Big( \summ_k \big\| dx_k b_k
dx_k^* \big\|_{(v/2)'}^{(v/2)'} \Big)^{\frac{1}{4t-2}}.
\end{eqnarray}
Using $\frac{1}{(v/2)'} - \frac{1}{(s/2)'} = \frac{2}{s} -
\frac{2}{v} = \frac{2}{p'}$, we find $$\Big( \summ_k \big\| dx_k
b_k dx_k^* \big\|_{(v/2)'}^{(v/2)'} \Big)^{\frac{1}{(v/2)'}} \le
\|dx\|_{p'} \|b\|_{(s/2)'} \|dx^*\|_{p'} \le \Big( \summ_k
\|dx_k\|_{p'}^{p'} \Big)^{\frac{2}{p'}}.$$ In other words, we have
$$\Big( \summ_k \big\| dx_k b_k
dx_k^* \big\|_{(v/2)'}^{(v/2)'} \Big)^{\frac{1}{4t-2}} \le
\|x\|_{\Ha_{p'}^{p'}(\mathcal{M})}^{\frac{4t}{4t-2}}.$$ For the
first term on the right of \eqref{Eq-ts} we use Lemma
\ref{Lemma-BMO3} with $$\big( \eta,s,t \big) = \Big( 1 - \theta, 2
(v/2)', 2(s/2)' \Big).$$ This yields
\begin{eqnarray*}
\Big\| \summ_k \mathcal{E}_{k-1} (dx_k b_k dx_k^*) \Big\|_{(v/2)'}
& = & \Big\| \summ_k \delta_k \otimes dx_k b_k^{\frac12}
\Big\|_{L_{\mathrm{cond}}^{2(v/2)'}(\mathcal{M}; \ell_2^r)}^2 \\
& \le & \Big( \summ_k \|b_k\|_{(s/2)'}^{(s/2)'}
\Big)^{\frac{1}{{(s/2)'}}} \|x\|_{[\Ha_{p'}^{p'}(\mathcal{M}),
h_{p'}^r(\mathcal{M})]_\theta}^2.
\end{eqnarray*}
Hence, since $(b_k)$ is in the unit ball of
$L_{(s/2)'}(\mathcal{M}; \ell_{(s/2)'})$, we conclude $$\Big\|
\summ_k \mathcal{E}_{k-1} (dx_k b_k dx_k^*)
\Big\|_{(v/2)'}^{\frac{t-1}{2t-1}} \le
\|x\|_{[\Ha_{p'}^{p'}(\mathcal{M}),
h_{p'}^r(\mathcal{M})]_\theta}^{\frac{4t-4}{4t-2}}$$ The estimates
above and \eqref{Eq-ts} give rise to $$\Big\| \big( \summ_k
\mathcal{E}_{k-1} (dx_k b_k^* dx_k^* dx_k b_k dx_k^* )
\big)^{\frac12} \Big\|_{(v/2)'} \!\! \le \!\! \max \Big\{
\|x\|_{\Ha_{p'}^{p'}(\mathcal{M})}^2,
\|x\|_{[\Ha_{p'}^{p'}(\mathcal{M}),
h_{p'}^r(\mathcal{M})]_\theta}^2 \Big\}.$$ Taking square roots as
imposed by \eqref{Eq-RaizCuadrada} and keeping track of the
constants, we obtain the assertion for
$\mathcal{Z}_p^r(\mathcal{M},\theta)$. The column case follows by
taking adjoints.
\end{proof}

At the beginning of this paragraph, we assumed that the von
Neumann algebra $\mathcal{M}$ was finite and equipped with a
\emph{n.f.} state $\varphi$ with respect to which the associated
density satisfied $c_1 1_\mathcal{M} \le \mathrm{D} \le c_2
1_\mathcal{M}$. This assumption was only needed in Lemma
\ref{Lemma-BMO1} (and its conditional version) in order to apply a
variation of Devinatz's theorem. On the other hand, this result
has been only applied to prove Lemma \ref{Lemma-BMO2}. Therefore,
if we are able to show that Lemma \ref{Lemma-BMO2} (and its
conditional version) holds for arbitrary $\sigma$-finite von
Neumann algebras, the same will hold for Proposition
\ref{Proposition-BMO}. As we shall see, this is relevant for our
aims since we shall use these results in the context of free
products. Let us indicate how to derive Lemma \ref{Lemma-BMO2} for
$\sigma$-finite von Neumann algebras. As expected, we apply
Haagerup's construction and consider $\mathcal{R} = \mathcal{M}
\rtimes_{\sigma} \mathrm{G}$ for the discrete group
$$\mathrm{G} = \bigcup_{n \in \N} 2^{-n} \Z.$$ The crossed product
$\mathcal{R}$ is a direct limit of a family of finite von Neumann
algebras $\mathcal{R}_1, \mathcal{R}_2, \ldots$ (we change our
usual notation here since in this chapter $\mathcal{M}_1,
\mathcal{M}_2, \ldots$ stand for a filtration of $\mathcal{M}$)
which are obtained as centralizers constructed from the modular
action for $\varphi \circ \mathsf{E}_\mathcal{M}$, where
$$\mathsf{E}_\mathcal{M}: \sum_{g \in \mathrm{G}} x_g \lambda(g)
\in \mathcal{R} \mapsto x_0 \in \mathcal{M}$$ denotes the natural
conditional expectation onto $\mathcal{M}$. The trace in
$\mathcal{R}_n$ is given by $\tau_n(x) = \varphi \circ
\mathsf{E}_\mathcal{M}(d_n x)$ where $c_{1n} 1_{\mathcal{R}_n} \le
d_n \le c_{2n} 1_{\mathcal{R}_n}$. It is then easily checked that
\begin{eqnarray*}
\hat{\mathcal{M}}_n & = & \mathcal{M}_n \rtimes_\sigma \mathrm{G},
\\ \hat{\mathcal{M}}_n(m) & = & \hat{\mathcal{M}}_n \cap
\mathcal{R}_m,
\end{eqnarray*}
are increasing filtrations in $\mathcal{R}$ and $\mathcal{R}_m$
respectively. Thus, Lemma \ref{Lemma-BMO2} (and its conditional
version) holds for $\mathcal{R}_m$ and the filtration
$\hat{\mathcal{M}}_1(m), \hat{\mathcal{M}}_2(m), \ldots$ for fixed
$m \ge 1$. Moreover, according to
\eqref{Equation-Limit-Expectation} and a simple density argument,
Lemma \ref{Lemma-BMO2} remains valid for $\mathcal{R}$ and the
filtration $\hat{\mathcal{M}}_1, \hat{\mathcal{M}}_2, \ldots$
Finally, it remains to see that $\mathsf{E}_\mathcal{M}$ extends
to a contraction on $\Ha_p^s(\mathcal{R})$ and
$h_p^s(\mathcal{R})$ where $s \in \{r,c,p\}$. This is obvious for
$s=p$, see e.g. \cite{JX} for the convention $h_p^p(\mathcal{R}) =
\Ha_p^p(\mathcal{R})$. For $s=r,c$ we recall that
$\mathsf{E}_{\mathcal{M}}$ and $\mathcal{E}_{k-1}$ commute. In the
case $2 \le p \le \infty$, this implies
\begin{eqnarray*}
\Big\| \summ_k \mathcal{E}_{k-1} \big(
d_k(\mathsf{E}_{\mathcal{M}}(x)) d_k(\mathsf{E}_\mathcal{M}(x))^*
\big) \Big\|_{p/2} & = & \Big\| \summ_k \mathcal{E}_{k-1} \big(
\mathsf{E}_\mathcal{M}(dx_k) \mathsf{E}_\mathcal{M}(dx_k^*) \big)
\Big\|_{p/2} \\ & \le & \Big\| \summ_k \mathcal{E}_{k-1} \big(
\mathsf{E}_\mathcal{M}(dx_k dx_k^*) \big) \Big\|_{p/2} \\ & = &
\Big\| \mathsf{E}_\mathcal{M} \big( \summ_k \mathcal{E}_{k-1}
(dx_k dx_k^*) \big) \Big\|_{p/2} \\
& \le & \Big\| \summ_k \mathcal{E}_{k-1}(dx_k dx_k^*)
\Big\|_{p/2}.
\end{eqnarray*}
The arguments for $\Ha_p^r(\mathcal{R})$, $p \ge 2$ (as well as
the analogues for the column spaces) are the same. For $1 \le p <
2$ we have to argue differently. By duality, it suffices to prove
the assertion for $L_q^r \mathcal{MO}$ with $2 < q \le \infty$,
see Theorem 4.1 of \cite{JX} for this duality result. Indeed, the
norm in that space is given by
$$\|x\|_{L_q^r \mathcal{MO}} = \Big\| \sup_{n \ge 1} \sum_{k \ge n}
\mathcal{E}_n (dx_k dx_k^*) \Big\|_{q/2}^{\frac12}.$$ Therefore,
using the inequality $$\big\| \sup_{n \ge 1}
\mathsf{E}_\mathcal{M}(z_n z_n^*) \big\|_{q/2} \le \big\| \sup_{n
\ge 1} z_n z_n^* \big\|_{q/2}$$ which follows easily from
\eqref{Eq-Dual-Norm-Sup}, we see that the same argument above
applies. On the other hand, it is easily checked (as in \cite{JX},
Theorem 4.1) that the corresponding dual in the conditional case
is given by $$\|x\|_{L_q^r mo} = \Big\| \sup_{n \ge 1} \sum_{k >
n} \mathcal{E}_n (dx_k dx_k^*) \Big\|_{q/2}^{\frac12}.$$

\subsection{Two-sided estimates}

Now we will perform a similar task considering two-sided terms.
Since most of the arguments are the same, we shall only sketch the
main ideas in the proofs. Again, we begin by assuming that
$\mathcal{M}$ is finite and the density $\mathrm{D}$ associated to
the state $\varphi$ satisfies $c_1 1_\mathcal{M} \le \mathrm{D}
\le c_2 1_\mathcal{M}$. The above argument via Haagerup's
construction leads to the $\sigma$-finite case.

\begin{lemma} \label{Lemma-BMO4}
Let $1 \le p \le 2 \le q \le \infty$ be such that $1/p = 1/2 +
1/q$. Given $0 < \theta < 1$, let $x$ be a norm one element in
$[\Ha_p^r(\mathcal{M}), \Ha_p^c(\mathcal{M})]_\theta$ and let us
consider the indices $2 \le w_r, w_c \le \infty$ defined as
follows $$1/w_r = (1-\theta)/q \quad \mbox{and} \quad 1/w_c =
\theta/q.$$ Then there exists $(a_r,a_c) \in
L_{w_r/2}(\mathcal{M})_+ \times L_{w_c/2}(\mathcal{M})_+$ and $b_k
\in L_2(\mathcal{M}_k)$ such that
$$dx_k = \mathcal{E}_k(a_r)^{\frac12} b_k
\mathcal{E}_k(a_c)^{\frac12} \quad \mbox{and} \quad \max \Big\{
\|a_r\|_{w_r/2}, \Big( \sum_{k \ge 1} \|b_k\|_2^2 \Big)^{\frac12},
\|a_c\|_{w_c/2} \Big\} \le 2.$$
\end{lemma}

\begin{proof}
Assume by approximation that $x$ is a finite martingale in
$L_p(\mathcal{M}_m)$ for some integer $m \ge 1$. Therefore, since
$x$ is of norm $1$ there exists an analytic function $f:
\mathcal{S} \to L_p(\mathcal{M}_m)$ of the form
$$f(z) = \sum_{k=1}^m d_k(z),$$ which satisfies $f(\theta) = x$
and the estimate $$\max \left\{ \sup_{z \in
\partial_0} \Big\| \big( \sum_{k=1}^m d_k(z) d_k(z)^*
\big)^{\frac12} \Big\|_p, \, \sup_{z \in
\partial_1} \Big\| \big( \sum_{k=1}^m d_k(z)^* d_k(z)
\big)^{\frac12} \Big\|_p \right\} \le 1.$$ Now we define
$$g_k^r(z) = \left\{ \begin{array}{ll} \big( \sum_{j=1}^k
d_j(z) d_j(z)^* + \delta \mathrm{D}^{\frac2p} \big)^{\frac{p}{q}},
& \mbox{if} \ z \in \partial_0, \\ 1, & \mbox{if} \ z \in
\partial_1. \end{array} \right.$$ $$g_k^c(z) = \left\{
\begin{array}{ll} 1, & \mbox{if} \ z \in
\partial_0, \\ \big( \sum_{j=1}^k d_j(z)^* d_j(z) + \delta
\mathrm{D}^{\frac2p} \big)^{\frac{p}{q}}, & \mbox{if} \ z \in
\partial_1.  \end{array} \right.$$ As in Lemma \ref{Lemma-BMO1},
we are now in position to apply Devinatz's theorem. This provides
us with two analytic functions $h_k^r$ and $h_k^c$ with analytic
inverses and satisfying the following relations
\begin{equation} \label{Eq-Devinatz-FactorRC}
\begin{array}{ccc}
h_k^r(z) h_k^r(z)^* = g_k^r(z) & \mbox{for all} & z \in
\partial \mathcal{S}, \\ [5pt] h_k^c(z)^* h_k^c(z) =
g_k^c(z) & \mbox{for all} & z \in \partial \mathcal{S}.
\end{array}
\end{equation}
According to the argument in Step 1 of the proof of Lemma
\ref{Lemma-BMO1}, we have
\begin{equation} \label{Eq-RC-h}
\begin{array}{rcl}
\displaystyle \Big\| \sum_{k=1}^m \delta_k \otimes h_k^r(\theta)
\Big\|_{L_{w_r}^r(\mathcal{M}; \ell_\infty)} & \le & \big( 1 +
\delta^{\frac{p}{2}} \big)^{\frac{1-\theta}{q}}, \\ [7pt]
\displaystyle \Big\| \sum_{k=1}^m \delta_k \otimes h_k^c(\theta)
\Big\|_{L_{w_c}^c(\mathcal{M}; \ell_\infty)} & \le & \big( 1 +
\delta^{\frac{p}{2}} \big)^{\frac{\theta}{q}}.
\end{array}
\end{equation}
Now we consider the functions $$w_k(z) = h_k^r(z)^{-1} d_k(z)
h_k^c(z)^{-1}.$$ Note that we have unitaries $u_k^r(z)$ and
$u_k^c(z)$ for which $$h_k^r(z) = g_k^r(z)^{\frac12} u_k^r(z)
\quad \mbox{and} \quad h_k^r(z) = u_k^c(z) g_k^c(z)^{\frac12}.$$
Thus, $w_k$ can be rewritten as follows on $\partial \mathcal{S}$
\begin{eqnarray*}
w_k(z) & = & u_k^r(z)^* g_k^r(z)^{-\frac12} d_k(z) u_k^c(z)^*
\quad \mbox{on} \quad \partial_0, \\ w_k(z) & = & u_k^r(z)^*
d_k(z) g_k^c(z)^{-\frac12} u_k^c(z)^* \quad \mbox{on} \quad
\partial_1.
\end{eqnarray*}
In particular, the same argument as in Lemma \ref{Lemma-BMO1}
(second part of Step 2) yields $$\sup_{z \in \partial \mathcal{S}}
\Big( \sum_{k=1}^m \|w_k(z)\|_2^2 \Big)^{\frac12} \le \sqrt{2
\big(1 + \delta^{\frac{p}{2}} \big)}.$$ Therefore, the same bound
holds for $z = \theta$ by the three lines lemma. For the moment,
we have seen that $dx_k = h_k^r(\theta) w_k(\theta)
h_k^c(\theta)$. Now, recalling that $g_k^r$ and $g_k^c$ take
values in $L_{q/2}(\mathcal{M}_k)$, we deduce that the sequence
$h_1^r(\theta), h_2^r(\theta), \ldots$ as well as $h_1^c(\theta),
h_2^c(\theta), \ldots$ are adapted. In particular, the argument in
Step 3 of Lemma \ref{Lemma-BMO1} gives rise to
$$h_k^r(\theta) = \mathcal{E}_k(a_r)^{\frac12} \gamma_k^r \quad
\mbox{and} \quad h_k^c(\theta) = \gamma_k^c
\mathcal{E}_k(a_c)^{\frac12}$$ for some contractions $\gamma_k^r,
\gamma_k^c \in \mathcal{M}_k$ and some positive elements $a_r,a_c$
satisfying $$\max \Big\{ \|a_r\|_{w_r/2}, \|a_c\|_{w_c/2} \Big\}
\le \big( 1 + \delta^{\frac{p}{2}} \big)^{1 + \frac1q}.$$ The
proof is completed by taking the elements $b_k$ to be $\gamma_k^r
w_k \gamma_k^c$ for $1 \le k \le m$.
\end{proof}

Exactly as we did in Lemma \ref{Lemma-BMO2}, the following result
is a direct application of Lemma \ref{Lemma-BMO4} for the dual
space $$\mathcal{Z}_p(\mathcal{M},\theta) = \big[
\Ha_p^r(\mathcal{M}), \Ha_p^c(\mathcal{M}) \big]_{\theta}^*.$$

\begin{lemma} \label{Lemma-BMO5}
If $p,q,w_r,w_c$ are as above, we have
$$\|x\|_{\mathcal{Z}_p(\mathcal{M},\theta)}
\le \ 8 \sup_{\begin{subarray}{c} \|a_r\|_{w_r/2}, \|a_c\|_{w_c/2}
\le 1 \\ a_r, a_c \ge 0 \end{subarray}} \Big( \summ_k \big\|
\mathcal{E}_k(a_r)^{\frac12} dx_k \mathcal{E}_k(a_c)^{\frac12}
\big\|_2^2 \Big)^{1/2}.$$
\end{lemma}

In Proposition \ref{Proposition-BMO} we found the constant
$c(p,\theta)$. Now we define
$$c'(p,\theta) = \max \Big\{ c(p,\theta), c(p,1-\theta) \Big\}.$$

\begin{proposition} \label{Proposition-BMO2}
If $p,q,w_r, w_c$ are as above, we have
$$\ \|x\|_{\mathcal{Z}_p(\mathcal{M},\theta)} \le c'(p,\theta)
\max \Big\{ \|x\|_{\Ha_{p'}^{p'}}, \|x\|_{[h_{p'}^r,
\Ha_{p'}^{p'}]_\theta}, \|x\|_{[\Ha_{p'}^{p'}, h_{p'}^c]_\theta},
\|x\|_{[h_{p'}^r, h_{p'}^c]_\theta} \Big\}.$$
\end{proposition}

\begin{proof}
According to Lemma \ref{Lemma-BMO5} we have to estimate the term
$$\mathrm {A} = \Big( \summ_k \big\|
\mathcal{E}_k(a_r)^{\frac12} dx_k \mathcal{E}_k(a_c)^{\frac12}
\big\|_2^2 \Big)^{1/2}$$ for any pair $(a_r,a_c)$ of positive
elements satisfying $\|a_r\|_{w_r/2} \le 1$ and $\|a_c\|_{w_c/2}
\le 1$. To that aim, we decompose it into the following three
terms
\begin{eqnarray}
\label{EqUse1RC} \mathrm{A}^2 & = & \summ_k \mbox{tr} \big( dx_k
\mathcal{E}_k(a_c) dx_k^* \mathcal{E}_k(a_r) \big) \\ \nonumber &
\le & \Big| \summ_k \mbox{tr} \big( dx_k d_k(a_c) dx_k^* d_k(a_r)
\big) \Big| \\ \nonumber & + & \Big| \summ_k \mbox{tr} \big( dx_k
\mathcal{E}_{k-1}(a_c) dx_k^* d_k(a_r) \big) \Big| \\ \nonumber &
+ & \Big| \summ_k \mbox{tr} \big( dx_k d_k(a_c) dx_k^*
\mathcal{E}_{k-1}(a_r) \big) \Big| \\ \nonumber & + & \Big|
\summ_k \mbox{tr} \big( dx_k \mathcal{E}_{k-1}(a_c) dx_k^*
\mathcal{E}_{k-1}(a_r) \big) \Big| = \mathrm{A}_1^2 +
\mathrm{A}_2^2 + \mathrm{A}_3^2 + \mathrm{A}_4^2.
\end{eqnarray}
In particular, we have $\mathrm{A} \le \mathrm{A}_1 + \mathrm{A}_2
+ \mathrm{A}_3 + \mathrm{A}_4$. The estimate for $\mathrm{A}_4$ is
rather simple. Indeed, arguing as in Remark
\ref{Remark-Results-CondHardy}, the conditional version of Lemma
\ref{Lemma-BMO5} follows after replacing $\mathcal{E}_k$ by
$\mathcal{E}_{k-1}$. Moreover, as in \eqref{Eq-Cond-Char-Hardy}
the argument in Lemma \ref{Lemma-BMO6} gives
$$\sup_{\begin{subarray}{c} \|a_r\|_{w_r/2}, \|a_c\|_{w_c/2} \le 1
\\ a_r, a_c \ge 0 \end{subarray}} \Big( \summ_k \big\|
\mathcal{E}_{k-1}(a_r)^{\frac12} dx_k
\mathcal{E}_{k-1}(a_c)^{\frac12} \big\|_2^2 \Big)^{1/2} \lesssim
\|x\|_{[h_{p'}^r(\mathcal{M}), h_{p'}^c(\mathcal{M})]_\theta},$$
with absolute constants. Therefore, we find
$$\mathrm{A}_4 = \Big( \summ_k \big\|
\mathcal{E}_{k-1}(a_r)^{\frac12} dx_k
\mathcal{E}_{k-1}(a_c)^{\frac12} \big\|_2^2 \Big)^{1/2} \lesssim
\, \|x\|_{[h_{p'}^r(\mathcal{M}),
h_{p'}^c(\mathcal{M})]_\theta}.$$ It remains to estimate the terms
$\mathrm{A}_1$, $\mathrm{A}_2$ and $\mathrm{A}_3$.

\vskip5pt

\noindent \textsc{Step 1}. We first estimate the term
$\mathrm{A}_1$. Recalling that $1/w_r + 1/w_c = 1/q \le 1/2$, we
must have $4 \le \max(w_r,w_c) \le \infty$. Moreover, since both
cases can be argued in the same way, we assume without lost of
generality that $4 \le w_r \le \infty$. In this case we have
\begin{eqnarray*}
\mathrm{A}_1^2 & = & \Big| \summ_k \mbox{tr} \big( dx_k d_k(a_c)
dx_k^* d_k(a_r) \big) \Big| \\ & \le & \Big( \summ_k
\|d_k(a_r)\|_{w_r/2}^{w_r/2} \Big)^{2/w_r} \Big( \summ_k \big\|
dx_k d_k(a_c) dx_k^* \big\|_{(w_r/2)'}^{(w_r/2)'}
\Big)^{1/(w_r/2)'}.
\end{eqnarray*}
The first term on the right is controlled by $\|a_r\|_{w_r/2}$
since we are assuming that $w_r \ge 4$. For the second term on the
right, we observe that $2 (w_r/2)'$ is determined by the following
relation $$\frac{1}{2(w_r/2)'} = \frac{1}{2} - \frac{1}{w_r} =
\frac{1}{2} - \frac{1-\theta}{q} = 1 - \frac{1}{p} +
\frac{\theta}{q} = 1 - \frac{1}{u} = \frac{1}{s}.$$ Thus, this
term is estimated as follows
\begin{eqnarray*}
\Big( \summ_k \big\| dx_k d_k(a_c) dx_k^* \big\|_{s/2}^{s/2}
\Big)^{2/s} & \le & \Big( \summ_k \big\| dx_k
\mathcal{E}_k(a_c)^{\frac12} \big\|_s^s \Big)^{2/s} \\ & + & \Big(
\summ_k \big\| dx_k \mathcal{E}_{k-1}(a_c)^{\frac12} \big\|_s^s
\Big)^{2/s} = \mathrm{A}_{11}^2 + \mathrm{A}_{12}^2.
\end{eqnarray*}
By (the proof of) Proposition \ref{Proposition-BMO} we conclude
that $$\mathrm{A}_{11} \le c(p,\theta) \, \max \Big\{
\|x\|_{\Ha_{p'}^{p'}(\mathcal{M})},
\|x\|_{[\Ha_{p'}^{p'}(\mathcal{M}), h_{p'}^c(\mathcal{M})]_\theta}
\Big\}.$$ On the other hand, \eqref{Eq-Cond-Char-Hardy} yields the
estimate $$\mathrm{A}_{12} \lesssim
\|x\|_{[\Ha_{p'}^{p'}(\mathcal{M}),
h_{p'}^c(\mathcal{M})]_\theta}.$$ The case $4 \le w_c \le \infty$
is similar and yields $$\mathrm{A}_{1} \le c(p,\theta) \, \max
\Big\{ \|x\|_{\Ha_{p'}^{p'}(\mathcal{M})},
\|x\|_{[h_{p'}^r(\mathcal{M}), \Ha_{p'}^{p'}(\mathcal{M})]_\theta}
\Big\}.$$ Therefore, in the general case we conclude
$$\mathrm{A}_{1} \le c(p,\theta) \, \max \Big\{
\|x\|_{\Ha_{p'}^{p'}(\mathcal{M})}, \|x\|_{[h_{p'}^r(\mathcal{M}),
\Ha_{p'}^{p'}(\mathcal{M})]_\theta},
\|x\|_{[\Ha_{p'}^{p'}(\mathcal{M}), h_{p'}^c(\mathcal{M})]_\theta}
\Big\}.$$

\noindent \textsc{Step 2}. The same arguments as in Step 1 yield
the right estimate for $\mathrm{A}_2$ in the case $4 \le w_r \le
\infty$ and for $\mathrm{A}_3$ in the case $4 \le w_c \le \infty$.
Of course, there is an obvious symmetry between both cases so that
we only prove the estimate for $\mathrm{A}_3$ in the case $4 \le
w_c \le \infty$. We have $$\mathrm{A}_3^2 \le \Big( \summ_k
\|d_k(a_c)\|_{w_c/2}^{w_c/2} \Big)^{2/w_c} \Big( \summ_k \big\|
\mathcal{E}_{k-1}(a_r)^{\frac12} dx_k
\big\|_{2(w_c/2)'}^{2(w_c/2)'} \Big)^{1/(w_c/2)'}.$$ The first
term on the right is controlled by $\|a_c\|_{w_c/2}$ while
$$\frac{1}{2(w_c/2)'} = \frac{1}{2} - \frac{1}{w_c} = \frac{1}{2}
- \frac{\theta}{q} = 1 - \frac{1}{p} + \frac{1-\theta}{q}.$$ That
is, the roles of $\theta$ and $1-\theta$ have exchanged with
respect to the situation in Step 1 above. Therefore, according to
the equivalence \eqref{Eq-Cond-Char-Hardy} we easily conclude that
$$\mathrm{A}_3 \lesssim \|x\|_{[h_{p'}^r(\mathcal{M}),
\Ha_{p'}^{p'}(\mathcal{M})]_\theta}.$$ When $4 \le w_r \le \infty$
we obtain the estimate $$\mathrm{A}_2 \lesssim
\|x\|_{[\Ha_{p'}^{p'}(\mathcal{M}),
h_{p'}^c(\mathcal{M})]_\theta}.$$

\noindent \textsc{Step 3}. Now we estimate $\mathrm{A}_2$ for $2 <
w_r < 4$ and $\mathrm{A}_3$ for $2 < w_c < 4$. Again by symmetry,
we only prove the estimate for $\mathrm{A}_2$. The proof of this
estimate follows the argument given in Step 2 of Proposition
\ref{Proposition-BMO}. By the noncommutative Burkholder
inequality, we may find a decomposition $d_k(a_r) = d_k(\alpha_r)
+ d_k(\beta_r) + d_k(\gamma_r)$ into three martingales satisfying
\eqref{Eq-NCBurk-Est} with $(\alpha,\beta,\gamma) \rightsquigarrow
(\alpha_r, \beta_r, \gamma_r)$ and $v \rightsquigarrow w_r$. Then
we have $\mathrm{A}_2^2 \le \mathrm{A}_2(\alpha)^2 +
\mathrm{A}_2(\beta)^2 + \mathrm{A}_2(\gamma)^2$ with
\begin{eqnarray*}
\mathrm{A}_2(\alpha)^2 & = & \Big| \summ_k \mbox{tr} \big( dx_k
\mathcal{E}_{k-1}(a_c) dx_k^* d_k(\alpha_r) \big) \Big|, \\
\mathrm{A}_2(\beta)^2 & = & \Big| \summ_k \mbox{tr} \big( dx_k
\mathcal{E}_{k-1}(a_c) dx_k^* d_k(\beta_r) \big) \Big|, \\
\mathrm{A}_2(\gamma)^2 & = & \Big| \summ_k \mbox{tr} \big( dx_k
\mathcal{E}_{k-1}(a_c) dx_k^* d_k(\gamma_r) \big) \Big|.
\end{eqnarray*}
The term $\mathrm{A}_2(\alpha)$ is estimated as in Step 2, due to
the first inequality in \eqref{Eq-NCBurk-Est}. The terms
$\mathrm{A}_2(\beta)$ and $\mathrm{A}_2(\gamma)$ are estimated in
the same way so that we only show how to estimate
$\mathrm{A}_2(\beta)$. We proceed as in Proposition
\ref{Proposition-BMO} again and obtain
\begin{eqnarray*}
\mathrm{A}_2(\beta)^2 & = & \Big| \summ_k \mbox{tr} \big(
\mathcal{E}_{k-1} (d_k(\beta_r) dx_k \mathcal{E}_{k-1}(a_c)
dx_k^*) \big) \Big| \\ & \le & \Big\| \big( \summ_k
\mathcal{E}_{k-1}
(d_k(\beta_r) d_k(\beta_r)^*) \big)^{\frac12} \Big\|_{w_r/2} \\
& \times & \Big\| \big( \summ_k \mathcal{E}_{k-1} (dx_k
\mathcal{E}_{k-1}(a_c) dx_k^*)^2 \big)^{\frac12} \Big\|_{(w_r/2)'} \\
& \le & c_{w_r} \|a_r\|_{w_r/2} \Big\| \big( \summ_k
\mathcal{E}_{k-1} (dx_k \mathcal{E}_{k-1}(a_c) dx_k^*)^2
\big)^{\frac12} \Big\|_{(w_r/2)'}.
\end{eqnarray*}
To estimate the last term on the right we define $z_k = (dx_k
\mathcal{E}_{k-1}(a_c) dx_k^*)^{\frac12}$. Recalling, our
assumption $2 < w_r < 4$, we have $2 < (w_r/2)' < \infty$ and
Lemma 5.2 of \cite{JX} gives for $t = \frac12 (w_r/2)'$ $$\Big\|
\summ_k \mathcal{E}_{k-1} (z_k^4) \Big\|_t \le \Big\| \summ_k
\mathcal{E}_{k-1} (z_k^2) \Big\|_{2t}^{\frac{2(t-1)}{2t-1}} \Big(
\summ_k \|z_k\|_{4t}^{4t} \Big)^{\frac{1}{2t-1}}.$$ In our
situation, this implies
\begin{eqnarray*}\Big\| \big( \summ_k \mathcal{E}_{k-1} (dx_k
\mathcal{E}_{k-1}(a_c) dx_k^*)^2 \big)^{\frac12}
\Big\|_{(\frac{w_r}{2})'} \nonumber \!\!\! & \le & \!\!\! \Big\|
\summ_k \mathcal{E}_{k-1} \big( dx_k \mathcal{E}_{k-1}(a_c) dx_k^*
\big) \Big\|_{(w_r/2)'}^{\frac{t-1}{2t-1}}
\\ \!\!\! & \times & \!\!\! \Big( \summ_k \big\| dx_k
\mathcal{E}_{k-1}(a_c) dx_k^* \big\|_{(w_r/2)'}^{(w_r/2)'}
\Big)^{\frac{1}{4t-2}}.
\end{eqnarray*}
If we take $$\mathrm{M}_\theta = \max \Big\{
\|x\|_{\Ha_{p'}^{p'}}, \|x\|_{[h_{p'}^r, \Ha_{p'}^{p'}]_\theta},
\|x\|_{[\Ha_{p'}^{p'}, h_{p'}^c]_\theta}, \|x\|_{[h_{p'}^r,
h_{p'}^c]_\theta} \Big\},$$ we have already seen in Step 1 above
that
$$\Big( \summ_k \big\| dx_k \mathcal{E}_{k-1}(a_c) dx_k^*
\big\|_{(w_r/2)'}^{(w_r/2)'} \Big)^{\frac{1}{4t}} \lesssim
\mathrm{M}_\theta.$$ We claim that
\begin{equation} \label{Eq-Thelclaim}
\Big\| \summ_k \mathcal{E}_{k-1} \big( dx_k \mathcal{E}_{k-1}(a_c)
dx_k^* \big) \Big\|_{(w_r/2)'} \le \|x\|_{[h_{p'}^r,
h_{p'}^c]_\theta}^2 \le \mathrm{M}_\theta^2.
\end{equation}
If we prove \eqref{Eq-Thelclaim} then it is easy to see that
$\mathrm{A}_2(\beta) \lesssim \mathrm{M}$ and the estimate for
$\mathrm{A}_2$ will be completed. Arguing as in Lemma
\ref{Lemma-BMO6}, it is not difficult to check that the left hand
side of \eqref{Eq-Thelclaim} does interpolate. Hence, it suffices
to estimate the extremal cases. When $\theta = 0$, we have
$(w_r,w_c) = (q,\infty)$ and $(w_r/2)' = p'/2$. Consequently, we
find $$\Big\| \summ_k \mathcal{E}_{k-1} \big( dx_k
\mathcal{E}_{k-1}(a_c) dx_k^* \big) \Big\|_{p'/2} \le
\|a_c\|_\infty \Big\| \summ_k \mathcal{E}_{k-1} (dx_k dx_k^*)
\Big\|_{p'/2} \le \|x\|_{h_{p'}^r(\mathcal{M})}^2.$$ When $\theta
= 1$ we have $(w_r,w_c) = (\infty,q)$ and $(w_r/2)' = 1$ so that
\begin{eqnarray*}
\Big\| \summ_k \mathcal{E}_{k-1} \big( dx_k \mathcal{E}_{k-1}(a_c)
dx_k^* \big) \Big\|_1 & = & \summ_k \mbox{tr} \big(
\mathcal{E}_{k-1}(a_c) dx_k^* dx_k \big) \\ & = & \summ_k
\mbox{tr} \big( a_c \mathcal{E}_{k-1}(dx_k^* dx_k) \big) \\ & \le
& \|a_c\|_{\frac{q}{2}} \Big\| \summ_k \mathcal{E}_{k-1} (dx_k^*
dx_k) \Big\|_{\frac{p'}{2}} \le \|x\|_{h_{p'}^c(\mathcal{M})}^2.
\end{eqnarray*}
Note that the first identity assumes that $a_c$ is positive and
this is not necessarily true on the boundary. However, decomposing
into a linear combination of four positive elements and allowing
an additional constant $2$, we may and do assume positivity.
Therefore, \eqref{Eq-Thelclaim} follows from the tree lines lemma.
A detailed reading of the proof gives now the constant
$c'(p,\theta)$ stated above. This completes the proof.
\end{proof}

\section{Interpolation of $2$-term intersections}

Let us fix some notation which will be used in the sequel. As
usual, we begin by fixing a von Neumann algebra $\mathcal{M}$ and
a von Neumann subalgebra $\mathcal{N}$ with conditional
expectation $\mathsf{E}: \mathcal{M} \to \mathcal{N}$. Given $1
\le q \le p \le \infty$ and a positive integer $n \ge 1$, the main
spaces in this paragraph will be the following
\begin{eqnarray*}
\mathcal{R}_{2p,q}^n(\mathcal{M},\mathsf{E}) & = &
n^{\frac{1}{2p}} L_{2p}(\mathcal{M}) \cap n^{\frac{1}{2q}}
L_{(\frac{2pq}{p-q},\infty)}^{2p}(\mathcal{M}, \mathsf{E}), \\
\mathcal{C}_{2p,q}^n \, (\mathcal{M},\mathsf{E}) & = &
n^{\frac{1}{2p}} L_{2p}(\mathcal{M}) \cap n^{\frac{1}{2q}}
L_{(\infty,\frac{2pq}{p-q})}^{2p}(\mathcal{M}, \mathsf{E}).
\end{eqnarray*}
In order to study these spaces we need to introduce some
terminology. We set $\mathsf{A}_k$ to be $\mathcal{M} \oplus
\mathcal{M}$ for $1 \le k \le n$. Then we consider the reduced
amalgamated free product $\mathcal{A} = *_{\mathcal{N}}
\mathsf{A}_k$ where the conditional expectation
$\mathsf{E}_{\mathcal{N}}: \mathcal{A} \rightarrow \mathcal{N}$,
defined in Paragraph \ref{Subsubsection5.1.1} as
$\mathsf{E}_{\mathcal{N}}(a) = \mathcal{Q}_{\emptyset} a
\mathcal{Q}_{\emptyset}$, has the following form when restricted
to $\mathsf{A}_k$
$$\mathsf{E}_{\mathcal{N}}(x_1,x_2) = \frac{1}{2} \big(
\mathsf{E}(x_1) + \mathsf{E}(x_2) \big).$$ Given a \emph{n.f.}
state $\varphi: \mathcal{N} \rightarrow \C$, let $\varphi_2:
\mathcal{M} \oplus \mathcal{M} \rightarrow \C$ be the \emph{n.f.}
state $$\varphi_2(x_1,x_2) = \frac{1}{2} \big(
\varphi(\mathsf{E}(x_1)) + \varphi(\mathsf{E}(x_2)) \big) =
\varphi \big( \mathsf{E}_{\mathcal{N}}(x_1,x_2) \big).$$ We shall
write $\mathcal{A}_{\oplus n}$ for the direct sum $\mathcal{A}
\oplus \mathcal{A} \oplus \ldots \oplus \mathcal{A}$ with $n$
terms. If $\phi$ stands for the free product state on
$\mathcal{A}$, we consider the \emph{n.f.} state $\phi_n:
\mathcal{A}_{\oplus n} \rightarrow \C$ and the conditional
expectation $\mathcal{E}_n: \mathcal{A}_{\oplus n} \rightarrow
\mathcal{A}$ given by $$\phi_n \big( \sum_{k=1}^n x_k \otimes
\delta_k \big) = \frac{1}{n} \sum_{k=1}^n \phi(x_k) \qquad
\mbox{and} \qquad \mathcal{E}_{\oplus_n} \big( \sum_{k=1}^n x_k
\otimes \delta_k \big) = \frac{1}{n} \sum_{k=1}^n x_k.$$ Let
$\pi_k: \mathsf{A}_k \to \mathcal{A}$ denote the embedding of
$\mathsf{A}_k$ into $\mathcal{A}$ as defined at the beginning of
this chapter. Moreover, given $x \in \mathcal{M}$ we shall write
$x_k$ as an abbreviation of $\pi_k(x,-x)$. Note that $x_k$ is a
mean-zero element for $1 \le k \le n$. In the following we shall
use with no further comment the identities
$$\mathsf{E}_{\mathcal{N}}(x_kx_k^*) = \mathsf{E}(xx^*) \quad
\mbox{and} \quad \mathsf{E}_{\mathcal{N}}(x_k^*x_k) =
\mathsf{E}(x^*x).$$ Let us consider the following map
\begin{equation} \label{Equation-Map-u}
u: x \in \mathcal{M} \mapsto \sum_{k=1}^n x_k \otimes \delta_k \in
\mathcal{A}_{\oplus n}.
\end{equation}
Moreover, if $d_{\widehat{\varphi}}$ denotes the density
associated to the \emph{n.f.} state $\widehat{\varphi} = \varphi
\circ \mathsf{E}$ on $\mathcal{M}$ and $d_{\phi_n}$ stands for the
density associated to the \emph{n.f.} state $\phi_n$ on
$\mathcal{A}_{\oplus_n}$, we may extend the definition of $u$ to
other indices by taking $$u \big( d_{\widehat{\varphi}}^{\frac1p}
x \big) = d_{\phi_n}^{\frac1p} u(x) \quad \mbox{and} \quad u \big(
x d_{\widehat{\varphi}}^{\frac1p} \big) = u(x)
d_{\phi_n}^{\frac1p}.$$ In the following we shall consider the
filtration on the von Neumann algebra $\mathcal{A}$ given by
$\mathcal{A}_k = \mathsf{A}_1 *_\mathcal{N} \mathsf{A}_2
*_\mathcal{N} \cdots *_\mathcal{N} \mathsf{A}_k$. In particular,
for any $x \in L_p(\mathcal{M})$ we obtain that $u(x) =
(x_1,x_2,\ldots, x_n)$ is the sequence of martingale differences
of $\summ_k x_k$.

\begin{lemma} \label{Lemma-Complemented-Isomorphism}
If $1 \le p < \infty$, the following mappings are isomorphisms
onto complemented subspaces
$$\begin{array}{rrcl} u: & \mathcal{R}_{2p,1}^n(\mathcal{M},
\mathsf{E}) & \rightarrow & \Ha_{2p}^r(\mathcal{A}), \\ [5pt] u: &
\mathcal{C}_{2p,1}^n(\mathcal{M}, \mathsf{E}) & \rightarrow &
\Ha_{2p}^c(\mathcal{A}). \end{array}$$ Moreover, the constants are
independent of $n$ and remain bounded as $p \to 1$.
\end{lemma}

\begin{proof} Let us observe that
\begin{eqnarray*}
\mathcal{R}_{2p,1}^n(\mathcal{M},\mathsf{E}) & = &
n^{\frac{1}{2p}} L_{2p}(\mathcal{M}) \cap \sqrt{n}
L_{2p}^r(\mathcal{M}, \mathsf{E}), \\
\mathcal{C}_{2p,1}^n \, (\mathcal{M},\mathsf{E}) & = &
n^{\frac{1}{2p}} L_{2p}(\mathcal{M}) \cap \sqrt{n}
L_{2p}^c(\mathcal{M}, \mathsf{E}).
\end{eqnarray*}
Given $x \in \mathcal{R}_{2p,1}^n(\mathcal{M}, \mathsf{E})$,
Corollary \ref{Corollary-Voiculescu} gives $$\big\| u(x)
\big\|_{\Ha_{2p}^r(\mathcal{A})} = \Big\| \big( \sum_{k=1}^n x_k
x_k^* \big)^{1/2} \Big\|_{2p} \sim \Big( \sum_{k=1}^n
\|x_k\|_{2p}^{2p} \Big)^{2p} + \Big\| \big( \sum_{k=1}^n
\mathsf{E}_{\mathcal{N}} (x_k x_k^*) \big)^{\frac12}
\Big\|_{2p}.$$ In other words, we have
$$\big\| u(x) \big\|_{\Ha_{2p}^r(\mathcal{A})} \sim
n^{\frac{1}{2p}} \|x\|_{L_{2p}(\mathcal{M})} + \sqrt{n} \,
\|x\|_{L_{2p}^r(\mathcal{M}, \mathsf{E})}.$$ This proves that $u:
\mathcal{R}_{2p,1}^n(\mathcal{M},\mathsf{E}) \to
\Ha_{2p}^r(\mathcal{A})$ is an isomorphism onto its image with
relevant constants independent of $p,n$. A similar argument yields
to the same conclusion for the column case. To prove the
complementation, we recall that $$\Ha_{2p}^r(\mathcal{A})^* \simeq
\Ha_{(2p)'}^r(\mathcal{A}) \quad \mbox{for} \quad 1 \le p <
\infty$$ (with relevant constants which remain bounded as $p \to
1$) and consider the map
$$\omega: x \in \frac{1}{n^{\frac{1}{2p}}}L_{(2p)'}(\mathcal{M}) +
\frac{1}{\sqrt{n}} \, L_{(2p)'}^r(\mathcal{M}, \mathsf{E})
\longmapsto \frac{1}{n} \sum_{k=1}^n x_k \otimes \delta_k \in
\Ha_{(2p)'}^r(\mathcal{A}).$$ Let $d_{\widehat{\varphi}}$ be the
density associated to $\varphi \circ \mathsf{E}$. Assume by
approximation that $$x = \alpha d_{\widehat{\varphi}}^{1/(2p)'}
a$$ for some $(\alpha,a) \in \mathcal{N} \times \mathcal{M}$.
Then, taking $d_\phi$ to be the density associated to the state
$\phi$ on $\mathcal{A}$ and defining $a_k = \pi_k(a,-a)$, the
following estimate holds by Theorem 7.1 in \cite{JX}
\begin{eqnarray*}
\big\| \omega(x) \big\|_{\Ha_{(2p)'}^r(\mathcal{A})} & = &
\frac{1}{n} \, \Big\| \alpha d_\phi^{1/(2p)'} \Big( \sum_{k=1}^n
a_k a_k^* \Big) d_\phi^{1/(2p)'} \alpha^*
\Big\|_{L_{(2p)'/2}(\mathcal{A})}^{1/2}
\\ & \le & \frac{1}{n} \, \Big\| \alpha d_{\varphi}^{1/(2p)'}
\Big( \sum_{k=1}^n \mathsf{E}_{\mathcal{N}} \big( a_k a_k^* \big)
\Big) d_{\varphi}^{1/(2p)'} \alpha^*
\Big\|_{L_{(2p)'/2}(\mathcal{N})}^{1/2}.
\end{eqnarray*}
This gives
\begin{equation} \label{Equation-Expectation-p<1}
\big\| \omega(x) \big\|_{\Ha_{(2p)'}^r(\mathcal{A})} \le
\frac{1}{\sqrt{n}} \, \|x\|_{L_{(2p)'}^r(\mathcal{M},
\mathsf{E})}.
\end{equation}
On the other hand, the inequality
\begin{equation} \label{Equation-Estimate-ComplexInt}
\big\| \omega(x) \big\|_{\Ha_{(2p)'}^r(\mathcal{A})} \le
\frac{1}{n^{\frac{1}{2p}}} \, \|x\|_{L_{(2p)'}(\mathcal{M})}
\end{equation}
follows by the complex interpolation method between the (trivial)
extremal cases for $p=1$ and $p=\infty$. The estimates
(\ref{Equation-Expectation-p<1}) and
(\ref{Equation-Estimate-ComplexInt}) show that the map $\omega$ is
a contraction. Note also that $$\big\langle u(x_1), \omega(x_2)
\big\rangle = \frac{1}{n} \sum_{k=1}^n \mbox{tr}_{\mathcal{A}}
\big( x_{1k}^* x_{2k}^{} \big) = \frac{1}{n} \sum_{k=1}^n
\mbox{tr}_{\mathcal{M}} (x_1^* x_2) = \langle x_1, x_2 \rangle.$$
In particular, since we have
$$\mathcal{R}_{2p,1}^n(\mathcal{M},\mathsf{E}) = \Big(
\frac{1}{n^{\frac{1}{2p}}}L_{(2p)'}(\mathcal{M}) +
\frac{1}{\sqrt{n}} \, L_{(2p)'}^r(\mathcal{M}, \mathsf{E})
\Big)^*,$$ it turns out that the map $\omega^* u$ is the identity
on $\mathcal{R}_{2p,1}^n(\mathcal{M},\mathsf{E})$ and $u\omega^*$
is a bounded projection onto the image of $u$ with constants
independent of $n$ and bounded for $p \sim 1$. This completes the
proof in the row case. The column case is the same.
\end{proof}

Before proving our interpolation result, we need to consider a
variation of Lemma \ref{Lemma-Complemented-Isomorphism}. Namely,
we know from \cite{JX} that $\Ha_{2p}^r(\mathcal{A}) \simeq
L_{2p}^r\mathcal{MO}(\mathcal{A})$ for $1 < p < \infty$ and with
constants depending on $p$ which diverge as $p \to \infty$. We
claim however that Lemma \ref{Lemma-Complemented-Isomorphism}
still holds in this setting with bounded constants as $p \to
\infty$.

\begin{lemma} \label{Lemma-LpMO}
If $1 < p \le \infty$, the following mappings are isomorphisms
onto complemented subspaces
$$\begin{array}{rrcl} u: & \mathcal{R}_{2p,1}^n(\mathcal{M},
\mathsf{E}) & \rightarrow & L_{2p}^r \mathcal{MO}(\mathcal{A}), \\
[5pt] u: & \mathcal{C}_{2p,1}^n(\mathcal{M}, \mathsf{E}) &
\rightarrow & L_{2p}^c\mathcal{MO}(\mathcal{A}).
\end{array}$$ Moreover, the constants are independent of $n$ and
remain bounded as $p \to \infty$.
\end{lemma}

\begin{proof}
The noncommutative Doob's inequality \cite{J1} gives
$$\big\| u(x) \big\|_{L_{2p}^r\mathcal{MO}} = \Big\| \sup_{1 \le m
\le n} \mathcal{E}_m \big( \sum_{k = m}^n x_k x_k^* \big)
\Big\|_p^{\frac12} \le \gamma_p \Big\| \sum_{k=1}^n x_k x_k^*
\Big\|_p^{\frac12}.$$ Note that $\gamma_p \to \infty$ as $p \to 1$
but $\gamma_p \le 2$ for $p \ge 2$. On the other hand, we may
estimate the term on the right by using the free Rosenthal
inequality (see Corollary \ref{Corollary-Voiculescu}) one more
time $$\Big\| \sum_{k=1}^n x_k x_k^* \Big\|_p^{\frac12} \sim \Big(
\sum_{k=1}^n \|x_k\|_{2p}^{2p} \Big)^{2p} + \Big\| \big(
\sum_{k=1}^n \mathsf{E}_{\mathcal{N}} (x_k x_k^*) \big)^{\frac12}
\Big\|_{2p} = \|x\|_{\mathcal{R}_{2p,1}^n(\mathcal{M},
\mathsf{E})}.$$ This shows that $u:
\mathcal{R}_{2p,1}^n(\mathcal{M}, \mathsf{E}) \to L_{2p}^r
\mathcal{MO}(\mathcal{A})$ is bounded with constant $\gamma_p$. To
prove complementation and the boundedness of the inverse we
proceed by duality as in Lemma
\ref{Lemma-Complemented-Isomorphism}. Indeed, using the map $w$
one more time and recalling that $$L_{2p}^r \mathcal{MO}
(\mathcal{A}) \simeq \Ha_{(2p)'}^r(\mathcal{A})^*$$ with constants
which remain bounded as $ p \to \infty$ (\emph{c.f.} Theorem 4.1
in \cite{JX}), we may follows verbatim the proof of Lemma
\ref{Lemma-Complemented-Isomorphism} to conclude that the map
$\omega^* u$ is the identity on
$\mathcal{R}_{2p,1}^n(\mathcal{M},\mathsf{E})$ and $u\omega^*$ is
a bounded projection onto the image of $u$ with constants
independent of $n$ and bounded for $p \sim \infty$. This completes
the proof in the row case. The column case follows in the same
way.
\end{proof}

\begin{theorem} \label{Theorem-Intersection1}
If $1 \le p \le \infty$ and $1/q = 1-\theta + \theta/p$, we have
\begin{eqnarray*}
\big[ \mathcal{R}_{2p,1}^n(\mathcal{M},\mathsf{E}),
\mathcal{R}_{2p,p}^n(\mathcal{M},\mathsf{E}) \big]_{\theta} &
\simeq & \mathcal{R}_{2p,q}^n(\mathcal{M},\mathsf{E}), \\ \big[
\hskip1.5pt \mathcal{C}_{2p,1}^n \, (\mathcal{M},\mathsf{E}), \,
\mathcal{C}_{2p,p}^n \, (\mathcal{M},\mathsf{E}) \big]_{\theta} &
\simeq & \, \mathcal{C}_{2p,q}^n \, (\mathcal{M},\mathsf{E}),
\end{eqnarray*}
isomorphically with relevant constant $c(p,q)$ independent of $n$
and such that $$c(p,q) \lesssim \sqrt{\frac{p-q}{pq+q-p}} \quad
\mbox{as} \quad (p,q) \to (\infty,1).$$
\end{theorem}

\begin{proof} By Corollary
\ref{Corollary-Interpolation-Conditional2} i), we have contractive
inclusions
\begin{eqnarray*}
\big[ \mathcal{R}_{2p,1}^n(\mathcal{M},\mathsf{E}),
\mathcal{R}_{2p,p}^n(\mathcal{M},\mathsf{E}) \big]_{\theta} &
\subset & \mathcal{R}_{2p,q}^n(\mathcal{M},\mathsf{E}), \\
\big[ \hskip1.5pt \mathcal{C}_{2p,1}^n \,
(\mathcal{M},\mathsf{E}), \, \mathcal{C}_{2p,p}^n \,
(\mathcal{M},\mathsf{E}) \big]_{\theta} & \subset & \,
\mathcal{C}_{2p,q}^n \, (\mathcal{M},\mathsf{E}).
\end{eqnarray*}
To prove the converse, we consider again the map given in
(\ref{Equation-Map-u}). It is clear that $$\big\| u(x)
\big\|_{\Ha_{2p}^{2p}(\mathcal{A})} = \Big( \sum_{k=1}^n
\|x_k\|_{2p}^{2p} \Big)^{\frac{1}{2p}} = n^{\frac{1}{2p}}
\|x\|_{2p}.$$ This shows that $u: n^{\frac{1}{2p}}
L_{2p}(\mathcal{M}) \to \Ha_{2p}^{2p}(\mathcal{A})$ is an
isometric isomorphism. Moreover, arguing as in the proof of Lemma
\ref{Lemma-Complemented-Isomorphism} we easily obtain that the
image of $u$ is contractively complemented. This observation
together with Lemmas \ref{Lemma-Complemented-Isomorphism} and
\ref{Lemma-LpMO} give rise to the following equivalences
$$\begin{array}{rcll} \|x\|_{[\mathcal{R}_{2p,1}^n (\mathcal{M},
\mathsf{E}), \mathcal{R}_{2p,p}^n (\mathcal{M},
\mathsf{E})]_\theta} & \sim & \big\|u(x) \big\|_{[\Ha_{2p}^r
(\mathcal{A}), \Ha_{2p}^{2p}(\mathcal{A}) ]_\theta} & \mbox{for
small} \ p,
\\ [5pt] \|x\|_{[\mathcal{R}_{2p,1}^n (\mathcal{M}, \mathsf{E}),
\mathcal{R}_{2p,p}^n (\mathcal{M}, \mathsf{E})]_\theta} & \sim &
\big\|u(x) \big\|_{[L_{2p}^r \mathcal{MO} (\mathcal{A}),
\Ha_{2p}^{2p}(\mathcal{A}) ]_\theta} & \mbox{for large} \, \ p,
\end{array}$$
with constants independent of $p,q,n$. On the other hand, Berg's
theorem gives isometric inclusions
\begin{eqnarray*}
\big[ \Ha_{2p}^r (\mathcal{A}), \Ha_{2p}^{2p}(\mathcal{A})
\big]_\theta & \subset & \big[ \Ha_{2p}^r (\mathcal{A}),
\Ha_{2p}^{2p}(\mathcal{A}) \big]^\theta, \\ \big[ L_{2p}^r
\mathcal{MO} (\mathcal{A}), \Ha_{2p}^{2p}(\mathcal{A})
\big]_\theta & \subset & \big[ L_{2p}^r \mathcal{MO}
(\mathcal{A}), \Ha_{2p}^{2p}(\mathcal{A}) \big]^\theta.
\end{eqnarray*}
Now we can use duality and obtain
\begin{eqnarray*}
\big[ \Ha_{2p}^r (\mathcal{A}), \Ha_{2p}^{2p}(\mathcal{A})
\big]^\theta & \simeq & \big[ \Ha_{(2p)'}^r (\mathcal{A}),
\Ha_{(2p)'}^{(2p)'}(\mathcal{A}) \big]_\theta^* \quad \mbox{for} \
1 \le p < \infty, \\ \big[ L_{2p}^r \mathcal{MO} (\mathcal{A}),
\Ha_{2p}^{2p}(\mathcal{A}) \big]^\theta & \simeq & \big[
\Ha_{(2p)'}^r (\mathcal{A}), \Ha_{(2p)'}^{(2p)'}(\mathcal{A})
\big]_\theta^* \quad \mbox{for} \ 1 < p \le \infty,
\end{eqnarray*}
where the constants in the first isomorphism remain bounded as $p
\to 1$ and the constants in the second one remain bounded as $p
\to \infty$. Therefore, recalling the terminology used in the
previous section $$\big[ \Ha_{(2p)'}^r (\mathcal{A}),
\Ha_{(2p)'}^{(2p)'}(\mathcal{A}) \big]_\theta^* =
\mathcal{Z}_{(2p)'}^r (\mathcal{A}, 1-\theta)$$ and taking
adjoints, we obtain
\begin{eqnarray*}
\|x\|_{[\mathcal{R}_{2p,1}^n (\mathcal{M}, \mathsf{E}),
\mathcal{R}_{2p,p}^n (\mathcal{M}, \mathsf{E})]_\theta} & \sim &
\|u(x)\|_{\mathcal{Z}_{(2p)'}^r (\mathcal{A}, 1-\theta)} =
\mathrm{A}_p^r(\mathcal{A},\theta),
\\ [5pt] \|x\|_{[ \, \mathcal{C}_{2p,1}^n (\mathcal{M},
\mathsf{E}) \, , \, \mathcal{C}_{2p,p}^n (\mathcal{M}, \mathsf{E})
\, ]_\theta} & \sim & \|u(x)\|_{\mathcal{Z}_{(2p)'}^c(\mathcal{A},
1-\theta)} = \mathrm{A}_p^c(\mathcal{A},\theta),
\end{eqnarray*}
with constants independent of $p,q,n$. According to Proposition
\ref{Proposition-BMO} we have
$$\mathrm{A}_p^r(\mathcal{A},\theta) \le c(p,q) \, \max \Big\{
\|u(x)\|_{\Ha_{2p}^{2p}(\mathcal{A})},
\|u(x)\|_{[\Ha_{2p}^{2p}(\mathcal{A}),
h_{2p}^r(\mathcal{A})]_{1-\theta}} \Big\}.$$ Let us estimate the
two terms on the right $$\|u(x)\|_{\Ha_{2p}^{2p}(\mathcal{A})} =
\Big( \sum_{k=1}^n \|x_k\|_{2p}^{2p} \Big)^{\frac{1}{2p}} =
n^{\frac{1}{2p}} \|x\|_{2p} \le
\|x\|_{\mathcal{R}_{2p,q}^n(\mathcal{M}, \mathsf{E})}.$$ For the
second term we observe that
\begin{eqnarray*}
\|u(x)\|_{h_{2p}^r(\mathcal{A})} & = & \Big\| \Big( \sum_{k=1}^n
\mathcal{E}_{k-1}(x_k x_k^*) \Big)^{\frac12} \Big\|_{2p} \\ & = &
\Big\| \Big( \sum_{k=1}^n \mathsf{E}_{\mathcal{N}}(x_k x_k^*)
\Big)^{\frac12} \Big\|_{2p} \\ & = & \Big\| \Big( \sum_{k=1}^n
\mathsf{E}(x x^*) \Big)^{\frac12} \Big\|_{2p} = \sqrt{n} \, \big\|
\mathsf{E}(x x^*) \big\|_{2p},
\end{eqnarray*}
where the second inequality follows by freeness. By complex
interpolation and Corollary
\ref{Corollary-Interpolation-Conditional2} we conclude
\begin{eqnarray*}
\|u(x)\|_{[h_{2p}^r(\mathcal{A}),
\Ha_{2p}^{2p}(\mathcal{A})]_{\theta}} & \le &
n^{\frac{1-\theta}{2} + \frac{\theta}{2p}} \,
\|x\|_{L_{(\frac{2p'}{1-\theta},\infty)}^{2p}(\mathcal{M},
\mathsf{E})} \\ & = & n^{\frac{1}{2q}} \,
\|x\|_{L_{(\frac{2pq}{p-q},\infty)}^{2p}(\mathcal{M}, \mathsf{E})}
\le \|x\|_{\mathcal{R}_{2p,q}^n(\mathcal{M}, \mathsf{E})}.
\end{eqnarray*}
In summary, we have proved that
$$\|x\|_{\mathcal{R}_{2p,q}^n(\mathcal{M}, \mathsf{E})} \le
\|x\|_{[\mathcal{R}_{2p,1}^n(\mathcal{M}, \mathsf{E}),
\mathcal{R}_{2p,p}^n(\mathcal{M}, \mathsf{E})]_\theta} \le c(p,q)
\, \|x\|_{\mathcal{R}_{2p,q}^n(\mathcal{M}, \mathsf{E})}$$ where
(recalling that $p \rightsquigarrow (2p)'$ and $\theta
\rightsquigarrow 1-\theta$), it follows from Proposition
\ref{Proposition-BMO} that $$c(p,q) \sim \sqrt{\frac{p-q}{pq+q-p}}
\quad \mbox{as} \quad (p,q) \to (\infty,1).$$ This proves the
assertion for rows. The column case follows by taking adjoints.
\end{proof}

\begin{remark}
\emph{At the time of this writing, we do not know whether or not
the relevant constants in Theorem \ref{Theorem-Intersection1} are
uniformly bounded in $p$ and $q$. Our constants are not uniformly
bounded due to fact that we use the noncommutative Burkholder
inequality from \cite{JX} in our approach. We take this
opportunity to pose this question as a problem for the interested
reader. The same question applies to Theorems
\ref{Theorem-Intersection2} and \ref{Theorem-Interpolation-J}
below.}
\end{remark}

\section{Interpolation of $4$-term intersections}

In this section we study the interpolation spaces between
$\mathcal{R}_{2p,1}^n(\mathcal{M}, \mathsf{E})$ and
$\mathcal{C}_{2p,1}^n(\mathcal{M}, \mathsf{E})$. Of course, as it
is to be expected, our main tools will be the free Rosenthal
inequalities and the two-sided estimates for BMO type norms. We
shall use below the constant $c'(p,\theta)$ in Proposition
\ref{Proposition-BMO2}.

\begin{theorem} \label{Theorem-Intersection2}
If $1 \le p \le \infty$, we have $$\big[
\mathcal{R}_{2p,1}^n(\mathcal{M},\mathsf{E}),
\mathcal{C}_{2p,1}^n(\mathcal{M},\mathsf{E}) \big]_{\theta} \
\simeq \bigcap_{u,v \in \{ 2p', \infty \}}^{\null}
n^{\frac{1-\theta}{u} + \frac{1}{2p} + \frac{\theta}{v}} \,
L_{(\frac{u}{1-\theta},\frac{v}{\theta})}^{2p} (\mathcal{M},
\mathsf{E})$$ isomorphically with relevant constant controlled by
$c'(p,\theta)$ and independent of $n$.
\end{theorem}

\begin{proof}
According to Corollary \ref{Corollary-Interpolation-Conditional2},
we have
\begin{eqnarray*} \big[ \sqrt{n} L_{2p}^r(\mathcal{M},
\mathsf{E}), n^{\frac{1}{2p}} L_{2p}(\mathcal{M}) \big]_{\theta} &
= & n^{\frac{1-\theta}{2} + \frac{\theta}{2p}}
L_{(\frac{2p'}{1-\theta}, \infty)}^{2p}(\mathcal{M}, \mathsf{E}),
\\ \big[ n^{\frac{1}{2p}} L_{2p}(\mathcal{M}), \sqrt{n}
L_{2p}^c(\mathcal{M}, \mathsf{E}) \big]_{\theta} & = &
n^{\frac{1-\theta}{2p} + \frac{\theta}{2}} \ L_{(\infty,
\frac{2p'}{\theta})}^{2p} \ (\mathcal{M}, \mathsf{E}),
\end{eqnarray*}
for $1 \le p \le \infty$. Moreover, if $1 \le p < \infty$ the same
result gives $$\big[ \sqrt{n} \, L_p^r(\mathcal{M}, \mathsf{E}),
\sqrt{n} \, L_p^c(\mathcal{M},\mathsf{E}) \big]_{\theta} =
\sqrt{n} \, L_{(\frac{2p'}{1-\theta},
\frac{2p'}{\theta})}^{2p}(\mathcal{M}, \mathsf{E}).$$ In the
extremal case we claim that we have a contractive inclusion
$$\big[ L_{\infty}^r(\mathcal{M}, \mathsf{E}),
L_{\infty}^c(\mathcal{M}, \mathsf{E}) \big]_{\theta} \subset
L_{(\frac{2}{1-\theta}, \frac{2}{\theta})}^{\infty}(\mathcal{M},
\mathsf{E}).$$ Indeed, let us consider the multi-linear mappings
$$\begin{array}{rrcl} \mathsf{T}_1: & (\alpha,x,\beta) \in
L_2(\mathcal{M}) \times L_{\infty}^r(\mathcal{M}, \mathsf{E})
\times L_{\infty}(\mathcal{M}) & \mapsto & \alpha x \beta \in
L_2(\mathcal{M}), \\ \mathsf{T}_2: & (\alpha,x,\beta) \in
L_{\infty}(\mathcal{M}) \times L_{\infty}^c(\mathcal{M},
\mathsf{E}) \times L_2(\mathcal{M}) & \mapsto & \alpha x \beta \in
L_2(\mathcal{M}). \end{array}$$ By the definition of
$L_{\infty}^r(\mathcal{M}, \mathsf{E})$ and
$L_{\infty}^c(\mathcal{M}, \mathsf{E})$ it is clear that both
$\mathsf{T}_1$ and $\mathsf{T}_2$ are contractions. In particular,
it is easily checked that our claim follows by multi-linear
interpolation, details are left to the reader. Therefore,
according to the observation above and Corollary
\ref{Corollary-Interpolation-Conditional2}, we obtain the lower
estimate with constant $1$. In other words, there exists a
contractive inclusion $$\big[
\mathcal{R}_{2p,1}^n(\mathcal{M},\mathsf{E}),
\mathcal{C}_{2p,1}^n(\mathcal{M},\mathsf{E}) \big]_{\theta} \
\subset \bigcap_{u,v \in \{ 2p', \infty \}}^{\null}
n^{\frac{1-\theta}{u} + \frac{1}{2p} + \frac{\theta}{v}} \,
L_{(\frac{u}{1-\theta},\frac{v}{\theta})}^{2p} (\mathcal{M},
\mathsf{E}).$$ To prove the converse, we consider again the map
given in (\ref{Equation-Map-u}). Arguing as in the proof of
Theorem \ref{Theorem-Intersection1} and according to Lemmas
\ref{Lemma-Complemented-Isomorphism} and \ref{Lemma-LpMO}, we
obtain the following equivalences
$$\begin{array}{rcll} \|x\|_{[\mathcal{R}_{2p,1}^n (\mathcal{M},
\mathsf{E}), \mathcal{C}_{2p,1}^n (\mathcal{M},
\mathsf{E})]_\theta} & \sim & \big\|u(x) \big\|_{[\Ha_{2p}^r
(\mathcal{A}), \Ha_{2p}^c(\mathcal{A}) ]_\theta} & \mbox{for} \ p
\le 2,
\\ [5pt] \|x\|_{[\mathcal{R}_{2p,1}^n (\mathcal{M}, \mathsf{E}),
\mathcal{C}_{2p,1}^n (\mathcal{M}, \mathsf{E})]_\theta} & \sim &
\big\|u(x) \big\|_{[L_{2p}^r \mathcal{MO} (\mathcal{A}), L_{2p}^c
\mathcal{MO} (\mathcal{A}) ]_\theta} & \mbox{for} \, \ p \ge 2,
\end{array}$$
with constants independent of $p,q,n$. On the other hand, Berg's
theorem gives isometric inclusions
\begin{eqnarray*}
\big[ \Ha_{2p}^r (\mathcal{A}), \Ha_{2p}^c(\mathcal{A})
\big]_\theta & \subset & \big[ \Ha_{2p}^r (\mathcal{A}),
\Ha_{2p}^c(\mathcal{A}) \big]^\theta, \\ \big[ L_{2p}^r
\mathcal{MO} (\mathcal{A}), L_{2p}^c \mathcal{MO} (\mathcal{A})
\big]_\theta & \subset & \big[ L_{2p}^r \mathcal{MO}
(\mathcal{A}), L_{2p}^c \mathcal{MO} (\mathcal{A}) \big]^\theta.
\end{eqnarray*}
Now we can use duality an obtain
\begin{eqnarray*}
\big[ \Ha_{2p}^r (\mathcal{A}), \Ha_{2p}^c(\mathcal{A})
\big]^\theta & \simeq & \big[ \Ha_{(2p)'}^r (\mathcal{A}),
\Ha_{(2p)'}^c(\mathcal{A}) \big]_\theta^* \quad \mbox{for} \ 1 \le
p < \infty, \\ \big[ L_{2p}^r \mathcal{MO} (\mathcal{A}), L_{2p}^c
\mathcal{MO} (\mathcal{A}) \big]^\theta & \simeq & \big[
\Ha_{(2p)'}^r (\mathcal{A}), \Ha_{(2p)'}^c(\mathcal{A})
\big]_\theta^* \quad \mbox{for} \ 1 < p \le \infty,
\end{eqnarray*}
where the constants in the first isomorphism remain bounded as $p
\to 1$ and the constants in the second one remain bounded as $p
\to \infty$. Therefore, recalling the terminology used above
$$\big[ \Ha_{(2p)'}^r (\mathcal{A}),
\Ha_{(2p)'}^c(\mathcal{A}) \big]_\theta^* = \mathcal{Z}_{(2p)'}
(\mathcal{A}, \theta),$$ we deduce the following equivalence
$$\|x\|_{[\mathcal{R}_{2p,1}^n (\mathcal{M},
\mathsf{E}), \mathcal{C}_{2p,1}^n (\mathcal{M},
\mathsf{E})]_\theta} \sim \|u(x)\|_{\mathcal{Z}_{(2p)'}
(\mathcal{A}, \theta)}.$$ According to Proposition
\ref{Proposition-BMO2}, the right hand side is controlled by
$$c'(p,\theta) \, \max \Big\{ \|u(x)\|_{\Ha_{2p}^{2p}},
\|u(x)\|_{[h_{2p}^r, \Ha_{2p}^{2p}]_\theta},
\|u(x)\|_{[\Ha_{2p}^{2p}, h_{2p}^c]_\theta}, \|u(x)\|_{[h_{2p}^r,
h_{2p}^c]_\theta} \Big\}.$$ The first two terms are estimated as
in Theorem \ref{Theorem-Intersection1}
\begin{eqnarray*}
\|u(x)\|_{\mathcal{H}_{2p}^{2p}(\mathcal{A})} & = &
n^{\frac{1}{2p}} \|x\|_{2p}, \\
\|u(x)\|_{[h_{2p}^r(\mathcal{A}),
\Ha_{2p}^{2p}(\mathcal{A})]_\theta} & \le & n^{\frac{1}{2q}}
\|x\|_{L_{(\frac{2pq}{p-q}, \infty)}^{2p}(\mathcal{M},
\mathsf{E})}.
\end{eqnarray*}
Note that the latter term is the norm of $x$ in
$$n^{\frac{1-\theta}{u} + \frac{1}{2p} + \frac{\theta}{v}} \,
L_{(\frac{u}{1-\theta},\frac{v}{\theta})}^{2p} (\mathcal{M},
\mathsf{E}) \quad \mbox{with} \quad (u,v) = (2p', \infty).$$
Taking adjoints and replacing $\theta$ by $1 - \theta$, we obtain
$$\|u(x)\|_{[\Ha_{2p}^{2p}(\mathcal{A}),
h_{2p}^c(\mathcal{A})]_\theta} \le n^{\frac{1}{2p} +
\frac{\theta}{2p'}} \|x\|_{L_{(\infty,
\frac{2p'}{\theta})}^{2p}(\mathcal{M}, \mathsf{E})}.$$ It remains
to estimate the last term in the maximum. As in Theorem
\ref{Theorem-Intersection1}
\begin{eqnarray*}
\|u(x)\|_{h_{2p}^r(\mathcal{A})} & = & \sqrt{n} \, \big\|
\mathsf{E}(x x^*) \big\|_{2p}, \\ \|u(x)\|_{h_{2p}^c(\mathcal{A})}
& = & \sqrt{n} \, \big\| \mathsf{E}(x^* x) \big\|_{2p}.
\end{eqnarray*}
Thus, by complex interpolation we conclude $$\|u(x)\|_{[h_{2p}^r,
h_{2p}^c]_\theta} \le \sqrt{n} \, \|x\|_{L_{(\frac{2p'}{1-\theta},
\frac{2p'}{\theta})}^{2p}(\mathcal{M}, \mathsf{E})}.$$ Combining
the estimates obtained above, the assertion follows.
\end{proof}

\begin{remark}
\emph{As we shall see in Chapter \ref{Section7} below, another
useful way to write the intersection space appearing on the right
side of Theorem \ref{Theorem-Intersection2} is in the following
form $$\bigcap_{\alpha, \beta \in \{ 2p, 2q \}}
n^{\frac{1-\eta}{\alpha} + \frac{\eta}{\beta}} \,
L_{(\frac{4pq}{(1-\eta)(2p-\alpha)},\frac{4pq}{\eta(2p -
\beta)})}^{2p} (\mathcal{M}, \mathsf{E}).$$ }
\end{remark}

\begin{remark}
\emph{It might be of independent interest to mention that our
methods immediately imply that the spaces
$\mathcal{R}_{2p,1}^n(\mathcal{M}, \mathsf{E})$ and
$\mathcal{C}_{2p,1}^n(\mathcal{M}, \mathsf{E})$ form interpolation
families with respect to the index $p$. In other words, for any $1
\le p \le \infty$}
\begin{eqnarray*}
\big[ \mathcal{R}_{\infty,1}^n(\mathcal{M}, \mathsf{E}),
\mathcal{R}_{2,1}^n(\mathcal{M}, \mathsf{E}) \big]_{1/p} & \simeq
& \mathcal{R}_{2p,1}^n(\mathcal{M}, \mathsf{E}), \\ \big[ \,
\mathcal{C}_{\infty,1}^n(\mathcal{M}, \mathsf{E}) \, , \,
\mathcal{C}_{2,1}^n(\mathcal{M}, \mathsf{E}) \, \big]_{1/p} &
\simeq & \, \mathcal{C}_{2p,1}^n \, (\mathcal{M}, \mathsf{E}).
\end{eqnarray*}
\emph{Moreover, using anti-linear duality we may replace
intersections by sums and extend our results to the whole range $1
\le 2p \le \infty$. These generalizations of Theorems
\ref{Theorem-Intersection1} and \ref{Theorem-Intersection2} are
out of the scope of this paper.}
\end{remark}

\chapter{Factorization of $\mathcal{J}_{p,q}^n(\mathcal{M},
\mathsf{E})$}
\label{Section6}

Let $(\mathrm{X}_1, \mathrm{X}_2)$ be a pair of operator spaces
containing a von Neumann algebra $\mathcal{M}$ as a common
two-sided ideal. We define the \emph{amalgamated} Haagerup tensor
product $\mathrm{X}_1 \otimes_{\mathcal{M}, h} \mathrm{X}_2$ as
the quotient of the Haagerup tensor product $\mathrm{X}_1
\otimes_h \mathrm{X}_2$ \label{AmalHaa} by the closed subspace
$\mathcal{I}$ generated by the differences $$x_1 \gamma \otimes
x_2 - x_1 \otimes \gamma x_2 \quad \mbox{with} \quad \gamma \in
\mathcal{M}.$$ We shall be interested only in the Banach space
structure of the operator spaces $\mathrm{X}_1
\otimes_{\mathcal{M}, h} \mathrm{X}_2$. In particular, we shall
write $\mathrm{X}_1 \oM \mathrm{X}_2$ to denote the underlying
Banach space of $\mathrm{X}_1 \otimes_{\mathcal{M},h}
\mathrm{X}_2$. According to the definition of the Haagerup tensor
product and recalling the isometric embeddings $\mathrm{X}_j
\subset \mathcal{B}(\mathcal{H}_j)$, we have
\begin{equation} \label{Equation-Norm-Circ}
\|x\|_{\mathrm{X}_1 \oM \mathrm{X}_2} = \inf \left\{ \Big\| \big(
\summ_k x_{1k} x_{1k}^* \big)^{1/2}
\Big\|_{\mathcal{B}(\mathcal{H}_1)} \Big\| \big( \summ_k x_{2k}^*
x_{2k} \big)^{1/2} \Big\|_{\mathcal{B}(\mathcal{H}_2)} \right\},
\end{equation}
where the infimum runs over all possible decompositions of $x +
\mathcal{I}$ into a finite sum $$x = \summ_k x_{1k} \otimes x_{2k}
+ \mathcal{I}.$$

\begin{remark}
\emph{Our definition \eqref{Equation-Norm-Circ} of the norm in
$\mathrm{X}_1 \otimes_{\mathcal{M}} \mathrm{X}_2$ uses the
operator space structure of $\mathrm{X}_1$ and $\mathrm{X}_2$
since the row and column square functions live in
$\mathcal{B}(\mathcal{H}_j)$ but not necessarily in
$\mathrm{X}_j$. However, in the sequel it will be important to
note that much less structure on $(\mathrm{X}_1, \mathrm{X}_2)$ is
needed to define the norm in $\mathrm{X}_1 \otimes_{\mathcal{M}}
\mathrm{X}_2$. Indeed, we just need to impose conditions under
which the row and column square functions become closed operations
in $\mathrm{X}_1$ and $\mathrm{X}_2$ respectively. In particular,
this is guaranteed if $\mathrm{X}_1$ is a right
$\mathcal{M}$-module and $\mathrm{X}_2$ is a left
$\mathcal{M}$-module. On the other hand, note that this structure
does not provide us with a natural operator space structure on
$\mathrm{X}_1 \otimes_\mathcal{M} \mathrm{X}_2$, as we did above
with $\mathrm{X}_1 \otimes_{\mathcal{M}, h} \mathrm{X}_2$.}
\end{remark}

Let us consider any pair of indices $1 \le p,q \le \infty$
satisfying $q \le p$ and let us define $1/r = 1/q - 1/p$. Then,
given a positive integer $n$, the rest of this paper will be
devoted to study the following intersection spaces
$$\mathcal{J}_{p,q}^n(\mathcal{M}, \mathsf{E}) = \bigcap_{u,v \in
\{ 2r, \infty \}} n^{\frac{1}{u} + \frac{1}{p} + \frac{1}{v}}
L_{(u,v)}^p(\mathcal{M}, \mathsf{E}).$$ The aim of this chapter is
the following \emph{factorization} result for the spaces
$\mathcal{J}_{pq}^n(\mathcal{M}, \mathsf{E})$.

\begin{theorem} \label{Theorem-LaFactorizacion}
If $1 \le q \le p \le \infty$, we have
$$\mathcal{J}_{p,q}^n(\mathcal{M}, \mathsf{E}) \simeq
\mathcal{R}_{2p,q}^n(\mathcal{M}, \mathsf{E}) \oM
\mathcal{C}_{2p,q}^n(\mathcal{M}, \mathsf{E}),$$ isomorphically.
Moreover, the constants are independent of $p,q,n$.
\end{theorem}

\section{Amalgamated tensors}

Before any other consideration, let us simplify the expression
\eqref{Equation-Norm-Circ} for the amalgamated Haagerup tensor
product in Theorem \ref{Theorem-LaFactorizacion} above. Since
$\mathcal{R}_{2p,q}^n(\mathcal{M}, \mathsf{E})$ and
$\mathcal{C}_{2p,q}^n(\mathcal{M}, \mathsf{E})$ coincide with
$L_{2p}(\mathcal{M})$ as a set, the product $ab$ of any two
elements $(a,b) \in \mathcal{R}_{2p,q}^n(\mathcal{M}, \mathsf{E})
\times \mathcal{C}_{2p,q}^n(\mathcal{M}, \mathsf{E})$ is
well-defined and the amalgamation over $\mathcal{M}$ allows us to
identify finite sums
$$\summ_k a_k b_k \simeq \summ_k a_k \otimes b_k + \mathcal{I}.$$
Moreover, it is easily seen that
\begin{eqnarray*}
a_1, a_2, \ldots, a_n \in L_{(u,\infty)}^{2p}(\mathcal{M},
\mathsf{E}) & \Rightarrow & \big( \summ_k a_k a_k^*
\big)^{\frac12} \in L_{(u,\infty)}^{2p}(\mathcal{M}, \mathsf{E}), \\
b_1 \hskip1pt, b_2 \hskip1pt, \ldots, b_n \hskip1pt \in
L_{(\infty,v)}^{2p}(\mathcal{M}, \mathsf{E}) & \Rightarrow & \big(
\summ_k \hskip1pt b_k^* b_k \hskip1pt \big)^{\frac12} \in
L_{(\infty,v)}^{2p}(\mathcal{M}, \mathsf{E}).
\end{eqnarray*}
In particular, it turns out that (\ref{Equation-Norm-Circ})
simplifies in this case as follows
$$\|x\|_{\mathcal{R}_{2p,q}^n(\mathcal{M}, \mathsf{E}) \oM
\mathcal{C}_{2p,q}^n(\mathcal{M}, \mathsf{E})} = \inf \left\{
\Big\| \big( \summ_k a_k a_k^* \big)^{1/2}
\Big\|_{\mathcal{R}_{2p,q}^n} \Big\| \big( \summ_k b_k^* b_k
\big)^{1/2} \Big\|_{\mathcal{C}_{2p,q}^n} \right\},$$ where the
infimum runs over all possible decompositions $$x = \summ_k a_k
b_k.$$ Of course, this argument holds in a more general context.
Indeed, arguing as in Proposition \ref{Proposition-Cond-Linfty},
we deduce that the conditional $L_p$ space
$L_{(u,v)}^p(\mathcal{M}, \mathsf{E})$ embeds contractively into
$L_s(\mathcal{M})$ with $$1/s = 1/u + 1/p + 1/v.$$ Thus, the same
arguments lead to the same simplification for the spaces
\begin{eqnarray*}
\mathsf{A} & = & L_{2p}(\mathcal{M}) \oM L_{2p}(\mathcal{M}),
\\ [9pt] \mathsf{B} & = & L_{2p}(\mathcal{M}) \oM
L_{(\infty, \frac{2pq}{p-q})}^{2p}(\mathcal{M}, \mathsf{E}), \\
[5pt] \mathsf{C} & = &
L_{(\frac{2pq}{p-q},\infty)}^{2p}(\mathcal{M}, \mathsf{E}) \oM
L_{2p}(\mathcal{M}), \\ [5pt] \mathsf{D} & = &
L_{(\frac{2pq}{p-q},\infty)}^{2p}(\mathcal{M}, \mathsf{E}) \oM
L_{(\infty,\frac{2pq}{p-q})}^{2p}(\mathcal{M}, \mathsf{E}).
\end{eqnarray*}

Our first step in the proof of Theorem
\ref{Theorem-LaFactorizacion} is the following.

\begin{lemma} \label{Lemma-LaFactorizacion}
We have $$\mathcal{R}_{2p,q}^n(\mathcal{M}, \mathsf{E}) \oM
\mathcal{C}_{2p,q}^n(\mathcal{M}, \mathsf{E}) \simeq
n^{\frac{1}{p}} \mathsf{A} \cap n^{\frac{1}{2p} + \frac{1}{2q}}
\mathsf{B} \cap n^{\frac{1}{2p} + \frac{1}{2q}} \mathsf{C} \cap
n^{\frac{1}{q}} \mathsf{D},$$ where the relevant constants in the
isomorphism above are independent of $p,q,n$.
\end{lemma}

\begin{proof}
It is clear that $$\mathcal{R}_{2p,q}^n(\mathcal{M}, \mathsf{E})
\oM \mathcal{C}_{2p,q}^n(\mathcal{M}, \mathsf{E}) \subset
n^{\frac{1}{p}} \mathsf{A} \cap n^{\frac{1}{2p} + \frac{1}{2q}}
\mathsf{B} \cap n^{\frac{1}{2p} + \frac{1}{2q}} \mathsf{C} \cap
n^{\frac{1}{q}} \mathsf{D}$$ contractively. To prove the reverse
inclusion it suffices to see that
\begin{eqnarray}
\label{eqsub1} n^{\frac{1}{p}} \mathsf{A} \cap n^{\frac{1}{2p} +
\frac{1}{2q}} \mathsf{C} & \subset &
\mathcal{R}_{2p,q}^n(\mathcal{M}, \mathsf{E}) \oM
n^{\frac{1}{2p}} L_{2p}(\mathcal{M}) = \mathcal{X}_1, \\
\label{eqsub2} n^{\frac{1}{2p} + \frac{1}{2q}} \mathsf{B} \cap
n^{\frac{1}{q}} \mathsf{D} & \subset &
\mathcal{R}_{2p,q}^n(\mathcal{M}, \mathsf{E}) \oM n^{\frac{1}{2q}}
L_{(\infty, \frac{2pq}{p-q})}^{2p}(\mathcal{M}, \mathsf{E}) =
\mathcal{X}_2,
\end{eqnarray}
with constants independent of $p,q,n$ and also
\begin{equation} \label{eqsub3}
\mathcal{X}_1 \cap \mathcal{X}_2 \subset
\mathcal{R}_{2p,q}^n(\mathcal{M}, \mathsf{E}) \oM
\mathcal{C}_{2p,q}^n(\mathcal{M}, \mathsf{E})
\end{equation}
with absolute constants. In fact, the three inclusions can be
proved using the same principle, which obviously works in a much
more general setting. Indeed, let us prove \eqref{eqsub1}. If $x$
is a norm one element in $n^{\frac{1}{p}} \mathsf{A} \cap
n^{\frac{1}{2p} + \frac{1}{2q}} \mathsf{C}$ and $\delta
> 0$, we may find decompositions $x = \summ_k a_{1k} a_{2k}$ and
$x = \summ_k c_{1k} c_{2k}$ satisfying the following estimates
\begin{eqnarray*} \max \Big\{ n^{\frac{1}{2p}} \Big\| \big(
\summ_k a_{1k} a_{1k}^* \big)^{\frac12} \Big\|_{2p},
n^{\frac{1}{2p}} \Big\| \big( \summ_k a_{2k}^* a_{2k}
\big)^{\frac12} \Big\|_{2p} \Big\} & \le & 1 + \delta, \\
\max \Big\{ n^{\frac{1}{2q}} \Big\| \big( \summ_k c_{1k} c_{1k}^*
\big)^{\frac12} \Big\|_{L_{(\frac{2pq}{p-q},
\infty)}^{2p}(\mathcal{M}, \mathsf{E})}, n^{\frac{1}{2p}} \Big\|
\big( \summ_k c_{2k}^* c_{2k} \big)^{\frac12} \Big\|_{2p} \Big\} &
\le & 1 + \delta.
\end{eqnarray*}
Let us consider the following element in $L_{2p}(\mathcal{M})$
$$\gamma = \Big( \summ_k a_{2k}^* a_{2k} + \summ_k c_{2k}^* c_{2k}
+ \delta \mathrm{D}^{1/p} \Big)^{1/2}.$$ Since $\gamma$ is
invertible as a measurable operator, we may define $\xi$ by $x =
\xi \gamma$. Moreover, we also define $(a_{2k}', c_{2k}')$ by
$a_{2k} = a_{2k}' \gamma$ and $c_{2k} = c_{2k}' \gamma$. This
gives rise to the following expressions for $x$ that will be used
below
\begin{eqnarray*}
x = \xi \gamma & = & \big( \summ_k a_{1k} a_{2k}' \big) \gamma, \\
x = \xi \gamma & = & \big( \summ_k c_{1k} \hskip1pt c_{2k}'
\hskip1pt \big) \gamma.
\end{eqnarray*}
The norm of $x$ in $\mathcal{X}_1$ is estimated as follows
$$\|x\|_{\mathcal{X}_1} \le \|\xi\|_{\mathcal{R}_{2p,q}^n
(\mathcal{M}, \mathsf{E})} \|\gamma\|_{n^{\frac{1}{2p}}
L_{2p}(\mathcal{M})},$$ where the norm of $\xi$ in the space
$\mathcal{R}_{2p,q}^n(\mathcal{M}, \mathsf{E})$ is given by
$$\max \Big\{ n^{\frac{1}{2p}} \Big\| \summ_k a_{1k} a_{2k}'
\Big\|_{2p}, n^{\frac{1}{2q}} \Big\| \summ_k c_{1k} c_{2k}'
\Big\|_{L_{(\frac{2pq}{p-q}, \infty)}^{2p}(\mathcal{M},
\mathsf{E})} \Big\} = \max \big\{ \mathrm{A}, \mathrm{C} \big\}.$$
However, since $$\summ_k a_{2k}^{{}_{'} *} a_{2k}' = \gamma^{-1} (
\summ_k a_{2k}^* a_{2k} ) \gamma^{-1} \le 1,$$ we obtain
$$\mathrm{A} = n^{\frac{1}{2p}} \Big\| \summ_k a_{1k} a_{2k}' \Big\|_{2p}
\le n^{\frac{1}{2p}} \Big\| \big( \summ_k a_{1k}
a_{1k}^*\big)^{\frac12} \Big\|_{2p} \Big\| \big( \summ_k
a_{2k}^{{}_{'} *} a_{2k}' \big)^{\frac12} \Big\|_{\infty} \le 1 +
\delta.$$ Similarly, we have $$\summ_k c_{2k}^{{}_{'} *} c_{2k}' =
\gamma^{-1} ( \summ_k c_{2k}^* c_{2k} ) \gamma^{-1} \le 1$$ and
therefore we deduce
\begin{eqnarray*}
\mathrm{C} & = & n^{\frac{1}{2q}} \sup \Big\{ \Big\| \alpha
\summ_k c_{1k} c_{2k}' \Big\|_{L_{2q}(\mathcal{M})} \, \big| \
\|\alpha\|_{L_{\frac{2pq}{p-q}}(\mathcal{N})} \le 1 \Big\} \\ &
\le & n^{\frac{1}{2q}} \sup \Big\{ \Big\| \summ_k \alpha c_{1k}
c_{1k}^* \alpha^* \Big\|_{L_q(\mathcal{M})}^{1/2} \, \big| \
\|\alpha\|_{L_{\frac{2pq}{p-q}}
(\mathcal{N})} \le 1 \Big\} \\
[4pt] & = & n^{\frac{1}{2q}} \sup \Big\{ \Big\| \alpha \big(
\summ_k c_{1k} c_{1k}^* \big)^{\frac12}
\Big\|_{L_{2q}(\mathcal{M})} \, \big| \
\|\alpha\|_{L_{\frac{2pq}{p-q}}(\mathcal{N})} \le 1 \Big\} \le 1 +
\delta.
\end{eqnarray*}
Thus, we have proved $$\|\xi\|_{\mathcal{R}_{2p,q}^n(\mathcal{M},
\mathsf{E})} \le 1 + \delta.$$ It remains to estimate the norm of
$\gamma$ in $n^{\frac{1}{2p}} L_{2p}(\mathcal{M})$
$$n^{\frac{1}{2p}} \|\gamma\|_{2p} \le n^{\frac{1}{2p}} \Big(
\Big\| \summ_k a_{2k}^* a_{2k} \Big\|_p + \Big\| \summ_k c_{2k}^*
c_{2k} \Big\|_p + \delta \Big)^{\frac12} \le \sqrt{2} (1+\delta) +
\delta^{\frac12} n^{\frac{1}{2p}}.$$ In conclusion, letting
$\delta \to 0^+$ we obtain $$\|x\|_{\mathcal{X}_1} \le \sqrt{2}.$$
This proves \eqref{eqsub1} and inclusion \eqref{eqsub2} is proved
in a similar way. To prove \eqref{eqsub3} we just need to observe
that the \emph{common factor} $\mathcal{R}_{2p,q}^n(\mathcal{M},
\mathsf{E})$ in $\mathcal{X}_1 \cap \mathcal{X}_2$ is on the left,
in contrast with the previous situations where the common factor
was on the right. The only consequence is that the roles of $\xi$
and $\gamma$ above must be exchanged.
\end{proof}

\begin{lemma} \label{Lemma-LaFactorizacion2}
We have
\begin{eqnarray*}
\mathsf{A} & = & L_p(\mathcal{M}),
\\ [9pt] \mathsf{B} & = & L_{(\infty, \frac{2pq}{p-q})}^p
(\mathcal{M}, \mathsf{E}), \\ [5pt] \mathsf{C} & = &
L_{(\frac{2pq}{p-q},\infty)}^p(\mathcal{M}, \mathsf{E}), \\ [5pt]
\mathsf{D} & = & L_{(\frac{2pq}{p-q},
\frac{2pq}{p-q})}^p(\mathcal{M}, \mathsf{E}),
\end{eqnarray*}
isometrically whenever the indices $p$ and $q$ satisfy $1 \le q
\le p < \infty$ or $1 < q \le p \le \infty$.
\end{lemma}

\begin{proof}
Let us define
$$1/s = 1/u + 1/p + 1/v = 1/s_1 + 1/s_2$$ with the indices $s_1$ and $s_2$
given by
\begin{eqnarray*}
1/s_1 & = & 1/u + 1/2p, \\ 1/s_2 & = & 1/2p + 1/v.
\end{eqnarray*}
Then we obtain the following estimate
\begin{eqnarray*}
\Big\| \summ_k a_k b_k \Big\|_{L_{(u,v)}^p(\mathcal{M},
\mathsf{E})} & = & \sup_{\alpha, \beta} \Big\| \summ_k \alpha a_k
b_k \beta \Big\|_s \\ & \le & \sup_{\alpha, \beta} \Big\| \alpha
\big( \summ_k a_k a_k^* \big)^{\frac12} \Big\|_{s_1} \Big\| \big(
\summ_k b_k^* b_k \big)^{\frac12} \beta \Big\|_{s_2} \\ & = &
\Big\| \big( \summ_k a_k a_k^* \big)^{\frac12}
\Big\|_{L_{(u,\infty)}^{2p}(\mathcal{M}, \mathsf{E})} \Big\| \big(
\summ_k b_k^* b_k \big)^{\frac12}
\Big\|_{L_{(\infty,v)}^{2p}(\mathcal{M}, \mathsf{E})},
\end{eqnarray*}
where the supremum runs over all $\alpha$ in the unit ball in
$L_u(\mathcal{N})$ and all $\beta$ in the unit ball of
$L_v(\mathcal{N})$. Therefore, taking infima on the right we
obtain a contractive inclusion
$$L_{(u,\infty)}^{2p}(\mathcal{M}, \mathsf{E}) \oM
L_{(\infty,v)}^{2p}(\mathcal{M}, \mathsf{E}) \subset
L_{(u,v)}^p(\mathcal{M}, \mathsf{E}).$$ This shows at once the
lower estimate for all isometries and for $1 \le q \le p \le
\infty$. That is, with no restriction on the indices. In order to
show the reverse inequalities, we restrict the indices $u$ and $v$
to be either $\frac{2pq}{p-q}$ or $\infty$. We shall obtain
$$L_{(u,v)}^p(\mathcal{M}, \mathsf{E}) \subset
L_{(u,\infty)}^{2p}(\mathcal{M}, \mathsf{E}) \oM
L_{(\infty,v)}^{2p}(\mathcal{M}, \mathsf{E})$$ contractively. To
do so we begin with some remarks. First, the isometry $\mathsf{A}
= L_p(\mathcal{M})$ is very well-known and there is nothing to
prove. Thus, the case $q=p$ is trivial since the spaces on the
left collapse into $L_{2p}(\mathcal{M}) \oM L_{2p}(\mathcal{M})$
while the spaces on the right coincide with $L_p(\mathcal{M})$.
Therefore, we just need to consider the cases $1 \le q < p <
\infty$ and $1 < q < p \le \infty$. In both cases we have $$2 <
\frac{2pq}{p-q} < \infty.$$ In other words, we may assume that $2
< \min(u,v) < \infty$. In particular, we are in position to apply
the standard Grothendieck-Pietsch separation argument as in
Theorem \ref{Theorem-Grothendieck} and Observation
\ref{Observation-max(u,v)}. This will be our main tool in the
proof. Let us consider the following norm on $L_p(\mathcal{M})$
$$|||x||| = \inf \left\{ \Big\| \big( \summ_k a_k a_k^*
\big)^{\frac12} \Big\|_{L_{(u,\infty)}^{2p}(\mathcal{M},
\mathsf{E})} \Big\| \big( \summ_k b_k^* b_k \big)^{\frac12}
\Big\|_{L_{(\infty,v)}^{2p}(\mathcal{M}, \mathsf{E})}\right\}$$
where the infimum runs over all decompositions of $x$ into a
finite sum $\summ_k a_k b_k$ with $a_k, b_k \in
L_{2p}(\mathcal{M})$. Since this norm majorizes that of
$L_{(u,\infty)}^{2p}(\mathcal{M}, \mathsf{E}) \oM
L_{(\infty,v)}^{2p}(\mathcal{M}, \mathsf{E})$ in
$L_p(\mathcal{M})$, it suffices by density to see that the norm in
$L_{(u,v)}^p(\mathcal{M}, \mathsf{E})$ controls $||| \ |||$ from
above. To that aim, given $x \in L_p(\mathcal{M})$, we consider a
norm one functional $$\phi: (L_p(\mathcal{M}), ||| \ |||) \to \C
\quad \mbox{satisfying} \quad |\phi(x)| = |||x|||.$$ Then, the
reverse inequalities will follow from
\begin{equation} \label{Equation-Loqfalta}
|\phi(x)| \le \sup \Big\{ \|\alpha x \beta\|_{L_s(\mathcal{M})} \,
\big| \ \alpha \in \mathsf{B}_{L_u(\mathcal{N})}, \beta \in
\mathsf{B}_{L_v(\mathcal{N})} \Big\}.
\end{equation}
Applying the standard Grothendieck-Pietsch separation argument, we
may find positive elements $\alpha \in
\mathsf{B}_{L_u(\mathcal{N})}$ and $\beta \in
\mathsf{B}_{L_v(\mathcal{N})}$ satisfying
$$|\phi(ab)| \le \|\alpha a \|_{L_{s_1}(\mathcal{M})} \|b
\beta\|_{L_{s_2}(\mathcal{M})}.$$ Then, taking $q_\alpha$ and
$q_\beta$ to be the support projections of $\alpha$ and $\beta$
respectively, we obtain after the usual arguments (see e.g.
Theorem \ref{Theorem-Grothendieck}) a right $\mathcal{M}$-module
map $\Psi: q_\alpha L_{s_1}(\mathcal{M}) \to q_\beta
L_{s_2'}(\mathcal{M})$ determined by $$\phi(ab) =
\mbox{tr}_{\mathcal{M}} \Big( \Psi(\alpha a) b \beta \Big).$$ In
particular, there exists $m_{\Psi} \in
\mathsf{B}_{L_{s'}(\mathcal{M})}$ satisfying $\Psi(\alpha a) =
m_{\Psi} \alpha a$ so that
$$|\phi(x)| = \big| \mbox{tr}_{\mathcal{M}} \big( m_{\Psi}
\alpha x \beta \big) \big| \le \sup \Big\{ \|\alpha x
\beta\|_{L_s(\mathcal{M})} \, \big| \ \alpha \in
\mathsf{B}_{L_u(\mathcal{N})}, \beta \in
\mathsf{B}_{L_v(\mathcal{N})} \Big\}.$$ This completes the proof
of inequality (\ref{Equation-Loqfalta}) and thus the proof is
concluded.
\end{proof}

Let us observe that Lemmas \ref{Lemma-LaFactorizacion} and
\ref{Lemma-LaFactorizacion2} give Theorem
\ref{Theorem-LaFactorizacion} for every pair of indices $(p,q)$
except for the case of $\mathcal{J}_{\infty,1}^n(\mathcal{M},
\mathsf{E})$. However, this is for several reasons one of the most
important factorization results that we need. In order to
factorize the space $\mathcal{J}_{\infty,1}^n(\mathcal{M},
\mathsf{E})$, we note that Lemma \ref{Lemma-LaFactorizacion} is
still valid. Moreover, due to the obvious isometries
$$\mathrm{X}_1 \oM \mathcal{M} = \mathrm{X}_1 \quad \mbox{and} \quad
\mathcal{M} \oM \mathrm{X}_2 = \mathrm{X}_2$$ for any right (resp.
left) $\mathcal{M}$-module, we deduce that the first three
isometries in Lemma \ref{Lemma-LaFactorizacion2} are trivial in
the limit case $p=\infty$. Hence, we just need to show that the
isometry for $\mathsf{D}$ still holds in the case $(p,q) =
(\infty,1)$. Again, we recall that the lower estimate can be
proved as in Lemma \ref{Lemma-LaFactorizacion2} so that it
suffices to prove the upper estimate. Unfortunately, the
application of Grothendieck-Pietsch separation argument in this
case is much more delicate and we need some preparation. After
some auxiliary results in the next paragraph, we will go back to
this question.

\section{Conditional expectations and ultraproducts}

We study certain ultraproduct von Neumann algebras and the
corresponding conditional expectations. These auxiliary results
will be used to factorize the norm of
$\mathcal{J}_{\infty,1}^n(\mathcal{M}, \mathsf{E})$ as explained
above.

\begin{lemma} \label{Lemma-Basis}
Let $\mathsf{F}$ be a finite dimensional subspace of
$\mathcal{M}^*$ and let $\mathsf{G}$ be a finite dimensional
subspace of $\mathcal{M}$. Then, given any $\delta > 0$ there
exists a linear mapping $$\omega: \mathsf{F} \to \mathcal{M}_*$$
satisfying the following properties:
\begin{itemize}
\item[i)] $\|\omega\|_{cb} \le 1 + \delta$. \item[ii)] The space
$\omega(\mathsf{F} \cap \mathcal{N}^*)$ is contained in
$\mathcal{N}_*$. \item[iii)] The following estimate holds for any
$f \in \mathsf{F}$ and $g \in \mathsf{G}$ $$\big| g(\omega(f))-
f(g) \big| \le \delta \, \|f\|_{\mathsf{F}} \,
\|g\|_{\mathsf{G}}.$$
\end{itemize}
\end{lemma}

\begin{proof} Let $(f_1, f_2, \ldots f_k;
\mathsf{f}_1^*, \mathsf{f}_2^*, \ldots, \mathsf{f}_k^*)$ be an
Auerbach basis \label{Auerbach} of $\mathsf{F} \cap
\mathcal{N}^*$. That is, $(f_1,f_2, \ldots, f_k)$ is a basis of
$\mathsf{F} \cap \mathcal{N}^*$ with $\|f_j\|=1$ for $1 \le j \le
k$ and the $\mathsf{f}_j^*$'s are functionals on $\mathsf{F} \cap
\mathcal{N}^*$ satisfying $\mathsf{f}_i^* (f_j) = \delta_{ij}$.
Let us take Hahn-Banach extensions $f_1^*, \ldots, f_k^*:
\mathsf{F} \to \C$ of $\mathsf{f}_1^*, \ldots, \mathsf{f}_k^*$
respectively. Then we may define the projection
$$\mathsf{P}: f \in \mathsf{F} \mapsto \sum_{j=1}^k f_j^*(f) f_j
\in \mathsf{F} \cap \mathcal{N}^*.$$ Now, using $\mathsf{P}$ we
may also consider an Auerbach basis $(f_{k+1}, \ldots, f_n;
\mathsf{f}_{k+1}^*, \ldots, \mathsf{f}_n^*)$ of $(1_{\mathsf{F}} -
\mathsf{P})(\mathsf{F})$. Then, given any $k+1 \le j \le n$ we
consider the linear functional $f_j^*: \mathsf{F} \to \C$ defined
by the relation $$f_j^*(f) = \mathsf{f}_j^* \Big( f -
\mathsf{P}(f) \Big).$$ Finally we consider an Auerbach basis
$(g_1, g_2, \ldots, g_m; g_1^*, g_2^*, \ldots, g_m^*)$ of
$\mathsf{G}$. This allows us to define the following set for any
$\varepsilon > 0$ $$\mathsf{C}(\varepsilon) = \mathrm{conv} \Big\{
\omega: \mathsf{F} \to \mathcal{M}_* \, \big| \ \omega(\mathsf{F}
\cap \mathcal{N}^*) \subset \mathcal{N}_*, \ \big|
g_k(\omega(f_j))-f_j(g_k) \big| \le \varepsilon \Big\}.$$ Let us
assume that $$(1+\delta) \,
\mathsf{B}_{\mathcal{CB}(\mathsf{F},\mathcal{M}_*)} \cap
\overline{\mathsf{C}(\varepsilon)} = \emptyset,$$ where
$\overline{\mathsf{C}(\varepsilon)}$ denotes the closure of
$\mathsf{C}(\varepsilon)$ in the $\sigma(\mathcal{CB}(\mathsf{F},
\mathcal{M}_*), \mathsf{F} \otimes \mathcal{M})$ topology. We will
show below that $\mathsf{C}(\varepsilon)$ is not empty. Thus, by
the Hahn-Banach theorem we may find a linear functional $\xi:
\mathcal{CB}(\mathsf{F}, \mathcal{M}_*) \to \C$ such that
$$\mathrm{Re} \big( \xi(\omega_1) \big) \le 1 \le \mathrm{Re}
\big( \xi(\omega_2) \big)$$ for all $\omega_1 \in (1+\delta) \,
\mathsf{B}_{\mathcal{CB}(\mathsf{F}, \mathcal{M}_*)}$ and
$\omega_2 \in \mathsf{C}(\varepsilon)$. This implies
$\|\xi\|_{\mathcal{CB}(\mathsf{F}, \mathcal{M}_*)^*} \le
(1+\delta)^{-1}$. After identifying the space
$\mathcal{CB}(\mathsf{F}, \mathcal{M}_*)$ with the minimal tensor
product $\mathsf{F}^* \otimes_{\mathrm{min}} \mathcal{M}_*$, we
consider the associated linear map $T_{\xi}: \mathcal{M}_* \to
\mathsf{F}$ defined by $$f^*(T_{\xi}(m_*)) = \xi(f^*\otimes
m_*).$$ Then, taking $m_j= f_j^*\circ T_{\xi} \in \mathcal{M}$ it
turns out that $$\xi = \sum_{j=1}^n m_j \otimes f_j.$$ We claim
that $$\|\xi\|_{\mathrm{N}(\mathcal{M}_*,\mathsf{F})} =
\|\xi\|_{\mathrm{I}(\mathcal{M}_*,\mathsf{F})} \le
(1+\delta)^{-1},$$ where $\mathrm{N}(\mathcal{M}_*, \mathsf{F})$
(resp. $\mathrm{I}(\mathcal{M}_*,\mathsf{F})$) denotes the space
of completely nuclear maps (resp. completely integral maps) from
$\mathcal{M}_*$ to $\mathsf{F}$. Indeed, according to the main
result of \cite{EJR}, $\mathcal{M}_*$ is locally reflexive and so
the first identity follows from Proposition 4.4 in \cite{EJR}. The
inequality following it holds by Corollary 12.3.4 of \cite{ER} and
the fact that $\|\xi\|_{\mathcal{CB}(\mathsf{F}, \mathcal{M}_*)^*}
\le (1+\delta)^{-1}$. Moreover, since $\mathsf{F}$ is finite
dimensional $$\mathrm{N}(\mathcal{M}_*, \mathsf{F}) = \mathsf{F}
\hat{\otimes} \mathcal{M} \quad \mbox{and} \quad
\mathrm{N}(\mathcal{M}_*, \mathsf{F})^* = \mathcal{CB}(\mathsf{F},
\mathcal{M}^*).$$ This means that given any linear map $\omega:
\mathsf{F} \to \mathcal{M}^*$, we have
$$|\xi(\omega)| = \Big| \sum_{j=1}^n \omega(f_j)(m_j) \Big| \le
\|\xi\|_{\mathrm{N}(\mathcal{M}_*, \mathsf{F})} \|\omega\|_{cb}
\le (1+\delta)^{-1} \|\omega\|_{cb}.$$ In particular, the
inclusion map $j: \mathsf{F} \to \mathcal{M}^*$ satisfies
$$|\xi(j)| \le (1+\delta)^{-1}.$$ On the other hand we can write
$$j = \sum_{j=1}^n f_j^* \otimes f_j.$$ Given $1 \le j \le k$, let
$(f_{j,\alpha})_{\alpha \in \Lambda} \subset \mathcal{N}_*$ be a
net converging to $f_j$ in the $\sigma(\mathcal{N}^*,\mathcal{N})$
topology. Thus $f_{j,\alpha} \circ \mathsf{E}$ converges to $f_j
\circ \mathsf{E}$ in the $\sigma(\mathcal{M}^*,\mathcal{M})$
topology. When $j=k+1, k+2, \ldots, n$ we may also fix a net
$(f_{j,\alpha})_{\alpha \in \Lambda} \subset \mathcal{M}_*$
converging to $f_j$ in the $\sigma(\mathcal{M}^*,\mathcal{M})$
topology. This implies that the maps $\omega_{\alpha}: \mathsf{F}
\to \mathcal{M}_*$ defined by $$\omega_{\alpha}(f)= \sum_{j=1}^n
f_j^*(f)f_{j,\alpha} \quad \mbox{satisfy} \quad
\omega_{\alpha}(\mathsf{F} \cap \mathcal{N}^*) \subset
\mathcal{N}_*$$ since $f_{j,\alpha} \in \mathcal{N}_*$ for $j = 1,
2, \ldots, k$ and $f_j^*(\mathsf{P}(\mathsf{F})) = 0$ for $j =
k+1, k+2, \ldots, n$. Moreover, for large enough $\alpha$ we
clearly find $$\big| g_k(\omega_{\alpha}(f_j))-f_j(g_k) \big| \le
\varepsilon.$$ This shows that there exists $\alpha$ for which
$\omega_{\alpha} \in \mathsf{C}(\varepsilon)$. Finally, we have
$$\lim_{\alpha} \xi(\omega_{\alpha}) = \lim_{\alpha} \sum_{j=1}^n
m_j \big( \omega_{\alpha}(f_j) \big) = \lim_{\alpha} \sum_{j=1}^n
m_j \big( f_{j,\alpha} \big) = \xi(j).$$ Therefore we may find
$\alpha$ such that $$\omega_{\alpha} \in \mathsf{C}(\varepsilon)
\quad \mbox{and} \quad |\xi(\omega_{\alpha})|<1.$$ However, any
$\omega \in \mathsf{C}(\varepsilon)$ satisfies $\mbox{Re} \big(
\xi(\omega) \big) \ge 1$. Therefore, we have a contradiction so
that we can find $$\omega \in (1+\delta)\,
\mathsf{B}_{\mathcal{CB}(\mathsf{F},\mathcal{M}_*)} \cap
\mathsf{C}(\varepsilon).$$ Such a mapping clearly satisfies
$$\big| g(\omega(f))-g(f) \big| = \Big| \sum_{j,k} f_j^*(f)
g_k^*(g) \big( g_k(\omega(f_j)) - f_j(g_k) \big) \Big| \le
\varepsilon \sum_{j=1}^n |f_j^*(f)| \sum_{k=1}^m |g_k^*(g)|.$$
Then, recalling the definition of $f_j^*$ and $g_k^*$ we easily
obtain $$\big| g(\omega(f))-g(f) \big| \le \varepsilon m n
\|1_{\mathsf{F}} - \mathsf{P}\| \, \|f\|_{\mathsf{F}}
\|g\|_{\mathsf{G}} \le 2 \varepsilon mn^2 \|f\|_{\mathsf{F}}
\|g\|_{\mathsf{G}}$$ for all $f \in \mathsf{F}$ and $g \in
\mathsf{G}$. Therefore, taking $\varepsilon < \delta/2mn^2$  the
assertion follows. \end{proof}

\noindent In the following we use the notation $$(x_i)^{\bullet} =
\ \mbox{Equivalence class of} \ (x_i) \ \mbox{in} \
\prodd_{i,\mathcal{U}} \mathrm{X}_i.$$

\begin{lemma} \label{Lemma-Basis2} There exist an
ultrafilter $\mathcal{U}$ on an index set $I$ and a linear map
$$\alpha: \mathcal{M}^* \to \prodd_{i,\mathcal{U}} \mathcal{M}_*$$
satisfying the following properties:
\begin{itemize}
\item[i)] The map $\alpha$ is a complete contraction. \item[ii)]
The space $\alpha(\mathcal{N}^*)$ is contained in
$\prod_{i,\mathcal{U}} \mathcal{N}_*$. \item[iii)] The following
identity holds for all $\psi\in \mathcal{M}^*$ and $m \in
\mathcal{M}$
$$\limm_{i,\mathcal{U}} \, m(\alpha(\psi)_i) = \psi(m).$$
\end{itemize}
\end{lemma}

\begin{proof} Let $I$ be the set of tuples $(\mathsf{F},\mathsf{G})$
with $\mathsf{F}$ a finite dimensional subspace of $\mathcal{M}^*$
and $\mathsf{G}$ a finite dimensional subspace of $\mathcal{M}$.
Let $\mathcal{U}$ be an ultrafilter containing all the subsets of
$I$ of the form $$I_{\mathsf{F},\mathsf{G}} = \Big\{ (\mathsf{F}',
\mathsf{G}') \, \big| \ \mathsf{F} \subset \mathsf{F}', \
\mathsf{G} \subset \mathsf{G}' \Big\}.$$ Note that this can be
done since these sets have the finite intersection property. For
fixed $\mathsf{F}, \mathsf{G}$ we choose $\omega_{\mathsf{F},
\mathsf{G}}: \mathsf{F} \to \mathcal{M}_*$ satisfying the
assumptions of Lemma \ref{Lemma-Basis} for $\delta = (\dim
\mathsf{F} \dim \mathsf{G})^{-1}$. Then we define $$\alpha(\psi) =
\big( \omega_{\mathsf{F}, \mathsf{G}}(\psi))^{\bullet}.$$ Note
that for $(\mathsf{F}, \mathsf{G})\in I_{\langle \psi \rangle,
\langle 0 \rangle}$ this is well-defined. Hence $\alpha$ is
well-defined on $\mathcal{M}^*$. It is easily checked $\alpha$ is
linear and completely contractive. By construction we have
$\alpha(\mathcal{N}^*) \subset \prod_{i,\mathcal{U}}
\mathcal{N}_*$ and $\alpha(\psi)(m) = \lim_{i,\mathcal{U}}
m(\alpha(\psi)_i)$ for all $m \in \mathcal{M}$. \end{proof}

\begin{lemma} \label{Lemma-Basis3}
There exists normal conditional expectations $$\mathbf{E}: \Big(
\prodd_{i, \mathcal{U}} \mathcal{M}_* \Big)^* \to \Big( \prodd_{i,
\mathcal{U}} \mathcal{N}_* \Big)^* \quad \mbox{and} \quad
\mathcal{E}: \Big( \prodd_{i, \mathcal{U}} \mathcal{M}_* \Big)^*
\to \mathcal{M}^{**}.$$ Moreover, they are related to each other
by the identity $\mathsf{E}^{**} \circ \mathcal{E} = \mathcal{E}
\circ \mathbf{E}$.
\end{lemma}

\begin{proof} Let us consider the map $$\mathsf{E}^{\mathcal{U}}:
(x_i)^{\bullet} \in \prodd_{i, \mathcal{U}} \mathcal{M} \mapsto
(\mathsf{E}(x_i))^{\bullet} \in \prodd_{i, \mathcal{U}}
\mathcal{N}.$$ By strong density of $\prod_{i,\mathcal{U}}
\mathcal{M}$ in $\big( \prod_{i,\mathcal{U}} \mathcal{M}_*
\big)^*$, we may define $$\mathbf{E}: \Big( \prodd_{i,
\mathcal{U}} \mathcal{M}_* \Big)^* \to \Big( \prodd_{i,
\mathcal{U}} \mathcal{N}_* \Big)^*$$ with predual
$\mathsf{E}_*^{\mathcal{U}}$. On the other hand, the map $$\gamma:
(\varphi_i)^{\bullet} \in \prodd_{i, \mathcal{U}} \mathcal{M}_*
\mapsto \limm_{i, \mathcal{U}} \varphi_i(\cdot) \in
\mathcal{M}^*$$ clearly satisfies $\alpha \circ \gamma =
1_{\prod_{i,\mathcal{U}} \mathcal{M}_*}$ with $\alpha$ being the
map constructed in Lemma \ref{Lemma-Basis2}. Taking adjoints we
find that $$\gamma^*: \mathcal{M}^{**} \to \Big( \prodd_{i,
\mathcal{U}} \mathcal{M}_* \Big)^*$$ is an injective
$*$-homomorphism so that the adjoint $\mathcal{E}$ of $\alpha$
$$\mathcal{E}: \Big( \prodd_{i, \mathcal{U}} \mathcal{M}_* \Big)^*
\to \mathcal{M}^{**}$$ is a normal conditional expectation. We are
interested in proving the relation $\mathsf{E}^{**} \circ
\mathcal{E} = \mathcal{E} \circ \mathbf{E}$. Note that since
$\alpha (\mathcal{N}^*)$ is contained in $\prod_{i, \mathcal{U}}
\mathcal{N}_*$, we have $$\mathcal{E} \Big( \prodd_{i,
\mathcal{U}} \mathcal{N}_* \Big)^* \subset \mathcal{N}^{**}.$$ In
particular, both maps $\mathsf{E}^{**} \circ \mathcal{E}$ and
$\mathcal{E} \circ \mathbf{E}$ end in $\mathcal{N}^{**}$. However,
if we \emph{predualize} this identity it turns out that it
suffices to see that $$\alpha \circ \mathsf{E}^* =
\mathsf{E}_*^{\mathcal{U}} \circ \alpha_{\mid_{\mathcal{N}^*}}.$$
This is a reformulation of Lemma \ref{Lemma-Basis2} iii). The
proof is complete.
\end{proof}

\section{Factorization of the space
$\mathcal{J}_{\infty,1}^n(\mathcal{M}, \mathsf{E})$}

According to \cite{Ra}, $$L_p(\mathcal{M}_{\mathcal{U}}) =
\prodd_{i, \mathcal{U}} L_p(\mathcal{M}) \quad \mbox{with} \quad
\mathcal{M}_{\mathcal{U}} = \Big( \prodd_{i, \mathcal{U}}
\mathcal{M}_* \Big)^* \quad \mbox{and} \quad 1 \le p < \infty.$$
Then, by Lemma \ref{Lemma-Basis3} the inclusion below is a
complete contraction $$\xi_p: L_p(\mathcal{M}^{**}) \to
L_p(\mathcal{M}_{\mathcal{U}}).$$

\begin{lemma} \label{Lemma-Basis4} The map $\xi_1:
L_1(\mathcal{M}^{**}) \to L_1(\mathcal{M}_{\mathcal{U}})$
satisfies $$\xi_1 \big( (ax)(yb) \big) = \xi_1(a^2)^{\frac12} \,
xy \, \xi_1(b^2)^{\frac12}$$ for all $(x,y) \in
L_{\infty}^r(\mathcal{M}, \mathsf{E}) \times
L_{\infty}^c(\mathcal{M}, \mathsf{E})$ and all positive elements
$a,b \in L_2^+(\mathcal{N}^{**})$.
\end{lemma}

\begin{proof} Accordingly to the terminology used in Lemma
\ref{Lemma-Basis3}, we may identify $\mathcal{M}^{**}$ with its
image $\gamma^*(\mathcal{M}^{**})$ in $\mathcal{M}_{\mathcal{U}}$.
In particular, the map $\xi_p$ can be regarded as an
$\mathcal{M}^{**}$ bimodule map satisfying
$$\xi_p(\mathrm{D}_{\psi}^{1/p}) = \mathrm{D}_{\psi \circ
\mathcal{E}}^{1/p}$$ for every functional $\psi \in \mathcal{M}^*$
with associated density $\mathrm{D}_{\psi} \in
L_1(\mathcal{M}^{**})$ so that $\psi \circ \mathcal{E}$ is a
functional on $\mathcal{M}_{\mathcal{U}}$ with associated density
$\mathrm{D}_{\psi \circ \mathcal{E}} \in
L_1(\mathcal{M}_{\mathcal{U}})$. Therefore, we have
$$\xi_1(\mathrm{D}_{\psi}) = \mathrm{D}_{\psi \circ
\mathcal{E}}^{1/2} \mathrm{D}_{\psi \circ \mathcal{E}}^{1/2} =
\xi_2(\mathrm{D}_{\psi}^{1/2}) \xi_2(\mathrm{D}_{\psi}^{1/2}).$$
Now, let us consider $a,b \in L_2(\mathcal{M}^{**})$ and define
$\mathrm{D}_{\psi} = aa^* + b^*b$ with corresponding positive
functional $\psi$. Let $m_1, m_2 \in \mathcal{M}^{**}$ such that
$a = \mathrm{D}_{\psi}^{1/2} m_1$ and $b = m_2
\mathrm{D}_{\psi}^{1/2}$. Assuming that $m_1$ and $m_2$ are
$\psi$-analytic and using the bimodule property of $\xi_1$ we
obtain
\begin{eqnarray*}
\xi_1(ab) & = & \xi_1 \big( \mathrm{D}_{\psi}^{1/2} m_1m_2
\mathrm{D}_{\psi}^{1/2} \big) \\ & = & \xi_1 \big(
\sigma_{-i/2}^{\psi}(m_1) \mathrm{D}_{\psi}
\sigma_{i/2}^{\psi}(m_2) \big) \\ & = & \sigma_{-i/2}^{\psi \circ
\mathcal{E}}(m_1) \mathrm{D}_{\psi \circ \mathcal{E}}^{1/2}
\mathrm{D}_{\psi \circ \mathcal{E}}^{1/2} \sigma_{i/2}^{\psi \circ
\mathcal{E}}(m_2) \\ & = & \mathrm{D}_{\psi \circ
\mathcal{E}}^{1/2} m_1 m_2 \mathrm{D}_{\psi \circ
\mathcal{E}}^{1/2} \\ & = & \xi_2(\mathrm{D}_{\psi}^{1/2}) m_1 m_2
\xi_2 (\mathrm{D}_{\psi}^{1/2}).
\end{eqnarray*}
Thus, by approximation with $\psi$-analytic elements, we conclude
$\xi_1(ab) = \xi_2(a) \xi_2(b)$. Assuming in addition that $a,b$
are positive and $x,y \in \mathcal{M}^{**}$
\begin{equation} \label{Equation-M-Identity}
\xi_1 \big( (ax)(yb) \big) = \xi_2(ax) \xi_2(yb) = \xi_2(a) xy
\xi_2(b) = \xi_1(a^2)^{\frac12} \, xy \, \xi_1(b^2)^{\frac12}.
\end{equation}
This proves the assertion with $a,b \in L_2^+(\mathcal{M}^{**})$
and $x,y \in \mathcal{M}$. Before considering elements $(x,y) \in
L_{\infty}^r(\mathcal{M}, \mathsf{E}) \times
L_{\infty}^c(\mathcal{M}, \mathsf{E})$, we observe that
$\xi_p(L_p(\mathcal{N}^{**}))$ is contained in the space
$$L_p(\mathcal{N}_{\mathcal{U}}) = \prodd_{i, \mathcal{U}}
L_p(\mathcal{N}) \quad \mbox{with} \quad \mathcal{N}_{\mathcal{U}}
= \Big( \prodd_{i, \mathcal{U}} \mathcal{N}_* \Big)^*.$$ Indeed,
according to Lemma \ref{Lemma-Basis3} we know that
$\mathsf{E}^{**} \circ \mathcal{E} = \mathcal{E} \circ \mathbf{E}$
so that $\mathcal{E}(L_p(\mathcal{N}_{\mathcal{U}}))$ is contained
in $L_p(\mathcal{N}^{**})$. According to
(\ref{Equation-M-Identity}) and the density of $\mathcal{M}$ in
$L_{\infty}^r(\mathcal{M}, \mathsf{E})$ and
$L_{\infty}^c(\mathcal{M}, \mathsf{E})$, to prove the assertion it
suffices to see that given $a,b \in L_2^+(\mathcal{N}^{**})$ and
$x,y\in \mathcal{M}$ we have
\begin{equation} \label{Equation-Continuity}
\big\| \xi_1 \big( (ax)(yb)\big)
\big\|_{L_1(\mathcal{M}_{\mathcal{U}})} \le
\|a\|_{L_2(\mathcal{N}^{**})} \big\| \mathsf{E}(xx^*)
\big\|_{\mathcal{M}}^{1/2} \big\| \mathsf{E}(y^*y)
\big\|_{\mathcal{M}}^{1/2} \|b\|_{L_2(\mathcal{N}^{**})}.
\end{equation}
By (\ref{Equation-M-Identity}) we have $$\big\| \xi_1 \big(
(ax)(yb)\big) \big\|_{L_1(\mathcal{M}_{\mathcal{U}})} \le
\mbox{tr}_{\mathcal{M}_{\mathcal{U}}} \Big( \xi_2(a) x x^*
\xi_2(a) \Big)^{1/2} \mbox{tr}_{\mathcal{M}_{\mathcal{U}}} \Big(
\xi_2(b) y^* y \xi_2(b) \Big)^{1/2}.$$ Now, using $\mathsf{E}^{**}
\circ \mathcal{E} = \mathcal{E} \circ \mathbf{E}$ and $\xi_2:
L_2(\mathcal{N}^{**}) \hookrightarrow
L_2(\mathcal{N}_{\mathcal{U}})$, we have
\begin{eqnarray*}
\mbox{tr}_{\mathcal{M}_{\mathcal{U}}} \Big( \xi_2(a) x x^*
\xi_2(a) \Big)^{1/2} & = & \mbox{tr}_{\mathcal{M}_{\mathcal{U}}}
\Big( \mathcal{E} \circ \mathbf{E} \big( \xi_2(a) x x^* \xi_2(a)
\big) \Big)^{1/2} \\ & = & \mbox{tr}_{\mathcal{M}_{\mathcal{U}}}
\Big( \xi_2(a) \mathsf{E}(x x^*) \xi_2(a) \Big)^{1/2} \\ & \le &
\|\xi_2(a)\|_{L_2(\mathcal{N}_{\mathcal{U}})} \big\|
\mathsf{E}(xx^*) \big\|_{\mathcal{M}}^{1/2}.
\end{eqnarray*}
Similarly, we have $$\mbox{tr}_{\mathcal{M}_{\mathcal{U}}} \Big(
\xi_2(b) y^* y \xi_2(b) \Big)^{1/2} \le \big\| \mathsf{E}(y^*y)
\big\|_{\mathcal{M}}^{1/2}
\|\xi_2(b)\|_{L_2(\mathcal{N}_{\mathcal{U}})}.$$ Therefore,
(\ref{Equation-Continuity}) follows since $\xi_2$ is a
contraction. This concludes the proof. \end{proof}

\begin{proposition} \label{Proposition-Factor-Norm}
Any $x \in \mathcal{M}$ satisfies
$$\|x\|_{L_{(2,\infty)}^{\infty}(\mathcal{M}, \mathsf{E}) \oM
L_{(\infty,2)}^{\infty}(\mathcal{M}, \mathsf{E})} =
\|x\|_{L_{(2,2)}^{\infty}(\mathcal{M}, \mathsf{E})}.$$
\end{proposition}

\begin{proof}
The lower estimate can be proved as in Lemma
\ref{Lemma-LaFactorizacion2}. To prove the upper estimate, we
consider the following norm on $\mathcal{M}$
$$|||x||| = \inf \left\{ \Big\| \big( \summ_k \mathsf{E}(a_ka_k^*)
\big)^{1/2} \Big\|_{L_{\infty}(\mathcal{N})} \Big\| \big( \summ_k
\mathsf{E}(b_k^*b_k) \big)^{1/2} \Big\|_{L_{\infty}(\mathcal{N})}
\right\}$$ where the infimum runs over all decompositions of $x$
into a finite sum $\summ_k a_kb_k$ with $a_k,b_k \in \mathcal{M}$.
Let $x \in \mathcal{M}$ and take a norm one functional $\phi:
\big( \mathcal{M}, ||| \ ||| \big) \to \C$ satisfying $|||x||| =
|\phi(x)|.$  Then it is clear that it suffices to see that
\begin{equation} \label{Equation-Sufficient}
|\phi(x)| \le \sup \Big\{ \|\alpha x \beta\|_{L_1(\mathcal{M})} \,
\big| \ \alpha, \beta \in \mathsf{B}_{L_2(\mathcal{N})} \Big\}.
\end{equation}
Applying the Grothendieck-Pietsch separation argument, we may find
states $\psi_1$ and $\psi_2$ in $\mathcal{N}^*$ with associated
densities $\mathrm{D}_1$ and $\mathrm{D}_2$ in
$L_1(\mathcal{N}^{**})$ satisfying the following inequality
$$|\phi(ab)| \le \psi_1(aa^*)^{\frac12} \psi_2(b^*b)^{\frac12}.$$
By Kaplansky's density theorem, $\mathrm{D}_j^{1/2} \mathcal{M}$
is norm dense in $\mathrm{D}_j^{1/2} \mathcal{M}^{**}$ for
$j=1,2$. Therefore, taking $e_j$ to be the support of
$\mathrm{D}_j$ for $j=1,2$, we may consider the map $$\Psi: e_1
L_2(\mathcal{M}^{**}) \to e_2 L_2(\mathcal{M}^{**})$$ determined
by the relation $$\phi(ab) = \big\langle \mathrm{D}_2^{\frac12}
b^*, \Psi(\mathrm{D}_1^{\frac12} a) \big\rangle =
\mbox{tr}_{\mathcal{M}^{**}} \Big( \Psi(\mathrm{D}_1^{\frac12} a)
b \mathrm{D}_2^{\frac12} \Big).$$ From this it easily follows that
$\Psi(\mathrm{D}_1^{\frac12}am) = \Psi(\mathrm{D}_1^{\frac12}a)m$
for all $m \in \mathcal{M}$. Hence, by density we deduce that
$\Psi$ commutes on $L_2(\mathcal{M}^{**})$ with the right action
on $\mathcal{M}^{**}$. In particular, we may find a contraction
$m_{\Psi} \in \mathcal{M}^{**}$ such that
$$\Psi(\mathrm{D}_1^{\frac12} a) = m_{\Psi} \mathrm{D}_1^{\frac12}
a.$$ This gives $$\phi(ab) = \mbox{tr}_{\mathcal{M}^{**}} \Big(
m_{\Psi} \mathrm{D}_1^{\frac12} a b \mathrm{D}_2^{\frac12}
\Big).$$ On the other hand, by Goldstein's theorem there exists a
net $(m_\lambda)$ in the unit ball of $\mathcal{M}$ which
converges to $m_{\Psi}$ in the $\sigma (\mathcal{M}^{**},
\mathcal{M}^*)$ topology. Therefore, since $\mathrm{D}_1^{1/2} a b
\mathrm{D}_2^{1/2} \in L_1(\mathcal{M}^{**})$ we deduce that
$$\phi(ab) = \lim_{\lambda} \mbox{tr}_{\mathcal{M}^{**}} \Big(
m_\lambda \mathrm{D}_1^{\frac12} a b \mathrm{D}_2^{\frac12}
\Big).$$ By Lemma \ref{Lemma-Basis2} we have
$$\mbox{tr}_{\mathcal{M}^{**}} \Big( m_\lambda
\mathrm{D}_1^{\frac12} a b \mathrm{D}_2^{\frac12} \Big) =
\psi_{\mathrm{D}_1^{\frac12} a b \mathrm{D}_2^{\frac12}}
(m_\lambda) = \mbox{tr}_{\mathcal{M}_{\mathcal{U}}} \Big(
m_\lambda \alpha \big( \psi_{\mathrm{D}_1^{\frac12} a b
\mathrm{D}_2^{\frac12}} \big) \Big).$$ Observing that $\xi_1 =
\alpha$ we may apply Lemma \ref{Lemma-Basis4} to obtain
$$\mbox{tr}_{\mathcal{M}_{\mathcal{U}}} \Big( m_\lambda \alpha
\big( \psi_{\mathrm{D}_1^{\frac12} a b \mathrm{D}_2^{\frac12}}
\big) \Big) = \limm_{i, \mathcal{U}} \mbox{tr}_{\mathcal{M}} \Big(
m_\lambda \xi_2(\mathrm{D}_1^{\frac12})_i \, ab \,
\xi_2(\mathrm{D}_2^{\frac12})_i \Big).$$ Therefore we conclude
$$\phi(x) = \lim_{\lambda} \limm_{i, \mathcal{U}}
\mbox{tr}_{\mathcal{M}} \Big( m_\lambda
\xi_2(\mathrm{D}_1^{\frac12})_i \, x \,
\xi_2(\mathrm{D}_2^{\frac12})_i \Big).$$ Moreover, since there
exist nets $(\alpha_i)$ and $(\beta_i)$ in $L_2^+(\mathcal{N})$
such that $$(\alpha_i)^{\bullet} = \xi_2(\mathrm{D}_1^{\frac12})
\quad \mbox{and} \quad (\beta_i)^{\bullet} = \xi_2(b^{\frac12}),$$
we obtain the following expression for $\phi(x)$ $$\phi(x) =
\lim_{\lambda} \limm_{i, \mathcal{U}} \mbox{tr}_{\mathcal{M}}
\big( m_\lambda \alpha_i \, x \, \beta_i \big).$$ Recalling that
$\lim_{i, \mathcal{U}} \|\alpha_i\|_2 \le 1$ and $\lim_{i,
\mathcal{U}} \|\beta_i\|_2 \le 1$, we obtain $$|\phi(x)| \le \sup
\Big\{ \|\alpha x \beta\|_{L_1(\mathcal{M})} \, \big| \ \alpha,
\beta \in \mathsf{B}_{L_2(\mathcal{N})} \Big\}.$$ This proves
(\ref{Equation-Sufficient}) and implies the assertion. The proof
is completed. \end{proof}

\chapter{Main results}
\label{Section7}

In this chapter we construct an isomorphic embedding
$$\mathcal{J}_{p,q}^n(\mathcal{M}, \mathsf{E}) \hookrightarrow
L_p(\mathcal{A}; \ell_q^n),$$ with $\mathcal{A}$ being a
sufficiently large von Neumann algebra. After that, we conclude by
studying the \emph{right} analogue of inequality $(\Sigma_{pq})$
(see the Introduction) in the noncommutative and operator space
levels. As we shall see, this appears as a particular case of our
embedding of $\mathcal{J}_{p,q}^n(\mathcal{M}, \mathsf{E})$ into
$L_p(\mathcal{A}; \ell_q^n)$ when considering the so-called
\emph{asymmetric} $L_p$ spaces, a particular case of conditional
noncommutative $L_p$ spaces. These results complete the general
method constructed in this paper and are closely related (as
explained in the Introduction) with several problems in the theory
of noncommutative $L_p$ spaces.

\section{Embedding into $L_p(\mathcal{A}; \ell_q^n)$}

Now we use the factorization results proved in the previous
chapter to identify the spaces $\mathcal{J}_{p,q}^n(\mathcal{M},
\mathsf{E})$ as an interpolation scale in $q$. This will give rise
to an isomorphic embedding of $\mathcal{J}_{p,q}^n(\mathcal{M},
\mathsf{E})$ into the space $L_p(\mathcal{A}; \ell_q^n)$, with
$\mathcal{A}$ determined by the map (\ref{Equation-Map-u})
$$u: x \in \mathcal{M} \mapsto \sum_{k=1}^n x_k \otimes \delta_k
\in \mathcal{A}_{\oplus n}.$$ Here $x_k = \pi_k(x,-x)$ in the
terminology of Chapter \ref{Section5}. That is,
$\mathcal{A}_{\oplus n}$ denotes the direct sum $\mathcal{A}
\oplus \mathcal{A} \oplus \ldots \oplus \mathcal{A}$ (with $n$
terms) of the $\mathcal{N}$-amalgamated free product
$*_{\mathcal{N}} \mathsf{A}_k$ with $\mathsf{A}_k = \mathcal{M}
\oplus \mathcal{M}$ for $1 \le k \le n$ and $\pi_k$ denotes the
natural embedding of $\mathsf{A}_k$ into $\mathcal{A}$. The
following lemma generalizes the main result in \cite{JP}.

\begin{lemma} \label{Lemma-Generalization-JP}
If $1 \le p \le \infty$, the map $$u: x \in
\mathcal{J}_{p,1}^n(\mathcal{M}, \mathsf{E}) \mapsto \sum_{k=1}^n
x_k \otimes \delta_k \in L_p(\mathcal{A}; \ell_1^n)$$ is an
isomorphism with complemented image and constants independent of
$p,n$.
\end{lemma}

\begin{proof} Since the result is clear for $p=1$, we shall assume
in what follows that $1 < p \le \infty$. According to Theorem
\ref{Theorem-LaFactorizacion}, we identify the intersection space
$\mathcal{J}_{p,1}^n(\mathcal{M}, \mathsf{E})$ with the
amalgamated tensor $\mathcal{R}_{2p,1}^n(\mathcal{M}, \mathsf{E})
\oM \mathcal{C}_{2p,1}^n(\mathcal{M}, \mathsf{E})$. Now, using the
characterization of $L_p(\mathcal{A}; \ell_1^n)$ given in
\cite{J1}, we have $$\big\| u(x) \big\|_{L_p(\mathcal{A};
\ell_1^n)} = \inf \left\{ \Big\| \big( \sum_{j,k} a_{kj} a_{kj}^*
\big)^{1/2} \Big\|_{L_{2p}(\mathcal{A})} \Big\| \big( \sum_{j,k}
b_{kj}^* b_{kj} \big)^{1/2} \Big\|_{L_{2p}(\mathcal{A})}
\right\},$$ with the infimum running over all possible
decompositions
$$x_k = \summ_j a_{kj} b_{kj}$$ and where sum above is required to
converge in the norm topology for $1 \le p < \infty$ and in the
weak operator topology otherwise. Then, given any decomposition of
$x$ into a finite sum $x = \sum_j \alpha_j \beta_j$ with $\alpha_j
\in \mathcal{R}_{2p,1}^n(\mathcal{M}, \mathsf{E})$ and $\beta_j
\in \mathcal{C}_{2p,1}^n(\mathcal{M}, \mathsf{E})$, we have $$x_k
= \summ_j a_{kj} b_{kj} \quad \mbox{with} \quad a_{kj} =
\pi_k(\alpha_j, \alpha_j) \quad \mbox{and} \quad b_{kj} =
\pi_k(\beta_j, -\beta_j).$$ This observation provides the
following estimate
\begin{eqnarray}
\label{Eq-SinSqFunct} \big\| u(x) \big\|_{L_p(\mathcal{A};
\ell_1^n)} & \le & \Big\| \summ_k \Lambda_{r,k} (\alpha) \otimes
e_{1k} \Big\|_{L_{2p}(\mathcal{A}; R_{2p}^n)} \\ \nonumber &
\times & \Big\| \summ_k \Lambda_{c,k}(\beta) \otimes e_{k1}
\Big\|_{L_{2p}(\mathcal{A}; C_{2p}^n)}
\end{eqnarray}
for any possible decomposition $x = \sum_j \alpha_j \beta_j$ and
where
\begin{eqnarray*}
\Lambda_{r,k}(\alpha) & = & \pi_k \Big[ \big( \summ_j \alpha_j
\alpha_j^* \big)^{1/2}, - \big( \summ_j \alpha_j \alpha_j^*
\big)^{1/2} \Big], \\ \Lambda_{c,k}(\beta) & = & \pi_k \Big[ \big(
\summ_j \beta_j^* \beta_j \big)^{1/2}, - \big( \summ_j \beta_j^*
\beta_j \big)^{1/2} \Big].
\end{eqnarray*}
We want to reformulate \eqref{Eq-SinSqFunct} in terms of the
square functions $$\mathcal{S}_r(\alpha) = \big( \summ_j \alpha_j
\alpha_j^* \big)^{1/2} \quad \mbox{and} \quad \mathcal{S}_c(\beta)
= \big( \summ_j \beta_j^* \beta_j \big)^{1/2}.$$ By the definition
of the map $u$, we find $$\big\| u(x) \big\|_{L_p(\mathcal{A};
\ell_1^n)} \le \inf \left\{ \big\| u \big( \mathcal{S}_r(\alpha)
\big) \big\|_{L_{2p}(\mathcal{A}; R_{2p}^n)} \big\| u \big(
\mathcal{S}_c(\beta) \big) \big\|_{L_{2p}(\mathcal{A}; C_{2p}^n)}
\right\},$$ where the infimum runs over decompositions $x = \sum_k
\alpha_j \beta_j$. However, applying Corollary
\ref{Corollary-Voiculescu} as we did in the proof of Lemma
\ref{Lemma-Complemented-Isomorphism}, the spaces
$\mathcal{R}_{2p,1}^n(\mathcal{M}, \mathsf{E})$ and
$\mathcal{C}_{2p,1}^n(\mathcal{M}, \mathsf{E})$ are isomorphic to
their images via $u$ in $L_{2p}(\mathcal{A}; R_{2p}^n)$ and
$L_{2p}(\mathcal{A}; C_{2p}^n)$ respectively. In particular, we
obtain the following inequality up to a constant independent of
$p,n$ $$\big\| u(x) \big\|_{L_p(\mathcal{A}; \ell_1^n)} \lesssim
\inf \left\{ \big\| \mathcal{S}_r(\alpha)
\big\|_{\mathcal{R}_{2p,1}^n(\mathcal{M}, \mathsf{E})} \big\|
\mathcal{S}_c(\beta) \big\|_{\mathcal{C}_{2p,1}^n(\mathcal{M},
\mathsf{E})} \right\}.$$ The right hand side is the norm of $x$ in
$\mathcal{R}_{2p,1}^n(\mathcal{M}, \mathsf{E}) \oM
\mathcal{C}_{2p,1}^n(\mathcal{M}, \mathsf{E})$. Thus, we have
proved that $u$ is a bounded map. Now we prove the reverse
estimate arguing by duality. First we note that it easily follows
from Theorem \ref{Theorem-Duality-Conditional} and Remark
\ref{Remark-Isometric-Duality} that \label{K4}
$$\mathcal{J}_{p,1}^n(\mathcal{M}, \mathsf{E}) =
\mathcal{K}_{p',\infty}^n(\mathcal{M}, \mathsf{E})^* \quad
\mbox{for} \quad 1 < p \le \infty$$ where, following the notation
of \cite{JP}, the latter space is given by
$$\mathcal{K}_{p',\infty}^n(\mathcal{M}, \mathsf{E}) = \sum_{u,v \in
\{2p',\infty\}} n^{-(\frac{1}{u} + \frac{1}{p} + \frac{1}{v})}
L_u(\mathcal{N}) L_{\rho_{uv}}(\mathcal{M}) L_v(\mathcal{N}),$$
with $\rho_{uv}$ determined by $1/u + 1/\rho_{uv} + 1/v = 1/p'$.
We claim that
\begin{equation} \label{Equation-Claim-Dual}
w: x \in \mathcal{K}_{p',\infty}^n(\mathcal{M}, \mathsf{E})
\mapsto \frac{1}{n} \sum_{k=1}^n x_k \otimes \delta_k \in
L_{p'}(\mathcal{A}; \ell_{\infty}^n)
\end{equation}
is a contraction. It is not difficult to see that the result
follows from our claim. Indeed, using the anti-linear duality
bracket we find that $$\langle u(x), w(y) \rangle = \frac{1}{n}
\sum_{k=1}^n\mbox{tr}_{\mathcal{A}} \big( \pi_k(x,-x)^*
\pi_k(y,-y) \big) = \frac{1}{n} \sum_{k=1}^n
\mbox{tr}_{\mathcal{M}} (x^*y) = \langle x,y \rangle.$$ In
particular, it turns out that $w^*u$ is the identity map on
$\mathcal{J}_{p,1}^n(\mathcal{M}, \mathsf{E})$ and $uw^*$ is a
bounded projection from $L_p(\mathcal{A}; \ell_1^n)$ onto
$u(\mathcal{J}_{p,1}^n(\mathcal{M}, \mathsf{E}))$. Thus, $u$ is an
isomorphism and its image is complemented with constants
independent of $p,n$. Therefore, it remains to prove our claim.
Since the initial space of $w$ is a sum of four Banach spaces
indexed by the pairs $(u,v)$ with $u,v \in \{2p',\infty\}$, it
suffices to see that the restriction of $w$ to each of these
spaces is a contraction. Let us start with the case $(u,v) =
(\infty,\infty)$. In this case the associated space is simply
$n^{-1/p} L_{p'}(\mathcal{M})$ and we have by complex
interpolation in $1 \le p'\le \infty$ that
$$\frac{1}{n} \Big\| \sum_{k=1}^n x_k \otimes \delta_k
\Big\|_{L_{p'}(\mathcal{A}; \ell_{\infty}^n)} \le n^{-\frac1p}
\|x\|_{L_{p'}(\mathcal{M})}.$$ Indeed, when $p=1$ we have an
identity while for $p=\infty$ the assertion follows by the
triangle inequality. Therefore, recalling from \cite{JX2} that the
spaces $L_{p'}(\mathcal{M}, \ell_\infty^n)$ form an interpolation
scale in $p'$, our claim follows. The space associated to $u = 2p'
= v$ is
$$\frac{1}{n} \, L_{2p'}(\mathcal{N}) L_{\infty}(\mathcal{M})
L_{2p'}(\mathcal{N}).$$ Here we use the characterization of the
norm of $L_{p'}(\mathcal{A}; \ell_{\infty}^n)$ given in \cite{J1}
$$\Big\| \sum_{k=1}^n x_k \otimes \delta_k
\Big\|_{L_{p'}(\mathcal{A}; \ell_{\infty}^n)} = \inf_{x_k = a d_k
b} \left\{ \|a\|_{L_{2p'}(\mathcal{A})} \Big( \sup_{1 \le k \le n}
\|d_k\|_{L_{\infty}(\mathcal{A})} \Big)
\|b\|_{L_{2p'}(\mathcal{A})} \right\}.$$ Then we consider the
decompositions of $x_k$ that come from decompositions of $x$. In
other words, any decomposition $x = \alpha y \beta$ with $\alpha,
\beta \in L_{2p'}(\mathcal{N})$ and $y \in
L_{\infty}(\mathcal{M})$ gives rise to the decompositions $x_k =
\alpha \pi_k(y,-y) \beta$. Hence, since
$$\|\pi_k(y,-y)\|_{L_{\infty}(\mathcal{A})} =
\|y\|_{L_{\infty}(\mathcal{M})} \quad \mbox{for} \quad 1 \le k \le
n,$$ we obtain the desired estimate $$\Big\| \sum_{k=1}^n x_k
\otimes \delta_k \Big\|_{L_{p'}(\mathcal{A}; \ell_{\infty}^n)} \le
\inf_{x = \alpha y \beta} \Big\{ \|\alpha\|_{L_{2p'}(\mathcal{N})}
\|y\|_{L_{\infty}(\mathcal{M})} \|\beta\|_{L_{2p'}(\mathcal{N})}
\Big\} = \|x\|_{2p' \cdot \infty \cdot 2p'}.$$ Finally, it remains
to consider the cross terms associated to $(u,v) = (2p',\infty)$
and $(u,v) = (\infty,2p')$. Since both can be treated using the
same arguments, we only consider the case $(u,v) = (2p',\infty)$
whose associated space is $$\frac{1}{n^{\frac{1}{2} +
\frac{1}{2p}}} \, L_{2p'}(\mathcal{N}) L_{2p'}(\mathcal{M})
L_{\infty}(\mathcal{N}).$$ Here we observe that given any $\alpha
\in L_{2p'}(\mathcal{N})$, the left multiplication mapping
$$\mathsf{L}_{\alpha}: \sum_{k=1}^n a_k \otimes \delta_k \in
L_{2p'}(\mathcal{A};\ell_{2p'}^n) \mapsto \alpha \sum_{k=1}^n a_k
\otimes \delta_k \in L_{p'}(\mathcal{A};\ell_{\infty}^n)$$ is a
bounded map with norm $\le \|\alpha\|_{L_{2p'}(\mathcal{A})}$.
Indeed, by complex interpolation we only study the extremal cases.
Noting that the result is trivial for $p'= \infty$, it suffices to
see it for $p'= 1$. Let us assume that
$$\sum_{k=1}^n \|a_k\|_{L_2(\mathcal{A})}^2 \le 1.$$ Then we
define the invertible operator $\beta = \big( \summ_k a_k^*a_k +
\delta \mathrm{D}_{\phi} \big)^{1/2}$ ($\mathrm{D}_{\phi}$ being
the density associated to the state $\phi$ of $\mathcal{A}$) so
that $a_k = b_k \beta$ for $1 \le k \le n$ where the operators
$b_1, b_2, \ldots, b_n$ satisfy
$$\sum_{k=1}^n b_k^* b_k \le 1.$$ This provides a factorization
$\alpha a_k = \alpha b_k \beta$ from which we deduce $$\Big\|
\alpha \sum_{k=1}^n a_k \otimes \delta_k \Big\|_{L_1(\mathcal{A};
\ell_{\infty}^n)} \le \|\alpha\|_{L_2(\mathcal{A})} \Big( \sup_{1
\le k \le n} \|b_k\|_{L_{\infty}(\mathcal{A})} \Big)
\|\beta\|_{L_2(\mathcal{A})} \le (1+\delta) \,
\|\alpha\|_{L_2(\mathcal{A})}.$$ Thus, we conclude letting $\delta
\to 0$. Now assume that $x$ is a norm 1 element of
$$\frac{1}{n^{\frac{1}{2} + \frac{1}{2p}}} \, L_{2p'}(\mathcal{N})
L_{2p'}(\mathcal{M}) L_{\infty}(\mathcal{N}),$$ so that for any
$\delta
> 0$ we can find a factorization $x = \alpha y$ such that
$$\|\alpha\|_{L_{2p'}(\mathcal{N})} \le 1 \quad \mbox{and} \quad
\|y\|_{L_{2p'}(\mathcal{M})} \le (1+\delta) \, n^{\frac12 +
\frac{1}{2p}}.$$ Then, since $\alpha \in L_{2p'}(\mathcal{N})$ we
have
$$\sum_{k=1}^n x_k \otimes \delta_k = \alpha \sum_{k=1}^n
\pi_k(y,-y) \otimes \delta_k = \mathsf{L}_{\alpha} \Big(
\sum_{k=1}^n \pi_k(y,-y) \otimes \delta_k \Big).$$ In particular,
\begin{eqnarray*}
\frac{1}{n} \, \Big\| \sum_{k=1}^n x_k \otimes \delta_k
\Big\|_{L_{p'}(\mathcal{A}; \ell_{\infty}^n)} & \le & \frac{1}{n}
\, \|\alpha\|_{L_{2p'}(\mathcal{N})} \Big\| \sum_{k=1}^n
\pi_k(y,-y) \otimes \delta_k \Big\|_{L_{2p'}(\mathcal{A};
\ell_{2p'}^n)} \\ & \le & \frac{1}{n} \, \Big( \sum_{k=1}^n
\|\pi_k(y,-y)\|_{L_{2p'}(\mathcal{A})}^{2p'} \Big)^{1/2p'} \le
1+\delta.
\end{eqnarray*}
Letting $\delta \to 0$, we conclude that the mapping defined in
(\ref{Equation-Claim-Dual}) is a contraction. \end{proof}

\begin{theorem} \label{Theorem-Interpolation-J}
If $1 \le p \le \infty$, then $$\big[
\mathcal{J}_{p,1}^n(\mathcal{M}, \mathsf{E}),
\mathcal{J}_{p,p}^n(\mathcal{M}, \mathsf{E}) \big]_{\theta} \simeq
\mathcal{J}_{p,q}^n(\mathcal{M}, \mathsf{E})$$ with $1/q =
1-\theta + \theta/p$ and with relevant constants independent of
$n$.
\end{theorem}

\begin{proof} By Theorem
\ref{Theorem-Interpolation-Conditional1} and Observation
\ref{Observation-Extension}, we have $$\big[
\mathcal{J}_{p,1}^n(\mathcal{M}, \mathsf{E}),
\mathcal{J}_{p,p}^n(\mathcal{M}, \mathsf{E}) \big]_{\theta}
\subset \mathcal{J}_{p,q}^n(\mathcal{M}, \mathsf{E})$$
contractively. To prove the reverse inclusion we consider the map
$u: \mathcal{M} \to \mathcal{A}_{\oplus n}$.

\vskip5pt

\noindent \textsc{Step 1}. Since $\mathcal{J}_{p,p}^n(\mathcal{M},
\mathsf{E}) = n^{\frac1p} L_p(\mathcal{M})$, we clearly have
$$\|u(x)\|_{L_p(\mathcal{A}; \ell_p^n)} =
\|x\|_{\mathcal{J}_{p,p}^n(\mathcal{M}, \mathsf{E})}.$$ Moreover,
$u(\mathcal{J}_{p,p}^n(\mathcal{M}, \mathsf{E}))$ is clearly
contractively complemented in $L_p(\mathcal{A}; \ell_p^n)$. On the
other hand, according to Lemma \ref{Lemma-Generalization-JP},
$\mathcal{J}_{p,1}^n (\mathcal{M}, \mathsf{E})$ is isomorphic to
its image via $u$, which is complemented in
$L_p(\mathcal{A};\ell_1^n)$ with constants not depending on $p,n$.
Therefore, the norm of $x$ in the space
$[\mathcal{J}_{p,1}^n(\mathcal{M}, \mathsf{E}),
\mathcal{J}_{p,p}^n(\mathcal{M}, \mathsf{E})]_{\theta}$ turns out
to be equivalent to the norm of
$$\sum_{k=1}^n x_k \otimes \delta_k \quad \mbox{in} \quad
L_p(\mathcal{A}; \ell_q^n).$$ Thus, it remains to see that
\begin{equation} \label{Equ-Last1}
\Big\| \sum_{k=1}^n x_k \otimes \delta_k
\Big\|_{L_p(\mathcal{A}:\ell_q^n)} \lesssim
\|x\|_{\mathcal{J}_{p,q}^n(\mathcal{M}, \mathsf{E})}.
\end{equation}

\vskip3pt

\noindent \textsc{Step 2}. We need an auxiliary result. Let us
consider the bilinear map $$\Lambda: (a,b) \in
\mathcal{R}_{2p,q}^n(\mathcal{M}, \mathsf{E}) \times
\mathcal{C}_{2p,q}^n(\mathcal{M}, \mathsf{E}) \mapsto \sum_{k=1}^n
\pi_k(ab,-ab) \otimes \delta_k \in L_p(\mathcal{A};\ell_q^n).$$ We
claim that $\Lambda$ is bounded by $c(p,q) \sim (p-q)/(pq+q-p)$.
By Theorem \ref{Theorem-Intersection1}, we may use bilinear
interpolation and it suffices to show our claim for the extremal
values of $q$. Let us note that, since we are applying Theorem
\ref{Theorem-Intersection1} for $\mathcal{R}_{2p,q}^n(\mathcal{M},
\mathsf{E})$ and $\mathcal{C}_{2p,q}^n(\mathcal{M}, \mathsf{E})$,
the operator norm of $\Lambda$ is controlled by the square of the
relevant constant in Theorem \ref{Theorem-Intersection1}, as we
have claimed. The boundedness of $\Lambda$ for $q=1$ follows from
Theorem \ref{Theorem-LaFactorizacion} and Lemma
\ref{Lemma-Generalization-JP}. The estimate for $q=p$ is much
easier $$\big\| \Lambda(a,b) \big\|_{L_p(\mathcal{A};\ell_p^n)} =
\Big( \sum_{k=1}^n \big\| \pi_k(ab,-ab) \big\|_p^p \Big)^{\frac1p}
= n^{\frac1p} \, \|ab\|_p \le \|a\|_{\mathcal{R}_{2p,p}}
\|b\|_{\mathcal{C}_{2p,p}}.$$

\vskip3pt

\noindent \textsc{Step 3}. Now it is easy to deduce inequality
\eqref{Equ-Last1} from the boundedness of the mapping $\Lambda$.
Indeed, according to Theorem \ref{Theorem-LaFactorizacion} we know
that \eqref{Equ-Last1} is equivalent to
\begin{equation} \label{Equ-Last2}
\Big\| \sum_{k=1}^n x_k \otimes \delta_k
\Big\|_{L_p(\mathcal{A}:\ell_q^n)} \lesssim \Big\| \Big( \summ_j
a_j a_j^* \Big)^{\frac12} \Big\|_{\mathcal{R}_{2p,q}^n} \Big\|
\Big( \summ_j b_j^* b_j \Big)^{\frac12}
\Big\|_{\mathcal{C}_{2p,q}^n},
\end{equation}
for any decomposition of $x$ into a finite sum $\sum_j a_j b_j$.
Let us assume that the index $j$ runs from $1$ to $m$ for some
finite $m$. Then we consider the matrix amplifications
$\hat{\mathcal{N}} = \mathrm{M}_m \otimes \mathcal{N}$ and
$\hat{\mathcal{M}} = \mathrm{M}_m \otimes \mathcal{M}$. Similarly,
we consider $\hat{\mathcal{A}} = \mathrm{M}_m \otimes
\mathcal{A}$. According to \eqref{Equation-CB-Amalgamated}, we
have
$$\hat{\mathcal{A}} = \mathrm{M}_m \otimes
\mathsf{A}_1 \, *_{\hat{\mathcal{N}}} \, \mathrm{M}_m \otimes
\mathsf{A}_2 \, *_{\hat{\mathcal{N}}} \, \cdots
*_{\hat{\mathcal{N}}} \, \mathrm{M}_m \otimes \mathsf{A}_n = \big(
\hat{\mathcal{M}} \oplus \hat{\mathcal{M}}
\big)^{*_{\hat{\mathcal{N}}}^n}.$$ By Step2, we deduce
\begin{equation} \label{Equ-Last3}
\Big\| \sum_{k=1}^n \mathbf{x}_k \otimes \delta_k
\Big\|_{L_p(\hat{\mathcal{A}}; \ell_q^n)} \lesssim
\|\mathbf{a}\|_{\hat{\mathcal{R}}_{2p,q}^n}
\|\mathbf{b}\|_{\hat{\mathcal{C}}_{2p,q}^n}
\end{equation}
where the elements $\mathbf{x}_k$ are given by $\mathbf{x}_k =
\pi_k(\mathbf{ab}, - \mathbf{ab})$ and
\begin{eqnarray*}
\hat{\mathcal{R}}_{2p,q}^n & = & \mathcal{R}_{2p,q}^n
\big( \hat{\mathcal{M}}, 1_{\mathrm{M}_m} \otimes \mathsf{E} \big), \\
\hat{\mathcal{C}}_{2p,q}^n & = & \ \mathcal{C}_{2p,q}^n \big(
\hat{\mathcal{M}}, 1_{\mathrm{M}_m} \otimes \mathsf{E} \big).
\end{eqnarray*}
To prove the remaining estimate \eqref{Equ-Last2} we fix $x$ in
$\mathcal{J}_{p,q}^n(\mathcal{M}, \mathsf{E})$ and decompose it
into a finite sum $\sum_j a_j b_j$ with $m$ terms. Then, we define
the row and column matrices
$$\mathbf{a} = \sum_{j=1}^m a_j \otimes e_{1j} \quad \mbox{and}
\quad \mathbf{b} = \sum_{j=1}^m b_j \otimes e_{j1}.$$ According to
\eqref{Equ-Last3}, it just remains to show that
\begin{eqnarray*}
\Big\| \sum_{k=1}^n \mathbf{x}_k \otimes \delta_k
\Big\|_{L_p(\hat{\mathcal{A}}; \ell_q^n)} & = & \Big\|
\sum_{k=1}^n x_k \otimes \delta_k \Big\|_{L_p(\mathcal{A};
\ell_q^n)}, \\ \|\mathbf{a}\|_{\hat{\mathcal{R}}_{2p,q}^n} & = &
\Big\| \Big( \summ_j a_j a_j^* \Big)^{\frac12}
\Big\|_{\mathcal{R}_{2p,q}^n}, \\
\|\mathbf{b}\|_{\hat{\mathcal{C}}_{2p,q}^n} \hskip3pt & = & \Big\|
\Big( \summ_j b_j^* \, b_j \Big)^{\frac12}
\Big\|_{\mathcal{C}_{2p,q}^n}.
\end{eqnarray*}
The first identity follows easily from
\begin{eqnarray*}
\mathbf{x}_k & = & \pi_k(\mathbf{ab}, - \mathbf{ab}) \\ & = &
\summ_j \pi_k \big( a_j b_j \otimes e_{11}, -a_j b_j \otimes
e_{11} \big) \\ & = & \pi_k \Big( \summ_j a_j b_j, - \summ_j a_j
b_j \Big) \otimes e_{11} \\ & = & x_k \otimes e_{11},
\end{eqnarray*}
where the third identity holds because $\pi_k$ is an
$\hat{\mathcal{N}}$-bimodule map. The relation $\mathbf{x}_k = x_k
\otimes e_{11}$ shows that $\sum_k \mathbf{x}_k \otimes \delta_k$
lives in fact in the $(1,1)$ corner, which by \cite{JX4} is
isometrically isomorphic to $L_p(\mathcal{A}; \ell_q^n)$. The
identities for $\mathbf{a}$ and $\mathbf{b}$ are proved in the
same way, so that we only consider the first of them. We have
$$\|\mathbf{a}\|_{\hat{\mathcal{R}}_{2p,q}^n} = \max \Big\{
n^{\frac{1}{2p}} \, \|\mathbf{a}\|_{L_{2p}(\hat{\mathcal{M}})},
n^{\frac{1}{2q}} \, \|\mathbf{a}\|_{L_{(\frac{2pq}{p-q},
\infty)}^{2p}(\hat{\mathcal{M}}, 1_{\mathrm{M}_m} \otimes
\mathsf{E})} \Big\}.$$ It is clear that
$$\|\mathbf{a}\|_{L_{2p}(\hat{\mathcal{M}})} = \Big\| \Big(
\summ_j a_j a_j^* \Big)^{\frac12} \Big\|_{L_{2p}(\mathcal{M})}.$$
Therefore, we just need to show that
$$\|\mathbf{a}\|_{L_{(\frac{2pq}{p-q}, \infty)}^{2p}
(\hat{\mathcal{M}}, 1_{\mathrm{M}_m} \otimes \mathsf{E})} = \Big\|
\Big( \summ_j a_j a_j^* \Big)^{\frac12}
\Big\|_{L_{(\frac{2pq}{p-q}, \infty)}^{2p}(\mathcal{M},
\mathsf{E})}.$$ By definition, we have
\begin{eqnarray*}
\|\mathbf{a}\|_{L_{(\frac{2pq}{p-q}, \infty)}^{2p}
(\hat{\mathcal{M}}, 1_{\mathrm{M}_m} \otimes \mathsf{E})} & = &
\sup_{\alpha} \, \big\| \alpha \mathbf{a} \mathbf{a}^* \alpha^*
\big\|_q^{\frac12} \\ & = & \sup_{\alpha} \, \Big\| \alpha \Big[
\Big( \summ_j a_j a_j^* \Big) \otimes e_{11} \Big] \alpha^*
\Big\|_q^{\frac12} \\ & = & \sup_{\alpha} \, \Big\| \alpha \Big[
\Big( \summ_j a_j a_j^* \Big)^{\frac12} \otimes e_{11} \Big]
\Big\|_{2q},
\end{eqnarray*}
where the supremum runs over all $\alpha$ such that
$$\|\alpha\|_{L_{\frac{2pq}{p-q}}(\hat{\mathcal{N}})} \le 1.$$ If
$(\alpha_{ij})$ are the $m \times m$ entries of $\alpha$, we
clearly have $$\alpha \Big[ \Big( \summ_j a_j a_j^*
\Big)^{\frac12} \otimes e_{11} \Big] = \sum_{s=1}^m \alpha_{s1}
\Big( \summ_j a_j a_j^* \Big)^{\frac12} \otimes e_{s1}.$$
Therefore, we obtain the following estimate
\begin{eqnarray*}
\lefteqn{\Big\| \alpha \Big[ \Big( \summ_j a_j a_j^*
\Big)^{\frac12} \otimes e_{11} \Big]
\Big\|_{L_{2q}(\hat{\mathcal{M}})}} \\ & = & \Big\| \Big( \summ_s
\Big[ \summ_j a_j a_j^* \Big]^{\frac12} \alpha_{s1}^* \alpha_{s1}
\Big[ \summ_j a_j a_j^* \Big]^{\frac12}
\Big)^{\frac12}\Big\|_{L_{2q}(\mathcal{M})}
\\ & = & \Big\| \Big( \summ_s \alpha_{s1}^* \alpha_{s1}
\Big)^{\frac12} \Big( \summ_j a_j a_j^* \Big)^{\frac12}
\Big\|_{L_{2q}(\mathcal{M})} \\ & \le & \Big\| \summ_s \alpha_{s1}
\otimes e_{s1} \Big\|_{L_{\frac{2pq}{p-q}}(\hat{\mathcal{N}})}
\Big\| \Big( \summ_j a_j a_j^* \Big)^{\frac12}
\Big\|_{L_{(\frac{2pq}{p-q}, \infty)}^{2p}(\mathcal{M},
\mathsf{E})}.
\end{eqnarray*}
Hence, the result follows since the column projection is
contractive on $L_{\frac{2pq}{p-q}}(\hat{\mathcal{N}})$.
\end{proof}

The following result is the main result of this paper.

\begin{theorem} \label{Theorem-Embed-OH}
If $1 \le q \le p \le \infty$, the map
$$u: x \in \mathcal{J}_{p,q}^n (\mathcal{M}, \mathsf{E}) \mapsto
\sum_{k=1}^n x_k \otimes \delta_k \in L_p(\mathcal{A}; \ell_q^n)$$
is an isomorphism with complemented image and constants
independent of $n$.
\end{theorem}

\begin{proof}
It is an immediate consequence of Lemma
\ref{Lemma-Generalization-JP} and Theorem
\ref{Theorem-Interpolation-J}
\end{proof}

\begin{remark}
\emph{If $1 \le q \le p \le \infty$, let us define
$$\mathcal{K}_{p',q'}^n (\mathcal{M}, \mathsf{E}) = \sum_{u,v \in
\{2r,\infty\}}^{\null} n^{-1 + \frac{1}{\rho_{uv}}} \,
L_u(\mathcal{N}) L_{\rho_{uv}}(\mathcal{M}) L_v(\mathcal{N})$$
where $1/r = 1/q - 1/p$ and $\rho_{uv}$ is determined by $1/u +
1/\rho_{uv} + 1/v = 1/p'$. Note that this definition is consistent
with the space $\mathcal{K}_{p',\infty}^n(\mathcal{M},
\mathsf{E})$ introduced in the proof of Lemma
\ref{Lemma-Generalization-JP}. Arguing by duality as in \cite{JP},
we easily conclude from Theorem \ref{Theorem-Embed-OH} that we
have an isomorphism with complemented image
\begin{equation} \label{Eq-Sum-Nocaps}
\xi_{free}: x \in \mathcal{K}_{p',q'}^n (\mathcal{M}, \mathsf{E})
\mapsto \sum_{k=1}^n x_k \otimes \delta_k \in
L_{p'}(\mathcal{A}_{free}; \ell_{q'}^n)
\end{equation}
where $\mathcal{A}_{free}$ denotes the usual free product algebra
$\mathcal{A}$ from Theorem \ref{Theorem-Embed-OH}. Moreover, it is
important to note that replacing free products by tensors and
freeness by noncommutative independence as in \cite{JP}, Theorem
\ref{Theorem-Embed-OH} and \eqref{Eq-Sum-Nocaps} hold in the range
$1 < p' \le q' \le \infty$. In other words, we replace
$\mathcal{A}_{free}$ by the tensor product algebra
$$\mathcal{A}_{ind} = \mathcal{M} \otimes \mathcal{M} \otimes
\cdots \otimes \mathcal{M}$$ with $n$ terms and $x_k$ is now given
by $$x_k = 1 \otimes \cdots 1 \otimes x \otimes 1 \otimes \cdots
\otimes 1$$ where $x$ is placed in the $k$-th position. The
question now is whether or not \eqref{Eq-Sum-Nocaps} holds in the
non-free setting for $p'= 1$. Using recent results from \cite{J3},
we shall see in a forthcoming paper that we have an isomorphism
with complemented image $$\xi_{ind}: x \in \mathcal{K}_{1,2}^n
(\mathcal{M}, \mathsf{E}) \mapsto \sum_{k=1}^n x_k \otimes
\delta_k \in L_{p'}(\mathcal{A}_{ind}; \mathrm{OH}_n).$$ Let us
note that the same question for $\mathcal{K}_{1,q'}^n
(\mathcal{M}, \mathsf{E})$ with $q' \neq 2$ is still open.}
\end{remark}

\section{Asymmetric $L_p$ spaces and $(\Sigma_{pq})$}

Let $\mathcal{M}$ be a von Neumann algebra equipped with a
\emph{n.f.} state $\varphi$. Now let $2 \le u,v \le \infty$ be
such that $1/p = 1/u + 1/v$ for some $1 \le p \le \infty$. Then,
we define the \emph{asymmetric $L_p$ space} \label{LpAsimetrico}
associated to the pair $(u,v)$ as the $\mathcal{M}$-amalgamated
Haagerup tensor product
\begin{equation} \label{Equation-Asymmetric}
L_{(u,v)}(\mathcal{M}) = L_u^r(\mathcal{M})
\otimes_{\mathcal{M},h} L_v^c(\mathcal{M}),
\end{equation}
where we recall that $L_q^r(\mathcal{M})$ and $L_q^c(\mathcal{M})$
were defined in \eqref{Et-LpAsimet} for $2 \le q \le \infty$. That
is, we consider the quotient of $L_u^r(\mathcal{M}) \otimes_h
L_v^c(\mathcal{M})$ by the closed subspace $\mathcal{I}$ generated
by the differences $x_1 \gamma \otimes x_2 - x_1 \otimes \gamma
x_2$ with $\gamma \in \mathcal{M}$. Recall that we are using the
notation $\otimes_{\mathcal{M},h}$ instead of $\oM$ because, in
contrast with the previous chapters, we shall be interested here
in the operator space structure rather than in the Banach space
one. By a well known factorization argument (see e.g. Lemma 3.5 in
\cite{P2}), the norm of an element $x$ in $L_{(u,v)}(\mathcal{M})$
is given by $$\|x\|_{(u,v)} = \inf_{x = \alpha \beta}
\|\alpha\|_{L_u(\mathcal{M})} \|\beta\|_{L_v(\mathcal{M})}.$$
According to this observation, it turns out that asymmetric $L_p$
spaces arise as a particular case of the amalgamated
noncommutative $L_p$ spaces defined in Chapter \ref{Section2} when
we take $q=\infty$ and $\mathcal{N}$ to be $\mathcal{M}$ itself.
Asymmetric $L_p$ spaces were introduced in \cite{JP} for finite
matrix algebras. There the amalgamated Haagerup tensor product
used in (\ref{Equation-Asymmetric}) was not needed in \cite{JP} to
define the asymmetric Schatten classes. In fact, if $\mathcal{M}$
is the algebra $\mathrm{M}_m$ of $m \times m$ matrices and
$\mathrm{X}$ is an operator space, we can define the
\emph{vector-valued asymmetric Schatten class} \label{VVAsymSp}
$$S_{(u,v)}^m(\mathrm{X}) = C_{u/2}^m \otimes_h \mathrm{X}
\otimes_h R_{v/2}^m.$$ Note that this definition is consistent
with (\ref{Equation-Asymmetric}). Indeed, recalling that
$$L_u^r(\mathrm{M}_m) = C_{u/2}^m \otimes_h R_m \quad \mbox{and}
\quad L_v^c(\mathrm{M}_m) = C_m \otimes_h R_{v/2}^m,$$ it can be
easily checked from the definition of the Haagerup tensor product
that we have $L_u^r(\mathrm{M}_m) \otimes_{\mathcal{M},h}
L_v^c(\mathrm{M}_m) = S_{(u,v)}^m(\C)$ isometrically. Moreover,
according to \cite{JP} any linear map $u: \mathrm{X}_1 \rightarrow
\mathrm{X}_2$ satisfies
\begin{equation} \label{Equation-cb-Asymmetric}
\|u\|_{cb} = \sup_{n \ge 1} \, \Big\| 1_{\mathrm{M}_n} \otimes u:
S_{(u,v)}^n(\mathrm{X}_1) \rightarrow S_{(u,v)}^n(\mathrm{X}_2)
\Big\|.
\end{equation}
In particular, since it is clear that
$$S_{(u,v)}^n \big( L_u^r(\mathrm{M}_m) \otimes_{\mathcal{M},h}
L_v^c(\mathrm{M}_m) \big) = S_{(u,v)}^{mn}(\C) \quad \mbox{for
all} \quad m \ge 1,$$ we have the following completely isometric
isomorphism $$L_{(u,v)}(\mathrm{M}_m) = S_{(u,v)}^m(\C).$$

\begin{remark} \label{Rem-Lp-2p2p}
\emph{$L_p(\mathcal{M})$ is completely isometric to
$L_{2p}^r(\mathcal{M}) \otimes_{\mathcal{M},h}
L_{2p}^c(\mathcal{M})$.}
\end{remark}

We conclude by generalizing the inequalities $(\Sigma_{pq})$
stated in the Introduction to the noncommutative setting.
Moreover, we shall seek for a completely isomorphic embedding
rather than a Banach space one. As we shall see, this appears as a
particular case of Theorem \ref{Theorem-Embed-OH} module the
corresponding identifications. Indeed, let $\mathcal{M}$ be a von
Neumann algebra equipped with a \emph{n.f.} state $\varphi$ and
let us consider the particular situation in which $(\mathcal{M},
\mathcal{N}, \mathsf{E})$ above are replaced by
\begin{equation} \label{Equation-MNE}
(\mathcal{M}_m, \mathcal{N}_m, \mathsf{E}_m) = \Big( \mathrm{M}_m
\otimes \mathcal{M}, \mathrm{M}_m, 1_{\mathrm{M}_m} \otimes
\varphi \Big).
\end{equation}
If $\mathrm{D}_{\varphi}$ is the associated density, we consider
the densely defined maps
\begin{eqnarray*}
\rho_r: \mathrm{D}_{\varphi}^{1/2} \big( x_{ij} \big) \in
S_{\infty}^m \big( L_2^r(\mathcal{M}) \big) & \mapsto & \big(
x_{ij} \big) \in L_{\infty}^r (\mathcal{M}_m, \mathsf{E}_m),
\\ \rho_c: \big( x_{ij} \big)
\mathrm{D}_{\varphi}^{1/2} \in S_{\infty}^m \big(
L_2^c(\mathcal{M}) \big) & \mapsto & \big( x_{ij} \big) \in
L_{\infty}^c (\mathcal{M}_m, \mathsf{E}_m).
\end{eqnarray*}

\begin{lemma} \label{Lemma-Isometry-Rho}
The maps $\rho_r$ and $\rho_c$ extend to isometric isomorphisms.
\end{lemma}

\begin{proof} By \cite[p.56]{ER} we have
\begin{eqnarray*}
\quad \big\| \mathrm{D}_{\varphi}^{1/2} \big( z_{ij} \big)
\big\|_{S_{\infty}^m(L_2^r(\mathcal{M}))} & = & \Big\| \Big(
\sum_{k=1}^m \mbox{tr}_{\mathcal{M}} \big(
\mathrm{D}_{\varphi}^{1/2} z_{ik}^{}z_{jk}^*
\mathrm{D}_{\varphi}^{1/2} \big) \Big) \Big\|_{\mathrm{M}_m}^{1/2}
\\ & = & \Big\| \Big( \sum_{k=1}^m \varphi
(z_{ik}^{}z_{jk}^*) \Big) \Big\|_{\mathrm{M}_m}^{1/2} \\ & = &
\Big\| 1_{\mathrm{M}_m} \otimes \varphi \Big[ \big( z_{ij} \big)
\big( z_{ij} \big)^* \Big] \Big\|_{\mathrm{M}_m}^{1/2} = \big\|
\big( z_{ij} \big) \big\|_{L_{\infty}^r(\mathcal{M}_m,
\mathsf{E}_m)}.
\end{eqnarray*}
The proof of the second isometric isomorphism follows from
\cite[p.54]{ER} instead. \end{proof}

Now, assuming that $\mathcal{M}$, $\mathcal{N}$ and $\mathsf{E}$
are given as in (\ref{Equation-MNE}), our aim is to identify the
intersection space $\mathcal{J}_{p,q}^n(\mathcal{M}_m,
\mathsf{E}_m)$ in terms of asymmetric $L_p$ spaces. More
concretely, let us define the following intersection of asymmetric
spaces
$$\mathcal{J}_{p,q}^n(\mathcal{M}) = \bigcap_{u,v \in\{2p,2q\}}
n^{\frac{1}{u} + \frac{1}{v}} \, L_{(u,v)}(\mathcal{M}).$$

\begin{lemma} \label{Lemma-Assym-Cond}
If $m \ge 1$, we have an isometry
$$S_p^m \big( \mathcal{J}_{p,q}^n (\mathcal{M}) \big) =
\mathcal{J}_{p,q}^n(\mathcal{M}_m, \mathsf{E}_m).$$
\end{lemma}

\begin{proof}
By definition we have $$\mathcal{J}_{p,q}^n(\mathcal{M}_m,
\mathsf{E}_m) = \bigcap_{u,v \in \{\frac{2pq}{p-q}, \infty\}}
n^{\frac1u + \frac1p + \frac1v} \, L_{(u,v)}^p(\mathcal{M}_m,
\mathsf{E}_m).$$ Since the powers of $n$ fit, it clearly suffices
to show that
\begin{eqnarray*}
S_p^m \big( L_{(2p,2p)}(\mathcal{M}) \big) & = &
L_{(\infty,\infty)}^p(\mathcal{M}_m, \mathsf{E}_m), \\ [5pt] S_p^m
\big( L_{(2p,2q)}(\mathcal{M}) \big) & = & L_{(\infty,
\frac{2pq}{p-q})}^p (\mathcal{M}_m, \mathsf{E}_m), \\ [3pt] S_p^m
\big( L_{(2q,2p)}(\mathcal{M}) \big) & = &
L_{(\frac{2pq}{p-q},\infty)}^p (\mathcal{M}_m, \mathsf{E}_m), \\
[3pt] S_p^m \big( L_{(2q,2q)}(\mathcal{M}) \big) & = &
L_{(\frac{2pq}{p-q}, \frac{2pq}{p-q})}^p (\mathcal{M}_m,
\mathsf{E}_m).
\end{eqnarray*}
The first isometry follows from Remark \ref{Rem-Lp-2p2p}, we have
$L_p(\mathcal{M}_m)$ at both sides. The last one follows from
Example \ref{Remark-Conditional-I-II} (b), which uses one of
Pisier's identities stated in Chapter \ref{Section1}. It remains
to see the second and third isometries. By complex interpolation
on $q$, it suffices to assume $q=1$ since the case $q=p$ has been
already considered. In that case, we have to show that
\begin{eqnarray*}
S_p^m \big( L_{(2p,2)}(\mathcal{M}) \big) & = & L_{(\infty,
2p')}^p(\mathcal{M}_m, \mathsf{E}_m), \\ [3pt] S_p^m \big(
L_{(2,2p)}(\mathcal{M}) \big) & = &
L_{(2p',\infty)}^p(\mathcal{M}_m, \mathsf{E}_m).
\end{eqnarray*}
When $p=1$, the isometry follows again from Remark
\ref{Rem-Lp-2p2p}. Therefore, by complex interpolation one more
time, it suffices to assume that $(p,q) = (\infty,1)$. In that
case, we note that $$S_\infty^m \big( L_{(\infty,2)}(\mathcal{M})
\big) = S_\infty^m \big( L_2^c(\mathcal{M}) \big) =
L_\infty^c(\mathcal{M}_m, \mathsf{E}_m) = L_{(\infty,
2)}^\infty(\mathcal{M}_m, \mathsf{E}_m),$$
$$S_\infty^m \big( L_{(2,\infty)}(\mathcal{M}) \big) = S_\infty^m
\big( L_2^r(\mathcal{M}) \big) = L_\infty^r(\mathcal{M}_m,
\mathsf{E}_m) = L_{(2,\infty)}^\infty(\mathcal{M}_m,
\mathsf{E}_m).$$ We have used Lemma \ref{Lemma-Isometry-Rho} in
the second identities. This completes the proof.
\end{proof}

Let $\mathcal{A}_m$ stand for the usual amalgamated free product
von Neumann algebra constructed out of $\mathcal{M}_m =
\mathrm{M}_m \otimes \mathcal{M}$ with $\mathcal{N}_m =
\mathrm{M}_m$. According to (\ref{Equation-CB-Amalgamated}), we
have $\mathcal{A}_m = \mathrm{M}_ m \otimes \mathcal{A}$ with
$\mathcal{A} = (\mathcal{M} \oplus \mathcal{M})^{* n}$. In
particular, $L_p(\mathcal{A}_m; \ell_q^n) = S_p^m \big(
L_p(\mathcal{A}; \ell_q^n) \big)$ and we deduce the result below
from Theorem \ref{Theorem-Embed-OH} and Lemma
\ref{Lemma-Assym-Cond}.

\begin{theorem} \label{Theorem-Sigma-2Infty}
If $1 \le q \le p \le \infty$, the map
$$u: x \in \mathcal{J}_{p,q}^n(\mathcal{M}) \mapsto \sum_{k=1}^n
x_k \otimes \delta_k \in L_p(\mathcal{A}; \ell_q^n)$$ is a
cb-isomorphism with cb-complemented image and constants
independent of $n$.
\end{theorem}

\begin{remark} \label{Seacabo}
\emph{Let $\mathcal{M}$ be the matrix algebra $\mathrm{M}_m$
equipped with its normalized trace $\tau$ and let us consider the
direct sum $\mathcal{M}_{\oplus n} = \mathcal{M} \oplus
\mathcal{M} \oplus \cdots \oplus \mathcal{M}$ with $n$ terms. We
equip this algebra with the (non-normalized) trace $\tau_n = \tau
\oplus \tau \oplus \cdots \oplus \tau$. In this case, given an
operator space $\mathrm{X}$, we could also define the space
$$\mathcal{J}_{p,q}(\mathcal{M}_{\oplus n}; \mathrm{X}) =
\bigcap_{u,v \in \{2p,2q\}} L_u^r(\mathcal{M}_{\oplus n})
\otimes_{\mathcal{M}_{\oplus n},h} L_{\infty}(\mathcal{M}_{\oplus
n}; \mathrm{X}) \otimes_{\mathcal{M}_{\oplus n},h}
L_v^c(\mathcal{M}_{\oplus n}).$$ The analogue of Theorem
\ref{Theorem-Sigma-2Infty} with
\begin{itemize}
\item[i)] $(p,q)$ being $(p,1)$, \item[ii)] freeness replaced by
noncommutative independence, \item[iii)] vector-values in some
operator space $\mathrm{X}$ as explained above,
\end{itemize}
was the main result in \cite{JP}. In our situation, a
vector-valued analogue of Theorem \ref{Theorem-Sigma-2Infty} also
holds in the context of finite matrix algebras. However, note that
the use of free probability requires to define vector-valued $L_p$
spaces for free product von Neumann algebras \cite{J4}, which is
beyond the scope of this paper.}
\end{remark}

\bibliographystyle{amsplain}

\end{document}